\documentclass[12pt,b5paper,notitlepage]{report}
\usepackage[b5paper, margin={0.5in,0.65in}]{geometry}

\usepackage{amsmath,amscd,amssymb,amsthm,mathrsfs,amsfonts,layout,indentfirst,graphicx,caption,mathabx, stmaryrd,appendix,calc,imakeidx,upgreek} 
\usepackage{palatino}  
\usepackage{slashed} 
\usepackage{mathrsfs} 
\usepackage{extarrows} 
\usepackage{enumitem} 

\usepackage{fancyhdr} 

\usepackage[nottoc]{tocbibind}   

\makeindex

\usepackage{lipsum}
\let\OLDthebibliography\thebibliography
\renewcommand\thebibliography[1]{
	\OLDthebibliography{#1}
	\setlength{\parskip}{0pt}
	\setlength{\itemsep}{2pt} 
}

\allowdisplaybreaks  
\usepackage{latexsym}
\usepackage{chngcntr}
\usepackage[colorlinks,linkcolor=blue,anchorcolor=blue,linktocpage
]{hyperref}
\hypersetup{citecolor=[rgb]{0,0.5,0}}

\counterwithin{figure}{section}

\theoremstyle{definition}
\newtheorem{df}{Definition}[section]
\newtheorem{eg}[df]{Example}
\newtheorem{rem}[df]{Remark}

\newtheorem{cv}[df]{Convention}

\theoremstyle{plain}
\newtheorem{thm}[df]{Theorem}

\newtheorem{pp}[df]{Proposition}
\newtheorem{co}[df]{Corollary}
\newtheorem{lm}[df]{Lemma}

\newcommand{\fk}{\mathfrak}
\newcommand{\mc}{\mathcal}
\newcommand{\wtd}{\widetilde}
\newcommand{\wht}{\widehat}
\newcommand{\wch}{\widecheck}
\newcommand{\ovl}{\overline}
\newcommand{\tr}{\mathrm{t}} 

\newcommand{\End}{\mathrm{End}} 
\newcommand{\id}{\mathbf{1}}
\newcommand{\Hom}{\mathrm{Hom}}
\newcommand{\Conf}{\mathrm{Conf}}
\newcommand{\Res}{\mathrm{Res}}

\newcommand{\Dom}{\scr D}

\newcommand{\bk}[1]{\langle {#1}\rangle}

\newcommand{\scr}{\mathscr}

\newcommand{\xk}{\mathfrak x}
\newcommand{\yk}{\mathfrak y}
\newcommand{\zk}{\mathfrak z}

\newcommand{\im}{\mathbf{i}}

\newcommand{\shom}{\underline{\Hom}}

\newcommand{\sgm}{\varsigma}
\newcommand{\SX}{S_{\fk X}}
\newcommand{\DX}{D_{\fk X}}
\newcommand{\mbb}{\mathbb}
\newcommand{\mbf}{\mathbf}
\newcommand{\blt}{\bullet}
\newcommand{\coker}{\mathrm{coker}}
\newcommand{\Vbb}{\mathbb V}
\newcommand{\Ubb}{\mathbb U}
\newcommand{\Wbb}{\mathbb W}
\newcommand{\Mbb}{\mathbb M}
\newcommand{\Gbb}{\mathbb G}
\newcommand{\Cbb}{\mathbb C}
\newcommand{\Nbb}{\mathbb N}
\newcommand{\Zbb}{\mathbb Z}
\newcommand{\Pbb}{\mathbb P}
\newcommand{\Rbb}{\mathbb R}
\newcommand{\Ebb}{\mathbb E}
\newcommand{\cbf}{\mathbf c}
\newcommand{\wt}{\mathrm{wt}}
\newcommand{\Lie}{\mathrm{Lie}}
\newcommand{\btl}{\blacktriangleleft}
\newcommand{\btr}{\blacktriangleright}
\newcommand{\svir}{\mathcal V\!\mathit{ir}}
\newcommand{\Ker}{\mathrm{Ker}}

\newcommand{\Sbf}{\mathbf{S}}
\newcommand{\low}{\mathrm{low}}
\newcommand{\Crm}{\mathrm C}
\newcommand{\Brm}{\mathrm B}

\numberwithin{equation}{section}

\title{Conformal Blocks:\\[0.5ex] \large Vector bundle structures, Sewing, and Factorization}
\author{{\sc Bin Gui}
}
\date{Last update: June 15, 2026}
\begin{document}\sloppy 
	\pagenumbering{arabic}

	\maketitle

\noindent\textbf{Note.}
The main body of this monograph was completed and made available on the author's website in November 2020. Since then, only minor revisions and corrections have been made. 

This monograph has been uploaded to arXiv because it has been cited by others. However, it was originally intended as an informal note on the theory of conformal blocks, and a large part of its content has since been adapted and expanded in works of the author and his collaborator, notably
\href{https://arxiv.org/abs/2011.07450}{arXiv:2011.07450},
\href{https://arxiv.org/abs/2110.04774}{arXiv:2110.04774},
\href{https://arxiv.org/abs/2305.10180}{arXiv:2305.10180}, and
\href{https://arxiv.org/abs/2411.07707}{arXiv:2411.07707}.
In particular, readers should refer to these official works for formal mathematical treatments and references.

This monograph has not been fully polished or elaborated, and may contain immature or even incorrect mathematical ideas and results. The references are also imprecise and incomplete. Use at your own risk!\\




\makeatletter
\newcommand*{\toccontents}{\@starttoc{toc}}
\makeatother
\toccontents
	







\chapter*{Introduction}

Conformal blocks are central objects in  $2$-dimensional conformal field theory (CFT). Mathematically rigorous definitions of conformal blocks were first introduced for special examples (e.g. minimal models, Weiss-Zumino-Witten models, see \cite{BFM91,TUY89}), and were later given for general vertex operator algebras (VOAs) (see \cite{Zhu94,FB04}). Conformal blocks  are the building blocks of correlation functions of CFT. They are mainly discussed in the literature of algebraic geometry, although they have important applications also to many other areas related to CFT, such as low-dimensional topology, tensor categories, VOAs, von Neumann algebras and subfactors, etc.. Thus, we believe it is worthwhile to introduce the beautiful theory of conformal blocks to the people working in these areas without assuming they have previous knowledge in algebraic geometry. This is a main goal of the present monograph.

We shall give a comprehensive exploration of the theory of conformal blocks in the framework of VOAs and complex \emph{analytic} geometry. Unlike most approaches, we work in the complex analytic setting rather than algebraic one, using the language of complex manifolds and complex spaces\footnote{Nodal curves and their open subsets are the only singular complex spaces (i.e. complex spaces which are not complex manifolds) we will consider in this monograph.} rather than schemes or stacks. This is partly due to the author's own taste, but also due to the following reason: Despite that an analytic theory of conformal blocks is necessary for application to many areas,  some  results of conformal blocks (e.g. proving that spaces of conformal blocks form a holomorphic vector bundle) cannot be directly translated without suitable adaption from algebraic to  analytic setting, and certain results (e.g. convergence of sewing conformal blocks) can only be proved using analytic methods. Let me explain this in more details.

\subsection*{Vector bundle structures}

Given a (CFT-type) VOA $\Vbb$,   an $N$-pointed compact Riemann surface $\fk X=(C;x_1,\dots,x_N)$ (i.e. $C$ has $N$ distinct marked points $x_1,\dots,x_N$) and $\Vbb$-modules $\Wbb_1,\dots,\Wbb_N$, a conformal block $\upphi$ is a linear functional on $\Wbb_\blt:=\Wbb_1\otimes\cdots\otimes\Wbb_N$ invariant under the action of the sections of \textbf{sheaf of VOAs} $\scr V_C$ defined using $\Vbb$ and $C$. The vector space of all conformal blocks associated to $C$ and these $\Vbb$-modules is called the \textbf{space of conformal blocks} $\scr T_{\fk X}^*(\Wbb_\blt)$. Assume that $\Vbb$ satisfies $C_2$-cofinite property, a natural finiteness condition introduced by Zhu \cite{Zhu96}. Then  $\scr W_{\fk X}(\Wbb_\blt)$ is expected to be finite-dimensional. Moreover, one should expect that its dimension is independent of the complex structure of $C$ and the positions of marked points, and these vector spaces should form a holomorphic vector bundle over the moduli space of $N$-pointed compact Riemann surfaces (with possibly extra data). In other words, suppose we have a holomorphic family of $N$-pointed compact Riemann surfaces $\fk X=(\pi:\mc C\rightarrow\mc B;\sgm_1,\dots,\sgm_N)$ where $\mc C,\mc B$ are complex manifolds, $\sgm_1,\dots,\sgm_N:\mc B\rightarrow\mc C$ are sections (i.e. families of mark points). Let $\mc C_b$ be the fiber $\pi^{-1}(b)$ for each $b\in\mc B$, which is a compact Riemann surface. Let $\fk X_b$ be $\mc C_b$ with marked points $\sgm_1(b),\dots,\sgm_N(b)$. Then $b\in\mc B\mapsto\dim\scr T_{\fk X_b}^*(\Wbb_\blt)$ should be locally constant, and the vector spaces $\scr T_{\fk X_b}^*(\Wbb_\blt)$ (for all $b\in\mc B$) should form a holomorphic vector bundle over $\mc B$.

A usual way of constructing vector bundle structures for spaces of conformal blocks is to first define \textbf{sheaf of covacua} $\scr T_{\fk X}(\Wbb_\blt)$. This is an $\scr O_{\mc B}$-module (where $\scr O_{\mc B}$ is the structure sheaf of $\mc B$), which is locally a quotient of the sheaf $\scr W_{\fk X}(\Wbb_\blt)$ of $\Wbb_\blt=\Wbb_1\otimes\cdots\otimes\Wbb_N$-valued holomorphic functions. Its dual module $\scr T_{\fk X}^*(\Wbb_\blt)$ is called \textbf{sheaf of conformal blocks}. Using some basic results in complex analytic or algebraic geometry, one can identify the fibers of $\scr T_{\fk X}(\Wbb_\blt)$ with the spaces of covacua (the dual spaces of the spaces of conformal blocks) in a natural way (Thm. \ref{lb69}). Thus, once we have proved that $\scr T_{\fk X}(\Wbb_\blt)$ is locally free (of finite rank), i.e., $\scr T_{\fk X}(\Wbb_\blt)$ is a vector bundle, then $\scr T_{\fk X}^*(\Wbb_\blt)$ is also locally free, which is the vector bundle structure we are looking for. 

In the algebraic setting, one proves  that $\scr T_{\fk X}(\Wbb_\blt)$ admits  local connections \cite{FB04,DGT19a}, and that $\scr T_{\fk X}(\Wbb_\blt)$ is coherent \cite{DGT19b}. Then a standard argument shows that $\scr T_{\fk X}(\Wbb_\blt)$ is locally free. It is obvious that $\scr T_{\fk X}(\Wbb_\blt)$ is quasi-coherent, i.e., is the cokernal of a morphism between two possibly infinite-rank locally free sheaves. Thus, once we can show that $\scr T_{\fk X}(\Wbb_\blt)$ is a finite-type $\scr O_{\mc B}$-module, we can conclude that $\scr T_{\fk X}(\Wbb_\blt)$ is coherent, thanks to the fact that $\scr O_{\mc B}$ is Noetherian.  In the analytic setting, $\scr O_{\mc B}$ is not Noetherian, and $\scr T_{\fk X}(\Wbb_\blt)$ is not quasi-coherent (in the sense of \cite{EP96}). Thus, although one can still define connections and show that $\scr T_{\fk X}(\Wbb_\blt)$ is finite-type, one cannot conclude that $\scr T_{\fk X}(\Wbb_\blt)$ is coherent or locally free. In this monograph, we fix this issue by proving a stronger finiteness theorem for $\scr T_{\fk X}(\Wbb_\blt)$ (Thm. \ref{lb67}), and show that this result together with the existence of connections imply the local freeness (Thm. \ref{lb66}).

\subsection*{Convergence of sewing}

Suppose  we have an $(N+2)$-pointed compact Riemann surface
\begin{align*}
\wtd{\fk X}=(\wtd C;x_1,\dots,x_N,x',x'').
\end{align*}
Then we can sew $\wtd{\fk X}$ along the pair of points $x',x''$ to obtain another Riemann surface with possibly higher genus. More precisely, we choose $\xi,\varpi$ to be local coordinates of $\wtd C$ at $x',x''$. Namely, they are univalent  (i.e. holomorphic and injective) functions defined respectively in neighborhoods $U'\ni x',U''\ni x''$ satisfying $\xi(x')=0,\varpi(x'')=0$. For each $r>0$ we let $\mc D_r=\{z\in\Cbb:|z|<r\}$ and $\mc D_r^\times=\mc D_r-\{0\}$. We choose $r,\rho>0$ so that the neighborhoods $U',U''$ can be chosen to satisfy that $\xi(U')=\mc D_r$ and $\varpi(U'')=\mc D_\rho$, that  $U'\cap U''=\emptyset$, and that none of $x_1,\dots,x_N$ is in $U'$ or $U''$. Then, for each $q\in\mc D_{r\rho}^\times$, we remove the closed subdiscs of $U',U''$ determined respectively by $|\xi|\leq \frac {|q|}\rho$ and $|\varpi|\leq \frac{|q|}r$, and glue the remaining part using the relation $\xi\varpi=q$. Then we obtain an $N$-pointed compact Riemann surface
\begin{align*}
\fk X_q=(\mc C_q;x_1,\dots,x_N)
\end{align*}
which clearly depends on $\xi$ and $\varpi$. By varying $q$, we obtain a family of $N$-pointed compact Riemann surfaces $\fk X=(\pi:\mc C\rightarrow\mc D_{r\rho}^\times;\sgm_1,\dots,\sgm_N)$.

Now, if we associate $\Vbb$-modules $\Wbb_1,\dots,\Wbb_N,\Mbb,\Mbb'$ (where $\Mbb'$ is the contragredient (i.e. dual) module of $\Mbb$) to $x_1,\dots,x_N,x',x''$, and choose a conformal block $\uppsi$ associated to $\wtd{\fk X}$ and these $\Vbb$-modules, then its \textbf{sewing} $\mc S\uppsi$ is an $\Wbb_\blt^*=(\Wbb_1\otimes\cdots\Wbb_N)^*$-valued formal series of $q$ defined by sending each $w_\blt=w_1\otimes\cdots\otimes w_N\in\Wbb_\blt$ to
\begin{align*}
\mc S\uppsi(w_\blt)=\uppsi(w_\blt\otimes q^{L_0}\btr\otimes\btl)\quad\in\Cbb\{q\}
\end{align*}
where $\btr\otimes\btl$ is the element of the ``algebraic completion" of $\Mbb\otimes\Mbb'$ corresponding to the identity element of $\End_\Cbb(\Mbb)$, and $L_0$ is the zero mode of the Virasoro operators $\{L_n:n\in\Zbb\}$. The sewing problem is about proving that $\mc S\uppsi(w_\blt)$ converges absolutely to a (possibly) multivalued function on $\mc D_{r\rho}^\times$. Moreover, for each $q\in\mc D_{r\rho}^\times$, $\mc S\uppsi(\cdot,q)$ defines a conformal block associated to $\fk X_q$ and $\Wbb_1,\dots,\Wbb_N$. If we sew $\wtd C$ along $n$ pairs of points, and if we let  $x_1,\dots,x_N$ and  $\wtd C$ and $\uppsi$    vary and be parametrized holomorphically by variables $\tau_\blt=(\tau_1,\dots,\tau_m)$ (see Sec. \ref{lb24} for details), then  sewing  conformal blocks is also absolutely  convergent with respect to $q_1,\dots,q_n$ and (locally) uniform with respect to $\tau_\blt$. 	The sewing problem is analytic by nature. It cannot be proved using purely algebro-geometric method since $\fk X$ is not an algebraic family; in particular, $\mc D_{r\rho}^\times$ is not an algebraic variety or scheme, and is not considered in algebraic geometry. In the monograph, we will give a detailed proof of the  sewing problem using analytic methods. See Section \ref{lb147}.

\section*{Factorization}

In the above setting, the factorization property of conformal blocks says that if $\Vbb$ is $C_2$-cofinite and rational (which means that certain classes of generalized $\Vbb$-modules are completely reducible), then  for any $q\in\mc D_{r\rho}^\times$, any conformal block associated to $\fk X_q$ and $\Wbb_1,\dots,\Wbb_N$ is a sum of $\mc S\uppsi(q)$ where $\uppsi\in\scr T_{\wtd{\fk X}}^*(\Wbb_\blt\otimes\Mbb\otimes\Mbb')$ (i.e. $\uppsi$ is a conformal block associated to $\wtd{\fk X}$ and $\Wbb_1,\dots,\Wbb_N,\Mbb,\Mbb'$) and $\Mbb$ is simple.  One can phrase the factorization as a relation between the dimensions of spaces of conformal blocks associated to $\wtd{\fk X}$ and to $\fk X_q$ as follows. For each equivalence class of simple $\Vbb$-module we choose a representative and let them form a set $\mc E$. Then for each $q\in\mc D_{r\rho}^\times$,
\begin{align*}
\dim\scr T_{\fk X_q}^*(\Wbb_\blt)=\sum_{\Mbb\in\mc E}\dim \scr T_{\wtd{\fk X}}^*(\Wbb_\blt\otimes\Mbb\otimes\Mbb').
\end{align*}
Factorization in this form was proved in \cite{TUY89} for WZW models and in \cite{DGT19b} for any $C_2$-cofinite rational VOA. Note that the left hand side is independent of $q$, and is a priori no less than the right hand side since the linear map $\uppsi\mapsto\mc S\uppsi(q)$ is indeed injective (Thm. \ref{lb114}).

Let us call the right hand side of this equation to be $D$. So $D\leq \dim\scr T_{\fk X_q}^*(\Wbb_\blt)$. To prove the factorization, we add $\fk X_0$, which is an $N$-pointed nodal curve, to the family $\fk X$. Then $\fk X$ is a family of $N$-pointed complex curves with base manifold $\mc D_{r\rho}$. We can still define the sheaf of covacua $\scr T_{\fk X}(\Wbb_\blt)$ and show that it is finitely generated (Thm. \ref{lb67}). Then by Nakayama's lemma (Prop. \ref{lb60}), the dimension of fibers of $\scr T_{\fk X}(\Wbb_\blt)$ (which can be identified with spaces of conformal blocks (Thm. \ref{lb69})) is upper-semicontinuous with respect to $q\in \mc D_{r\rho}$. Thus, for each $q\neq 0$, $\dim\scr T_{\fk X_q}^*(\Wbb_\blt)\leq \dim\scr T_{\fk X_0}^*(\Wbb_\blt)$. Therefore, it suffices to prove the \textbf{nodal factorization}:  that $\dim\scr T_{\fk X_0}^*(\Wbb_\blt)\leq D$.

To prove the nodal factorization, we should realize each $\upphi\in\scr T_{\fk X_0}^*(\Wbb_\blt)$ as an element of $\bigoplus_{\Mbb\in\mc E}\scr T_{\wtd{\fk X}}^*(\Wbb_\blt\otimes\Mbb\otimes\Mbb')$. For that purpose, we consider $\wtd{\fk X}$ as a Riemann surface with input points $x_1,\dots,x_N$ and output points $x',x''$.  Then we can define the \textbf{dual $\wtd{\fk X}$-tensor product} of $\Wbb_1,\dots,\Wbb_N$ (associated to $x_1,\dots,x_N$ respectively) to be a vector space $\boxbackslash_{\wtd{\fk X}}(\Wbb_\blt)$, which is the subspace of linear functional on $\Wbb_\blt$ satisfying certain properties that $\upphi$ satisfy. The difficulty of this approach is to define an action of $\Vbb\times\Vbb$ on $\boxbackslash_{\wtd{\fk X}}(\Wbb_\blt)$ which makes the latter a weak $\Vbb\times\Vbb$-module. Then using the $C_2$-cofiniteness and rationality, it is not hard to show that  $\boxbackslash_{\wtd{\fk X}}(\Wbb_\blt)$ (or a suitable $\Vbb\times\Vbb$-submodule containing $\upphi$) is a direct sum of irreducible $\Vbb\times\Vbb$ modules which must be of the form $\Mbb\otimes\Mbb'$ for some $\Mbb\in\mc E$. Nodal factorization follows.

\cite{DGT19b} defines the weak $\Vbb\times\Vbb$-module structure on $\boxbackslash_{\wtd{\fk X}}(\Wbb_\blt)$ using Zhu's algebras. Our approach does not use Zhu's algebra, but relies heavily on the \textbf{propagations of conformal blocks} (Sec.\ref{lb149}) and more generally, \textbf{propagations  of dual tensor product elements} (Sec. \ref{lb110}). Moreover, our proofs of propagations and some related properties are analyic and different from those in \cite{FB04} or \cite{DGT19b}: we use the sewing of an $N$-pointed compact Riemann surface and a $3$-pointed $\Pbb^1$; the convergence of the corresponding sewing of conformal blocks is due to  the strong residue theorem for \emph{families of} compact Riemann surfaces (Thm. \ref{lb18}).

\subsection*{Prerequisite and outline}

We assume the readers know some basic properties of complex manifolds; sheaves,  sheaves of modules, and their morphisms; sheaf (\v Cech) cohomology. See for instance \cite[Sec. 0.2, 0.3]{GH78}, \cite[Chapter B]{Huy06}, \cite[Annex]{GR84}. Some familiarity with computations in VOA (cf. for instance \cite{FHL93}) is helpful but not necessary. No knowledge in algebraic geometry is required.

More advanced topics in complex geometry will be discussed in Chapter 1. In particular, we review the basic properties of compact Riemann surfaces, which will be generalized to nodal curves. Since we do not assume the readers have any previous knowledge on complex spaces or nodal curves, we give complete and self-contained account of these properties. We also introduce the necessary tools for studying families of compact Riemann surfaces and, more generally, families of complex curves. We  give a detailed description of how to sew a family of compact Riemann surfaces along several pairs of points to obtain a family of complex curves. We prove strong residue theorem for families of compact Riemann surfaces, which is necessary for proving the propagations of conformal blocks. Basic properties of holomorphic differential equations are recalled, which will be used to prove the convergence of sewing conformal blocks. We also give criteria on local freeness of sheaves.

In Chapter 2 we discuss sheaves of VOAs for families of complex curves introduced in \cite{FB04} (for smooth curves) and \cite{DGT19a} (for nodal curves). These sheaves are infinite rank holomorphic vector bundles whose transition functions are discovered in \cite{Hua97}. We also give a formula for Lie derivatives of sheaves of VOAs, generalizing those of tangent fields and tensor fields. This formula is due to the author, and is used in the next chapter to define connections on sheaves of conformal blocks.

Following \cite{FB04,DGT19a}, we define in Chapter 3 sheaves of conformal blocks for families of complex curves. We define the sewing of conformal blocks which corresponds to the geometric construction of sewing families of compact Riemann surfaces. We prove that the sewing of conformal blocks are also conformal blocks in the formal sense, which is due to \cite{DGT19b}. We then prove  propagation of conformal blocks. Single propagation is due to \cite{FB04}, and multiple propagations are due to the author. Double propagations play an important role of defining weak $\Vbb\times\Vbb$-module structures on dual tensor products in Chapter 4. As mentioned previously, our treatment of single propagation is new and relies on sewing. We then prove that sheaves of conformal blocks support locally logarithmic connections. This result is due to \cite{FB04} and \cite{DGT19a} for smooth families and general families respectively. Our treatment is different from theirs and uses the result on Lie derivatives in Chapter 2. With the help of connections, we then prove that for $C_2$-cofinite VOAs, the sheaves of conformal blocks are locally free.

The first three sections of Chapter 4 are devoted to the proof that sewing conformal blocks is convergent when $\Vbb$ is $C_2$-cofinite. Our treatment of projective structures is motivated by \cite[Sec. 8.2]{FB04}. We also prove that the sewing map is injective. We then define the vector space of dual tensor products, and use (single and double) propagations of dual tensor product elements to define weak $\Vbb\times\Vbb$-module structures.  For an approach using Zhu's algebra, see \cite{DGT19b}. We then prove factorization for conformal blocks associated to $C_2$-cofinite rational VOAs, which is originally due to \cite{DGT19b}.

The connections defined locally in Chapter 3 are in general not flat but only projectively flat. In chapter 5, we explain how to slightly modify the definition and obtain flat connections on sheaves of conformal blocks. Typically, the construction of flat connections uses determinant line bundles. It turns out that tensor products of these line bundles are equivalent to sheaves of conformal blocks associated to holomorphic VOAs. We use the latter sheaves instead of determinant line bundles. This treatment is motivated by \cite{AU07a,AU07b}. We also provide in this chapter all the necessary results for constructing modular functors from conformal blocks. In the last section, we explain how the famous and mysterious factor $q^{-\frac c{24}}$ appears in genus $1$ CFT.

We remark that before \cite{FB04,DGT19a,DGT19b}, the definition of conformal blocks, and the proof of propagation, local freeness, and factorization of conformal blocks were given in \cite{BFM91,TUY89} for minimal models and WZW-models respectively, and in \cite{NT05} for general VOAs (satisfying $C_2$-cofiniteness, rationality, and some other small conditions) but only genus $0$ curves. Proofs of convergence of sewing and factorization were also given in \cite{Hua95,Hua98,Hua05a} for general VOAs as above for genus $0$ curves, and in \cite{Zhu96,Hua05b} for genus $1$ curves.

\chapter{Basics of complex geometry}

\section{Sheaves of modules}
Let us  fix some notations.  \index{N@$\mathbb N=\{0,1,2,\dots\}$}  \index{Z+@$\mathbb Z_+=\{1,2,3,\dots\}$}
\begin{gather*}
\im=\sqrt{-1}.\\
\Nbb=\{0,1,2,\dots\},\qquad \Zbb_+=\{1,2,3,\dots\}.
\end{gather*}
Throughout this monograph, $\scr O_X$ of a given complex manifold $X$ (or more generally, a complex space) always denotes the sheaf (of germs) of holomorphic functions on $X$. Thus, $\scr O_X(X)$ is the space of holomorphic functions on $X$. We will sometimes write $\scr O_X(X)$ as $\scr O(X)$ for short. For any $x\in X$, $\scr O_{X,x}$ denotes the stalk of $\scr O_X$ at $x$, and $\fk m_x$ denotes the ideal of all germs $f\in\scr O_{X,x}$ satisfying $f(x)=0$.

In general, if $\scr E$ is a sheaf on $X$ then $\scr E_x$ denotes the stalk at $x$. \index{Ex@$\scr E_x$} If $U$ is an open set containing $x$, and if $s\in\scr E(U)$, then $s_x\in\scr E_x$ denotes the germ of $s$ at $x$. \index{Ex@$\scr E_x$!$s_x$} If $\scr E$ is a (sheaf of) $\scr O_X$-module, then  $\scr E_x/\fk m_x\scr E_x$ is a complex vector space, called the \textbf{fiber} of $\scr E$ at $x$. (Fibers can also be defined using pull backs of sheaves; see Section \ref{lb1}. It will be denoted by $\scr E|x$ in the future.) \index{EX@ ${\scr E\lvert X},\scr E\lvert x$} It's dimension $r_x$ is called the \textbf{rank} of $\scr E$ at $x$. The function $r:x\in X\mapsto r_x$ is called the \textbf{rank function}.

A homomorphism of $\scr O_X$-modules $\scr E\rightarrow\scr F$ is an isomorphism (i.e., the induced homomorphism of $\scr O_X(U)$-modules $\scr F(U)\rightarrow\scr F(U)$ is an isomorphism for any open subset $U$) if any only if the corresponding stalk map $\scr F_x\rightarrow\scr G_x$ is an isomorphism for each $x\in X$. $\scr E$ is called \textbf{locally free} if  each $x\in X$ has a neighborhood $U$ such that $\scr E|_U\simeq \scr O_U^n$ for some natural number $n$. $\scr E$ is locally free if and only if it is the sheaf of germs of a holomorphic (finite rank) vector bundle. Thus, locally free sheaves and vector bundles are regarded as the same things. It is clear that the rank function is locally constant for any locally free sheaf. Unless otherwise sated, we assume that locally free sheaves  have (locally) finite rank.

If $U$ is an open subset of $X$, and $s_1,\dots,s_n\in\scr E(U)$, we say that $s_1,\dots,s_n$ \textbf{generate} $\scr E_U$, if for each $x\in U$, the stalk $\scr E_x$ is generated by (the germs of) $s_1,\dots, s_n$. This is equivalent to saying that the $\scr O_U$-module homomorphism $\scr O_U^n\rightarrow\scr E_U$ defined by 
\begin{align*}
\scr O_U^n(V)\rightarrow\scr E_U(V),\qquad (f_1,\dots,f_n)\mapsto f_1s_1+\cdots+f_ns_n
\end{align*}
(where $V$ is any open subset of $U$) is a surjective sheaf map. Also, it is equivalent to that for any $x\in U$, $V\subset U$ a neighborhood of $x$, and $s\in\scr E(V)$, there exists a neighborhood $W$ of $x$ inside $V$ such that $s=f_1s_1+\cdots+f_ns_n$ for some $f_1,\dots,f_n\in\scr O(W)$.  If the above homomorphism $\scr O^n_U\rightarrow\scr E_U$ is an isomorphism, then we say that $s_1,\dots,s_n$ \textbf{generate freely} $\scr E_U$. If $s_1,\dots,s_n$ generate $\scr E_U$, then they generate freely $\scr E_U$ if and only if for any open subset $V\subset U$ and any $f_1,\dots,f_n\in\scr O(V)$, $f_1s_1|_V+\cdots+f_ns_n|_V=0$ implies $f_1=\cdots=f_n=0$.  We say that $\scr E$ is a \textbf{finite-type} $\scr O_X$-module if each $x\in X$ is contained in a neighborhood $U$ such that $\scr E_U$ is generated by finitely many elements of $\scr E(U)$.

The above notion of generating sections can be generalized to any subset $E$ of $\scr E(U)$, i.e., that (the elements of) $E$ \textbf{generate} $\scr E_U$ if for each $x\in U$, the germs of the elements of $E$ at $x$ generate the $\scr O_{U,x}$-module $\scr E_x$. This is not the same as saying that $E$ generates (the $\scr O(U)$-module) $\scr E(U)$, which means that each element of $\scr E(U)$ is an $\scr O(U)$-linear combination of elements of $E$.

Most sheaves  we will encounter in this monograph are locally free. However, sometimes we need to consider quotients of locally free sheaves, which are not necessarily locally free. Here is the precise definition:  An $\scr O_X$-module $\scr E$ is called a \textbf{coherent}  $\scr O_X$-module (or coherent sheaf) if  each $x\in X$ is contained in a neighborhood $U$ such that the restriction $\scr E_U$ is isomorphic to $\coker(\varphi)$ where $\varphi:\scr O_U^m\rightarrow\scr O_U^n$ is a homomorphism of $\scr O_U$-modules and $m,n\in\Nbb$. A locally free $\scr O_X$-module is clearly coherent.

Let $\scr E$ and $\scr F$ be $\scr O_X$-modules. For any open $U$ in $X$, let $\scr E_U$ and $\scr F_U$ be respectively the restrictions of $\scr E$ and $\scr F$ to $U$. Let $\Hom_{\scr O_U}(\scr E_U,\scr F_U)$ \index{Hom@$\Hom_{\scr O_U}(\scr E_U,\scr F_U)$} be the set of $\scr O_U$-morphisms from $\scr E$ to $\scr F$. So any element  $\phi\in\Hom_{\scr O_U}(\scr E_U,\scr F_U)$ is
described as follows. For any open $V\subset U$, we have an $\scr O(V)$-module homomorphism $\phi=\phi_V:\scr E(V)\rightarrow \scr F(V)$. $\phi$ is compatible with the restriction of sections, i.e., for any open $W\subset V\subset U$ and $s\in\scr E(V)$, we have 
\begin{align*}
\phi(s)|_W=\phi(s|_W).
\end{align*}
For each $x\in V$, $\phi$ induces homomorphisms of $\scr O_{X,x}$-modules and vector spaces
\begin{gather}
\phi:\scr E_x\rightarrow\scr F_x,\qquad \phi:\scr E_x/\fk m_x\scr E_x\rightarrow \scr F_x/\fk m_x\scr F_x.\label{eq223}
\end{gather}
It is clear that
\begin{align}
\phi(\scr E)_x=\phi(\scr E_x)	
\end{align}
where $\phi(\scr E)$ is the image sheaf. For each section $s$ of $\scr E$ defined near $x$, if we let $s_x$ and $s(x)$ denote the values of $s$ in $\scr E_x$ and $\scr E_x/\fk m_x\scr E_x$, and adopt similar notations for $\phi(s)$, then
\begin{gather}
\phi(s_x)=\phi(s)_x,\qquad \phi(s(x))=\phi(s)(x).\label{eq224}
\end{gather}

Note that  $\Hom_{\scr O_U}(\scr E_U,\scr F_U)$ is clearly an $\scr O(U)$ module. Then we have the so called sheaf of $\scr O_X$-homomorphisms $\shom_{\scr O_X}(\scr E,\scr F)$, \index{Hom@$\shom_{\scr O_X}(\scr E,\scr F)$} where for any open  $U\subset X$, $\shom_{\scr O_X}(\scr E,\scr F)(U)=\Hom_{\scr O_U}(\scr E_U,\scr F_U)$ whose sections are all those $\phi_U$. $\shom_{\scr O_X}(\scr E,\scr F)$ is obviously an $\scr O_X$-module.  We call 
\begin{align*}
\scr E^*:=\shom_{\scr O_X}(\scr E,\scr O_X)
\end{align*}
the \textbf{dual sheaf} of $\scr E$. \index{E@$\scr E^*$} Choose $s\in\scr E(U)$ and $t\in\scr E^*(U)$, then $t(s)$, as a section in $\scr O_X(U)$, is also denoted by $\bk{s,t}$ or $\bk{t,s}$.  Note that when $\scr E$ is locally free, then $\scr E^*$ is also locally free, and $\scr E^*$ is dual to $\scr E$ as holomorphic vector bundles. In particular, $\scr E^{**}$ can be naturally identified with $\scr E$. 

The collection $\{\scr E(U)\otimes_{\scr O(U)}\scr F(U)\}$ over all open $U\subset X$ forms a presheaf of $\scr O_X$-modules. The restriction of sections of this presheaf is defined in an obvious way. Its sheafification  $\scr E\otimes_{\scr O_X}\scr F$ \index{EF@$\scr E\otimes_{\scr O_X}\scr F=\scr E\otimes\scr F$}, which is clearly an $\scr O_X$-module, is called the tensor product of $\scr E$ and $\scr F$ over $\scr O_X$. Unless otherwise stated, we will write $\scr E\otimes_{\scr O_X}\scr F$ as  $\scr E\otimes\scr F$ for short. If we consider the tensor product over $\Cbb$, we will write $\scr E\otimes_\Cbb\scr F$ instead. Note that when the two $\scr O_X$-modules are locally free, their tensor product is nothing but the tensor product of holomorphic vector bundles. Another easy fact is the isomorphms $\scr E\otimes\scr O_X\simeq\scr E\simeq\scr O_X\otimes\scr E$ in a natural way. We write $\scr E^{\otimes n}$ as the $n$-th tensor power of $\scr E$ \index{En@$\scr E^{\otimes n}$} for any $n\in\mathbb N$.  Note also that $\scr L\otimes \scr L^*\simeq\scr O_X$ for a (holomorphic) line bundle $\scr L$  (i.e., $\scr L$ is a  rank $1$ locally free $\scr O_C$-module). In this case we write $\scr L^{-1}=\scr L^*$, \index{L-1@$\scr L^{-1}$} and more generally, $\scr L^{\otimes(-n)}=(\scr L^*)^{\otimes n}$.

A useful method of constructing sheaves is called \textbf{gluing}. Let $(U_\alpha)_{\alpha\in\fk A}$ be an open cover of $X$. Suppose that for each $\alpha\in\fk A$, we have an $\scr O_{U_\alpha}$-module $\scr E^\alpha$, that for any $\alpha,\beta\in\fk A$, we have an $\scr O_{U_\alpha\cap U_\beta}$-module isomorphism $\phi_{\beta,\alpha}:\scr E^\alpha_{U_\alpha\cap U_\beta}\xrightarrow{\simeq}\scr E^\beta_{U_\alpha\cap U_\beta}$, that $\phi_{\alpha,\alpha}=\id$, and that $\phi_{\gamma,\alpha}=\phi_{\gamma,\beta}\phi_{\beta,\alpha}$ when restricted to $U_\alpha\cap U_\beta\cap U_\gamma$. Then we can define a sheaf $\scr E$ on $X$ as follows. For any open $V\subset X$, $\scr E(V)$ is the set of all $(s_\alpha)_{\alpha\in\fk A}\in\prod_{\alpha\in\fk A}\scr E^\alpha(V\cap U_\alpha)$ (where each component $s_\alpha$ is in $\scr E^\alpha(V\cap U_\alpha)$) satisfying that $s_\beta|_{U_\alpha\cap U_\beta}=\phi_{\beta,\alpha}(s_\alpha|_{U_\alpha\cap U_\beta})$ for any $\alpha,\beta\in\fk A$. If $W$ is an open subset of $V$, then the restriction $\scr E(V)\rightarrow\scr E(W)$ is defined by that of $\scr E^\alpha(V\cap U_\alpha)\rightarrow \scr E^\alpha(W\cap U_\alpha)$. The action of $\scr O(V)$ on $\scr E(V)$ is defined by the one of $\scr O(V\cap U_\alpha)$ on $\scr E^\alpha(V\cap U_\alpha)$. It is easy to see that $\scr E$ is a sheaf of $\scr O_X$-modules. Moreover, for each $\beta\in\fk A$, we have a canonical isomorphism (trivialization) $\phi_\beta:\scr E_{U_\beta}\xrightarrow{\simeq} \scr E^\beta_{U_\beta}$ defined by $(s_\alpha)_{\alpha\in\fk A}\mapsto s_\beta$. It is clear that for each $\alpha,\beta\in\fk A$, we have $\phi_\beta=\phi_{\beta,\alpha}\phi_\alpha$ when restricted to $U_\alpha\cap U_\beta$.

For instance, any locally free sheaf is obtained by gluing a collection of free sheaves associated to an open cover of $X$.

\section{Compact Riemann surfaces}\label{lb4}

\subsection*{Serre duality}
 
Let $C$ be a compact Riemann surface, and let $\scr E$ be a (sheaf of) locally free $\scr O_C$-module. We list some basic facts about the cohomology groups of $\scr E$. Choose $q\in\mathbb N$. The first important fact is that $H^q(C,\scr E)$ is finite dimensional. Moreover, $\dim H^q(C,\scr E)=0$ when $q$ is greater than $1$, the complex dimension of $C$. These facts can be proved by Hodge theory. Moreover, let $\omega_C$ be the \textbf{dualizing sheaf} (also called canonical line bundle) of $C$.\index{zz@$\omega_C$} In other words, $\omega_C$ is the line bundle of holomorphic $1$-forms on $C$. So the sections of $\omega_C$ look locally like $f(z)dz$ where $f$ is a holomorphic function. Then \textbf{Serre duality} (which can also follow from Hodge theory) says that for any $p\in\{0,1\}$ there is an isomorphism of vector spaces
\begin{align}
H^1(C,\scr E\otimes \omega_C^p)\simeq H^0(C,\scr E^*\otimes \omega_C^{1-p})^*,\label{eq2}
\end{align} 
where $\omega_C^0=\scr O_C$ and $\omega_C^1=\omega_C$. (Cf. \cite{Huy06} Proposition 4.1.16.) In other words, there is a perfect pairing
\begin{align}
\bk{\cdot,\cdot}:H^1(C,\scr E\otimes \omega_C^p)\otimes H^0(C,\scr E^*\otimes \omega_C^{1-p})\rightarrow \mathbb C.
\end{align}
We now describe such a pairing, called the \textbf{residue pairing}.  This explicit description will be used in the proof of strong residue theorem.

Recall that the \v{C}ech cohomology group $H^1(C,\scr E\otimes \omega_C^p)$ is the direct limit of $H^1(\fk U,\scr E\otimes \omega_C^p)={Z^1(\fk U,\scr E\otimes \omega_C^p)}/{B^1(\fk U,\scr E\otimes \omega_C^p)}$ over all open covers $\fk U$ of $C$. Now choose any $N\in\mathbb Z_+$.  The data $\fk X=(C;x_1,\dots,x_N)$ is a called an \textbf{$N$-pointed compact Riemann surface}, if $x_1,\dots,x_N$ are distinct points on $C$. Choose mutually disjoint connected open subsets $U_1,\dots,U_N\subset C$ containing $x_1,\dots,x_N$ respectively, and define $U_0=C-\{x_1,\dots,x_N\}$. Then $\fk U=\{U_0,U_1,\dots,U_N\}$ is an open cover of $C$.  We now construct some cocycles in $Z^{1}(\fk U,\scr E\otimes \omega_C^p)$. For any $1\leq n\leq N$, choose
\begin{align*}
\sigma_n\in\mc (\scr E\otimes\omega_C^p)(U_n-\{x_n\}).
\end{align*}
Note that $U_n-\{x_n\}=U_n\cap U_0$. We now define \v{C}ech $1$-cocycle $s=(s_{m,n})_{m,n=0,1,\dots,N }\in Z^{1}(\fk U,\scr E\otimes \omega_C^p)$ (where each $s_{m,n}\in (\scr E\otimes\omega_C^p)(U_m\cap U_n)$) in the following way. Set $s_{0,0}=0$; if $m,n> 0$ then $s_{m,n}$ is not defined since $U_m\cap U_n=\emptyset$; if $n>0$ then $s_{n,0}=-s_{0,n}=\sigma_n$. Then $s$ can also be regarded as an element in $H^1(C,\scr E\otimes \omega_C^p)$. (Indeed, any element in $H^1(C,\scr E\otimes \omega_C^p)$ arises in such way. One way to see this is to note that since each $U_j$ is Stein manifold, $H^p(U_j,\scr E\otimes\omega_C^p)=0$ when $p>0$ by Cartan's Theorem B. (See Sec. \ref{lb20}.) So by Leray's theorem, $H^1(C,\scr E\otimes \omega_C^p)=H^1(\fk U,\scr E\otimes \omega_C^p)$.)  Choose any $t\in H^0(C,\scr E^*\otimes\omega_C^{1-p})$, which is a global section of $\scr E^*\otimes\omega_C^{1-p}$ on $C$. Then for any $n$, the evaluation $\bk{\sigma_n,t}$ is an element of $\omega_C(U_n-\{x_n\})$. So we have the residue
\begin{align}\label{eq1}
\Res_{x_n}\bk{\sigma_n,t}=\frac{1}{2\im\pi}\oint_{\gamma_n}\bk{\sigma_n,t},
\end{align}
where $\gamma_n$ is an arbitrary loop around $x_n$ whose orientation is anticlockwise in any local coordinate at $x_n$.  Now, the residue pairing of Serre duality \eqref{eq2} is described by
\begin{align}
\bk{s,t}=\sum_{n=1}^N\Res_{x_n}\bk{\sigma_n,t}.
\end{align}

Let us explain, assuming the existence of an isomorphism \eqref{eq2}, why an explicit isomorphism can be realized by the above pairing. It is an easy exercise that $\bk{s,t}=0$ when $s$ is inside the coboundary $B^1(\fk U,\scr E\otimes \omega_C^p)$, hence when $s$ is zero when regarded as in $H^1(C,\scr E\otimes \omega_C^p)=H^1(\fk U,\scr E\otimes \omega_C^p)$. Use the residue pairing to define a linear map from the left to the right hand side of \eqref{eq2}. It suffices to prove that this map is surjective.  It is easy to see that this map is independent of the sizes of $U_1,\dots,U_N$. So we may assume the vector bundle $\scr E$ can be trivialized on each $U_1,\dots,U_N$. Now, for each linear functional $H^0(C,\scr E^*\otimes_C^{1-p})\rightarrow\Cbb$,  using linear algebra, it is easy to realize it as the residue pairing with some element of $Z^1(\fk U,\scr E\otimes\omega_C^p)$. This proves the surjectivity.

\subsection*{Vanishing theorems}

When studying a family of vector bundles $\{\scr E_b:b\in\mc B \}$ over a family of compact Riemann surfaces $\{\mc C_b:b\in\mc B \}$, it is important to know if the collection $\{H^0(C_b,\scr E_b):b\in\mc B\}$ forms a vector bundle on $\mc B$ in a natural way.  Clearly, a necessary condition is that  $\dim H^0(\mc C_b,\scr E_b)$ is locally independent of $b$. As we shall see, a theorem of Grauert implies that this is also a sufficient condition.

The constancy of $\dim H^0(\mc C_b,\scr E_b)$ is not always true in general. However, if we define the \textbf{character} of $\scr E$ to be \index{zz@$\chi(C,\scr E)$}
\begin{align}
\chi(C,\scr E)=\sum_{n\in\mathbb N}(-1)^n\dim H^n(C,\scr E)=\dim H^0(C,\scr E)-\dim H^1(C,\scr E),\label{eq13}
\end{align}
then $\chi(\mc C_b,\scr E_b)$ is indeed always constant over $b$. On the other hand, if one can use vanishing theorems to show that $\dim H^1(\mc C_b,\scr E_b)=0$ for all $b$, then the constancy of $\dim H^0$ immediately follows from that of the characters. In the following, we discuss several vanishing results which will be useful for the study of conformal blocks.

Let $D$ be a \textbf{divisor} of $C$. In other words, $D$ is a finite formal sum $D=\sum_i n_ix_i$, where $\{x_i\}$ are  points of $C$, and each $n_i\in\mathbb Z$. We say $D$ is \textbf{effective} and write $D\geq 0$ if any $n_i$ is non-negative.  Recall that the \textbf{degree} $\deg D$ of $D=\sum_i n_ix_i$ is $\sum_i n_i$.\index{deg@$\deg D$}    Regard $\scr E$ as a holomorphic vector bundle. For any open $U\subset C$, let $\scr E(D)(U)$ be the set of all $s\in \scr E(D)(U-\{x_i\})$  satisfying that for any $x_i$ and any local coordinate $\eta_i$ near $x_i$, $\eta_i^{n_i}s$ has removable singularity at $x_i$. Then the collection $\{\scr E(D)(U) \}$ over all open $U\subset C$ forms an $\scr O_C$-module $\scr E(D)$. \index{ED@$\scr E(D),\scr O_C(D)$}

Note that $\scr O_C(D)$ is a line bundle, and it is well known that any line bundle is isomorphic to some $\scr O_C(D)$. (A proof is sketched in Remark \ref{lb63}.) One has a natural isomorphism of $\scr O_C$-modules $\scr E(D)\simeq\scr E\otimes\scr O_C(D)$.  Therefore $\scr E(D)$ is locally free. We understand $\scr E(D)$ as $\scr E\otimes\scr O_C(D)$ even when $\scr E$ is not locally free. One also has $\scr O_C(-D)\simeq\scr O_C(D)^*$, and $\scr O_C(D_1+D_2)\simeq\scr O_C(D_1)\otimes\scr O_C(D_2)$ for two divisors $D_1,D_2$.


\begin{pp}\label{lb62}
Assume that $C$ is connected, $\scr E$ is a locally free $\scr O_C$-module, and  $D$ is a non-zero effective divisor of $C$. Then   there exits $N\in\mathbb N$ such that $H^0(C,\scr E(-nD))=0$ for any $n> N$.
\end{pp}

\begin{proof}
Consider $\scr E$ as a vector bundle. Write $D=\sum_i x_i$. For any $x=x_i$, $H^0(C,\scr E(-D))$ is a subspace of $H^0(C,\scr E(-x))$. Thus it suffices to prove that $H^0(C,\scr E(-nx))=0$ when $n$ is large enough. 

Since $H^0(C,\scr E)$ is finite dimensional, there exist finitely many global sections $\{s_k:k=1,2,\dots \}$ of $\scr E$ spanning $H^0(C,\scr E)$. Regard a neighborhood $U$ of $x$ as an open subset of $\Cbb$, and assume that $x=0\in\mathbb C$. Assume also that $\scr E_U$ has trivialization $\scr E_U\simeq E\otimes_{\mathbb C} \scr O_U$ where $E$ is a finite dimensional vector space. Then for any $k$, $s_k$ has expansion $s_k(z)=\sum_{j=0}^{\infty}v_{k,j}z^j$ near $z=0$, where each $v_{k,j}\in E$.

For any $n\in\mathbb N$ , let $\vec{v}_k|_n$ be $(v_{k,0},v_{k,1},\dots,v_{k,n})$ in $E\otimes\mathbb C^{n+1}$. Let $F_n$ be the subspace of $E\otimes\mathbb C^{n+1}$ spanned by $\{\vec{v}_k|_n:k=1,2,\dots\}$. Then $\dim F_n$ is an increasing function of $n$ whose values are bounded from above. Choose $N\in\mathbb N$ such that $\dim F_n$ is the constant $K$ for all $n\geq N$. Assume without loss of generality that $\vec{v}_1|_N,\dots,\vec{v}_K|_N$ are linearly independent. So for any $n\geq N$, $\vec{v}_1|_n,\dots,\vec{v}_K|_n$ are also linearly independent, which therefore form a basis of $F_n$. Choose any $k$. Then for any $n\geq N$, $\vec{v}_k|_n=c_{1,n}\vec{v}_1|_n+\cdots+c_{K,n}\vec{v}_K|_n$ for some unique $c_{1,n},\dots,c_{K,n}\in\mathbb C$. By such uniqueness, we conclude that $c_{1,n}=c_{1,N},\dots,c_{K,n}=c_{K,N}$ for all $n\geq N$. Therefore $s_k=c_{1,N}s_1+\cdots+c_{K,N}s_K$ near $x$. This equation holds globally since $C$ is connected. We thus conclude that $s_1,\dots, s_K$ form a basis of $H^0(C,\scr E)$. In particular, $K=\dim H^0(C,\scr E)$.

Now choose any $n>N$ and any $\sigma\in H^0(C,\scr E(-nx))\subset H^0(C,\scr E)$. Then there exist $c_1,\dots,c_K\in\mathbb C$ such that $\sigma=c_1s_1+\cdots+c_Ks_K$. Near $x$ one has series expansion $\sigma(z)=\sum_{j=0}^\infty \nu_jz^j$. Then $(\nu_0,\nu_1,\dots,\nu_N)=c_1\vec{v}_1|_N+\dots+c_K\vec{v}_K|_N$. Since $n>N$ and $z^{-n}\sigma(z)$ has removable singularity near $0$, we have $\nu_0=\dots=\nu_N=0$. Therefore $c_1=\dots=c_K=0$ by the linear independence of $\vec{v}_1|_N,\dots,\vec{v}_K|_N$. This proves that $\sigma=0$.
\end{proof}

Let $\Theta_C$ be the  \textbf{tangent sheaf} of $C$, \index{zz@$\Theta_C$} i.e., the sheaf of holomorphic tangent vectors on $C$. So $\Theta_C\simeq\omega_C^{-1}$. Consequently, $\Theta_C\otimes\omega_C\simeq\scr O_C$.

\begin{co}\label{lb6}
Assume that $C$ is connected,  $\scr E$ is a locally free $\scr O_C$-module, and  $D$ is a non-zero effective divisor of $C$. Then   there exits $N\in\mathbb N$ such that $H^1(C,\scr E(nD))=0$ for any $n> N$.
\end{co}

\begin{proof}
By the above proposition, $H^0(C,\scr E^*\otimes\omega_C(-nD))=0$ for any sufficiently large $n$. Note that $(\scr E^*\otimes\omega_C(-nD))^*\otimes\omega_C\simeq\scr E\otimes\Theta_C\otimes\scr O_C(nD)\otimes\omega_C\simeq \scr E(nD)$. Thus, by Serre duality, $H^1(C,\scr E(nD))=0$ for any sufficiently large $n$.
\end{proof}

The above corollary is also true when $D$ has positive degree. See \cite[Prop. 5.2.7]{Huy06}

\begin{rem}\label{lb128}
In Corollary \ref{lb6}, if we know that $H^1(C,\scr E(ND))=0$, then $H^1(C,\scr E(nD))=0$ for any $n\geq N$. Indeed, by Serre duality, we have $H^0(C,\scr E^*\otimes\omega_C(-ND))=0$. Then $H^0(C,\scr E^*\otimes\omega_C(-nD))$, which is naturally a subspace of $H^0(C,\scr E^*\otimes\omega_C(-ND))$, must also be trivial. By Serre duality again, we obtain $H^1(C,\scr E(nD))=0$.
\end{rem}

It will be important to know what $N$ precisely is in the above corollary. When $\scr E$ is a line bundle, we can find such $N$ with the help of Kodaira vanishing theorem and Riemann-Roch theorem. To begin with, recall:

\begin{thm}(Kodaira vanishing theorem)\label{lb64}
Assume that $C$ is connected. Let $D$ be a divisor on $C$ with $\deg D>0$. Then $H^0(C,\scr O_C(-D))=0$ and $H^1(C,\omega_C(D))=0$.
\end{thm}

\begin{proof}
Suppose $H^0(C,\scr O_C(-D))$ is nontrivial.  Choose any non-zero $f\in H^0(C,\scr O_C(-D))$. Then $f$ is a global meromorphic function on $C$ which is not always zero on any open subset of $C$. Thus its degree $\deg f$ must be $0$.\footnote{For any $x\in C$ with an arbitrary local coordinate $\eta_x$ (so $\eta_x(x)=0$), let $n_x$ be the smallest integer such that $f/\eta_x^{n_x}$ has removable singularity at $x$. Then $\deg f$ is defined to be $\sum_{x\in C} n_x$. Clearly $n_x=\Res_x f^{-1}df$. Thus we have $\deg f=0$ by residue theorem.} But $\deg f-\deg D$ must be non-negative by the definition of $\scr O_C(-D)$, which contradicts $\deg D>0$. So $H^0(C,\scr O_C(-D))=0$. By Serre duality, $H^1(C,\omega_C(D))=0$.
\end{proof}

For any line bundle $\scr L$, we choose a divisor $D$ such that $\scr L\simeq\scr O_C(D)$, and  define the degree $\deg \scr L $\index{deg@$\deg \scr L$} of $\scr L$ to be $\deg D$. This is well defined and independent of the choice of $D$. Now assume that $C$ is connected. Then $\dim H^0(C,\scr O_C)=1$. Since $C$ admits a  K\"ahler structure, we have the Hodge structure $\dim H^1(C,\mathbb C)=\dim H^1(C,\scr O_C)+\dim H^0(C,\omega_C)$ (\cite{Huy06} corollary 3.2.12), which together with Serre duality implies $\dim H^1(C,\scr O_C)=\frac 12\dim H^1(C,\mathbb C)$. Thus the \textbf{genus} $g:=\dim H^1(C,\scr O_C)$ of $C$ depends only on the topological structure but not the complex structure of $C$. (This fact also follows from the base change Theorem \ref{lb2}.) Now, the \textbf{Riemann-Roch theorem} tells us that
\begin{align}
\chi(C,\scr L)=1-g+\deg\scr L.\label{eq162}
\end{align}
(See also Remark \ref{lb63}.) Apply this formula to $\omega_C$ and use Serre duality, we obtain
\begin{align*}
&1-g+\deg\omega_C=\chi(C,\omega_C)=\dim H^0(C,\omega_C)-\dim H^1(C,\omega_C)\\
=&\dim H^1(C,\scr O_C)-\dim H^0(C,\scr O_C)=g-1.
\end{align*}
This shows that
\begin{align}
\deg\omega_C=2g-2,\qquad \deg\Theta_C=2-2g.\label{eq3}
\end{align}
 
\begin{thm}\label{lb3}
Let $C$ be a compact connected Riemann surface with genus $g$, let $D$ be a divisor of $C$, and let $\scr L$ be a line bundle on $C$. Then   $H^1(C,\scr L(D))=0$ when $\deg D>2g-2-\deg\scr L$.
\end{thm}

\begin{proof}
	Choose  divisors $T,L$ such that $\Theta_C\simeq\scr O_C(T),\scr L\simeq\scr O_C(L)$. Then $\scr L(D)\simeq\omega_C\otimes\Theta_C\otimes\scr L(D)\simeq\omega_C(T+L+D)$. By \eqref{eq3}, $\deg T=2-2g$. Therefore, when $\deg D>2g-2-\deg\scr L$, we have $\deg(T+L+D)>0$, which implies $H^1(C,\scr L(D))=0$ by Kodaira vanishing theorem.
\end{proof}

\begin{co}\label{lb5}
Let $C$ be a compact connected Riemann surface with genus $g$, let $D$ be a divisor of $C$, and choose $n\in\mathbb Z$. Then   $H^1(C,\Theta_C^{\otimes n}(D))=0$ when $\deg D>(n+1)(2g-2)$.
\end{co}

\begin{proof}
This follows immediately from the above theorem and \eqref{eq3}.
\end{proof}

This corollary will be the most useful vanishing result for our future study of VOA bundles. The most remarkable point   is that  the threshold $(n+1)(2g-2)$ for $H^1$ to vanish  is independent of the complex structure of $C$.

\section{Families of compact Riemann surfaces}\label{lb1}

\subsection*{Higher direct images}

Let $\pi:X\rightarrow Y$ be a holomorphic map of complex manifolds. Let $\scr E$ be an $\scr O_{X}$-module. For any open $ U\subset Y$, we let $\pi_*(\scr E)( U)=\scr E(\pi^{-1}( U))$. Then $\pi_*( U)$ is an $\scr O( U)$-module: if $s\in\scr E(\pi^{-1}( U))$ and $f\in\scr O( U)$, then the product $f\cdot s$ is defined to be the section $(f\circ\pi)s$. The collection $\{\pi_*(\scr E)( U)\}$ over all open $ U\subset Y$ forms a sheaf of $\scr O_{Y}$-module $\pi_*(\scr E)$, called the \textbf{direct image} of $\scr E$ (under $\pi$). 

More generally, for any $q\in\mathbb N$, the collection $\{H^q(\pi^{-1}( U),\scr E)\}$ over all open $ U\subset Y$ forms a presheaf of $\scr O_{Y}$-module, whose sheafification $R^q\pi_*(\scr E)$ \index{Rq@$R^q\pi_*(\scr E)$, $\pi_*(\scr E)$} is called the $q$-th order \textbf{higher direct image} of $\scr E$ under $\pi$. Then $R^0\pi_*(\scr E)$ is just $\pi_*(\scr E)$.  Note that when $Y$ is a single point and $\pi$ is surjective, higher direct images are nothing but the cohomology groups of $\scr E$. Similar to cohomology groups, whenever there is a short exact sequence of $\scr O_{X}$-modules $0\rightarrow\scr E\rightarrow\scr F\rightarrow\scr G\rightarrow 0$, there is a long exact sequence of $\scr O_{Y}$-modules
\begin{align}
0\rightarrow \pi_*(\scr E)\rightarrow\pi_*(\scr F)\rightarrow\pi_*(\scr G)\xrightarrow{\delta} R^1\pi_*(\scr E)\rightarrow R^1\pi_*(\scr F)\rightarrow R^1\pi_*(\scr G)\xrightarrow{\delta} R^2\pi_*(\scr E)\rightarrow\cdots.\label{eq4}
\end{align}

All the maps in this exact sequence, except the connecting homomorphism $\delta$,  are understood in an obvious way. We now describe the first $\delta$, i.e. $\delta:\pi_*(\scr G)\rightarrow R^1\pi_*(\scr E)$, which is sufficient for the purpose of this monograph. Choose any open $ U\subset Y$ and any $s\in\scr G(\pi^{-1}( U))$. Then by the surjectivity of the stalk map $\scr F_x\rightarrow\scr G_x$ for all $x\in X$, $\pi^{-1}( U)$ has an open cover $\fk V=(V_\alpha)_\alpha$ such that for each $\alpha$ there exists $\sigma_\alpha\in \scr F( V_\alpha)$ whose image in $\scr G( V_\alpha)$ equals $s|_{ V_\alpha}$. Define the $1$-cochain $\varsigma=(\varsigma_{\alpha,\beta})_{\alpha,\beta}$ such that for any $ V_\alpha, V_\beta$ in $\fk V$, $\varsigma_{\alpha,\beta}=\sigma_\alpha|_{ V_\alpha\cap V_\beta}-\sigma_\beta|_{ V_\alpha\cap V_\beta}$. Then $\varsigma_{\alpha,\beta}\in\scr E( V_\alpha\cap V_\beta)$, and hence $\varsigma\in Z^1(\fk V,\scr E)$. So $\varsigma$ can be regarded as an element of $H^1(\pi^{-1}( U),\scr E)$ and hence of $R^1\pi_*(\scr E)(U)$. Then $\delta (s)$ is just this $\varsigma$.

\subsection*{Pulling back  sheaves}

We now assume that $\scr E$ is an $\scr O_Y$-module. We define an $\scr O_X$-module $\pi^*(\scr E)$, called the \textbf{pull back} of $\scr E$ under $\pi$, \index{Epi@$\pi^*(\scr E),\pi^*s$} as follows.

Let $U$ be an open subset of $X$. For any open subset $V$ of $Y$ satisfying $\pi(U)\subset V$, the ring $\scr O_X(U)$ is an $\scr O_Y(V)$-module: if $f\in\scr O_X(U)$ and $g\in\scr O_Y(V)$, then the action of $g$ on $f$ is $(g\circ \pi)|_U\cdot f$, noticing that $g\circ \pi\in\scr O_X(\pi^{-1}(V))$ and hence $(g\circ \pi)|_U\in\scr O_X(U)$. Then the tensor product $\scr O_X(U)\otimes_{\scr O_Y(V)}\scr E(V)$ is an $\scr O_X(U)$-module. The collection over all open $U\subset X$ of the $\scr O_X(U)$-modules
\begin{align*}
\varinjlim_{V\supset \pi(U)}\scr O_X(U)\otimes_{\scr O_Y(V)}\scr E(V)
\end{align*}
form a presheaf of $\scr O_X$-modules, whose sheafification is $\pi^*(\scr E)$. We can thus define  for any open $V\subset Y$ an $\scr O_Y(V)$-module homomorphism
\begin{align}
\pi^*:\scr E(V)\rightarrow\pi^*\scr E(\pi^{-1}(V)),\qquad s\mapsto\pi^*s,\label{eq139}
\end{align}
such that  $1\otimes s\in\scr O_X(\pi^{-1}(V))\otimes_{\scr O_Y(V)}\scr E(V)$, regarded as an element of $\pi^*\scr E(\pi^{-1}(V))$, is our $\pi^*s$. That $\pi^*$ is an $\scr O_Y(V)$-module homomorphism means that for any $g\in\scr O_Y(V)$,
\begin{align*}
\pi^*(gs)=(g\circ\pi)\pi^*s.
\end{align*}
Then for any $q\in\mathbb N$ we have a natural linear map
\begin{align*}
\pi^*:H^q(Y,\scr E)\rightarrow H^q(X,\pi^*\scr E),
\end{align*}
called the pull back of cohomology groups.

When $\scr E$ is locally free, for any $x\in X$, choose a neighborhood $W\subset Y$ of $\pi(x)$ such that $\scr E_W\simeq\scr O_W\otimes_{\mathbb C}E$, where $E$ is a finite dimensional vector space. Then for any open $U\subset \pi^{-1}(W)$ and $V\subset W$, the $\scr O_X(U)$-module $\scr O_X(U)\otimes_{\scr O_Y(V)}\scr E(V)$ is naturally isomorphic to $\scr O_X(U)\otimes_{\mathbb C}E$. So $\pi^*(\scr E)$ is also locally free whose rank is the same as that of $\scr E$. Indeed, the pull back of locally free modules is just the pull back of holomorphic vector bundles. 

If $\iota:X\rightarrow Y$ is an embedding of complex manifolds, for any  $\scr O_Y$-module $\scr E$, we call  
\begin{align*}
\scr E|X\equiv\scr E|_X:=\iota^*\scr E
\end{align*}
\index{EX@ $\scr E\lvert X\equiv\scr E\lvert_X,\scr E\lvert x$}the \textbf{restriction} of $\scr E$ to $X$. For any open $V\subset Y$, we have the $\scr O_Y(V)$-module homomorphism
\begin{gather*}
\iota^*:\scr E(V)\rightarrow\scr E|_X(V\cap X),\qquad s\mapsto \iota^*s.
\end{gather*}
We  write 
\begin{align*}
s|X\equiv s|_X:=\iota^*s
\end{align*}
and call it the \textbf{restriction} of $s$ to $X$. \index{EX@ ${\scr E\lvert X},\scr E\lvert x$!${s\lvert X=s\lvert_X}$} Similar notation will be used for restrictions of elements in cohomology groups.

Note that $\scr E|V$ is just $\scr E_V$,  the usual restriction of $\scr E$ to the open subset $V$.  For any $y\in Y$, $\scr E|y$ is a vector space isomorphic to the fiber $\scr E_y/\fk m_y\scr E_y$.  Indeed, notice that
\begin{align*}
\scr E|y=\varinjlim_{V\ni y}\mathbb C_y\otimes_{\scr O_Y(V)}\scr E(V).
\end{align*}
Here $\mathbb C_y$ is the same as $\mathbb C$ as a vector space. But it is also an $\scr O_Y(V)$-module, where for any $\lambda\in\mathbb C_y$ and any $g\in\scr O_Y(V)$, the action of $g$ on $\lambda$ is $\lambda g(y)$. It is then not hard to check that the linear map $1\otimes s\in\mathbb C_y\otimes_{\scr O_Y(V)}\scr E(V)\mapsto s\in\scr E_y/\fk m_y\scr E_y$ induces an isomorphism
\begin{align*}
\scr E|y\xrightarrow{\simeq}\scr E_y/\fk m_y\scr E_y.
\end{align*}
We will not distinguish between $\scr E|y$ and $\scr E_y/\fk m_y\scr E_y$ in the future. $s|y$, the restriction of $s$ to the $0$-dimensional submanifold $y$, is precisely the value of $s$ in $\scr E_y/\fk m_y\scr E_y$. We also write
\begin{align*}
s(y):=s|y,
\end{align*}
\index{EX@ ${\scr E\lvert X},\scr E\lvert x$!${s\lvert x=s(x)}$}which should not be confused with the germ $s_y$ in $\scr E_y$.

\subsection*{A theorem of Grauert}

Let $\mc C$ and $\mc B$ be complex manifolds, and $\pi:\mc C\rightarrow \mc B$ be a surjective proper holomorphic map. (The word ``proper" means that the preimage of any compact subset of $\mc B$ is compact.)  Assume that $\pi$ is a (holomorphic) submersion, which means that the linear map $d\pi$ between tangent spaces is always surjective. We say that $\pi:\mc C\rightarrow \mc B$ is a (holomorphic) \textbf{family of compact Riemann surfaces} (or a smooth family of complex curves), if for each $b\in\mc B$, the fiber $\mc C_b:=\pi^{-1}(b)$ \index{Cb@$\mc C_b=\pi^{-1}(b)$} is a compact Riemann surface. As an obvious fact, the complex dimensions of $\mc C$ and $\mc B$ differ by $1$. To simplify discussions, we also assume that $\mc B$ has finitely many connected components.

The following theorem of Grauert \cite{Gra60} is extremely helpful in studying families of compact Riemann surfaces.
\begin{thm}\label{lb2}
Let $\scr E$ be a locally free $\scr O_{\mc C}$-module.

(a) The function
\begin{gather*}
\mc B\rightarrow \mathbb Z, \qquad b\mapsto\chi(\mc C_b,\scr E|\mc C_b)
\end{gather*}
is locally constant. 

(b) For any $q\in\mathbb N$, the function  $b\mapsto\dim H^q(\mc C_b,\scr E|\mc C_b)$ is upper semi-continuous. Moreover, if this function is locally constant, then $R^q\pi_*\scr E$ is locally free of rank $\dim H^q(\mc C_b,\scr E|\mc C_b)$, and for any $b\in\mc B$, the linear map $(R^q\pi_*\scr E)_b\rightarrow H^q(\mc C_b,\scr E|\mc C_b)$ defined by  restriction of sections $s\mapsto s|{\mc C_b}$ induces an isomorphism of vector spaces
\begin{align}
\frac{(R^q\pi_*\scr E)_b}{\fk m_b\cdot (R^q\pi_*\scr E)_b}\simeq H^q(\mc C_b,\scr E|\mc C_b).
\end{align}
\end{thm}
To be more precise, note that $(R^q\pi_*\scr E)_b$ is the direct limit of $H^q(\pi^{-1}(U),\scr E)$ over all open $U\ni b$. Then the  pull back $H^q(\pi^{-1}(U),\scr E)\rightarrow H^q(\mc C_b,\scr E|\mc C_b)$   induced by the inclusion $\mc C_b\hookrightarrow \pi^{-1}(U)$ provides the desired map $(R^q\pi_*\scr E)_b\rightarrow H^q(\mc C_b,\scr E|\mc C_b)$.

The above theorem is indeed true for any proper submersion of complex manifolds $\pi:X\rightarrow Y$, assuming that $\scr E$ is a coherent $\scr O_X$-module. Even more generally, \cite{Gra60} proves this theorem when $X$ and $Y$ are complex spaces, $Y$ is reduced, $\scr E$ is coherent, and $\scr E$ is $\pi$-flat. The actual meanings of these three terms are not so important for understanding this monograph. Roughly speaking, complex spaces are generalizations of complex manifolds which may have singularities. If $Y$ (together with an associated structure sheaf $\scr O_Y$) is a complex space, then $(Y,\scr O_Y)$ looks locally like the solution space of finitely many holomorphic functions. (See  \cite{GR84} chapter 1.) The only non-smooth complex spaces that we will encounter in this monograph are nodal curves. $Y$ is called reduced if any stalk $\scr O_y$ has no non-zero nilpotent elements. That $\scr E$ is $\pi$-flat means that for any $x\in X$, the natural action of $\scr O_{Y,\pi(x)}$ on $\scr E_x$ makes $\scr E_x$ a flat $\scr O_{Y,\pi(x)}$-module. We refer the reader to \cite{GPR94} section III.4.2 and the reference therein for the general form of the theorem and its proof. (See especially theorem 4.7 of that section.) A proof in English can be found in \cite[Thm. III.4.12]{BS76} as well as \cite[Thm. 9.4.8]{EP96}.

As an immediate consequence of  Theorem \ref{lb2}, for a family $\pi:\mc C\rightarrow \mc B$ and a locally free $\scr O_{\mc C}$-module $\scr E$, we have $R^q\pi_*\scr E=0$ when $q>1$, since $\dim H^q(\mc C_b,\scr E|\mc C_b)$ is constantly $0$.

\subsection*{Relative tangent sheaves}

In the remaining part of this section, we apply Theorem \ref{lb2} to the families of tangent sheaves and their tensor products. Again, we fix a family of compact Riemann surfaces $\pi:\mc C\rightarrow\mc B$.

The first obvious question is how  the tangent bundles over all fibers can be assembled into a vector bundle on $\mc C$. Consider the tangent sheaves $\Theta_{\mc C}$ and $\Theta_{\mc B}$ of $\mc C$ and $\mc B$ respectively. Then we have a natural $\scr O_{\mc C}$-module homomorphism $d\pi:\Theta_{\mc C}\rightarrow\pi^*\Theta_{\mc B}$ described as follows. Locally the map $\pi$ looks like the projection $U\times V\rightarrow V$ where $U$ is an open subset of $\mathbb C$ and $V$ is an open subset of $\mathbb C^m$. (So $m$ is the dimension of $\mc B$.)  Let $z$ be the standard coordinate of $U$, and $\tau_1,\dots,\tau_m$ be the standard coordinates of $V$. Then $\partial_z,\partial_{\tau_1},\dots,\partial_{\tau_m}$ are sections in $\Theta_{\mc C}(U\times V)$, and  $\partial_{\tau_1},\dots,\partial_{\tau_m}$ can also be regarded as sections in $\Theta_{\mc B}(V)$.  We thus define the $\scr O_{\mc C}(U\times V)$-module homomorphism $d\pi:\Theta_{\mc C}(U\times V)\rightarrow\pi^*\Theta_{\mc B}(U\times V)$ by sending $\partial_z,\partial_{\tau_1},\dots,\partial_{\tau_m}$ to  $0,\pi^*\partial_{\tau_1},\dots,\pi^*\partial_{\tau_m}$ respectively. We write each $\pi^*\partial_{\tau_j}$ simply as $\partial_{\tau_j}$. Then for any holomorphic functions $f,g_1,\dots,g_n$ on $U\times V$, we have
\begin{align}
d\pi(f\partial_z+g_1\partial_{\tau_1}+\cdots+g_m\partial_{\tau_m})=g_1\partial_{\tau_1}+\cdots+g_m\partial_{\tau_m}.
\end{align} 
So $d\pi$ is the projection onto  ``horizontal components". The map $d\pi$ is independent of local coordinates. Thus we have an exact sequence of $\scr O_{\mc C}$-modules
\begin{align}
0\rightarrow \Theta_{\mc C/\mc B}\rightarrow \Theta_{\mc C}\xrightarrow{d\pi}\pi^*\Theta_{\mc B}\rightarrow 0,\label{eq81}
\end{align}
where $\Theta_{\mc C/\mc B}$ is the kernel sheaf of $d\pi$, \index{zz@$\Theta_{\mc C/\mc B}$} called the \textbf{relative tangent sheaf} of the family $\pi:\mc C\rightarrow\mc B$. One checks easily that if $\zeta,t_1,\dots,t_m$ are also local  coordinates chosen in a similar way, then the two sections $\partial_z$ and $\partial_\zeta$ of $\Theta_{\mc C/\mc B}$ are related by $\partial_z=\frac{\partial\zeta}{\partial z}\partial_\zeta$. From this one sees that  the restriction $\Theta_{\mc C/\mc B}|\mc C_b$ of the relative tangent sheaf on each fiber $\mc C_b$ is isomorphic to the $\scr O_{\mc C_b}$-module $\Theta_{\mc C_b}$. Therefore one can regard $\Theta_{\mc C/\mc B}$ as the sheaf of ``vertical sections" of $\Theta_{\mc C}$, or the sheaf of vectors that are tangent to the fibers. Since $\Theta_{\mc C/\mc B}$ is a line bundle, we call its dual sheaf $\omega_{\mc C/\mc B}:=\Theta_{\mc C/\mc B}^{-1}$ \index{zz@$\omega_{\mc C/\mc B}$} the \textbf{relative dualizing sheaf} of the family $\pi:\mc C\rightarrow\mc B$. The restriction of $\omega_{\mc C/\mc B}$ to each fiber $\mc C_b$ is therefore $\omega_{\mc C_b}$.

We next discuss families of divisors. A family of points can be represented by a  section $\sgm:\mc B\rightarrow\mc C$, i.e., a holomorphic map such that $\pi\circ \sgm=\id_{\mc B}$. The image of $\sgm$ is a hypersurface of $\mc C$. Let $\scr E$ be a locally free $\scr O_{\mc C}$-module, let $U$ be  an open subset of $\mc C$, and choose  $s\in\scr E(U-\sgm(\mc B))$. We say that $s$ has \textbf{removable singularity} at $\sgm(\mc B)$ if $s$ can be extended to an element in $\scr E(U)$. 

A \textbf{local coordinate} $\eta$ of the family $\pi:\mc C\rightarrow\mc B$ at $\sgm$ is a holomorphic function on a neighborhood $U$ of $\sgm(\mc B)$  such that for any $b\in \mc B$,  $\eta(\sgm(b))=0$, and  $\eta$  restricts to a biholomorphic map from $U\cap\mc C_b$ to an open subset of $\mathbb C$. (The second condition is equivalent to saying that $\eta|_{U\cap\mc C_b}$ is injective, i.e. univalent.)  Then $(\eta,\pi)$ is a biholomorphic map from $U$ to an open subset of $\mathbb C\times\mc B$. It is obvious that  any $b\in\mc B$ is contained in a neighborhood $V\subset\mc B$ such that the subfamily $\pi:\pi^{-1}(V)\rightarrow V$ has a local coordinate at $\sgm|_V$.

If $\pi:\mc C\rightarrow\mc B$ is a family of compact Riemann surfaces, we say that $(\pi:\mc C\rightarrow\mc B;\sgm_1,\dots,\sgm_N)$ is a \textbf{family of $N$-pointed compact Riemann surfaces}, if: 
\begin{enumerate}[label=(\alph*)]
\item $\sgm_1,\dots\sgm_N:\mc B\rightarrow\mc C$ are sections.
\item For any $b\in\mc B$ and any $1\leq i<j\leq N$, $\sgm_i(b)\neq\sgm_j(b)$.
\item For any $b\in\mc B$, each connected component of $\mc C_b$ contains $\sgm_i(b)$ for some $1\leq i\leq N$.
\end{enumerate}
(The third condition is not essential. It is introduced only to simply the statement  of vanishing theorems in the future.) If, moreover, $\eta_1,\dots,\eta_N$ are local coordinates at $\sgm_1,\dots,\sgm_N$ respectively, we say that
\begin{align*}
\fk X=(\pi:\mc C\rightarrow\mc B;\sgm_1,\dots,\sgm_N;\eta_1,\dots,\eta_N)
\end{align*}
is a \textbf{family of $N$-pointed compact Riemann surfaces with local coordinates}. \index{XX@$\wtd {\fk X}$ and $\fk X$} (In the case that $\mc B$ is a single point, $\fk X$ is called an $N$-pointed compact Riemann surface with local coordinates, and is denoted by $(C;x_1,\dots,x_N;\eta_1,\dots,\eta_N)$, where $\eta_1,\dots,\eta_N$ are respectively local coordinates at the distinct points $x_1,\dots,x_N$ on a compact Riemann surface $C$.)

Now, if we have $(\pi:\mc C\rightarrow\mc B;\sgm_1,\dots,\sgm_N)$, each hypersurface $\sgm_i(\mc B)$ is also a divisor of $\mc C$. Set $D=\sum_{i=1}^Nn_i\sgm_i(\mc B)$ where each $n_1,\dots,n_N\in\mathbb Z$. Let $\scr E$ be a locally free $\scr O_{\mc C}$-module. Then the $\scr O_{\mc C}$-module $\scr E(D)$ can also be defined in a similar way as in the case of Riemann surfaces: For each open $U\subset \mc C$, $\scr E(D)(U)$ is the $\scr O_{\mc C}(U)$-module of all $s\in\scr E(U-\bigcup_{i=1}^N\sgm_i(\mc B))$ such that for any $1\leq i\leq N$, any open subset $V\subset U$, and any (holomorphic) submersion $\eta:V\rightarrow \mathbb C$ vanishing at $V\cap\sgm_i(\mc B)$, $\eta^{n_i}s$ has removable singularity at $V\cap\sgm_i(\mc B)$. \index{ED@$\scr E(D),\scr O_C(D)$} When there are local coordinates $\eta_1,\dots,\eta_N$ at $\sgm_1,\dots,\sgm_N$ respectively, this is equivalent to saying that $\eta_i^{n_i}s$ has removable singularity at $\sgm_i(\mc B)$ for any $i$. In particular we have defined $\scr O_{\mc C}(D)$. We then have natural isomorphism $\scr E(D)\simeq\scr E\otimes\scr O_{\mc C}(D)$. In the general case that  $\scr E$ is not necessarily locally free, we define $\scr E(D)$ simply to be $\scr E\otimes\scr O_{\mc C}(D)$. Notice that for any exact sequence of $\scr O_{\mc C}$-modules $0\rightarrow\scr E\rightarrow\scr F\rightarrow\scr G\rightarrow 0$, one also has exact sequence $0\rightarrow\scr E(D)\rightarrow\scr F(D)\rightarrow\scr G(D)\rightarrow 0$. Indeed,  the line bundle $\scr O_{\mc C}(D)$ is locally equivalent to $\scr O_{\mc C}$, thus tensoring by $\scr O_{\mc C}(D)$ preserves the exactness of $\scr O_{\mc C}$-module homomorphisms (i.e. $\scr O_{\mc C}(D)$ is a flat $\scr O_{\mc C}$-module).

For a family $\fk X=(\pi:\mc C\rightarrow\mc B;\sgm_1,\dots,\sgm_N)$ of $N$-pointed compact Riemann surfaces, we always define divisors $\SX$ of $\mc C$ and $\SX(b)$ of $\mc C_b$ \index{SX@$\SX$, $\SX(b)$} to be
\begin{gather*}
\SX=\sum_{i=1}^N\sgm_i(\mc B),\qquad \SX(b)=\sum_{i=1}^N\sgm_i(b).
\end{gather*}

\begin{pp}\label{lb129}
Let $\scr E$ be a locally free $\scr O_{\mc C}$-module. Then for any precompact open subset $V\subset\mc B$, there exists $k_0\in\Nbb$ such that for any $k\geq k_0$ and $b\in V$, we have $H^1(\mc C_b,\scr E(k\SX)|\mc C_b)=0$.
\end{pp}

\begin{proof}
We shall prove that  each $b\in\mc B$ is contained in a neighborhood $U_b\subset\mc B$ such that one can find $k_b\in\Nbb$ satisfying $H^1(\mc C_{\wtd b},\scr E(k_b\SX)|\mc C_{\wtd b})=0$ for any  $\wtd b\in U_b$. Then by Remark \ref{lb128}, the same is true when $k_b$ is replaced by any $k\geq k_b$.  Therefore, the claim of this proposition follows since we can cover $\ovl V$ by finitely many such neighborhoods.

Choose any $b\in\mc B$. Then, by Corollary \ref{lb6}, we can find $k_b\in\Nbb$ such that $H^1(\mc C_b,\scr E(k_b\SX)|\mc C_b)=0$. By the upper semi-continuity Theorem \ref{lb2}-(b), we can find a neighborhood $U_b$ of $b$ such that for each $\wtd b\in U_b$, $H^1(\mc C_{\wtd b},\scr E(k_b\SX)|\mc C_{\wtd b})=0$.
\end{proof}

If $\mc B$ is connected, then by a result of Ehresmann, all fibers are diffeomorphic (cf. \cite{Huy06} proposition 6.2.2). We have the following:

\begin{thm}\label{lb8}
Suppose that $\mc B$ is connected, and $g$ is the maximal genus of the connected components of each fiber of $\pi:\mc C\rightarrow\mc B$. Then for any integers $n\geq-1$ and $k>(n+1)(2g-2)$, $\pi_*\Theta_{\mc C/\mc B}^{\otimes n}(k\SX)$ is a locally free $\scr O_{\mc B}$-module, and for any $b\in \mc B$ there is a natural isomorphism of vector spaces
\begin{align}
\frac{\pi_*\Theta_{\mc C/\mc B}^{\otimes n}(k\SX)_b}{\fk m_b\cdot\pi_*\Theta_{\mc C/\mc B}^{\otimes n}(k\SX)_b}\simeq H^0\big(\mc C_b,\Theta_{\mc C_b}^{\otimes n}(k\SX(b))\big)\label{eq6}
\end{align} 
defined by restriction of sections. In particular,  $\dim H^0\big(\mc C_b,\Theta_{\mc C_b}^{\otimes n}(k\SX(b))\big)$ is constant over $b$.

When $n<-1$ and $g$ is the minimal genus of the connected components of any fiber, the same result holds.
\end{thm}

\begin{proof}
Note that the restriction of $\Theta_{\mc C/\mc B}$ to each fiber $\mc C_b$ is naturally equivalent to $\Theta_{\mc C_b}$. Thus the restriction of $\Theta_{\mc C/\mc B}^{\otimes n}(k\SX)$ to $\mc C_b$ is $\Theta_{\mc C_b}^{\otimes n}(k\SX(b))$. Therefore, by Theorem \ref{lb2}-(a), $\chi\big(\mc C_b,\Theta_{\mc C_b}^{\otimes n}(k\SX(b))\big)$ is constant over $b$. By corollary  \ref{lb5}, $\dim H^1\big(\mc C_b,\Theta_{\mc C_b}^{\otimes n}(k\SX(b))\big)$ is always $0$. So  $\dim H^0\big(\mc C_b,\Theta_{\mc C_b}^{\otimes n}(k\SX(b))\big)$ is constant over $b$. The remaining part of the theorem follows from Theorem \ref{lb2}-(b).
\end{proof}

An important  consequence of the above theorem is that any global section  of $\Theta_{\mc C_b}^{\otimes n}(k\SX(b))$ on $\mc C_b$ can be extended to a \emph{holomorphic} family of global sections over a neighborhood of $b\in\mc B$. More precisely, for any $s\in H^0\big(\mc C_b,\Theta_{\mc C_b}^{\otimes n}(k\SX(b))\big)$ there exists a neighborhood $V$ of $b$ and $\wtd s\in H^0\big(\pi^{-1}(V),\Theta_{\mc C/\mc B}^{\otimes n}(k\SX)\big)$ such that $\wtd s$ restricts to $s$ on $\mc C_b$. 

We will not use this theorem directly. However,  $\Theta_{\mc C/\mc B}^{\otimes n}$ is closely related to sheaves of VOAs. (Cf. Proposition \ref{lb56}.) We will prove a similar theorem for sheaves of VOAs in the next chapter; see Theorem \ref{lb27}.

\section{Strong residue theorem}

Suppose that $\fk X=(C;x_1,\dots,x_N;\eta_1,\dots,\eta_N)$ is an $N$-pointed compact Riemann surface, and $\scr E$ is a locally free $\scr O_C$-module. By our notation in the last section, $\SX=x_1+\cdots+x_N$. Set \index{ESX@$\scr E(\bullet\SX)$} 
\begin{align}
\scr E(\bullet\SX)=\varinjlim_{n\in\mathbb N}\scr E(n\SX).\label{eq9}
\end{align}
Then $\scr E(\bullet\SX)$ is the sheaf of meromorphic sections of $\scr E$ whose only possible poles are $x_1,\dots,x_N$. Let $E_i=\scr E|x_i$ be the fiber of $\scr E$ at $x_i$. In a neighborhood $U_i$ of $x_i$, $\scr E_{U_i}$ has a trivialization $\scr E_{U_i}\simeq\scr O_{U_i}\otimes_{\mathbb C}E_i$, and  $\scr E^*_{U_i}$ has the corresponding dual trivialization $\scr E^*_{U_i}\simeq\scr O_{U_i}\otimes_{\mathbb C}E_i^*$. Choose any $s\in H^0(C,\scr E(\bullet\SX))$.  Then at any $x_i$ the section $s$ has formal Laurent series expansion
$$s_i=\sum_n e_{i,n}\eta_i^n\in E_i((\eta_i)),$$
where $e_{i,n}\in\scr E|x_i$ equals $0
$ when $n$ is sufficiently small.  (Its meaning is obvious when we regard $\scr E$ as a vector bundle.) Suppose now that $\sigma\in H^0(C,\scr E^*\otimes\omega_C(\bullet\SX))$. Let
$$\sigma_i\in E_i^*((\eta_i))d\eta_i$$
be similarly the formal Laurent series expansion of $\sigma$ at $x_i$ (with respect to the variable $\eta_i$). Then $\bk{s,\sigma}\in H^0(C,\omega_C(\bullet\SX))$ and $\bk{s_i,\sigma_i}\in\mathbb C((\eta_i))d\eta_i$. Thus, by residue theorem,
\begin{align}
\sum_{i=1}^N\Res_{\eta_i=0}\bk{s_i,\sigma_i}=\sum_{i=1}^N\oint_{x_i}\bk{s,\sigma}=0.\label{eq7}
\end{align}
The strong residue theorem  (\cite{FB04} section 9.2.9) says that if $s_1\in E_1((\eta_i)),\dots,s_N\in E_N((\eta_i))$ satisfy \eqref{eq7} for any  $\sigma\in H^0(C,\scr E^*\otimes\omega_C(\bullet\SX))$, then $s_1,\dots,s_N$ are series expansions of an $s\in H^0(C,\scr E(\bullet\SX))$ at $x_1,\dots,x_N$ respectively. In particular, each  $s_i=\sum_n e_{i,n}\eta_i^n$ converges absolutely when $\eta_i\in U_i$. A proof of this theorem can be found in \cite{Ueno08} theorem 1.22. In the following, we prove a strong residue theorem to families of compact Riemann surfaces.

\subsection*{Strong residue theorem for families of compact Riemann surfaces}

Let $\fk X=(\pi:\mc C\rightarrow\mc B;\sgm_1,\dots,\sgm_N;\eta_1,\dots,\eta_N)$ be a family of $N$-pointed compact Riemann surfaces with local coordinates, and let $\scr E$ be a (holomorphic) vector bundle on $\fk X$. We assume that $\mc B$ is small enough such that for  each $i$, $\sgm_i(\mc B)$ is contained in a neighborhood $U_i$ such that $\scr E_{U_i}$ has trivialization $\scr E_{U_i}\simeq \scr O_{U_i}\otimes_{\mathbb C}E_i$, where $E_i$ is a finite dimensional complex vector space. The result for this subsection is local with respect to $\mc B$. So the reader can regard $\mc B$ as an open subset of $\mbb C^m$ for convenience.

By shrinking $U_i$, we assume that $\eta_i$ is defined on $U_i$, and $U_i\cap U_j=\emptyset$ when $i\neq j$. Fix trivialization $\scr E^*_{U_i}\simeq \scr O_{U_i}\otimes_{\mathbb C}E_i^*$ to be dual to $\scr E_{U_i}\simeq \scr O_{U_i}\otimes_{\mathbb C}E_i$. Identify $U_i$ with an open subset of $\mathbb C\times\mc B$ (containing $0\times\mc B$) via the embedding $(\eta_i,\pi):U_i\rightarrow \mathbb C\times\mc B$, and identify $\scr E_{U_i}$ and $\scr E^*_{U_i}$ with their respective trivializations. Then $\eta_i$ is identified with the standard coordinate $z$ of $\mbb C$.

For each $i$ we choose 
\begin{align}
s_i=\sum_n e_{i,n}\cdot z^n\quad\in \big(\scr O(\mc B)\otimes_{\mathbb C}E_i\big)((z)).\label{eq10}
\end{align}
So $e_{i,n}=0$ when $n$ is sufficiently small, and each $e_{i,n}\in\scr O(\mc B)\otimes_{\mathbb C}E_i$ can be viewed as an $E_i$-valued holomorphic function on $\mc B$. For each $b\in\mc B$, set
\begin{align}
s_i(b)=\sum_n e_{i,n}(b)z^n\quad\in E_i((z)).\label{eq8}
\end{align}
Define $\scr E(\bullet\SX)$ again using \eqref{eq9}. Now suppose that $s$ is a section of $\scr E(\bullet\SX)$ defined on an open set containing $U_i$.  Then $s|_{U_i}=s|_{U_i}(z,b)$ is an $E_i$-valued meromomorphic function on $U_i$ with poles at $z=0$. We say that $s$ \textbf{has series expansion $s_i$ near $\sgm_i(\mc B)$} if for each $b\in \mc B$, the meromorphic function $s|_{U_i}(z,b)$ of $z$ has series expansion \eqref{eq8} near $z=0$.

For any given $b\in\mc B$, choose $\sigma_b\in H^0(\mc C_b,\scr E^*|\mc C_b\otimes\omega_{\mc C_b}(\bullet\SX(b)))$. Then in $U_{i,b}:=U_i\cap\pi^{-1}(b)$, $\sigma_b$ can be regarded as an $E_i^*\otimes dz$-valued meromorphic function with pole at $z=0$. So it has series expansion at $z=0$:
\begin{align}
\sigma_b|_{U_i,b}(z)=\sum_n \phi_{i,n}z^ndz \quad\in E_i^*((z))dz
\end{align}
where $\phi_{i,n}\in E_i^*$. We define the residue pairing $\Res_i\bk{s_i,\sigma_b}\in\mathbb C$ to be
\begin{align}
\Res_i \bk{s_i,\sigma_b}=&\Res_{z=0}\langle s_i(b),\sigma_b|_{U_i,b}(z) \rangle\nonumber\\
=&\Res_{z=0}\bigg(\Big\langle \sum_n e_{i,n}(b)z^n, \sum_n \phi_{i,n}z^n  \Big\rangle dz\bigg).\label{eq217}
\end{align}

\begin{thm}\label{lb18}
For each $1\leq i\leq N$, choose $s_i$ as in \eqref{eq10}.  Then the following statements are equivalent.

(a) There exists $s\in H^0(\mc C,\scr E(\bullet\SX))$ whose series expansion near $\sgm_i(\mc B)$ (for each $1\leq i\leq N$) is $s_i$.

(b) For any $b\in\mc B$ and any $\sigma_b\in H^0\big(\mc C_b,\scr E^*|\mc C_b\otimes\omega_{\mc C_b}(\bullet\SX(b))\big)$,
\begin{align}
\sum_{i=1}^N\Res_i\bk{s_i,\sigma_b}=0.
\end{align}
Moreover, when these statements hold, there is only one $s\in H^0(\mc C,\scr E(\bullet\SX))$ satisfying (a).
\end{thm}

\begin{proof}
That (a) implies (b) follows from Residue theorem. If $s$ satisfies (a), then for each $b\in\mc B$, $s|\mc C_b$ is uniquely determined by its series expansions near $\sgm_1(b),\dots,\sgm_N(b)$ (since each component of $\mc C_b$ contains some $\sgm_i(b)$ by the definition of families of $N$-pointed compact Riemann surfaces). Therefore the sections satisfying (a) is unique.

Now assume (b) is true.  Suppose we can prove that any $b\in\mc B$ is contained in a neighborhood $V$ such that there exists a unique section $s$ of $\scr E(\blt\SX)$ on $\pi^{-1}(V)$ whose series expansion near each $\sgm_i(V)$ is the restriction of $s_i$ to $V$. Then (a) clearly follows. Thus, by replacing $\mc B$ with an arbitrary precompact open subset and applying Proposition \ref{lb129}, we can assume that there exists $k_0\in\mbb N$ such that for any $b\in\mc B$ and any $k\geq k_0$,
\begin{align*}
H^1\big(\mc C_b,\scr E^*|\mc C_b\otimes\omega_{\mc C_b}(k\SX(b))\big)=0.
\end{align*}

Choose $p\in\mathbb N$ such that for each $1\leq i\leq N$, the $e_{i,n}$ in \eqref{eq10} equals $0$ when $n<-p$. For any $k\geq k_0$, the exact sequence
\begin{align*}
0\rightarrow \scr E(-k\SX)\rightarrow \scr E(p\SX)\rightarrow \scr E(p\SX)/\scr E(-k\SX)\rightarrow 0
\end{align*}
induces a long exact sequence
\begin{align}
0\rightarrow\pi_*\scr E(-k\SX)\rightarrow\pi_*\scr E(p\SX)\rightarrow \pi_*\big(\scr E(p\SX)/\scr E(-k\SX)\big)\xrightarrow{\delta}R^1\pi_*\scr E(-k\SX).\label{eq11}
\end{align}
By Serre duality, $\dim H^0\big(\mc C_b,(\scr E|\mc C_b)(-k\SX(b))\big)=\dim H^1\big(\mc C_b,\scr E^*|\mc C_b\otimes\omega_{\mc C_b}(k\SX(b))\big)=0$ for any $b$. Note also that $(\scr E|\mc C_b)(-k\SX(b))$ is the restriction of $\scr E(-k\SX)$ to $\mc C_b$. Therefore $\pi_*\scr E(-k\SX)=0$ by Theorem \ref{lb2}-(b).

For each $1\leq i\leq N$, set $s_i|_k=\sum_{n<k}e_{i,n}\cdot z^n$, which can be regarded as a section in $\scr E(p\SX)(U_i)$ via the identification $\eta_i=z$. Let $U_0=\mc C-\sgm_1(\mc B)\cup\cdots\cup\sgm_N(\mc B)$. Then $\fk U=\{U_0,U_1,\dots,U_N \}$ is an open cover of $\mc C$. Define \v{C}ech $0$-cocycle $\psi=(\psi_i)_{0\leq i\leq N}\in Z^0(\fk U,\scr E(p\SX)/\scr E(-k\SX))$ by setting $\psi_0=0$ and $\psi_i=s_i|_k$ ($1\leq i\leq N$). Then $\delta\psi=\big((\delta\psi)_{i,j}\big)_{0\leq i,j\leq N}\in Z^1(\fk U,\scr E(-k\SX))$ is described as follows: $(\delta\psi)_{0,0}=0$; if $i,j>0$ then $(\delta\psi)_{i,j}$ is not defined since $U_i\cap U_j=\emptyset$; if $1\leq i\leq N$ then $(\delta\psi)_{i,0}=-(\delta\psi)_{0,i}$ equals $s_i|_k$ (considered as a section in $\scr E(-k\SX)(U_i\cap U_0)\simeq \scr E(U_i-\sgm_i(\mc B))$).

Consider $\delta\psi$ as an element in $R^1\pi_*\scr E(-k\SX)$ whose restriction to $\mc C_b$ is denoted by $\delta\psi|\mc C_b\in H^1(\mc C_b,\scr E|\mc C_b(-k \SX(b))$. Then the residue pairing for the Serre duality
\begin{align*}
H^1(\mc C_b,\scr E|\mc C_b(-k \SX(b)))\simeq H^0\big(\mc C_b,\scr E^*|\mc C_b\otimes\omega_{\mc C_b}(k\SX(b))\big)^*
\end{align*}
(see Section \ref{lb4}), applied to $\delta\psi|\mc C_b$ and any $\sigma_b\in H^0\big(\mc C_b,\scr E^*|\mc C_b\otimes\omega_{\mc C_b}(k\SX(b))\big)$, is given by
\begin{align*}
\bk{\delta\psi|\mc C_b,\sigma_b}=\sum_{i=1}^N\Res_i\bk{s_i|_k,\sigma_b}.
\end{align*}
Since for each $1\leq i\leq N$, $\bk{s_i-s_i|_k,\sigma_b}$ has removable singularity at $z=0$, we have $\Res_i\bk{s_i-s_i|_k,\sigma_b}=0$. Therefore, 
\begin{align*}
\bk{\delta\psi|\mc C_b,\sigma_b}=\sum_{i=1}^N\Res_i\bk{s_i,\sigma_b}=0.
\end{align*}
Thus $\delta\psi|\mc C_b=0$ for any $b$. By Theorem \ref{lb2}-(a), $\dim H^1\big(\mc C_b,(\scr E|\mc C_b)(-k\SX(b))\big)$ is locally constant over $b\in\mc B$. So by Theorem \ref{lb2}-(b), $\delta\psi$ is constantly $0$ in each fiber of the locally free $\scr O_{\mc B}$-module $R^1\pi_*\scr E(-k\SX)$. This proves that $\delta\psi=0$. 

By \eqref{eq11} and that $\pi_*\scr E(-k\SX)=0$, for any $b\in\mc B$ there exist an open $V_b\subset\mc B$ containing $b$  and a unique $s|_k\in\big(\pi_*\scr E(p\SX)\big)(V_b)= H^0(\pi^{-1}(V_b),\scr E(p\SX))$ which is sent to $\psi\in \pi_*\big(\scr E(p\SX)/\scr E(-k\SX)\big)(V_b)$. So near $\sgm_i(V_b)$, $s|_k$ has series expansion
\begin{align}
s|_k=s_i|_k+\bullet z^k+\bullet z^{k+1}+\cdots.\label{eq12}
\end{align}
By such uniqueness, these locally chosen $s|_k$ (over all $V_p$) are compatible with each other, which produce a global $s|_k\in H^0(\pi^{-1}(\mc B),\scr E(p\SX))$ which has series expansion \eqref{eq12} near $\sgm_i(\mc B)$. Note that this result holds for any $k\geq k_0$. So we have $s|_{k_0}=s|_{k_0+1}=s|_{k_0+2}=\cdots$, again by the uniqueness. Let $s=s|_{k_0}$. Then $s$ has series expansion $s_i$ near $\sgm_i(\mc B)$ for each $i$.
\end{proof}

\section{Nodal curves}\label{lb20}

\subsection*{The structure sheaf $\scr O_C$, and the invertible sheaves $\Theta_C$ and $\omega_C$}

Choose $M\in\Nbb$, and let $\wtd{\fk X}=(\wtd C;y_1',\dots,y_M';y_1'',\dots,y_M'')$ be a $2M$-pointed compact Riemann surface. (Here we do not require that each component of $\wtd C$ contains a point. In particular, $M$ could be $0$.) One can define a (possibly singular) \textbf{complex curve} $C$ in the following way. As a topological space, $C$ is the quotient space of $\wtd C$ defined by the identification $y_1'=y_1'',\dots,y_M'=y_M''$. The quotient map is denoted by $\nu:\wtd C\rightarrow C$. The points $x_1'=\nu(y_1'),\dots,x_M'=\nu(y_M')$ in $\wtd C$ are called \textbf{nodes}, and $\wtd{\fk X}$ (or just $\wtd C$) is called the \textbf{normalization} of $C$. To define the structure sheaf $\scr O_C$, we choose any open $U\subset C$. Then $\scr O(U)\equiv\scr O_C(U)$ is the set of all  $f\in\scr O_{\wtd C}(\nu^{-1}(U))$ such that $f(y_j')=f(y_j'')$ whenever $x_j'\in U$. When $M=0$, $C$ is just a compact Riemann surface. If $M>0$, we say that $C$ is a \textbf{nodal curve}.

We now describe locally free $\scr O_C$-modules in terms of vector bundles. Suppose that $\wtd {\scr E}$ is a (holomorphic) vector bundle on $\wtd C$. For any $j=1,\dots,M$, we fix an identification of fibers $\wtd{\scr E}|y_j'\simeq \wtd{\scr E}|y_j''$. Then the sheaf $\scr E$ is defined such that for any open $U\subset C$, $\scr E(U)$ is the set of all $s\in\wtd{\scr E}(\nu^{-1}(U))$ satisfying $s(y_j')=s(y_j'')$. Then $\scr E$ is obviously an $\scr O_C$-module, which is easily seen to be locally free.  ($\wtd{\scr E}$ is indeed equivalent to the pull back of $\scr E$.) Conversely, it is not hard to show that any locally free $\scr O_C$-module $\scr E$ arises in such a way.

Next we define an \textbf{invertible} $\scr O_C$-module $\Theta_C$ (i.e., a locally free $\scr O_C$-module with rank 1).   \index{zz@$\Theta_C$} Its dual sheaf $\omega_C:=\Theta_C^{-1}$ will be called the \textbf{dualizing sheaf of $C$}. \index{zz@$\omega_C$} For any $1\leq j\leq M$ we choose local coordinates $\xi_j,\varpi_j$ of $y_j',y_j''$ respectively. If an open subset $U\subset C$ contains no $x_1',\dots,x_M'$, then $\Theta_C(U)$ is defined to be $\Theta_{\wtd C}(\nu^{-1}(U))$. If $U$ is  neighborhood of some $x_j'$ such that $\nu^{-1}(U)$ is a disjoint union of neighborhoods $V_j'\ni y_j'$ and $V_j''\ni y_j''$, and that $\xi_j$ and $\varpi_j$ are defined on $V_j',V_j''$ respectively,  then $\Theta_C(U)$ is the $\scr O_C(U)$-submodule of $\Theta_{\wtd C}(V_j'\cup V_j'')$ generated by
\begin{align}
\xi_j\partial_{\xi_j}-\varpi_j\partial_{\varpi_j}.\label{eq57}
\end{align}
To be more precise, $\xi_j\partial_{\xi_j}-\varpi_j\partial_{\varpi_j}$ is the section in $\Theta_{\wtd C}(V_j'\cup V_j'')$ whose restrictions to $V_j'$ and $V_j''$ are $\xi_j\partial_{\xi_j}$ and $-\varpi_j\partial_{\varpi_j}$ respectively. One checks easily that $\Theta_C(U)$ is a free $\scr O_C(U)$-module.  Such definition is independent of the choice of local coordinates. For a general $U$, $\Theta_C(U)$ is defined by gluing as described in Section \ref{lb4}.

Alternatively, consider the line bundle on $\wtd C$ defined by
\begin{align}
\wtd{\Theta_C}=\Theta_{\wtd C}\big(-\sum_{j=1}^M(y_j'+y_j'')\big),\label{eq14}
\end{align}
which will be the pull back of $\Theta_C$. For any $1\leq j\leq M$, we need to choose identification of the lines (one dimensional fibers) $\wtd{\Theta_C}|y_j'\simeq\wtd{\Theta_C}|y_j''$. For this purpose we choose a local coordinate $\xi_j$ of $y_j'$. Note that the restriction $\xi_j\partial_{\xi_j}|y_j'$ of $\xi_j\partial_{\xi_j}$ to $y_j'$, which is an element of the fiber $\wtd{\Theta_C}|y_j'\simeq\Theta_{\wtd C}(-y_j')|y_j'$, is independent of the choice of the local coordinate. Choose similarly a local coordinate $\varpi_j$ of $y_j''$. Then the identification $\wtd{\Theta_C}|y_j'\simeq\wtd{\Theta_C}|y_j''$ is defined by
\begin{align}
\xi_j\partial_{\xi_j}|y_j'= -\varpi_j\partial_{\varpi_j}|y_j''.
\end{align}
One can then define $\Theta_C$ using $\wtd{\Theta_C}$ and the chosen identification. For the dualizing sheaf $\omega_C$, its pull back $\wtd{\omega_C}$ is
\begin{align}
\wtd{\omega_C}=\omega_{\wtd C}\big(\sum_{j=1}^M(y_j'+y_j'')\big),\label{eq218}
\end{align}
and the identification $\wtd{\omega_C}|y_j'\simeq\wtd{\omega_C}|y_j''$ is given by
\begin{align}
\frac{d\xi_j}{\xi_j}\Big|y_j'=-\frac{d\varpi_j}{\varpi_j}\Big|y_j''.\label{eq159}
\end{align}

\subsection*{Vanishing of higher order cohomology groups}

We shall generalize some results for compact Riemann surfaces to the nodal curve $C$. Our first result is that $H^q(C,\scr E)=0$ for any $q>1$ and any locally free $\scr E$. Again, we let $\wtd{\fk X}=(\wtd C;y_1',\dots,y_M';y_1'',\dots,y_M'')$ be the normalization of $C$. So $x_1'=\nu(y_1'),\dots,x_M'=\nu(y_M')$ are the nodes. For each $1\leq j\leq M$, we choose a neighborhood $U_j\subset C$ of $x_j'$ such that $\nu^{-1}(U_j)$ is the disjoint union of neighborhoods $V_j'\ni y_j'$ and $V_j''\ni y_j''$, and that the complex manifolds with marked points $(V_j',y_j')$ and $(V_j'',y_j'')$ are biholomorphic to the open unit disc $\mc D_1=\{z\in\mathbb C:|z|<1 \}$ with marked point $0$. Then $U_j$ (as a complex space) is biholomorphic to the complex subspace of $\mc D_1\times \mc D_1$ defined by $\{(z,w)\in \mc D_1\times \mc D_1:zw=0\}$. 

The complex manifold $\mc D_1\times \mc D_1$ belongs to a very important class of complex spaces, called \textbf{Stein spaces}, which satisfy Cartan's theorems A and (equivalently) B. If $X$ is a Stein space and $\scr F$ is a coherent $\scr O_X$-module (in particular, if $\scr F$ is  locally free), then \textbf{Theorem B} says that $H^q(X,\scr F)=0$ when $q>0$, and \textbf{Theorem A} says that the global sections of $\scr F$ generate every stalk $\scr F_x$ as an $\scr O_{X,x}$-module.  All non-compact connected Riemann surfaces are Stein spaces. A  product  of two Stein spaces is a Stein space. Closed complex subspaces of Stein spaces are Stein spaces. A finite intersection of Stein open subsets is Stein. References of these results can be found in \cite{GR84} section 1.4 or \cite{GPR94} section III.3. 

As $U_j$ ($1\leq j\leq M$) is biholomorphic to a (singular) hypersurface of $\mc D_1\times \mc D_1$, $U_j$ is a Stein space. Set $U_0=C-\{x_1',\dots,x_M'\}$. Then $\fk U=\{U_j:0\leq j\leq M \}$ is an open cover of $C$.

\begin{thm}\label{lb7}
Let $\scr E$ be locally free. Then $H^q(C,\scr E)=0$ for any $q>1$.
\end{thm}

\begin{proof}
Assume without loss of generality that $C$ is connected. We know  $H^q(C,\scr E)=0$ when $C$ is smooth (i.e., a compact Riemann surface). Thus it suffices to assume that $C$ contains at least one node. Then each connected component of $U_0$ is not compact. The same is true for $U_0\cap U_j$ for any $1\leq j\leq M$. We assume that $U_1,\dots,U_M$ are small enough so that they are mutually disjoint. Then for any $0\leq i,j\leq M$, $U_j$ and $U_i\cap U_j$ are Stein spaces, thus $H^p(U_i)=H^p(U_i\cap U_j)=0$ for any $p> 0$. Notice that the intersection of any three distinct open sets in $\fk U$ is $0$. Thus $H^p(U_{i_1}\cap\cdots \cap U_{i_n},\scr E)=0$ for any $1\leq i_1,\dots, i_n\leq M$ and $p>0$. Therefore, by Leray's theorem, $H^q(C,\scr E)=H^q(\fk U,\scr E)$. When $q\geq 2$, $U_{i_1}\cap U_{i_2}\cap\cdots\cap U_{i_{q+1}}=\emptyset $ for any $0\leq i_1<i_2<\cdots<i_{q+1}\leq M$. Thus any \v Cech $q$-cocycle is zero, which shows $Z^q(\fk U,\scr E)=0$ and hence $H^q(\fk U,\scr E)=0$.
\end{proof}

\subsection*{Serre duality}

Serre duality holds for the nodal curve $C$; see \cite{ACG11} section 10.2. In the following, we present a proof using the Serre duality for compact Riemann surfaces. The pairing for duality is also constructed in the proof. Let us first introduce the following notation. Let $\wtd{\fk X}=(\wtd C;y_1,\dots,y_N;y_1',\dots,y_M';y_1'',\dots,y_M'')$ be $(N+2M)$-pointed. We assume each component of $\wtd C$ contains $y_i$ for some $1\leq i\leq N$. Then the quotient map $\nu:\wtd C\rightarrow C$ identifying $y_1',\dots,y_M'$ with $y_1'',\dots,y_M''$ respectively defines an  \textbf{$N$-pointed complex curve} $\fk X=(C;x_1,\dots,x_N)$, where $x_1=\nu(y_1),\dots,x_N=\nu(y_N)$ are the marked points of $C$, and $x_1'=\nu(y_1'),\dots,x_M'=\nu(y_M')$ are the nodes. Again, we call $\wtd{\fk X}$ the normalization of $\fk X$. Set $\SX=x_1+x_2+\cdots+x_N$ to be a divisor either of $C$ or of $\wtd C$.

\begin{thm}[Serre duality]\label{lb150}
For the nodal curve $C$ and a locally free $\scr O_C$-module $\scr E$, the relation \eqref{eq2} holds.
\end{thm}

\begin{proof}
Assume without loss of generality that $C$ is connected. We shall prove \eqref{eq2} for $p=1$, i.e., 
\begin{align*}
H^1(C,\scr E\otimes \omega_C)\simeq H^0(C,\scr E^*)^*.
\end{align*}
Then, by replacing $\scr E$ with $\scr E\otimes\Theta_C$, we obtain the relation for $p=0$. For simplicity, we assume $M=1$ and write $y_1'=y_{-1},y_1''=y_{-2},x_1'=x_{-1}$. Thus $\fk X=(C;x_1,\dots,x_N)$ is a nodal curve with node $x_{-1}$, and it is obtained from the smooth curve $\wtd{\fk X}=(\wtd C;y_1,\dots,y_N;y_{-1};y_{-2})$ by gluing $y_{-1}$ and $y_{-2}$. We identify $C-\{x_{-1}\}$ with $\wtd C-\{y_{-1},y_{-2}\}$ in a natural way. So $x_1=y_1,\dots,x_N=y_N$. Choose mutually disjoint Stein neighborhoods $U_1,\dots,U_N$ of $x_1,\dots,x_N$ respectively disjoint from $x_{-1}$ (and hence from $y_{-1},y_{-2}$). Choose open discs $V_{-1}\ni y_{-1},V_{-2}\ni y_{-2}$ in $\wtd C$ and disjoint from $U_1,\dots,U_N$. Then $U_{-1}:=\nu(V_{-1}\cup V_{-2})$ is obtained from $V_{-1}$ and $V_{-2}$ by gluing $y_{-1},y_{-2}$. $U_{-1}$ is also a Stein space. Set $U_0=C-\{x_1,\dots,x_N,x_{-1}\}$, which has no compact connected component and is hence a Stein space. It follows that $\fk U:=\{U_{-1},U_0,U_1,\dots,U_N\}$ is a Stein cover of $C$. Hence, by Leray's theorem, we have $H^1(C,\scr E\otimes \omega_C)\simeq H^1(\fk U,\scr E\otimes \omega_C)$. Set $V_0=U_0,V_1=U_1,\dots,V_N=U_N$. Then $\fk V:=\{V_{-1},V_{-2},V_0,V_1,\dots,V_N\}$ is a Stein cover of $\wtd C$. 

Define a linear map
\begin{align}
\Psi :Z^1(\fk U,\scr E\otimes\omega_C)\rightarrow H^0(C,\scr E^*)^*\label{eq219}
\end{align}
as follows. Choose any \v{C}ech $1$-cochain $s=(s_{m,n})_{m,n=-1,0,1,\dots N}\in Z^1(\fk U,\scr E\otimes\omega_C)$. Then all the components of $s$, except possibly $s_{n,0}$ and $s_{0,n}$ (where $-1\leq n\leq N$), are zero. We set $\sigma_n=s_{n,0}=-s_{0,n}$ for $-1\leq n\leq N$ (in particular, $\sigma_0=0$), and let $\sigma_{-2}=\sigma_{-1}$. For each $-2\leq n\leq N$, we choose an anticlockwise circle $\gamma_n$ in $V_n-\{y_n \}$ surrounding $y_n$. The anticlockwiseness is understood under local coordinates of $\wtd C$ at $y_n$. Then $\Psi(s)$ is defined such that for any $t\in H^0(C,\scr E^*)$, the evaluation of $\Psi(s)$ with $t$ is
\begin{align}
\bk{s,t}=\sum_{n=-2}^N~\oint_{\gamma_n}\bk{\sigma_n,t}.\label{eq157}
\end{align}
The proof will be completed if we can show that $\Psi$ is surjective with kernel $B^1(\fk U,\scr E\otimes\omega_C)$. We divide the proof into several steps.

Step 1. We first show that $\Psi$ is surjective. Let $\SX=x_1+\cdots x_N$. Then for sufficiently large $k$, we have $H^0(C,\scr E^*(-k\SX))=0$. (See the proof of Proposition \ref{lb121}.) Therefore, from the exact sequence
\begin{align*}
0\rightarrow H^0(C,\scr E^*(-k\SX))\rightarrow H^0(C,\scr E^*)\rightarrow H^0(C,\scr G)
\end{align*} 
where $\scr G=\scr E^*/\scr E^*(-k\SX)$, we see that $H^0(C,\scr E^*)$ is naturally a subspace of $H^0(C,\scr G)$. For each $1\leq n\leq N$, choose a trivilization $\scr E_{U_n}\simeq E_n\otimes_\Cbb\scr O_{U_n}$ which yields the dual trivilization $\scr E_{U_n}^*\simeq E_n^*\otimes_\Cbb\scr O_{U_n}$. Let $\eta_n$ be a coordinate of $U_n$ satisfying $\eta_n(x_n)=0$. Then we have a natural equivalence $H^0(C,\scr G)\simeq\bigoplus_n E_n^*\otimes \Cbb^k$ such that each  $(v_n^0,\dots,v_n^{k-1})_{1\leq n\leq N}\in \bigoplus_n E_n^*\otimes \Cbb^k$   corresponds to the section of $\scr G$ whose restriction to each $U_n$ is $\sum_{0\leq l\leq k-1}v_n^l\eta_n^l$, and whose restriction to $C-\{x_1,\dots,x_N\}$ is $0$. Choose any $\psi\in H^0(C,\scr E^*)^*$ and extend it to a linear functional on $H^0(C,\scr G)$. Then, for each $n$, one can choose $\alpha_n^0,\dots,\alpha_n^{k-1}\in E_n$ such that the evaluation of $\psi$ with any $(v_n^0,\dots,v_n^{k-1})_{1\leq n\leq N}$ is $\sum_{1\leq n\leq N}\sum_{0\leq l\leq k-1}\alpha_n^l(v_n^l)$. Choose $s=(s_{m,n})_{m,n=-1,0,1,\dots N}\in Z^1(\fk U,\scr E\otimes\omega_C)$ whose only possibly non-zero components are $s_{n,0}=-s_{0,n}=\sum_{0\leq l\leq k-1}\alpha_n^l\eta_i^{-l-1}d\eta_i$ (where $1\leq n\leq N$). One checks easily that $\Psi(s)=\psi$.


Step 2. To finish the proof, we need to show that  $\Ker(\Psi)=B^1(\fk U,\scr E\otimes\omega_C)$. Note that $\scr E$ can be obtained from a vector bundle $\wtd{\scr E}$  on  $C$  by identifying the fibers $\wtd{\scr E}|y_{-1}$ and $\wtd{\scr E}|y_{-2}$ via an isomorphism. Under this viewpoint, $H^0(C,\scr E^*)$ is naturally a subspace of $H^0(\wtd C,\wtd{\scr E}^*)$ consisting of global sections whose values at $y_{-1}$ and at $y_{-2}$ agree. In particular, this is true when the two values are $0$. We thus have
\begin{align*}
H^0(\wtd C,\wtd{\scr E}^*(-y_{-1}-y_{-2}))\subset H^0(C,\scr E^*)\subset H^0(\wtd C,\wtd{\scr E}^*).
\end{align*}
Therefore, $\Ker(\Psi)$ vanishes when evaluating on $H^0(\wtd C,\wtd{\scr E}^*(-y_{-1}-y_{-2}))$. Also, since $U_m\cap U_n$ does not contain the node $x_{-1}$ whenever $m\neq n$, we have natural identifications
\begin{align*}
Z^1(\fk V,\wtd{\scr E}\otimes \omega_{\wtd C}(y_{-1}+y_{-2}))=Z^1(\fk U,\scr E\otimes \omega_C)=Z^1(\fk V,\wtd{\scr E}\otimes \omega_{\wtd C}).
\end{align*}
It is easy to see that the pairing of the smooth Serre duality
\begin{align}
H^1(\wtd C,\wtd{\scr E}\otimes\omega_{\wtd C}(y_{-1}+y_{-2}))\simeq H^0(\wtd C,\wtd{\scr E}^*(-y_{-1}-y_{-2}))^*\label{eq158}
\end{align}
is compatible with the one of \eqref{eq219} defined by \eqref{eq157}. Thus, if we regard any  $s\in\Ker(\Psi)$ as an element of $Z^1(\fk V,\wtd{\scr E}\otimes \omega_{\wtd C}(y_{-1}+y_{-2}))$, then it becomes zero in $H^1(\fk V,\wtd{\scr E}\otimes\omega_{\wtd C}(y_{-1}+y_{-2}))$. Therefore,
\begin{align*}
\Ker(\Psi)\subset B^1(\fk V,\wtd{\scr E}\otimes\omega_{\wtd C}(y_{-1}+y_{-2})).
\end{align*}

Choose any
\begin{align*}
\sgm\in B^1(\fk V,\wtd{\scr E}\otimes\omega_{\wtd C}(y_{-1}+y_{-2})).
\end{align*}
We shall show that $\sgm\in B^1(\fk U,\scr E\otimes\omega_C)$ if and only if $\sgm\in\Ker(\Psi)$. This will finish the proof.



Step 2-(a). Choose 
\begin{align*}
\sigma\in C^0(\fk V,\wtd{\scr E}\otimes \omega_{\wtd C}(y_{-1}+y_{-2}))
\end{align*}
(i.e., $\sigma=(\sigma_n)_{-2\leq n\leq N}$ is a $0$-cochain of $\wtd{\scr E}\otimes \omega_{\wtd C}(y_{-1}+y_{-2})$ over the cover $\fk V$) such that
\begin{align*}
\sgm=\delta(\sigma).
\end{align*}
Since $\sgm$ can be regarded as an element of $Z^1(\fk U,\scr E\otimes\omega_C)$, we can  calculate the pairing $\bk{\sgm,t}$ defined by \eqref{eq157} for any
\begin{align*}
t\in H^0(C,\scr E^*).
\end{align*}
Write $\sgm=(\sgm_{m,n})_{-2\leq m,n\leq N}$. The only possibly non-zero components of $\sgm$ are $\sgm_{n,0}=-\sgm_{0,n}$ where $-2\leq n\leq N$ and $n\neq 0$. It is clear that $\sgm_{n,0}=(\sigma_n-\sigma_0)|_{U_n\cap U_0}$ for any such $n$. Then
\begin{align*}
\bk{\sgm,t}=\sum_{n=-2}^N~\oint_{\gamma_n}\bk{\sigma_n-\sigma_0,t}.
\end{align*}
By residue theorem (or Stokes theorem), we have $\sum_{n=-2}^{N}\oint_{\gamma_n}\bk{\sigma_0,t}=0$. Moreover, when $n\geq 1$, since $\sigma_n$ can be defined as a section of $\scr E$ on $U_n$, $\bk{\sigma_n,t}$ is a holomorphic function on $U_n$. So $\bk{\sigma_n,t}=0$. We thus have
\begin{align*}
\bk{\sgm,t}=\oint_{\gamma_{-1}}\bk{\sigma_{-1},t}+\oint_{\gamma_{-2}}\bk{\sigma_{-2},t}.
\end{align*}

We set $E=\wtd{\scr E}|y_{-1}=\wtd{\scr E}|y_{-2}$. Then its dual space is $E^*=\wtd{\scr E}^*|y_{-1}=\wtd{\scr E}^*|y_{-2}$. We choose trivialization $\scr E_{V_n}\simeq E\otimes_\Cbb\scr O_{V_n}$ which yields the dual one $\scr E^*_{V_n}\simeq E^*\otimes_\Cbb\scr O_{V_n}$ for $n=-1,-2$. Let $\xi$ (resp. $\varpi$) be a coordinate of $V_{-1}$ (resp. $V_{-2}$) satisfying $\xi(y_{-1})=0$ (resp. $\varpi(y_{-2})=0$). Since $t\in H^0(C,\scr E^*)$, $t$ can  be viewed as an element of $H^0(\wtd C,\wtd{\scr E}^*)$ such that $t(y_{-1})=t(y_{-2})$ (which are elements of $E^*$). Set 
\begin{align*}
e_1=\oint_{\gamma_{-1}}\sigma_{-1},\qquad e_2=-\oint_{\gamma_{-2}}\sigma_{-2}
\end{align*}
(which are vectors in $E$), i.e.
\begin{align}
\sigma_{-1}=e_1d\xi/\xi+\sum_{l\geq 0}\blt \xi^ld\xi,\qquad \sigma_{-2}=-e_2d\varpi/\varpi+\sum_{l\geq 0}\blt \varpi^ld\varpi \label{eq160}
\end{align}
We thus have
\begin{align}
\bk{\sgm,t}=\bk{e_1-e_2,t(y_{-2})}.\label{eq161}
\end{align}
By \eqref{eq159}, we have $\sigma\in C^0(\fk U,\scr E\otimes\omega_C)$ (i.e., $\sigma$ is a $0$-cochain of $\scr E\otimes\omega_C$) if and only if $e_1-e_2=0$. If $\sgm\in B^1(\fk U,\scr E\otimes\omega_C)$, then one can choose $\sigma$ (which satisfies $\sgm=\delta(\sigma)$) to be an element of $C^0(\fk U,\scr E\otimes\omega_C)$, which implies $e_1-e_2=0$ and hence $\bk{s,t}=0$ for any $t\in H^0(C,\scr E^*)$. So $\sgm\in\Ker(\Psi)$. 

Step 2-(b). We now prove the other direction. Assume $\sgm\in\Ker(\Psi)$. We shall show that $\sgm\in B^1(\fk U,\scr E\otimes\omega_C)$. Note that \eqref{eq161} is zero for each $t\in H^0(C,\scr E^*)$. Consider the short exact sequence of $\scr O_{\wtd C}$-modules
\begin{align*}
0\rightarrow\wtd{\scr E}\otimes\omega_{\wtd C}\rightarrow \wtd{\scr E}\otimes\omega_{\wtd C}(y_{-1}+y_{-2})\rightarrow \scr F\rightarrow 0
\end{align*}
where $\scr F=\wtd{\scr E}\otimes\omega_{\wtd C}(y_{-1}+y_{-2})/\wtd{\scr E}\otimes\omega_{\wtd C}$. Then we have an exact sequence
\begin{align*}
H^0(\wtd C,\wtd{\scr E}\otimes\omega_{\wtd C}(y_{-1}+y_{-2}))\rightarrow H^0(\wtd C,\scr F)\xrightarrow{\delta} H^1(\wtd C,\wtd{\scr E}\otimes\omega_{\wtd C}).
\end{align*}
For any $e_0\in E$, we choose $\varepsilon\in H^0(\wtd C,\scr F)$ (which clearly depends on $e_0$) such that
\begin{gather*}
\varepsilon|_{V_{-1}}=e_0d\xi/\xi,\qquad \varepsilon|_{V_{-2}}=(-e_0+e_2-e_1)d\varpi/\varpi ,\qquad  \varepsilon|_{\wtd C-\{y_{-1},y_{-2}\}}=0.
\end{gather*}
We claim that there exists $e_0$ such that $\delta(\varepsilon)=0$. Suppose this has been proved. Then there exists $\lambda\in H^0(\wtd C,\wtd{\scr E}\otimes\omega_{\wtd C}(y_{-1}+y_{-2}))$ which is sent to $\varepsilon$. We treat $\lambda=(\lambda_n)_{-2\leq n\leq N}$ as an element of $Z^0(\fk V,\wtd{\scr E}\otimes\omega_{\wtd C}(y_{-1}+y_{-2}))$. Then
\begin{align*}
\lambda_{-1}=e_0d\xi/\xi+\sum_{l\geq 0}\blt \xi^ld\xi,\qquad \lambda_{-2}=(-e_0+e_2-e_1)d\varpi/\varpi+\sum_{l\geq 0}\blt \varpi^ld\varpi.
\end{align*}
If we compare this relation with \eqref{eq160},  we see that $(\sigma+\lambda)_{-1}=(e_1+e_0)d\xi/\xi+\cdots$ and $(\sigma+\lambda)_{-2}=-(e_1+e_0)d\varpi/\varpi+\cdots$. Thus $\sigma+\lambda$ can be viewed as an element of $C^0(\fk U,\scr E\otimes\omega_C)$. Since $\delta(\lambda)=0$, we have $\delta(\sigma+\lambda)=\delta(\sigma)=\sgm$, which finishes the proof.

Let us prove the existence of $e_0$ such that $\delta(\varepsilon)=0$. Notice that $H^1(\wtd C,\wtd{\scr E}\otimes\omega_{\wtd C})=H^1(\fk V,\wtd{\scr E}\otimes\omega_{\wtd C})$, and $\delta(\varepsilon)$ can be represented by an element  $(\delta(\varepsilon)_{m,n})_{-2\leq m,n\leq N}$ in $Z^1(\fk V,\wtd{\scr E}\otimes\omega_{\wtd C})$ whose only non-zero components are $\delta(\varepsilon)_{-1,0}=-\delta(\varepsilon)_{0,-1}=e_0d\xi/\xi$ and $\delta(\varepsilon)_{-2,0}=-\delta(\varepsilon)_{0,-2}=(-e_0+e_2-e_1)d\varpi/\varpi$. Notice the smooth Serre duality
\begin{align*}
H^1(\wtd C,\wtd{\scr E}\otimes\omega_{\wtd C})\simeq H^0(\wtd C,\wtd{\scr E}^*)^*
\end{align*}
defined by residue pairing. If we regard $\delta(\varepsilon)$ as a linear functional on $H^0(\wtd C,\wtd{\scr E}^*)$, then its evaluation with any $\tau\in H^0(\wtd C,\wtd{\scr E}^*)$ is
\begin{align*}
\bk{\delta(\varepsilon),\tau}=&\bk{e_0,\tau(y_{-1})}+\bk{-e_0+e_2-e_1,\tau(y_{-2})}\\
=&\bk{e_0,\tau(y_{-1})-\tau(y_{-2})}-\bk{e_1-e_2,\tau(y_{-2})}.
\end{align*}
Then $\delta(\varepsilon)=0$ will follow if we can find $e_0\in E$ satisfying
\begin{align*}
\bk{e_0,\tau(y_{-1})-\tau(y_{-2})}=\bk{e_1-e_2,\tau(y_{-2})}
\end{align*}
for any $\tau\in H^0(\wtd C,\wtd{\scr E}^*)$. Indeed, one first  defines $e_0$ as a linear functional on $T:=\{\tau(y_{-1})-\tau (y_{-2}):\tau\in H^0(\wtd C,\wtd {\scr E}^*)\}$ using the above relation. Then $e_0$ is well defined: if $\tau(y_{-1})-\tau(y_{-2})=0$, then $\tau\in H^0(C,\scr E^*)$; since we assume $\sgm\in\Ker(\Psi)$, according to \eqref{eq161}, we have $\bk{e_2-e_1,\tau(y_{-2})}=0$. Now, as $T$ is a subspace of $E^*$, we can extend $e_0$ to a linear functional on $E^*$. Then $e_0$ is in $E$ and satisfies the desired relation.
\end{proof}

\begin{co}
If $\scr E$ is a locally free $\scr O_C$-module, then $\dim H^q(C,\scr E)<+\infty$ for any $q\in\mbb N$.
\end{co}

Thus one can define the character $\chi(C,\scr E)$ using \eqref{eq13}.

\begin{proof}
We have proved that $H^q(C,\scr E)=0$ when $q>1$. Let $\wtd{\scr E}=\nu^*\scr E$ where $\nu:\wtd C\rightarrow C$ is the normalization. Since $H^0(C,\scr E)$ is naturally a subspace of $H^0(\wtd C,\wtd{\scr E})$ where the latter is finite dimensional, so is $H^0(C,\scr E)$. Similarly, $H^0(C,\scr E^*\otimes\omega_C)$ is finite dimensional. So is $H^1(C,\scr E)$ by Serre duality. 
\end{proof}

\subsection*{Vanishing theorems}

\begin{pp}\label{lb121}
Assume that each connected component of $\wtd C$ (equivalently, each irreducible component of $C$) contains one of $x_1,\dots,x_N$, and set $D=x_1+\cdots+x_N$. Then Proposition \ref{lb62} and Corollary \ref{lb6} hold verbatim.
\end{pp}

\begin{proof}
The pull back of $\scr E$ along the normalization $\nu:\wtd C\rightarrow C$ is denoted by $\wtd{\scr E}$. We know that $\scr E$ can be obtained by gluing the fibers of $\wtd{\scr E}$ at the double points.  So $H^0(C,\scr E(-nD))$ is naturally a subspace of $H^0(\wtd C,\wtd{\scr E}(-nD))$ (consisting of sections whose values at $y_j'$ and at $y_j''$ agree), which vanishes when $n$ is sufficiently large. This proves Proposition \ref{lb62}. Corollary \ref{lb6} follows as in the smooth case by Serre duality.
\end{proof}

As an application of Prop. \ref{lb121}, we now give a better description of the pairing in Serre duality. Let $x_1,\dots,x_N$ be distinct smooth points on $C$ as in the proof of Serre duality Thm. \ref{lb150}. (Namely, each connected component of the normalization $\wtd C$ of $C$ contains at least one element of the preimage of $x_1,\dots,x_N$.)  We choose $U_0^+=C-\{x_1,\dots,x_N\}$ and $U_1,\dots,U_N$ mutually disjoint disks around $x_1,\dots,x_N$ that do not intersect the nodes. Choose an anticlockwise circle $\gamma_i\subset U_i$ around $x_i$. 

\begin{thm}
Let $\fk U^+=\{U_0^+,U_1,\dots,U_N\}$. Then $H^1(C,\scr E\otimes\omega_C)=H^1(\fk U^+,\scr E\otimes\omega_C)$. Moreover, Serre duality
\begin{align}
H^1(C,\scr E\otimes\omega_C)\simeq H^0(C,\scr E^*)^*
\end{align}
holds, and the isomorphism is realized by
\begin{align}
\Phi :Z^1(\fk U^+,\scr E\otimes\omega_C)\rightarrow H^0(C,\scr E^*)^*
\end{align}
such that for each $s=(s_{m,n})_{m,n=0,1,\dots,N}\in Z^1(\fk U^+,\scr E\otimes\omega_C)$ and $t\in H^0(C,\scr E^*)$, by setting $\sigma_n=s_{n,0}=-s_{0,n}$, we have
\begin{align}
\bk{s,t}=\sum_{n=1}^N\oint_{\gamma_n}\bk{\sigma_n,t}.
\end{align}
\end{thm}

\begin{proof}
This theorem follows directly from the fact that $U_0^+$ is indeed Stein, and hence that $\fk U^+$ is a Stein cover. (A quick argument is to embed $C$ as a closed complex subspace of $\Pbb^n$ for some $n$ and show that $U_0^+$ is the intersection of $C$ and the complement of a hyperplane in $\Pbb^n$.) But here we give a different argument which avoids showing that $U_0^+$ is Stein.

Assume for simplicity that $C$ has only one node $x_{-1}$; the general case follows from a similar argument. Following the notations in the proof of Thm. \ref{lb150}, we let $\fk U=\{U_{-1},U_0,U_1,\cdots,U_N\}$ where $U_{-1}$ is a small Stein neighborhood of the node $x_{-1}$, and $U_0=C-\{x_1,\dots,x_N, x_{-1}\}$. (So $U_0^+=U_0\cup U_{-1}$.) Then  $Z^1(\fk U^+,\scr E\otimes\omega_C)$ is naturally a subspace of $Z^1(\fk U,\scr E\otimes\omega_C)$ consisting of $1$-cocycles $s=(s_{m,n})_{m,n=-1,0,1,\dots,N}$ vanishing on $U_0\cap U_{-1}$, i.e. $s_{-1,0}=s_{0,-1}=0$. Moreover, an element in $Z^1(\fk U^+,\scr E\otimes\omega_C)$ is a coboundary with respect to $\fk U^+$ if and only if it is so with respect to $\fk U$. Therefore, the theorem follows from Thm. \ref{lb150}  and the pairing \eqref{eq219} \eqref{eq157} if we can show that every element of $H^1(C,\scr E\otimes\omega_C)=H^1(\fk U,\scr E\otimes\omega_C)$ is represented by an element of $Z^1(\fk U^+,\scr E\otimes\omega_C)$, i.e., represented by an element  $s\in Z^1(\fk U,\scr E\otimes\omega_C)$ such that $s_{-1,0}=-s_{0,-1}=0$.

Let $D=x_1+\cdots+x_N$. By Prop. \ref{lb121}, there is $k\in\Nbb$ such that $H^1(C,\scr E\otimes\omega_C(kD))=0$. The short exact sequence
\begin{align*}
0\rightarrow\scr E\otimes\omega_C\rightarrow\scr E\otimes\omega_C(kD)\rightarrow\scr G\rightarrow 0
\end{align*} 
(where $\scr G$ is the quotient of the previous two sheaves) gives a long one
\begin{align*}
H^0(C,\scr G)\xrightarrow{\delta} H^1(\fk U,\scr E\otimes\omega_C)\rightarrow H^1(C,\scr E\otimes\omega_C(kD))=0.
\end{align*}
So $\delta$ is surjective. Since $\scr G$ has support in $x_1,\dots,x_N$, by the explicit description of $\delta$, it is clear that any element in the image of $\delta$ is represented by a cocycle satisfying  $s_{-1,0}=-s_{0,-1}=0$.
\end{proof}

\begin{rem}\label{lb63}
If $C$ is connected, the number $g=\dim H^1(C,\scr O_C)$ is called the (arithmetic) \textbf{genus} of $C$. Thus $\chi(C,\scr O_C)=1-g$. Again, any line bundle $\scr L$ is equivalent to $\scr O_C(D)$ for some divisor $D=k_1x_1+\cdots +k_Nx_N$. (For the nodal curve, we assume none of $x_1,\dots,x_N$ is a node.) The argument is the same as for smooth curves:  By Proposition \ref{lb121}, one may find a divisor $D_0$ such that $H^1(C,\scr L(-D_0-x))=0$, where $x$ is any smooth point of $C$. Thus, the short exact sequence $0\rightarrow \scr L(-D_0-x)\rightarrow\scr L(-D_0)\rightarrow k_x\rightarrow 0$ (where $k_x\simeq\scr L(-D_0)/\scr L(-D_0-x)$ is the skyscraper sheaf) yields a surjective $H^0(C,\scr L(-D_0))\rightarrow H^0(C,k_x)$. So $H^0(C,\scr L(-D_0))$ is nonzero. Thus we may find a non-zero global meromorphic section $s$ of $\scr L$. Let $D=-\sum n_x\cdot x$ where the sum is over all smooth points of $C$, and $n_x$ is the unique integer such that (by choosing any local coordinate $z_x$ at $x$) $z_x^{n_x}s$ can be extended to a section of $\scr L$ on a neighborhood of $x$ whose value at the fiber $\scr L|x$ is non-zero. Then  $D$ is a finite divisor, and $f\mapsto fs$ defines an isomorphism $\scr O_C(D)\rightarrow \scr L$.

Define $\deg(\scr L)=\deg(D)$. Then the Riemann-Roch theorem \eqref{eq162} holds for $C$ and can be proved in exactly the same way: Identify $\scr L$ with $\scr O_C(D)$. Notice $\chi(C,k_x)=\dim H^0(C,k_x)=1$. Then the we have a short exact sequence $0\rightarrow\scr O_C(D)\rightarrow\scr O_C(D+x)\rightarrow k_x\rightarrow 0$ and hence a long one
\begin{align*}
&0\rightarrow H^0(C,\scr O_C(D))\rightarrow H^0(C,\scr O_C(D+x))\rightarrow H^0(C,k_x)\\
&\rightarrow  H^1(C,\scr O_C(D))\rightarrow H^1(C,\scr O_C(D+x))\rightarrow H^1(C,k_x)=0,
\end{align*}
which gives $\chi(C,\scr O_C(D))-\chi(C,\scr O_C(D+x))=-\chi(C,k_x)=-1$. Thus, \eqref{eq162} follows from induction and the base case $\chi(C,\scr O_C)=1-g$. Therefore, \eqref{eq3} also holds in the nodal case. 

As a consequence, we obtain again $\deg\omega_C=-\deg\Theta_C=2g-2$.
\end{rem}

As usual, a divisor $D$ of $C$ is called \textbf{effective} if $D=\sum_{i=1}^kn_ix_i$ where each $x_i$ is a smooth point and each $n_i\in\Zbb$ is non-negative.

\begin{thm}\label{lb10}
Let $C$ be a complex curve with $M$ nodes and normalization $\wtd C$. Let $x_1,\dots,x_N$ be smooth points of $C$, and assume that any connected component of $\wtd C$ contains one of these points. Set $D=x_1+\cdots+x_N$. Let $n\in\Zbb$, let $D'$ be an effective divisor of $C$,  and let $\wtd g$ be the maximal (resp. minimal) genus of the connected components of $\wtd C$ if $n\geq -1$ (resp. $n<-1$). Then for any $k>(n+1)(2\wtd g-2)+2M+\deg D'$, we have $H^1(C,\Theta_C^{\otimes n}(kD-D'))=0$.
\end{thm}

\begin{proof}
Choose any such $k$. By Serre duality, it suffices to prove $H^0(C,\omega_C^{\otimes(n+1)}(-kD+D'))=0$. Let $x_1',\dots,x_M'$ be the nodes of $C$. For each $j$, let $y_j',y_j''$ be the double points of $\wtd C$ mapped to $x_j'$ by $\nu$. Let $D_0=D'+\sum_{j=1}^M(y_j'+y_j'')$. Then, by \eqref{eq218}, $H^0(C,\omega_C^{\otimes(n+1)}(-kD+D'))$ is naturally a subspace of $H^0(\wtd C,\omega_{\wtd C}^{\otimes(n+1)}(-kD+D_0))$. Thus, it suffices to prove that for each connected component $\wtd C_0$ of $\wtd C$,  $H^0(\wtd C_0,\omega_{\wtd C_0}^{\otimes(n+1)}(-kD+D_0))$ is trivial. This follows from Theorem \ref{lb64} and the computation $\deg(\omega_{\wtd C_0}^{\otimes(n+1)}(-kD+D_0))\leq (n+1)(2\wtd g-2)-k+2M+\deg D'<0$. 
\end{proof}

\section{Families of complex curves}\label{lb24}

\subsection*{Smoothing the nodes}

For any $r>0$, let  $\mc D_r=\{z\in\mbb C:|z|<r \}$ and $\mc D_r^\times=\mc D_r-\{0\}$.\index{Dr@$\mc D_r,\mc D_r^\times$} If $r,\rho>0$, we define \index{zz@$\pi_{r,\rho}:\mc D_r\times\mc D_\rho\rightarrow\mc D_{r\rho}$}
\begin{align}
\pi_{r,\rho}:\mc D_r\times\mc D_\rho\rightarrow\mc D_{r\rho},\qquad (\xi,\varpi)\mapsto \xi\varpi.
\end{align}
$d\pi_{r,\rho}$ is surjective at $(\xi,\varpi)$ whenever $\xi\neq 0$ or $\varpi\neq 0$. So if $q\in\mc D_{r\rho}$ is not $0$, the fiber $\pi^{-1}(q)$ is smooth. But $\pi_{r,\rho}^{-1}(0)$ is singular, which is just a neighborhood of a node of a nodal curve. Denote also by $\xi$ and $\varpi$ the standard coordinates of $\mc D_r$ and $\mc D_\rho$, and set $q=\pi_{r,\rho}$, i.e.,
\begin{align*}
q:\mc D_r\times\mc D_\rho\rightarrow\Cbb,\qquad q=\xi\varpi.
\end{align*}
Then $(\xi,\varpi),(\xi,q),(\varpi,q)$ are coordinates  of $\mc D_r\times \mc D_\rho,\mc D_r^\times\times \mc D_\rho,\mc D_r\times \mc D_\rho^\times$ respectively. The standard tangent vectors of the coordinates $(\xi,\varpi),(\xi,q)$ are related by
\begin{gather}
\left\{ \begin{array}{l}
\partial_\xi=\partial_\xi-\xi^{-1}\varpi\cdot\partial_\varpi\\
\partial_q=\xi^{-1}\partial_\varpi
\end{array} \right.
\qquad
\left\{ \begin{array}{l}
\partial_\xi=\partial_\xi+\xi^{-1}q\cdot\partial_q\\
\partial_\varpi=\xi\partial_q
\end{array} \right.\label{eq74}
\end{gather}
The formulae between $(\xi,\varpi),(\varpi,q)$ are similar.

It is easy to see that  $(\xi,q)(\mc D_r^\times\times\mc D_\rho)$ (resp. $(\varpi,q)(\mc D_r\times\mc D_\rho^\times)$) is precisely the subset of all $(\xi_0,q_0)\in\mc D_r\times\mc D_{r\rho}$ (resp. $(\varpi_0,q_0)\in\mc D_\rho\times\mc D_{r\rho}$) satisfying
\begin{align}
\frac{|q_0|}{\rho}<|\xi_0|<r\qquad\text{resp.}\qquad  \frac{|q_0|}{r}<|\varpi_0|<\rho.\label{eq181}
\end{align}
We choose closed subsets $E_{r,\rho}'\subset \mc D_r\times\mc D_{r\rho}$ and $E_{r,\rho}''\subset\mc D_\rho\times\mc D_{r\rho}$ such that
\begin{gather}
(\xi,q):\mc D_r^\times\times\mc D_\rho\xrightarrow{\simeq}\mc D_r\times\mc D_{r\rho}-E_{r,\rho}',\nonumber\\
(\varpi,q): \mc D_r\times\mc D_\rho^\times\xrightarrow{\simeq}\mc D_\rho\times\mc D_{r\rho}-E_{r,\rho}''.\label{eq18}
\end{gather}


\subsection*{Sewing families of compact Riemann surfaces}

Choose $M\in\mathbb N$. Consider a family of $2M$-pointed compact Riemann surfaces with local coordinates
\begin{align}
\wtd{\fk X}=(\wtd\pi:\wtd{\mc C}\rightarrow\wtd{\mc B};\sgm_1',\dots,\sgm_M';\sgm_1'',\dots,\sgm_M'';\xi_1,\dots,\xi_M;\varpi_1,\dots,\varpi_M).\label{eq61}
\end{align}
We do not assume that every component of any fiber contains a marked point. For each $1\leq j\leq M$ we choose $r_j,\rho_j>0$ and a neighborhood $U_j'$ (resp. $U_j''$) of $\sgm_j'(\wtd {\mc B})$ (resp. $\sgm_j''(\wtd {\mc B})$) such that
\begin{gather}
(\xi_j,\wtd\pi):U_j'\xrightarrow{\simeq} \mc D_{r_j}\times\wtd{\mc B}\qquad\text{resp.}\qquad (\varpi_j,\wtd\pi):U_j''\xrightarrow{\simeq} \mc D_{\rho_j}\times\wtd{\mc B}\label{eq113}
\end{gather}
is a biholomorphic map. We also assume that these $r_i$ and $\rho_j$ are small enough so that the neighborhoods $U_1',\dots,U_M',U_1'',\dots,U_M''$ are mutually disjoint. Identify
\begin{gather*}
U_j'=\mc D_{r_j}\times\wtd{\mc B}\qquad\text{resp.}\qquad U_j''=\mc D_{\rho_j}\times\wtd{\mc B}
\end{gather*}
via the above maps. Then $\xi_j,\varpi_j$ (when restricted to the first components) become the standard coordinates of $\mc D_{r_j},\mc D_{\rho_j}$ respectively, and $\wtd\pi$ is the projection onto the $\wtd{\mc B}$-component. Set $q_j=\xi_j\varpi_j:\mc D_{r_j}\times\mc D_{\rho_j}\rightarrow\mc D_{r_j\rho_j}$ as previously.  

We now construct a family of complex curves $\fk X=(\pi:\mc C\rightarrow\mc B)$ as follows. \index{XX@$\wtd {\fk X}$ and $\fk X$} Let
\begin{align}
\mc D_{r_\bullet\rho_\bullet}=\mc D_{r_1\rho_1}\times\cdots\times\mc D_{r_M\rho_M},\qquad\mc B=\mc D_{r_\blt\rho_\blt}\times\wtd{\mc B}.\label{eq65}
\end{align}
We shall freely switch the order of Cartesian product. For each $1\leq j\leq M$, we also define $\mc D_{r_\blt\rho_\blt\setminus j}$ to be the product of all $\mc D_{r_1\rho_1},\dots,\mc D_{r_M\rho_M}$ except $\mc D_{r_j\rho_j}$. So
\begin{align*}
\mc D_{r_\bullet\rho_\bullet}=\mc D_{r_j\rho_j}\times \mc D_{r_\blt\rho_\blt\setminus j}.
\end{align*}
Recall that by \eqref{eq18}, $E_{r_j,\rho_j}'\subset \mc D_{r_j}\times\mc D_{r_j\rho_j}$ and $E_{r_j,\rho_j}''\subset \mc D_{\rho_j}\times\mc D_{r_j\rho_j}$. \index{Drrho@$\mc D_{r_\blt\rho_\blt},\mc D_{r_\blt\rho_\blt\setminus j}$} So
\begin{gather*}
F_j':=E_{r_j,\rho_j}'\times\mc D_{r_\blt\rho_\blt\setminus j}\times\wtd{\mc B}\qquad\subset \mc D_{r_j}\times\mc D_{r_\blt\rho_\blt}\times\wtd{\mc B}\quad(= U_j'\times \mc D_{r_\blt\rho_\blt}),\\
F_j'':=E_{r_j,\rho_j}''\times\mc D_{r_\blt\rho_\blt\setminus j}\times\wtd{\mc B}\qquad\subset \mc D_{\rho_j}\times\mc D_{r_\blt\rho_\blt}\times\wtd{\mc B}\quad(= U_j''\times \mc D_{r_\blt\rho_\blt})
\end{gather*}
are subsets of $\wtd{\mc C}\times\mc D_{r_\blt\rho_\blt}$. They are the subsets  we should discard in the sewing process.  Set \begin{align}
W_j:=\mc D_{r_j}\times\mc D_{\rho_j}\times\mc D_{r_\bullet\rho_\bullet\setminus j}\times\wtd{\mc B}.\label{eq19}
\end{align}
We  glue $\wtd{\mc C}\times \mc D_{r_\blt\rho_\blt}$ (with $F_j',F_j''$ all removed) with all $W_j$ and obtain a complex manifold $\mc C$. 

To be more precise, we define
\begin{align}
\mc C=\bigg(\coprod_{j=1}^MW_j\bigg)\coprod\bigg(\wtd{\mc C}\times \mc D_{r_\blt\rho_\blt}-\bigcup_{j=1}^M F_j'-\bigcup_{j=1}^M F_j''\bigg)\bigg/\sim \label{eq112}
\end{align}
where the equivalence $\sim$ is described as follows. Consider the following subsets of $W_j$:
\begin{gather}
W_j'=\mc D_{r_j}^\times\times\mc D_{\rho_j}\times\mc D_{r_\bullet\rho_\bullet\setminus j}\times\wtd{\mc B},\label{eq66}\\
W_j''= \mc D_{r_j}\times\mc D_{\rho_j}^\times\times\mc D_{r_\bullet\rho_\bullet\setminus j}\times\wtd{\mc B}.\label{eq67}
\end{gather}
Then the relation $\sim$ identifies $W_j'$ and $W_j''$  respectively via $(\xi_j,q_j,\id)$ and $(\varpi_j,q_j,\id)$ (where $\id$ is the identity map of $\mc D_{r_\blt\rho_\blt\setminus j}\times\wtd{\mc B}$) to
\begin{gather}
\mc D_{r_j}\times\mc D_{r_\blt\rho_\blt}\times\wtd{\mc B}-F_j'\qquad (\subset U_j'\times\mc D_{r_\blt\rho_\blt}),\label{eq220}\\
\mc D_{\rho_j}\times\mc D_{r_\blt\rho_\blt}\times\wtd{\mc B}-F_j''\qquad (\subset U_j''\times\mc D_{r_\blt\rho_\blt})\label{eq221}
\end{gather}
(recall \eqref{eq18}), which are subsets of $\wtd{\mc C}\times \mc D_{r_\blt\rho_\blt}-\bigcup_j F_j'-\bigcup_j F_j''$. (In particular, certain open subsets of \eqref{eq220} and \eqref{eq221} are glued together and identified with $W_j'\cap W_j''$.)

It is easy to see that the projection
\begin{align}
\wtd\pi\times \id:\wtd{\mc C}\times \mc D_{r_\blt\rho_\blt}\rightarrow \wtd{\mc B}\times \mc D_{r_\bullet\rho_\bullet}=\mc B,\label{eq234}
\end{align}
agrees with
\begin{gather}
\pi_{r_j,\rho_j}\times \id:W_j=\mc D_{r_j}\times\mc D_{\rho_j}\times\mc D_{r_\bullet\rho_\bullet\setminus j}\times\wtd{\mc B}\rightarrow \mc D_{r_j\rho_j}\times\mc D_{r_\bullet\rho_\bullet\setminus j}\times\wtd{\mc B}=\mc B\label{eq227}
\end{gather}
for each $j$ when restricted to $W_j',W_j''$. (Indeed, recalling \eqref{eq220} and \eqref{eq221}, they are the standard projections $(\xi_j,q_\blt,\wtd b)\mapsto(q_\blt,\wtd b)$ and $(\varpi_j,q_\blt,\wtd b)\mapsto(q_\blt,\wtd b)$.) Thus, we have a well-defined surjective holomorphic map
\begin{align*}
\pi:\mc C\rightarrow\mc B
\end{align*}
whose  restrictions to $\wtd{\mc C}\times \mc D_{r_\blt\rho_\blt}-\bigcup_j F_j'-\bigcup_j F_j''$ and to each $W_j$ are $\wtd\pi\otimes\id$ and $\pi_{r_j,\rho_j}\otimes\id$ respectively. We say that $\fk X=(\pi:\mc C\rightarrow\mc B)$ is obtained from $\wtd{\fk X}$ via \textbf{sewing}.

\subsection*{The discriminant locus $\Delta$ and the critical locus $\Sigma$}

For each $j$ we set
\begin{align*}
&\Delta_j=\{0\}\times \mc D_{r_\blt\rho_\blt\setminus j}\times\wtd{\mc B}\\
& (\subset\mc D_{r_j\rho_j}\times\mc D_{r_\bullet\rho_\bullet\setminus j}\times\wtd{\mc B}=\mc B).
\end{align*}
Then \index{zz@$\Delta_j,\Delta=\bigcup\limits_{j=1}^M\Delta_j$}
\begin{align}
\Delta=\bigcup_{j=1}^M\Delta_j\label{eq132}
\end{align}
is the set of all points $b\in\mc B$ such that the fiber $\mc C_b$ is singular, called the \textbf{discriminant locus}. Roughly speaking, if $b$ is (for example) in $\Delta_1,\dots,\Delta_j$ but not in $\Delta_{j+1},\dots,\Delta_M$, then  $\mc C_b$ is obtained  by attaching $\sgm_1',\dots,\sgm_M'$ with $\sgm_1'',\dots,\sgm_M''$ respectively and smoothing the last $M-j$ nodes. Therefore $\mc C_b$ has $j$ nodes. If $b$ is outside $\Delta$ then $\mc C_b$ is smooth (i.e. a compact Riemann surface). Set also 
\begin{align*}
&\Sigma_j=\{0\}\times\{0\}\times\mc D_{r_\blt\rho_\blt\setminus j}\times\wtd{\mc B}\\
&(\subset\mc D_{r_j}\times\mc D_{\rho_j}\times\mc D_{r_\bullet\rho_\bullet\setminus j}\times\wtd{\mc B}=W_j\subset\mc C).
\end{align*}
Then \index{zz@$\Sigma_j,\Sigma=\bigsqcup_{j=1}^M\Sigma_j$}
\begin{align*}
\Sigma=\bigcup_{j=1}^M\Sigma_j
\end{align*}
is the set of nodes, called the \textbf{critical locus}. The linear map $d\pi$  is not surjective precisely at $\Sigma$. Clearly we have $\pi(\Sigma)=\Delta$.

\subsection*{Sewing families of $N$-pointed compact Riemann surfaces}

\begin{rem}\label{lb39}
Suppose, moreover, that we are sewing a family of $N$-pointed compact Riemann surfaces with local coordinates
\begin{align*}
\wtd{\fk X}=(\wtd\pi:\wtd{\mc C}\rightarrow\wtd{\mc B};\sgm_1,\dots,\sgm_N;\sgm_1',\dots,\sgm_M';\sgm_1'',\dots,\sgm_M'';\eta_1,\dots,\eta_N;\xi_1,\dots,\xi_M;\varpi_1,\dots,\varpi_M).
\end{align*}
Then, we assume that \emph{each connected component of each fiber $\wtd{\mc C}_{\wtd b}$ of $\wtd{\mc C}$ contains one of $\sgm_1(\wtd b),\dots,\sgm_N(\wtd b)$, and  $\sgm_1(\wtd{\mc B}),\dots,\sgm_N(\wtd{\mc B})$ are disjoint from the neighborhoods $U_1',\dots,U_M',U_1'',\dots,U_M''$}. One thus have a family of $N$-pointed complex curves with local coordinates
\begin{align*}
\fk X=(\pi:\mc C\rightarrow\mc B;\sgm_1,\dots,\sgm_N;\eta_1,\dots,\eta_N),
\end{align*}
where each section $\sgm_i$ is defined on $\mc B=\wtd{\mc B}\times\mc D_{r_\blt\rho_\blt}$, takes values in $\wtd {\mc C}\times\mc D_{r_\blt\rho_\blt}-\bigcup_j U_j'-\bigcup_j U_j''$, is constant over $\mc D_{r_\blt\rho_\blt}$, and equals the original one  over  $\wtd{\mc B}$. Similarly, the local coordinate $\eta_i$ of $\wtd{\fk X}$ at $\sgm_i(\wtd{\mc B})$ is extended constantly to one of $\fk X$ at $\sgm_i(\mc B)$. We say that the $N$-points $\sgm_1,\dots,\sgm_N$ and the local coordinates $\eta_1,\dots,\eta_N$ of $\fk X$ are \textbf{constant with respect to sewing}.

In the case that the local coordinates $\eta_1,\dots,\eta_N$ are not assigned to $\wtd{\fk X}$, then $\fk X=(\pi:\mc C\rightarrow\mc B;\sgm_1,\dots,\sgm_N)$ is a family of $N$-pointed complex curves.
\end{rem}

\subsection*{Families of complex curves}

We now give a general definition of families of complex curves.   Suppose that $\fk X=(\pi:\mc C\rightarrow \mc B)$ is a surjective holomorphic map of  complex manifolds, where $\mc B$ has finitely many connected components. We say that $\fk X$ is a \textbf{family of complex curves} if  $\fk X$  is  either smooth or obtained via sewing.  
We let $\Sigma$ be the set nodes, i.e., the set of all $x\in\mc C$ such that $\mc C_{\pi(x)}$ is nodal and $x$ is a node of the nodal curve $\mc C_{\pi(x)}$. \index{zz@$\Sigma_j,\Sigma=\bigsqcup_{j=1}^M\Sigma_j$} Equivalently, $\Sigma$ is the set of all $x\in\mc C$ such that $d\pi$ is not surjective at $x$. We \index{zz@$\Delta_j,\Delta=\bigcup\limits_{j=1}^M\Delta_j$} \index{CDelta@$\mc C_\Delta=\pi^{-1}(\Delta)$}  define
\begin{align*}
\Delta=\pi(\Sigma),\qquad \mc C_\Delta=\pi^{-1}(\Delta).
\end{align*}
Then $\mc C_\Delta$ is the union of all singular fibers. $\Delta$ and $\mc C_\Delta$ will be considered later as (normal crossing) divisors. 

\begin{thm}\label{lb11}
Grauert's Theorem \ref{lb2} holds verbatim for a family of complex curves $\pi:\mc C\rightarrow\mc B$.
\end{thm}

As mentioned in  the paragraphs after Theorem \ref{lb2}, Grauert's theorem holds in general when $\mc C$ and $\mc B$ are complex spaces, $\scr E$ is a coherent $\scr O_{\mc C}$-module, $\pi$ is proper, and $\scr E$ is $\pi$-flat. To apply that theorem to  the family of complex curves $\pi:\mc C\rightarrow\mc B$ and a locally free $\scr O_{\mc C}$-module $\scr E$, we need to check  that $\pi$ is proper and $\scr E$ is $\pi$-flat. One can check the properness of $\pi$ by checking, for example, the sequential compactness of the preimages of  compact subsets. Also, one checks easily that $\pi$ is an open map. Thus $\scr O_{\mc C}$ is $\pi$-flat by  \cite[Sec. 3.20]{Fis76} (see also \cite[Thm. II.2.13]{GPR94} and \cite[Thm. V.2.13]{BS76}). (One simply says that $\pi$ is flat.) Since $\scr E$ is locally free, $\scr E$ is $\pi$-flat.

\begin{rem}\label{lb124}
Apply Grauert's theorem to $\scr O_{\mc C}$, we see that $b\mapsto\chi(\mc C_b,\scr O_{\mc C_b})$ is constant on each connected component of $\mc B$. If we assume that all the fibers of $\mc C$ are  connected, then the genus of the fiber $g(\mc C_b)=1-\chi(\mc C,\scr O_{\mc C_b})$ is locally constant over $b$. We conclude that the genus of a complex curve is unchanged under deformation. Consequently, the genus of a nodal curve equals the genus of its ``smoothing". 
\end{rem}

Consider a family $\pi:\mc C\rightarrow\mc B$ of complex curves. Let $\sgm_1,\dots,\sgm_N:\mc B\rightarrow\mc C$ be (holomorphic) sections whose images are mutually disjoint and are also disjoint from $\Sigma$. For each $b\in\mc B$, we assume that $(\mc C_b;\sgm_1(b),\dots,\sgm_N(b))$ is an $N$-pointed complex curve. Thus every connected component of $\wtd{\mc C}_b$ (the normalization  of $\mc C_b$) contains one of $\sgm_1(b),\dots,\sgm_N(b)$. Equivalently, every irreducible component of $\mc C_b$ contains one of $\sgm_1(b),\dots,\sgm_N(b)$. Notice that  $\sgm_1(b),\dots,\sgm_N(b)$ are not nodes since they are not in $\Sigma$. Then we say that $(\pi:\mc C\rightarrow\mc B;\sgm_1,\dots,\sgm_N)$ is a \textbf{family of $N$-pointed complex curves}. Assume moreover that for each $1\leq j\leq N$, the family $\pi:\mc C\rightarrow\mc B$ has a \textbf{local coordinate} $\eta_j$ at $\sgm_j$. This means that $\sgm_j(\mc B)$ is contained in an open set $U_j\subset\mc C-\Sigma$ such that $\eta_j:U_j\rightarrow\Cbb$ is holomorphic and equals $0$ on $\eta_j(\mc B)$,  and that $\eta_j|_{U_j\cap\mc C_b}$ is a univalent function  for each $b\in\mc B$. So $(\eta_j,\pi)$ is a biholomorphic map from $U_j$ to an open subset of $\mbb C\times\mc B$ containing $\{0\}\times\mc B$. Then we say that 
\begin{align*}
\fk X=(\pi:\mc C\rightarrow\mc B;\sgm_1,\dots,\sgm_N;\eta_1,\dots,\eta_N)
\end{align*}
is a \textbf{family of $N$-pointed complex curves with local coordinates}. As usual, we define a divisor
\begin{align*}
\SX=\sum_{j=1}^N\sgm_j(\mc B).
\end{align*}
Note that even if $\pi:\mc C\rightarrow\mc B$ is obtained by sewing a smooth family $\wtd{\pi}:\wtd{\mc C}\rightarrow\wtd{\mc B}$, the $N$-points and their local coordinates are not assumed to be constant with respect to sewing.

\begin{pp}
Proposition \ref{lb129} holds verbatim if $\fk X$ is a family of $N$-pointed complex curves.
\end{pp}

\begin{proof}
The proof is the same as in Proposition \ref{lb129}.
\end{proof}

\subsection*{The invertible sheaves $\Theta_{\mc C/\mc B}$ and $\omega_{\mc C/\mc B}$}

Given a family of complex curves $\pi:\mc C\rightarrow\mc B$, we can define a homomorphism of $\scr O_{\mc C}$-modules $d\pi:\Theta_{\mc C}\rightarrow\pi^*\Theta_{\mc B}$ in a similar way as  for families of compact Riemann surfaces. However, this map is no necessarily a surjective sheaf map, which means that we do not have a long exact sequence from this map. To remedy this issue we consider $d\pi:\Theta_{\mc C}(-\log \mc C_\Delta)\rightarrow \pi^*\big(\Theta_{\mc B}(-\log\Delta)\big)$ defined as follows. (We shall write $\pi^*\big(\Theta_{\mc B}(-\log\Delta)\big)$ as $\pi^*\Theta_{\mc B}(-\log\Delta)$ for short.)

Let us first describe $\Theta_{\mc B}(-\log\Delta)$. \index{zz@$\Theta_{\mc B}(-\log\Delta)$}  Choose any $b\in\mc B$. Then one can always find a  neighborhood $V$ of $b$ such that
\begin{align*}
V\simeq \mc D_{r_\blt\rho_\blt}\times\wtd{\mc B}=\mc D_{r_1\rho_1}\times\cdots\times\mc D_{r_M\rho_M}\times\wtd{\mc B}
\end{align*}
where $\wtd {\mc B}$ is an open subset of $\mathbb C^n$,   and
\begin{align*}
\Delta\cap V=\{q_{\blt,0}\in\mc D_{r_\blt\rho_\blt}:q_{1,0}\cdots q_{M,0}=0\}\times\wtd{\mc B}.
\end{align*}
Let $q_1,\dots,q_M$ be the standard coordinates of $\mc D_{r_1\rho_1},\dots,\mc D_{r_M\rho_M}$. Let also $\tau_\blt=(\tau_1,\cdots,\tau_n)$ be the standard coordinate of $\wtd {\mc B}$ inside $\mathbb C^n$. Then $\Theta_{\mc B}(-\log\Delta)_V$ is defined to be the $\scr O_V$-submodule of $\Theta_V=\Theta_{\mc B}|_V$ generated (freely) by
\begin{align}
q_1\partial_{q_1},\dots,q_M\partial_{q_M},\partial_{\tau_1},\dots\partial_{\tau_n}.\label{eq133}
\end{align}
By gluing, we obtain a locally free $\scr O_{\mc B}$-module $\Theta_{\mc B}(-\log\Delta)$.

Next, we describe $\Theta_{\mc C}(-\log \mc C_\Delta)$ and the map $d\pi:\Theta_{\mc C}(-\log \mc C_\Delta)\rightarrow\pi^*\Theta_{\mc B}(-\log\Delta)$. Choose any $x\in\mc C$. 

Case I. $x\notin\Sigma$. Then one can find a neighborhood $U$ of $x$ disjoint from $\Sigma$, and a neighborhood $V$ of $b=\pi(x)$ described as above, such that
\begin{gather*}
U\simeq \mc D_{r_0}\times V=\mc D_{r_0}\times\mc D_{r_\blt\rho_\blt}\times\wtd{\mc B},
\end{gather*}
and that $\pi: \mc D_{r_0}\times V\rightarrow V$ is the projection on to the $V$-component. Let $z$ be the standard coordinates of $\mc D_{r_0}$. Notice that
\begin{gather*}
\mc C_\Delta\cap U=\mc D_{r_0}\times\{q_{\blt,0}\in\mc D_{r_\blt\rho_\blt}:q_{1,0}\cdots q_{M,0}=0\}\times\wtd{\mc B},
\end{gather*}
which suggests that we define $\Theta_{\mc C}(-\log \mc C_\Delta)_U$ to be the $\scr O_U$-submodule of $\Theta_U=\Theta_{\mc C}|_U$ generated (freely) by
\begin{align}
\partial_z,q_1\partial_{q_1},\dots,q_M\partial_{q_M},\partial_{\tau_1},\dots,\partial_{\tau_n}.\label{eq21}
\end{align}
The homomorphism
\begin{align}
d\pi:\Theta_{\mc C}(-\log \mc C_\Delta)_U\rightarrow \big(\pi^*\Theta_{\mc B}(-\log\Delta)\big)_U\label{eq23}
\end{align}
is defined by sending $\partial_z$ to $0$ and keeping all the other elements in \eqref{eq21}. (Here we do not differentiate between $\partial_{q_\blt},\partial_{\tau_\blt}$ and their pull backs.) It is clear that $d\pi$ is surjective. We leave it to the reader to check that the above definitions are independent of the choice of local coordinates.

Case II. $x\in\Sigma$. Then one can find  neighborhoods $U$ of $x$ and $V$ of $b=\pi(x)$ described as above, such that
\begin{gather*}
U\simeq \mc D_{r_1}\times \mc D_{\rho_1}\times \mc D_{r_\blt\rho_\blt\setminus 1}\times\wtd {\mc B},\\
\pi=\pi_{r_1,\rho_1}\otimes\id\otimes\id:\mc D_{r_1}\times \mc D_{\rho_1}\times \mc D_{r_\blt\rho_\blt\setminus 1}\times\wtd {\mc B}\rightarrow\mc D_{r_\blt\rho_\blt}\times\wtd{\mc B}=V.
\end{gather*}
From this we deduce that
\begin{align*}
&\mc C_\Delta\cap U\\
=&\{(\xi_{1,0},\varpi_{1,0},q_{2,0},\dots,q_{M,0})\in \mc D_{r_1}\times \mc D_{\rho_1}\times \mc D_{r_\blt\rho_\blt\setminus 1}:\xi_{1,0}\cdot\varpi_{1,0}\cdot q_{2,0}\cdots q_{M,0}=0\}\times\wtd{\mc B}.
\end{align*}
Let $\xi_1,\varpi_1$ be respectively the standard coordinates of $\mc D_{r_1},\mc D_{\rho_1}$. Again $q_\blt$ and $\tau_\blt$ are respectively the standard coordinates of $\mc D_{r_\blt\rho_\blt}$ and $\wtd {\mc B}$. Note also the relation $q_1=\xi_1\varpi_1$. The description of $\mc C_\Delta\cap U$ suggests that we define $\Theta_{\mc C}(-\log \mc C_\Delta)_U$ to be the $\scr O_U$-submodule of $\Theta_U=\Theta_{\mc C}|_U$ generated (freely) by
\begin{align}
\xi_1\partial_{\xi_1},\varpi_1\partial_{\varpi_1},q_2\partial_{q_2},\dots,q_M\partial_{q_M},\partial_{\tau_1},\dots,\partial_{\tau_n}.\label{eq22}
\end{align}
The homomorphism $d\pi$ of \eqref{eq23} is defined by setting
\begin{align}
d\pi(\xi_1\partial_{\xi_1})=d\pi(\varpi_1\partial_{\varpi_1})=q_1\partial_{q_1}\label{eq182}
\end{align}
and keeping all the other elements in \eqref{eq22}. Again $d\pi$ is surjective and independent of the choice of local coordinates. Also, using \eqref{eq74}, it is easy to see that the $d\pi$ constructed in case I and II are compatible.

We may now glue the two cases together and obtain the locally free $\scr O_{\mc C}$-module $\Theta_{\mc C}(-\log \mc C_\Delta)$  \index{zz@$\Theta_{\mc C}(-\log \mc C_\Delta)$} and a surjective  homomorphism of $\scr O_{\mc C}$-modules $d\pi:\Theta_{\mc C}(-\log \mc C_\Delta)\rightarrow\pi^*\Theta_{\mc B}(-\log\Delta)$, whose kernel is denoted by $\Theta_{\mc C/\mc B}$. \index{zz@$\Theta_{\mc C/\mc B}$} Thus, there is a short exact sequence
\begin{align}
\boxed{~~0\rightarrow \Theta_{\mc C/\mc B}\rightarrow \Theta_{\mc C}(-\log \mc C_\Delta)\xrightarrow{d\pi}\pi^*\Theta_{\mc B}(-\log \Delta)\rightarrow 0~~}\label{eq138}
\end{align}
Choose any $x\in\mc C$ and a small neighborhood $U$ of $\Theta_{\mc C/\mc B}$. Then in case I resp. in case II $\Theta_{\mc C/\mc B}|_U$ is generated freely by
\begin{align}
\partial_z\qquad\text{resp.}\qquad \xi_1\partial_{\xi_1}-\varpi_1\partial_{\varpi_1}.
\end{align}
Therefore, $\Theta_{\mc C/\mc B}$ is an invertible $\scr O_{\mc C}$-module, whose dual sheaf is denoted by $\omega_{\mc C/\mc B}$.\index{zz@$\omega_{\mc C/\mc B}$}    We leave it to the reader to check that there are natural equivalences
\begin{align*}
\Theta_{\mc C/\mc B}|\mc C_b\simeq \Theta_{\mc C_b},\qquad \omega_{\mc C/\mc B}|\mc C_b\simeq \omega_{\mc C_b}
\end{align*}
for any $b\in\mc B$. One may use the following fact:  in case II, by \eqref{eq74}, in the $(\xi_1,q_\blt,\tau_\blt)$- resp. $(\varpi_1,q_\blt,\tau_\blt)$-coordinate, the section $\xi_1\partial_{\xi_1}-\varpi_1\partial_{\varpi_1}$ equals
\begin{align}
\xi_1\partial_{\xi_1}\qquad\text{resp.}\qquad -\varpi_1\partial_{\varpi_1}.\label{eq80}
\end{align}

We close this section with the following generalization of Theorem \ref{lb8}.

\begin{thm}\label{lb9}
Let $\fk X=(\pi:\mc C\rightarrow\mc B;\sgm_1,\dots,\sgm_N)$ be a family of $N$-pointed complex curves. Let $n\in\Nbb$. Then  there exists $k_0\in\mbb N$ such that for any $k\geq k_0$,  the $\scr O_{\mc B}$-module  $\pi_*\Theta_{\mc C/\mc B}^{\otimes n}(k\SX)$ is locally free, and for any $b\in\mc B$ there is a natural isomorphism of vector spaces
\begin{align}
\frac{\pi_*\Theta_{\mc C/\mc B}^{\otimes n}(k\SX)_b}{~\fk m_b\cdot\pi_*\Theta_{\mc C/\mc B}^{\otimes n}(k\SX)_b~}\simeq H^0\big(\mc C_b,\Theta_{\mc C_b}^{\otimes n}(k\SX(b))\big)
\end{align} 
defined by  restriction of sections. In particular,  $\dim H^0\big(\mc C_b,\Theta_{\mc C_b}^{\otimes n}(k\SX(b))\big)$ is locally constant over $b$.
\end{thm}

\begin{proof}
This follows easily from theorems \ref{lb7} and \ref{lb10},  and Grauert's Theorem \ref{lb11}.
\end{proof}

\section{Linear differential equations with simple poles}

We first introduce the following notations. If $W$ is a vector space and $z$ is a (formal) variable, we \index{z@$[[z]],[[z^{\pm 1}]],((z)),\{z\}$} define 
\begin{gather*}
W[[z]]=\bigg\{\sum_{n\in\mathbb N}w_nz^n:\text{each }w_n\in W\bigg\},\\
W[[z^{\pm 1}]]=\bigg\{\sum_{n\in\mathbb Z}w_nz^n:\text{each }w_n\in W\bigg\},\\
W((z))=\Big\{f(z):z^kf(z)\in W[[z]]\text{ for some }k\in\mbb Z \Big\},\\
W\{z\}=\Big\{\sum_{n\in\mbb C}w_nz^n :\text{each $w_n\in W$}\Big\}.
\end{gather*}

In this section, we fix $m,N\in\mbb N,r>0$, and let $V$ be an open subset of $\mbb C^m$. Let $q$ be the standard coordinate of $\mc D_r\subset\mbb C$, and let $\tau_\blt=(\tau_1,\dots,\tau_m)$ be the standard coordinates of $V$. $A$ is an $\End(\Cbb^N)$-valued holomorphic function on $\mc D_r\times V$, i.e.
\begin{align*}
A\in\End(\Cbb^N)\otimes_\Cbb\scr O(\mc D_r\times V).
\end{align*}
First, recall the following well-known fact.

\begin{thm}\label{lb12}
For any  $\varphi\in\Cbb^N\otimes_\Cbb\scr O(V)$, there exists a unique  $\psi\in\Cbb^N\otimes_\Cbb\scr O(\mc D_r\times V)$ satisfying the differential equation $\partial_q\psi(q,\tau_\blt)=A(q,\tau_\blt)\psi(q,\tau_\blt)$ and the initial condition $\psi(0,\tau_\blt)=\varphi(\tau_\blt)$. 
\end{thm}

\begin{proof}
This is an easy consequence of Picard iteration. We provide the details for the readers' convenience. It suffices to prove the existence and the uniqueness of such $\psi$ on any precompact open subset of $\mc D_r\times V$. Thus, we may well assume that $\lVert A(q,\tau_\blt)\lVert$ (where $\lVert\cdot\lVert$ is the operator norm) is uniformly bounded by a positive (finite) number $C$, and $\lVert \varphi(\tau_\blt)\lVert$ is also uniformly bounded. Fix any $\delta\in(0,1)$. We claim that for any $q_0\in \mc D_r$, if we let $\mc D_{\delta/C}(q_0)$ be the open disc with center $q_0$ and radius $\delta/C$ and set $W(q_0)=\mc D_{\delta/C}(q_0)\cap\mc D_r$, then for any uniformly bounded $\varphi\in\Cbb^N\otimes_\Cbb\scr O(V)$, there exists a unique $\psi\in\Cbb^N\otimes_\Cbb\scr O(W(q_0)\times V)$ satisfying $\partial_q\psi=A\psi$ and $\psi(q_0,\tau_\blt)=\varphi(\tau_\blt)$. Then the theorem is proved by covering $\mc D_r$ by finitely many open discs with radius $\delta/C$.

Fix $q_0\in\mc D_r$ and a uniformly bounded $\varphi\in\Cbb^N\otimes_\Cbb\scr O(V)$. Let $\fk H$ be the Banach space of uniformly bounded elements of $\Cbb^N\otimes_\Cbb\scr O(W(q_0)\times V)$, whose norm is given by $\lVert \cdot\lVert_\infty$. Define a map $T:\fk H\rightarrow\fk H$ such that for any $\psi\in\fk H$,
\begin{align*}
T\psi(q,\tau_\blt)=\varphi(\tau_\blt)+\int_{\gamma_q} A(z,\tau_\blt)\psi(z,\tau_\blt)dz.
\end{align*}
where $\gamma_q$ is any path in $W(q_0)$ from $q_0$ to $q$. Then  $\psi$ is a solution satisfying $\psi(q_0,\tau_\blt)=\varphi(\tau_\blt)$ if and only if $T\psi=\psi$. It is easy to see that $\lVert T(\psi_1-\psi_2)\lVert\leq\delta \lVert \psi_1-\psi_2\lVert$. So $T$ is a contraction, which proves the existence and the uniqueness of the fixed point of $T$.
\end{proof}

The above theorem clearly holds also when $\mc D_r$ is replaced by a simply connected domain.

Consider the following  differential equation with simple pole
\begin{align}
q\partial_q\psi=A\psi\label{eq27}
\end{align}
where $\psi\in\Cbb^N\otimes_\Cbb\scr O(V)((q))$ is a formal solution of this equation. By our assumption on $\psi$, we can write 
\begin{align}
\psi(q,\tau_\blt)=\sum_{n\in\mbb Z}\wht\psi_n(\tau_\blt)q^n,\label{eq24}
\end{align}
where $\wht\psi_n\in\Cbb^N\otimes_\Cbb\scr O(V)$, and  $\wht\psi_n=0$ when $n$ smaller than some  negative integer.

\begin{thm}\label{lb14}
Suppose that the formal series $\psi$ is a formal solution of \eqref{eq27}. Then $\psi$ is an element of $\Cbb^N\otimes_\Cbb\scr O(\mc D_r^\times\times V)$.
\end{thm}

\begin{proof}
Suppose the mode $\wht\psi_n$ is zero when $n<-K$. Then $\phi:=q^K\psi$ has no negative modes and satisfies a similar differential equation $q\partial_q\phi=(K+A)\phi$. Thus,  we may well assume that $\wht\psi_n=0$ when $n<0$. Consider the series expansion of $A$:
\begin{align*}
A(q,\tau_\blt)=\sum_{n\in\mbb N}\wht A_n(\tau_\blt)q^n,
\end{align*}
where each $\wht A_n$ is  in $\End(\Cbb^N)\otimes_\Cbb\scr O(V)$.  Then for each $n\in\mbb N$,
\begin{align*}
n\wht\psi_n=\sum_{j=0}^n\wht A_{n-j}\wht\psi_j.
\end{align*}
Choose any open subset $U$ of $V$ with compact closure, and choose $M\in\mbb N$ such that $\lVert\wht A_0(\tau_\bullet)\lVert\leq M$ whenever $\tau_\blt\in U$. (Here $\lVert\cdot\lVert$ is the operator norm.) Then for any $n>M$, $n\id-\wht A_0(\tau_\blt)$ is invertible (with inverse $n^{-1}\sum_{j=0}^\infty(\wht A_0(\tau_\blt)/n)^j $). Thus, whenever $n>M$,
\begin{align}
\wht\psi_n=(n-\wht A_0)^{-1}\sum_{j=0}^{n-1}\wht A_{n-j}\wht\psi_j.\label{eq25}
\end{align}
Choose any $r_1<r$ and set
\begin{align}
\alpha=\sup_{(q,\tau_\blt)\in {\ovl{\mc D_{r_1}}\times U}}\lVert A(q,\tau_\blt)\lVert.\label{eq28}
\end{align}
Using the fact that $\wht A_n(\tau_\blt)=\oint_{\partial \mc D_{r_1}} A(q,\tau_\blt)q^{-n-1}\frac{dq}{2\im\pi}$, we have
\begin{align}
\lVert\wht A_n(\tau_\blt)\lVert \leq \alpha r_1^{-n}\label{eq222}
\end{align}
for all $n$ and all  $\tau_\blt$ in $U$.

Choose $\beta>0$ such that  $\lVert (n-\wht A_0(\tau_\blt))^{-1} \lVert\leq \beta n^{-1}$ for any $n> M$ and $\tau_\blt\in U$. (Such $\beta$ can be found using the explicit formula of inverse matrix given above.) Set $\gamma\geq\max\{1,\alpha\beta\}$. Then, from \eqref{eq25} and \eqref{eq222}, we see that for any $n>M$ and $\tau_\blt\in U$,
\begin{align}
r_1^n\lVert \wht\psi_n(\tau_\blt) \lVert\leq \gamma n^{-1}\sum_{j=0}^{n-1}r_1^j\lVert \wht\psi_j(\tau_\blt) \lVert.\label{eq26}
\end{align}
By induction, one can show that there exists $c>0$ such that
\begin{align*}
r_1^n\lVert\wht\psi_n(\tau_\blt)\lVert\leq c\gamma^n 
\end{align*}
for any $n\in\mbb N$ and $\tau_\blt\in U$. Indeed, if this is true for $0,1,2,\dots,n-1$ where $n>N$, then by \eqref{eq26},
\begin{align*}
r_1^n\lVert \wht\psi_n(\tau_\blt) \lVert\leq \gamma n^{-1}\sum_{j=0}^{n-1}c\gamma^j\leq \gamma n^{-1}\sum_{j=0}^{n-1}c\gamma^{n-1}=c\gamma^n.
\end{align*}
Thus $\lVert\wht\psi_n(\tau_\blt)\lVert\leq c\gamma^nr_1^{-n}$ for all $n$ and $\tau_\blt\in U$. Therefore, if we choose any $r_0\in(0,\gamma^{-1}r_1)$, then   the series $\sum_n\lVert \wht \psi_n(\tau_\blt) \lVert \cdot|q|^n$ is uniformly bounded by some positive number for all  $|q|\leq r_0$ and $\tau_\blt\in U$. Since each $\wht\psi_n(\tau_\blt)$ is holomorphic over $\tau_\blt$, the series \eqref{eq24} must converge uniformly to a holomorphic function on $\mc D_{r_0}\times U$. 

By Theorem \ref{lb12}, for any simply connected open subset $W\subset \mc D_r^\times$ which overlaps with $\mc D_{r_0}$, there exists a unique  holomorphic solution of \eqref{eq27} on $W\times U$ which agrees with $\psi$ on $(W\cap\mc D_{r_0})\times U$. Thus, $\psi$ can be extended to a holomorphic $\Cbb^N$-valued function on $\mc D_r\times U$. Since $U$ is an arbitrary precompact open subset of $V$, $\psi$ is holomorphic on $\mc D_r\times V$.
\end{proof}

\begin{rem}\label{lb13}
The proof of the above proposition shows  that if $M\in\mbb N$ and $\lVert\wht A_0(\tau_\blt)\lVert\leq M$ whenever $\tau_\blt\in V$, then $\psi$ is holomorphic  on $\mc D_r^\times \times V$ provided that the first $(M+1)$ non-zero modes of $\psi$ are holomorphic on $\tau_\blt\in V$. 
\end{rem}

\section{Criteria on local freeness}

Let $X$ be a complex manifold and $\scr E$ be an $\scr O_X$-module. We say that $\scr E$ is a \textbf{finite-type $\scr O_X$-module} if any $x\in X$ is contained in an open subset $U$ such that $\scr E_U$ is  generated by some $s_1,\dots,s_n\in\scr E(U)$. Note that this does not mean that the $\scr O(U)$-module $\scr E(U)$ is generated by $s_1,\dots,s_n$.

Recall that for each open subset $U$ of $X$, and for each  $x\in U$ and $s\in\scr E(U)$, $s(x)$ denotes the value of $s$ in $\scr E|x\simeq \scr E_x/\fk m_x\scr E_x$. It is clear that if $s_1,\dots,s_n\in\scr E(U)$ generate $\scr E_U$, then $s_1(x),\dots,s_n(x)$ span the vector space $\scr E|x$. Conversely, we have:

\begin{pp}[Nakayama's lemma]\label{lb60}
Suppose that $\scr E$ is a finite-type $\scr O_X$-module, $x\in X$, $U$ is an open set containing $x$, $s_1,\dots,s_n\in\scr E(U)$, and $s_1(x),\dots,s_n(x)$ span the vector space $\scr E|x$. Then there exists an open subset $V$ of $U$ containing $x$ such that $s_1|_V,\dots,s_n|_V$ generate $\scr E_V$. 
\end{pp}

Consequently, the rank function of a finite-type $\scr O_X$-module is upper semi-continuous.

\begin{proof}
Since $\scr E$ is  finite-type, we may assume that $U$ is small enough such that $\scr E_U$ is generated  by $\sigma_1,\dots,\sigma_m\in\scr E(U)$. Thus the germs $\sigma_{1,x},\dots,\sigma_{m,x}$ generate the $\scr O_{X,x}$-module $\scr E_x$. In particular, $\sigma_{1,x},\sigma_{2,x}\dots,\sigma_{m,x},s_{1,x},\dots,s_{n,x}$ generate $\scr E_x$. We now show that $\sigma_{2,x},\sigma_{3,x}\dots,\sigma_{m,x},s_{1,x},\dots,s_{n,x}$ generate $\scr E_x$. Since $s_1(x),\dots,s_n(x)$ span $\scr E_x/\fk m_x\scr E_x$, there exist complex numbers $c_1,\dots,c_n$ such that $\sigma_{1,x}\in\sum_{k=1}^n c_ks_{k,x}+\fk m_x\scr E_x$. Therefore, there exists $f_1,\dots,f_m,g_1,\dots,g_n\in\fk m_x$ such that
\begin{align*}
\sigma_{1,x}=\sum_{k=1}^n c_ks_{k,x}+\sum_{i=1}^m f_i\sigma_{i,x}+\sum_{j=1}^n g_j s_{j,x}.
\end{align*}
Since $f_1(x)=0$, the element $1-f_1$ has inverse  in $\scr O_{X,x}$. For each $j=1,\dots,n$, set $\wtd g_j=c_j+g_j$ which is an element in $\scr O_{X,x}$. Then we clearly have
\begin{align*}
\sigma_{1,x}=(1-f_1)^{-1}\bigg(\sum_{i=2}^m f_i\sigma_{i,x}+\sum_{j=1}^n \wtd g_j s_{j,x}  \bigg).
\end{align*}
This shows that $\sigma_{2,x},\sigma_{3,x}\dots,\sigma_{m,x},s_{1,x},\dots,s_{n,x}$ generate $\scr E_x$. A similar argument shows that $\sigma_{3,x},\sigma_{4,x}\dots,\sigma_{m,x},s_{1,x},\dots,s_{n,x}$ generate $\scr E_x$. If we repeat this argument several times, we arrive at the conclusion that $s_{1,x},\dots,s_{n,x}$ generate $\scr E_x$. Thus, there exists an open subset $V$ of $U$ containing $x$ such that for each $i=1,\dots m$, there exist $h_{i,1},\dots,h_{i,n}\in\scr O(V)$ satisfying $\sigma_i|_V=\sum_{j=1}^n h_{i,j} s_j|_V$. Therefore, as $\sigma_1|_V,\dots,\sigma_m|_V$ generate  $\scr E_V$, so do $s_1|_V,\dots,s_n|_V$.
\end{proof}

\begin{thm}\label{lb59}
Suppose that $\scr E$ is finite-type. Then $\scr E$ is locally free if and only if the rank function $x\in X\mapsto r_x=\dim_{\mbb C}(\scr E|x)$ is locally constant.
\end{thm}

\begin{proof}
The only if part is obvious. Let us prove the if part. Suppose that the rank function is locally constant. Choose any $x\in X$ and let $r=r_x$. There exists a neighborhood $U$ of $x$ and $s_1,\dots,s_n\in\scr E(U)$ generating $\scr E_U$. We may also assume that $U$ is small enough such that for any $y\in U$, $r_y=r$. Since $s_1(x),\dots,s_n(x)$ span $\scr E|x$, we must have $n\geq r$. By rearranging $s_1,\dots,s_n$, we may assume that $s_1(x),\dots,s_r(x)$ form a basis of the vector space $\scr E_x$. Thus, by Nakayama's lemma, there exists a neighborhood $V$ of $x$ contained in $U$ such that $s_1|_V,\dots,s_r|_V$ generate $\scr E|_V$. We prove that $\scr E|_V$ is a free $\scr O_V$-module with basis $s_1|_V,\dots,s_r|_V$. Choose any open subset $W\subset V$. We need to show that $s_1|_W,\dots,s_r|_W$ are $\scr O(W)$-linear independent. Suppose that $f_1,\dots,f_r\in\scr O(W)$ and $f_1s_1|_W+\cdots+f_rs_r|_W=0$. Then for any $y\in W$, $f_1(y)s_1(y)+\cdots+f_r(y)s_r(y)=0$. Clearly $s_1(y),\dots,s_r(y)$ span $\scr E|y$. Since $\scr E|y$ has dimension $r$, $s_1(y),\dots,s_r(y)$ are linearly independent. So $f_1(y)=\cdots=f_r(y)=0$. Since this is true for any $y\in W$, we conclude that $f_1=\cdots=f_r=0$. 
\end{proof}

\begin{df}\label{lb49}
A \textbf{connection} $\nabla$ on $\scr E$ associates to each open subset $U\subset X$ a bilinear map
\begin{gather*}
\nabla:\Theta_X(U)\times\scr E(U)\rightarrow\scr E(U),\qquad (\fk x,s)\mapsto \nabla_{\fk x} s
\end{gather*}
satisfying the following conditions.

\begin{enumerate}[label=(\alph*)]
\item If $V$ is an open subset of $U$, then $\nabla_{\fk x|_V}s|_V=(\nabla_{\fk x} s)|_V$.
\item If $f\in\scr O
(U)$, then
\begin{gather*}
\nabla_{f\fk x}s=f\nabla_{\fk x}s,\\
\nabla_{\fk x}(fs)=\fk x(f)s+f\nabla_{\fk x}s.
\end{gather*}
\end{enumerate}
\end{df}

\begin{lm}\label{lb61}
Let $\scr E$ be a finite-type $\scr O_X$-module equipped with a connection $\nabla$. Let $x\in X$ and $U\subset X$ a connected neighborhood of $x$.  Choose $s_1,\dots,s_n\in\scr E(U)$ and assume that $s_1(x),\dots,s_n(x)$ form a basis of the fiber $\scr E|x$. Then $s_1,\dots,s_n$ are $\scr O(U)$-linear independent elements of $\scr E(U)$, i.e., if $f_1,\dots,f_n\in\scr O(U)$ satisfy $f_1s_1+\cdots+f_ns_n=0$, then $f_1=\cdots=f_n=0$.
\end{lm}

\begin{proof}
Assume $U$ is open in $\Cbb^m$ and $x=0\in\Cbb^m$. Since $f_1(0)s_1(0)+\cdots +f_n(0)s_n(0)=0$ and $s_1(0),\dots,s_n(0)$ form a basis of $\scr E|0$, we obtain $f_1(0)=\cdots=f_n(0)=0$.	

Apply $\nabla_{\partial_1}$ to $\sum_j f_js_j=0$ and take value at $0$, we get $\sum_j\partial_1f_j(0)s_j(0)+\sum_j f_j(0)(\nabla_{\partial_1}s_j)(0)=0$, which shows $\partial_1 f_j(0)=0$ for all $j$. Similarly, apply $\nabla_{\partial_1},\dots,\nabla_{\partial_m}$ successively to $\sum_j f_js_j=0$ and take values at $0$. Then an induction on  $|k_\blt|=k_1\cdots+k_m$ shows that $\partial_1^{k_1}\cdots \partial_m^{k_m}f_j(0)=0$ for all $1\leq j\leq n$ and $k_1,\dots,k_n\in\Nbb$. This proves $f_1=\cdots=f_n=0$ on $U$ because $U$ is connected.
\end{proof}



\begin{thm}\label{lb66}
Let $\scr E$ be an $\scr O_X$-module equipped with a connection $\nabla$. Assume that  each $x\in X$ is contained in a neighborhood $U$ such that the following conditions hold. Then $\scr E$ is locally free.
\begin{enumerate}[label=(\alph*)]
\item $\scr E|_U$ is equivalent to the cokernel of a homomorphism of (possibly infinite rank) free $\scr O_U$-modules, i.e., there exist cardinalities $I,J$ and a homomorphism of $\scr O_U$-modules $\varphi:\scr O_U^I\rightarrow\scr O_U^J$  such that $\scr E|_U\simeq\mathrm{coker}(\varphi)$.
\item Write $\scr F=\scr O^I_U$ and $\scr G=\scr O^J_U$. Then $\scr G(U)/\varphi(\scr F)(U)$ is a finitely generated $\scr O(U)$-module. 
\end{enumerate}
\end{thm}

Note that  $\varphi(\scr F)$ is the image sheaf of $\varphi$, which is defined by sheafifying the collection $\{\varphi(\scr F(V))\}$ over all open subsets $V\subset U$. (So $\varphi(\scr F(V))$ is only a subset of $\varphi(\scr F)(V)$.) $\scr G(U)/\varphi(\scr F)(U)$ is naturally an $\scr O(U)$-submodule of $\coker(\varphi)(U)$.

Also, if $s_1,\dots s_n$ are generating elements of $\scr G(U)/\varphi(\scr F)(U)$, then for any open subset $V\subset U$, $s_1|_V,\dots,s_n|_V$ generate $\scr G(V)/\varphi(\scr F)(V)$. This is due to the obvious fact that sections of $\scr G(U)$, when restricted to $V$, generate the $\scr O(V)$-module $\scr G(V)$. We will use this property in the following proof.

\begin{proof}
Since (the stalks of) $\coker(\varphi)$ are generated by elements of $\scr G(U)$, it is clear that $\scr E$ is finite-type. Choose any $x\in X$. We shall show that the rank function $y\in X\mapsto r_y=\dim(\scr E|y)$ is  constant in a neighborhood $U$ of $x$. Then $\scr E|_U$ is locally free by Theorem \ref{lb59}, which will finish the proof.

Choose a connected $U$ as in the assumption of this theorem, and identify $\scr E|_U$ with the cokernel of $\varphi$.   We may assume that $U$ is small enough such that there exist $\sigma_1,\dots,\sigma_r\in\scr E(U)$ such that $\sigma_1(x),\dots,\sigma_r(x)$ form a basis of the fiber $\scr E|x$. We may shrink $U$ such that one can find $\wht\sigma_1,\dots,\wht\sigma_r\in\scr G(U)$ whose corresponding equivalence classes in $\scr G(U)/\varphi(\scr F)(U)$ (and hence in $\scr E(U)$) are $\sigma_1,\dots,\sigma_r$. Thus $\sigma_1,\dots,\sigma_n$ are in $\scr G(U)/\varphi(\scr F)(U)$. Suppose that $s_1,\dots,s_n$ are generating elements of $\scr G(U)/\varphi(\scr F)(U)$. By Proposition \ref{lb60} (Nakayama's lemma), we may shrink $U$ so that $s_1,\dots,s_n$ are $\scr O(U)$-linear combinations of $\sigma_1,\dots,\sigma_r$. Therefore, the $\scr O(U)$-module $\scr G(U)/\varphi(\scr F)(U)$ is generated by $\sigma_1,\dots,\sigma_r$, i.e., each element of $\scr G(U)/\varphi(\scr F)(U)$ is an $\scr O(U)$-linear combination of $\sigma_1,\dots,\sigma_r$.

Since $\sigma_1,\dots,\sigma_r$ generate   $\scr E|_U$, for each $y\in U$, we know that $\sigma_1(y),\dots,\sigma_r(y)$ span $\scr E|y$. If we can show that $\sigma_1(y),\dots,\sigma_r(y)$ are linearly independent, then they form a basis of $\scr E|y$, which implies that $r_y=r=r_x$. This will finish the proof. Choose any $c_1,\dots,c_r\in\Cbb$ satisfying  $c_1\sigma_1(y)+\cdots+c_r\sigma_r(y)=0$. Then the germ of $c_1\sigma_1+\cdots+c_r\sigma_r$ at $y$ belongs to $\fk m_y\scr E_y$. Thus, the germ of $c_1\wht \sigma_1+\cdots+c_r\wht \sigma_r$ at $y$ belongs to $\fk m_y\scr G_y+\varphi(\scr F)_y$. (Note that $\varphi:\scr F\rightarrow\scr G$ descends to $\scr F_y\rightarrow \scr G_y$ and furthermore to $\scr F|y\rightarrow\scr G|y$; see \eqref{eq223} and \eqref{eq224}.) This means precisely that $c_1\wht \sigma_1(y)+\cdots+c_r\wht \sigma_r(y)\in\varphi(\scr F_y)/\fk m_y\scr G_y=\varphi(\scr F|y)$. 

It is clear that each fiber of $\scr F=\scr O_U^I$ is spanned by the values of all global sections. Thus, we can choose $\wht t\in\scr F(U)$ satisfying $c_1\wht \sigma_1(y)+\cdots+c_r\wht \sigma_r(y)=\varphi(\wht t(y))$. Set $\wht u=c_1\wht \sigma_1+\cdots+c_r\wht \sigma_r-\varphi(\wht t)$. Then $\wht u(y)=0$. Regard $\wht u$ as a $\Cbb^J$-valued holomorphic function, and notice that its value at $y$ is $0$. If we let $\{\wht e_j\}|_{j\in J}$ be the standard basis of $\Cbb^J$, regarded as constant sections of $\scr O_U^J(U)$, then there exist a finite subset $\{\wht e_{j_1},\dots,\wht e_{j_m}\}$ of $\{\wht e_j\}|_{j\in J}$ and $f_1,\dots,f_m\in\scr O(U)$ satisfying $\wht u=f_1\wht e_{j_1}+\cdots+f_m\wht e_{j_m}$, and that $f_1(y)=\cdots=f_m(y)=0$. Let $e_{j_1},\cdots,e_{j_m}\in \scr E(U)$ be the corresponding equivalence classes, which are clearly in $\scr G(U)/\varphi(\scr F)(U)$. Then $c_1\sigma_1+\cdots+c_r\sigma_r=f_1e_{j_1}+\cdots+f_me_{j_m}$. Since  $e_{j_1},\dots,e_{j_m}$ are $\scr O(U)$-linear combinations of $\sigma_1,\dots,\sigma_r$, one can find $g_1,\dots,g_r\in\scr O(U)$ whose values at $y$ are all $0$, such that $c_1\sigma_1+\cdots+c_r\sigma_r=g_1\sigma_1+\cdots+g_r\sigma_r$. Let $h_1=c_1-g_1,\dots,h_r=c_r-g_r$. Then $h_1(y)=c_1,\dots,h_r(y)=c_r$, and $h_1\sigma_1+\cdots+h_r\sigma_r=0$. By Lemma \ref{lb61}, we have $h_1=\cdots=h_r=0$. Thus $c_1=\cdots=c_r=0$.
\end{proof}

\begin{co}
If $\scr E$ is a coherent $\scr O_X$-module equipped with a connection $\nabla$, then $\scr E$ is locally free.
\end{co}

\chapter{Sheaves of VOAs}\label{lb90}

\section{Vertex operator algebras}
Let $\mbb V$ be a  complex vector space with grading $\mbb V=\bigoplus_{n\in\mathbb N}\mbb V(n)$ satisfying $\dim \mbb V(n)<\infty$ for each $n$ and $\dim\Vbb(0)=1$ (the \textbf{CFT type} condition)\footnote{The only reason we assume $\dim\Vbb(0)=1$ is to apply Buhl's result in \cite{Buhl02}; see Theorem \ref{lb122}.}.   We assume that there is a linear map
\begin{gather}
\mbb V\rightarrow(\text{End }(\mbb V))[[z^{\pm1}]]\nonumber\\
u\mapsto Y(u,z)=\sum_{n\in\mathbb Z}Y(u)_n z^{-n-1}\label{eq50}
\end{gather}
where each  $Y(u)_n\in\End(\Vbb)$ is called a \textbf{mode} of the operator $Y(u,z)$. Note that we have for any $u\in \Vbb,n\in\mbb Z$ that
\begin{align}
\Res_{z=0} Y(u,z)\cdot z^ndz=Y(u)_n.\label{eq36}
\end{align}

\begin{df}
We say that $(\mbb V,Y)$ (or $\mbb V$ for short) is a (positive energy) \textbf{vertex operator algebra} (VOA), if for any $u\in\Vbb$  the following conditions are satisfied.

(a) (\textbf{Lower truncation}) For any $v\in \Vbb$, 
\begin{gather*}
Y(u,z)v\in\Vbb((z)).
\end{gather*}

(b) (\textbf{Jacobi identity}) For any $u,v\in \Vbb$ and $m,n,h\in\mathbb Z$, we have
\begin{align}
&\sum_{l\in\mathbb N}{m\choose l}Y(Y(u)_{n+l}\cdot v)_{m+h-l}\nonumber\\
=&\sum_{l\in\mathbb N}(-1)^l{n\choose l}Y(u)_{m+n-l}Y(v)_{h+l}-\sum_{l\in\mathbb N}(-1)^{l+n}{n\choose l}Y(v)_{n+h-l}Y(u)_{m+l}.\label{eq39}
\end{align}

(c) There exists a vector $\mbf 1\in \Vbb(0)$ (the \textbf{vacuum vector}) \index{1@$\id$} such that $Y(\mbf 1,z)=\id_\Vbb$.

(d) (\textbf{Creation property}) For any $v\in\Vbb$, we have
\begin{align*}
Y(v,z)\id-v\in z\Vbb[[z]].
\end{align*}
This is equivalent to that $Y(v)_{-1}\id=v$ and $Y(v)_n\id=0$ for any $n\in\mathbb N$.

(e) There exists a vector $\cbf\in \Vbb$ (the \textbf{conformal vector}) \index{c@$\cbf$} such that the operators $L_n:=Y(\cbf)_{n+1}$ ($n\in\mathbb Z$) satisfy the Virasoro relation:
\begin{gather}
[L_m,L_n]=(m-n)L_{m+n}+\frac 1 {12}(m^3-m)\delta_{m,-n}c.\label{eq40}
\end{gather}
Here the number $c\in\mathbb C$ is called the \textbf{central charge} of $V$.

(f) If $v\in \Vbb(n)$ then $L_0v=nv$. $n$ is called the \textbf{conformal weight} (or the \textbf{energy}) of $v$ and will be denoted by $\wt(v)$. \index{wt@$\wt(v),\wt(w),\wtd\wt(w)$} $L_0$ is called the \textbf{energy operator}. We say that a vector $v\in \mbb V$ is \textbf{homogeneous} if $v\in\mbb V(n)$ for some $n\in\mbb Z$.

(g)	(\textbf{$L_{-1}$-derivative}) $\frac d{dz} Y(v,z)=Y(L_{-1}v,z)$ for any $v\in\Vbb$.
\end{df}

Note that from the creation property, we have $\cbf=Y(\cbf)_{-1}\id=L_{-2}\id$. Since $L_0\id=0$, by the Virasoro relation, we have $L_0\cbf=L_0L_{-2}\id=L_{-2}L_0\id+2L_{-2}\id=2\cbf$. We conclude
\begin{align*}
\cbf\in\Vbb(2).
\end{align*}

We explain the meaning of Jacobi identity. Let $\Vbb'=\bigoplus_{n\in\mbb Z}\Vbb(n)^*$ where each $\Vbb(n)^*$ is the dual vector space of $\Vbb(n)$. \index{V@$\Vbb',\Vbb^*$} Then $\Vbb'$ is a subspace of the dual space $\Vbb^*$ of $\Vbb$. A vector $w'\in\Vbb^*$ is inside $\Vbb'$ if and only if there exits $N\in\mbb N$ such that  $\bk{u,w'}=0$ whenever $u\in\Vbb(n)$ and $n>N$. For each $v\in\Vbb$ and $n\in\mbb Z$, the transpose $Y(v)_n^\tr$ of $Y(v)_n$ is a linear map on $\Vbb^*$. We then define
\begin{align*}
Y(v,z)^\tr=\sum_{n\in\mbb Z}Y(v)_n^\tr z^{-n-1}\quad\End(\Vbb^*)[[z]]
\end{align*}
to be the transpose of $Y(v,z)$. Then $Y(v,z)^\tr dz=-\sum_{m\in\mbb Z}Y(v)_m^\tr (z^{-1})^{m-1}d(z^{-1})$, which shows that for each $n\in\mbb Z$,
\begin{align}
\Res_{z^{-1}=0}Y(v,z)^\tr\cdot z^n dz=-Y(v)_n^\tr.\label{eq37}
\end{align}
We assume that the lower truncation property is satisfied, and that for any $v,w\in \Vbb,w'\in\Vbb'$,
\begin{gather}
\qquad Y(v,z)^\tr w'\in\Vbb'((z^{-1})).\label{eq31}
\end{gather}
Then we have
\begin{align*}
\bk{Y(v,z)w,w'}\in\mbb C[z,z^{-1}],
\end{align*}
namely, $z\mapsto \bk{Y(v,z)w,w'}$ is a meromorphic function on $\mbb P^1$ with poles possibly at $0,\infty$, namely, it is an element of $H^0(\mbb P^1,\scr O_{\mbb P^1}(\blt(0+\infty)))$. We will see that if $\Vbb$ is a VOA then these conditions hold automatically. Let $\mbb C^\times=\mbb C-\{0\}$. \index{C@$\mbb C^\times$} For each $z\in\mbb C^\times$, $\bk{Y(v,z)w,w'}$ is a complex number. Choose any $u\in\Vbb$, and let $\zeta$ be also a standard coordinate of $\mbb C$. Then $\zeta-z$ (where $z$ is a fixed complex number) is the standard coordinate of a neighborhood of $z$. By lower truncation property and condition \eqref{eq31}, we have
\begin{gather}
\bk{Y(v,z)Y(u,\zeta)w,w'}\in\mbb C((\zeta)),\label{eq32}\\
\bk{Y(Y(u,\zeta-z)v,z)w,w'}\in\mbb C((\zeta-z)),\label{eq33}\\
\bk{Y(u,\zeta)Y(v,z)w,w'}\in\mbb C((\zeta^{-1})).\label{eq34}
\end{gather}

\begin{thm}\label{lb19}
Let $(\Vbb,Y)$ satisfy the lower truncation property, and assume that for any $v,w\in\Vbb,w'\in\Vbb'$, condition \eqref{eq31} holds. Then the Jacobi identity is equivalent to the requirement that for any $u,v,w\in\Vbb$, $w'\in\Vbb'$, $z\in\mbb C^\times$, there exists $f\in H^0(\mbb P^1,\scr O_{\mbb P^1}(\blt (0+z+\infty)))$ whose Laurent series expansions near $0,z,\infty$ are \eqref{eq32}, \eqref{eq33}, and \eqref{eq34} respectively. 
\end{thm}

\begin{proof}
We apply the strong residue Theorem \ref{lb18} to the  single $3$-pointed Riemann sphere $(\mbb P^1;0,z,\infty)$ and the sheaf $\scr O_{\mbb P^1}$. Then such $f$ exists if and only if for any $\lambda\in H^0(\mbb P^1,\omega_{\mbb P^1}(\blt(0+z+\infty)))$,
\begin{align}
\Res_{\zeta-z=0}f_z\lambda=-\Res_{\zeta^{-1}=0} f_\infty\lambda-\Res_{\zeta=0}f_0\lambda,\label{eq35}
\end{align}
where $f_0,f_z,f_\infty$ are defined by  \eqref{eq32}, \eqref{eq33}, and \eqref{eq34} respectively. It is easy to see that $H^0(\mbb P^1,\omega_{\mbb P^1}(\blt(0+z+\infty)))$ is spanned by $\zeta^m(\zeta-z)^nd\zeta$ (where $m,n\in\mbb Z$). Thus, $f$ exists if and only if \eqref{eq35} holds whenever $\lambda=\zeta^m(\zeta-z)^nd\zeta$. Assuming $\lambda$ is defined like this. Then, using \eqref{eq36}, we compute
\begin{align*}
&\Res_{\zeta-z=0}f_z\lambda=\Res_{\zeta-z=0}\bk{Y(Y(u,\zeta-z)v,z)w,w'}\zeta^m(\zeta-z)^nd(\zeta-z)\\
=&\sum_{l\in\Nbb}{m\choose l}\Res_{\zeta-z=0}\bk{Y(Y(u,\zeta-z)v,z)w,w'}(\zeta-z)^{n+l}z^{m-l}d(\zeta-z)\\
=&\sum_{l\in\Nbb}{m\choose l}\bk{Y(Y(u)_{n+l}v,z)w,w'}z^{m-l}.
\end{align*}
Similar computations using \eqref{eq36} and \eqref{eq37} give the explicit expression of the two terms on the right hand side of \eqref{eq35}, which show that \eqref{eq35} is equivalent to
\begin{align}
&\sum_{l\in\Nbb}{m\choose l}\bk{Y(Y(u)_{n+l}v,z)w,w'}z^{m-l}\nonumber\\
=&\sum_{l\in\Nbb}{n\choose l}(-1)^l\bk{Y(u)_{m+n-l}Y(v,z)w,w'}z^l\nonumber\\
&-\sum_{l\in\Nbb}{n\choose l}(-1)^{n-l}\bk{Y(v,z)Y(u)_{m+l}w,w'}z^{n-l}.\label{eq38}
\end{align}
Note that by lower truncation property and condition \eqref{eq31}, the three terms in the above equation are all finite sums. We conclude that the requirement in this theorem holds if and only if for any $u,v,w\in\mbb V$ and $m,n\in\mbb Z$, \eqref{eq38} holds where $z$ is considered as a variable. This means that  we are now considering \eqref{eq38} as an equation of elements in $\mbb C[[z^{\pm 1}]]$. For each $h\in\mbb Z$, apply $\Res_{z=0}(\cdot)z^hdz$ to both sides of \eqref{eq38}, we get precisely \eqref{eq39} evaluated between $w$ and $w'$. Since $w,w'$ are arbitrary, we see that our requirement is equivalent to the Jacobi identity.
\end{proof}

From now on, we shall always assume that $\Vbb$ is a VOA (of CFT type). From the Jacobi identity and the $L_{-1}$-derivative, one has for each $n$ that
\begin{align*}
[L_0,Y(v)_n]=Y(L_0 v)_n-(n+1)Y(v)_n,
\end{align*}
This shows that if $v,w\in\mbb V$ are homogeneous, then $Y(v)_nw$ is also homogeneous with conformal weight $\wt(v)+\wt(w)-(n+1)$. Thus condition \eqref{eq31} follows easily. Note also that $[L_0,L_n]=-nL_n$ implies for any $n\in\mbb Z,s\in\Cbb$ that
\begin{align}
L_n\Vbb(s)=\Vbb(s-n).
\end{align}
In particular, when $n\in\Zbb_+$, we have $L_n\id=0$ since it is inside the trivial subspace $\Vbb(-n)$.

\begin{rem}\label{lb58}
When $n>2$, we have $L_n\cbf=0$ since $\Vbb(2-n)$ is trivial. Using the fact that $\cbf=L_{-2}\id$, that $L_n\id=0$ when $n\geq0$,  and the Virasoro relation \eqref{eq40}, we compute $L_1\cbf=L_1L_{-2}\id=[L_1,L_{-2}]\id=3L_{-1}\id=0$ since $Y(L_{-1}\id,z)=\partial_zY(\id,z)=\partial_z\id=0$. Also, $L_2\cbf=L_2L_{-2}\id=[L_2,L_{-2}]\id=4L_0\id+\frac c 2\id=\frac c 2\id$ where $c$ is the central charge. We conclude
\begin{gather*}
L_0\cbf=2\cbf\qquad L_1\cbf=0,\qquad L_2\cbf=\frac c2 \id.
\end{gather*}
\end{rem}

\section{VOA modules}

Let $\Wbb$ be a vector space equipped with a linear map
\begin{gather}
\Vbb\rightarrow(\text{End }\Wbb)[[z^{\pm1}]]\nonumber\\
u\mapsto Y_\Wbb(u,z)=\sum_{n\in\mathbb Z}Y_\Wbb(u)_n z^{-n-1}.\label{eq51}
\end{gather}

\begin{df}
We say that $(\Wbb,Y_\Wbb)$ \index{Y@$Y,Y_\Wbb$} (or $\Wbb$ for short) is a \textbf{weak $\Vbb$-module} if the lower truncation property holds, i.e., for any $u\in\Vbb,w\in\Wbb$, $Y_\Wbb(u,z)w\in\Wbb((z))$, if $Y_\Wbb(\id,z)=\id_\Wbb$, and if for any $m,n,h\in\mbb Z$ and $u,v\in\Vbb$, the Jacobi identity \eqref{eq39} holds with $Y$ replaced by $Y_\Wbb$, i.e.,
\begin{align}
&\sum_{l\in\mathbb N}{m\choose l}Y_\Wbb(Y(u)_{n+l}\cdot v)_{m+h-l}\nonumber\\
=&\sum_{l\in\mathbb N}(-1)^l{n\choose l}Y_\Wbb(u)_{m+n-l}Y_\Wbb(v)_{h+l}-\sum_{l\in\mathbb N}(-1)^{l+n}{n\choose l}Y_\Wbb(v)_{n+h-l}Y_\Wbb(u)_{m+l}.\label{eq192}
\end{align}
Homomorphism and endomorphisms of weak $\Vbb$-modules are the linear maps commuting with the actions of vertex operators.
\end{df}

Set $n=0$ in \eqref{eq192}, we obtain the commutator formula
\begin{align}
[Y_\Wbb(u)_m,Y_\Wbb(v)_h]=\sum_{l\in\mathbb N}{m\choose l}Y_\Wbb(Y(u)_l\cdot v)_{m+h-l}.\label{eq242}
\end{align}

If $\Wbb_1,\Wbb_2$ are weak $\Vbb$-modules, we set $\Hom_\Vbb(\Wbb_1,\Wbb_2)$ to be the space of homomorphisms from $\Wbb_1$ to $\Wbb_2$. We set $\End_\Vbb(\Wbb)=\Hom_\Vbb(\Wbb,\Wbb)$.

For a weak $\Vbb$-module $\Wbb$, we set $L_n=Y_\Wbb(\cbf)_{n+1}$.
Then the Virasoro relation \eqref{eq40} holds for the same central charge $c$. Moreover, the $L_{-1}$-derivative property holds: 
\begin{align}
\frac d{dz} Y_\Wbb(v,z)=Y_\Wbb(L_{-1}v,z)\label{eq134}
\end{align}
for any $v\in\Vbb$. We refer the reader to \cite{DLM97} for the proof. The $L_{-1}$-derivative property and the Jacobi identity implies
\begin{align}
[L_0,Y_\Wbb(v)_n]=Y_\Wbb(L_0 v)_n-(n+1)Y_\Wbb(v)_n,\label{eq104}
\end{align}
for any $v\in\Vbb$ and $n\in\Zbb$. 

\index{L0Ln@$\wtd L_0,\wtd L_n$}
\begin{df}
A weak $\Vbb$-module $\Wbb$ is called an \textbf{admissible $\Vbb$-module} if there exists $A\in\End_\Vbb(\Wbb)$ such that $\wtd L_0:=L_0+A$ is diagonal (on $\Wbb$), and that the eigenvalues of $\wtd L_0$ are natural numbers. If, moreover, each eigenspace of $\wtd L_0$ is finite-dimensional, we say that $\Wbb$ is a \textbf{finitely admissible $\Vbb$-module}.
\end{df}

\begin{rem}
According to \eqref{eq104}, that $\Wbb$ is admissible is equivalent to that there is a diagonalizable $\wtd L_0\in\End(\Wbb)$ with spectrum in $\Nbb$ satisfying $[\wtd L_0,Y_\Wbb(v)_n]=Y_\Wbb(\wtd L_0 v)_n-(n+1)Y_\Wbb(v)_n$. Equivalently, $\Wbb$ has grading $\Wbb=\bigoplus_{n\in\Nbb}\Wbb(n)$ \index{Wn@$\Wbb(n),\Wbb_{(n)}$} such that 
\begin{align}
Y_\Wbb(v)_m\Wbb(n)\subset \Wbb(n+\wt(v)-m-1)\label{eq203}
\end{align}
for any $m,n\in\Nbb$ and homogeneous $v\in\Vbb$. Indeed, $\Wbb(n)$ is the $n$-eigenspace of $\wtd L_0$. Moreover, $\Wbb$ is finitely admissible if and only if each $\Wbb(n)$ is finite dimensional. We write $\wtd\wt(w)=n$ if $w\in\Wbb(n)$, equivalently, $\wtd L_0w=\wtd\wt(w)w$.\index{wt@$\wt(v),\wt(w),\wtd\wt(w)$} We say that $w$ is \textbf{$\wtd L_0$-homogeneous} (with weight $n$) if $w$ is an eigenvector of $\wtd L_0$ (with eigenvalue $n$). 
\end{rem}

\begin{df}
We say that $\Wbb$ is a \textbf{(grading-restricted ordinary) $\Vbb$-module}, if $\Wbb$ is a weak $\Vbb$-module, if there exists a finite subset $E\subset\Cbb$ such that  $\Wbb$ has grading $\Wbb=\bigoplus_{s\in E+\Nbb}\Wbb_{(s)}$,  and if for each $s\in\mbb C$ we have $\dim \Wbb_{(s)}<+\infty$ and $L_0|_{\Wbb_{(s)}}=s\id_{\Wbb_{(s)}}$. \index{Wn@$\Wbb(n),\Wbb_{(n)}$}  A vector $w\in\Wbb$ is called \textbf{homogeneous} if $w\in W(s)$ for some $s\in \mbb C$. In this case, we write $\wt(w)=s$ \index{wt@$\wt(v),\wt(w),\wtd\wt(w)$} and call it the \textbf{(conformal) weight} of $w$. 
\end{df}
Note that when $v$ is homogeneous, \eqref{eq104} is equivalent to
\begin{align}
Y_\Wbb(v)_m\Wbb_{(s)}\subset \Wbb_{(s+\wt(v)-m-1)}.\label{eq105}
\end{align}

\begin{rem}\label{lb137}
If $\Wbb$ is a $\Vbb$-module, then $\Wbb$ is finitely admissible. Indeed, one suffices to assume that any two elements in $E$ do not differ by an integer.   By \eqref{eq105}, for each $\alpha\in E$, $\Wbb_{\alpha+\Nbb}:=\bigoplus_{s\in\alpha+\Nbb}\Wbb_{(s)}$ is a weak $\Vbb$-submodule of $\Wbb$. Moreover, $\Wbb$ is the direct sum of all such $\Wbb_{\alpha+\Nbb}$. One then defines $\wtd L_0$ whose action on each $\Wbb_{\alpha+\Nbb}$ is $L_0-\alpha$. This makes $\Wbb$ finitely admissible. 
\end{rem}

Note that $\Vbb$ itself is a $\Vbb$-module, called the \textbf{vacuum $\Vbb$-module}. For a $\Vbb$-module $\Wbb$, we can give a similar interpretation of Jacobi identity as in Theorem \ref{lb19}. We leave the details to the reader. 

For a $\Vbb$-module $\Wbb$, consider the dual vector space $\Wbb^*$ and the graded dual $\Wbb':=\bigoplus_{s\in\mbb C}\Wbb_{(s)}^*$ where each $\Wbb_{(s)}^*$ is the dual vector space of $\Wbb_{(s)}$. \index{W@$\Wbb',\Wbb^*$} Then $\Wbb'$ is equipped with a natural $\Vbb$-module structure: the vertex operator $Y_{\Wbb'}$ is defined such that for any $v\in\Vbb,w\in\Wbb,w'\in\Wbb'$,
\begin{align}
\bk{Y_{\Wbb'}(v,z)w',w}=\bk{w',Y_{\Wbb}(e^{zL_1}(-z^{-2})^{L_0} v,z^{-1})w}.\label{eq41}
\end{align}
Here, if $v$ is homogeneous, then $(-z^{-2})^{L_0}v$ is understood as $(-z^{-2})^{\wt(v)}v$. In general, $(-z^{-2})^{L_0}v$ is defined by linearity. We briefly write
\begin{align}
Y_{\Wbb'}(v,z)=Y_{\Wbb}(e^{zL_1}(-z^{-2})^{L_0} v,z^{-1})^\tr.\label{eq72}
\end{align}
The meaning of $e^{zL_1}(-z^{-2})^{L_0}$ will be explained in example \ref{lb25}. We call $\Wbb'$ the \textbf{contragredient module} of $\Wbb$. We have $\Wbb''=\Wbb$. See \cite{FHL93} chapter 5 for more details. When $v$ is homogeneous, it is easy to check that for any $n\in\mbb Z$,
\begin{align}
Y_{\Wbb'}(v)_n=\sum_{m\in\mathbb N}\frac{(-1)^{\wt(v)}}{m!}Y_\Wbb(L_1^mv)_{-n-m-2+2\wt(v)}^\tr.\label{eq211}
\end{align}
As a consequence, we have
\begin{align}
L_n^\tr=L_{-n}.
\end{align}
Let $\wtd L_0$ act on $\Wbb'$ as the transpose of $\wtd L_0\curvearrowright\Wbb$, i.e., set for any $w\in\Wbb,w'\in\Wbb'$ that
\begin{align}
\bk{\wtd L_0 w,w'}=\bk{w,\wtd L_0 w'}.\label{eq106}
\end{align}
Then $\wtd L_0$ makes $\Wbb'$ admissible. $\wtd L_0$ yields the grading $\Wbb'=\bigoplus_{n\in\Nbb}\Wbb(n)^*$. Thus $\Wbb'$ is also finitely admissible under $\wtd L_0$.

\begin{cv}\label{lb36}
For a $\Vbb$-module $\Wbb$, we always assume $\wtd L_0$ is chosen such that $\Wbb$ is finitely admissible (see Remark \ref{lb137}). For its contragredient module $\Wbb'$, we always assume that the actions of $\wtd L_0$ on $\Wbb$ and on $\Wbb'$ satisfy \eqref{eq106}.  If $\Wbb$ is semi-simple, i.e., a finite direct some of irreducible $\Vbb$-modules, except when $\Wbb=\Vbb$, we assume  each irreducible submodule $\Mbb$ of $\Wbb$  is $\wtd L_0$-invariant, and the lowest eigenvalue of $\wtd L_0|_{\Mbb}$ (with non-trivial eigenspace) is $0$. These  assumptions are compatible. For the vacuum module $\Vbb$, we set $\wtd L_0=L_0$. Thus $\Vbb(n)=\Vbb_{(n)}$.
\end{cv}

The following fact will be used later without explicit mentioning.

\begin{pp}\label{lb131}
Let $\Wbb$ be an irreducible $\Vbb$-module, and let $T$ be an endomorphism of $\Wbb$. Then $T$ is a scalar multiplication. 
\end{pp}
In particular, if  $\wtd L_0$ makes $\Wbb$ admissible, then $\wtd L_0-L_0$ is a scalar multiplication.

\begin{proof}
Choose any $s\in\Cbb$ so that the $L_0$-weight space $\Wbb_{(s)}$ is nontrivial. Note that $\Wbb_{(s)}$ is also finite-dimensional by the definition of ordinary modules. Since $T$ commutes with $L_0$, $T$ preserves the eigenvalues of the eigenvectors of $L_0$. So $T\Wbb_{(s)}\subset\Wbb_{(s)}$. Thus, we can find an eigenvalue $\lambda$ of $T|\Wbb_{(s)}$. It follows that the kernel of $T-\lambda\id_\Wbb$ is a nontrivial $\Vbb$-invariant subspace of $\Wbb$, which must be $\Wbb$. So $T=\lambda\id_\Wbb$.
\end{proof}

\begin{df}\label{lb100}
Let $\Vbb_1,\Vbb_2$ be VOAs. If $\Wbb$ is a vector space which is both a weak $\Vbb_1$-module $(\Wbb,Y_+)$ and a weak $\Vbb_2$-module $(\Wbb,Y_-)$. We say that $(\Wbb,Y_+,Y_-)$ is a weak $\Vbb_1\times\Vbb_2$-module if $[Y_+(u)_m,Y_-(v)_n]=0$ for any $u\in\Vbb_1,v\in\Vbb_2,m,n\in\Zbb$.
\end{df}

\section{Change of coordinates}\label{lb52}

We define a group $(\Gbb,\circ)$ as follows. \index{G@$\Gbb,\Lie(\Gbb)$} As a set $\Gbb$ consists of all $\rho\in\scr O_{\Cbb,0}$ such that $\rho(0)=0$ and $\rho'(0)\neq 0$. If $\rho_1,\rho_2\in\Gbb$, then their multiplication is just the composition $\rho_1\circ\rho_2$. We should understand elements in $\Gbb$ as maps but not functions.  $\Gbb$ acts on $\scr O_{\Cbb,0}$ as $\rho\star f=f\circ \rho^{-1}$ if $\rho\in\Gbb,f\in\scr O_{\Cbb,0}$. The Lie algebra $\Lie(\Gbb)$ of $\Gbb$ is spanned by $L_0,L_1,L_2,\dots$, where for each $n\in\Nbb$, 
\begin{align}
L_n=z^{n+1}\partial_z,\label{eq145}
\end{align}
whose action on each $f\in\scr O_{\Cbb,0}$, denoted by $L_n\star f$, is defined by
\begin{align}
L_n\star f:=-L_nf=-z^{n+1}\partial_zf.\label{eq45}
\end{align}
The Lie bracket relation is defined using above relation, i.e., satisfying $[L_m,L_n]\star f=L_m\star L_n\star f-L_n\star L_m\star f$. It is the negative of the usual Lie bracket for vector fields. One easily checks that $[\cdot,\cdot]$ is compatible with the Virasoro relation \eqref{eq40}.

Let $\Wbb$ be a $\Vbb$-module. Given any $\rho\in\Gbb$, one can define $\mc U(\rho)\in\End(\Vbb)$ as follows. \index{U@$\mc U(\rho),\mc U(\eta_\blt)$} Choose $c_0,c_1,c_2,\dots\in\Cbb$ such that $c_0\neq 0$, and 
\begin{align}
\rho=c_0^{L_0}\circ \exp\Big(\sum_{n>0}c_n L_n\Big)\label{eq43}
\end{align}
when acting on any $f\in\scr O_{\Cbb,0}$ by $\star$. Then 
\begin{align}
\mc U(\rho)=c_0^{\wtd L_0}\cdot \exp\Big(\sum_{n>0}c_n L_n\Big),\label{eq46}
\end{align}
where the Virasoro operators $\wtd L_0,L_n$ are acting on $\Wbb$. (The reason we use $\wtd L_0$ but not $L_0$ is that $c_0^{L_0}$ might not be single-valued on $\Wbb$.) 

Let $z$ be the standard coordinate of $\mbb C$, regarded as an element in $\scr O_{\Cbb,0}$. It is the identity element of $\Gbb$. Then, for any $f\in\scr O_{\Cbb,0}$, we have $\rho\star f=f\circ\rho^{-1}=f(\rho^{-1}(z))=f(\rho\star z)$. We conclude
\begin{align}
\rho\star f=f(\rho\star z),\label{eq44}
\end{align}
which shows that the action of $\Gbb$ on $\scr O_{\Cbb,0}$ is determined by its action on $z$. For example, since $z\partial_z z=z$ and hence $(z\partial_z)^nz=z$ for each $n$, we have
\begin{align*}
c_0^{-L_0}\star z=c_0^{z\partial_z}z=\exp(\log c_0\cdot z\partial_z)z=\sum_{n\in\Nbb}\frac 1{n!}(\log c_0)^n(z\partial_z)^n z=c_0z.
\end{align*}
Therefore, we conclude that
\begin{align}
c_0^{-L_0}\star f=f(c_0z).
\end{align}

We now give a more direct relation between $\rho$ and the coefficients $c_0,c_1,c_2,\dots$ in \eqref{eq43}. It is clear that $\rho(z)$, as an element in $\scr O_{\Cbb,0}$, equals
\begin{align*}
&\rho(z)=z\circ\rho=\rho^{-1}\star z=\exp\Big(-\sum_{n>0}c_n L_n\Big)\star c_0^{-L_0}\star z\\
=&\exp\Big(-\sum_{n>0}c_n L_n\Big)\star (c_0z) \xlongequal{\eqref{eq44}} c_0\cdot \exp\Big(-\sum_{n>0}c_n L_n\Big)\star z,
\end{align*}
which, together with \eqref{eq45}, shows that
\begin{align}
\boxed{\rho(z)=c_0\cdot \exp\Big(\sum_{n>0}c_n z^{n+1}\partial_z \Big) z}
\end{align}
One can use the above equation to completely determine the coefficients $c_0,c_1,\dots$. For instance, if we write
\begin{align}
\rho(z)=a_1z+a_2z^2+a_3z^3+\cdots,\label{eq48}
\end{align}
then one has
\begin{gather}
c_0=a_1,\nonumber\\
c_1c_0=a_2,\nonumber\\
c_2c_0+c_1^2c_0=a_3.\nonumber
\end{gather}
In particular, one has $c_0=\rho'(0)$. Thus \eqref{eq46} could be rewritten as
\begin{align}
\boxed{\mc U(\rho)=\rho'(0)^{\wtd L_0}\cdot \exp\Big(\sum_{n>0}c_n L_n\Big)}\label{eq47}
\end{align}
Notice $a_n=\rho^{(n)}(0)/n!$, we have
\begin{gather}
c_0=\rho'(0),\nonumber\\
c_1=\frac 12\frac{\rho''(0)}{\rho'(0)},\nonumber\\
c_2=\frac 16 \frac{\rho'''(0)}{\rho'(0)}-\frac 14\Big(\frac{\rho''(0)}{\rho'(0)}\Big)^2.\label{eq165}
\end{gather}
The following is (essentially) proved in \cite{Hua97} section 4.2:
\begin{thm}
For each $\Vbb$-module $\Wbb$, $\mc U$ is a representation of $\Gbb$ on $\Wbb$. Namely, we have $\mc U(\rho_1\circ\rho_2)=\mc U(\rho_1)\mc U(\rho_2)$ for each $\rho_1,\rho_2\in\Gbb$.
\end{thm}

\begin{eg}\label{lb25}
We have seen that
\begin{align}
c_0^{z\partial_z}z=c_0z.\label{eq70}
\end{align}
We now calculate $\exp(c_1z^2\partial_z)z$. It is easy to see that $(c_1z^2\partial_z)^nz=n!c_1^nz^{n+1}$. Thus $\exp(c_1z^2\partial_z)z=\sum_{n=0}^{\infty}c_1^nz^{n+1}=z/(1-c_1z)$. We conclude
\begin{align}
\exp(c_1z^2\partial_z)z=\frac z{1-c_1z}.\label{eq71}
\end{align}
Then it is easy to see that $\exp(c_1z^2\partial_z)c_0^{z\partial_z}z=c_0z/(1-c_0c_1z)=c_0^{z\partial_z}\exp(c_0c_1z^2\partial_z)z$. We conclude
\begin{align}
e^{c_1L_1}c_0^{\wtd L_0}=c_0^{\wtd L_0}e^{c_0c_1L_1}.
\end{align}
Set \index{zz@$\upgamma_\xi,\upgamma_1$}
\begin{align}
\boxed{~~\upgamma_\xi(z)=\frac{1}{\xi+z}-\frac 1\xi.~~}
\end{align}
Then the inverse of $\upgamma_\xi$ is $\upgamma_{\xi^{-1}}$. By \eqref{eq70} and \eqref{eq71}, it is easy to see that the following identity holds when acting on any $\Wbb$:
\begin{align}
\boxed{~~\mc U(\upgamma_\xi)=e^{\xi L_1}(-\xi^{-2})^{\wtd L_0}.~~}\label{eq73}
\end{align}
Thus \eqref{eq72} could be rewritten as
\begin{align}
Y_{\Wbb'}(v,z)=Y_{\Wbb}(\mc U(\upgamma_z)v,z^{-1})^\tr.\label{eq92}
\end{align}
It is easy to see that $\upgamma_\xi(\xi z)=\xi^{-1}\upgamma_1(z)$.
Therefore,
\begin{align}
\mc U(\upgamma_\xi)\xi^{\wtd L_0}=\xi^{-\wtd L_0}\mc U(\upgamma_1).\label{eq78}
\end{align}
\end{eg}

Let $n\in\Zbb$. Then $\mc U(\rho)$ does not preserve the vector space $\Wbb(n)$ (the $n$-eigenspace of $\wtd L_0$). However, $\Wbb$ has filtration $\emptyset=\Wbb^{\leq -1}\subset \Wbb^{\leq 0}\subset\Wbb^{\leq 1}\subset\Wbb^{\leq 2}\subset\cdots$, where \index{Vn@$\Vbb^{\leq n},\Wbb^{\leq n}$}
\begin{align}
\Wbb^{\leq n}=\bigoplus_{k\leq n}\Wbb(k).
\end{align}
Then, by \eqref{eq203}, $L_m\Wbb^{\leq n}\subset\Wbb^{\leq n-m}$ for any $m\in\Zbb$. Thus, by \eqref{eq47}, we conclude that $\mc U$ restricts to a representation of $\Gbb$ on $\Wbb^{\leq n}$. Moreover, if $w\in\Wbb$ is $\wtd L_0$-homogeneous, then
\begin{align}
\mc U(\rho)w=\rho'(0)^{\wtd\wt(w)}w~~\mod~~ \Wbb^{\leq \wtd\wt(w)-1}.\label{eq56}
\end{align}
In other words, the action of $\mc U(\rho)$ on $\Wbb^{\leq n}/\Wbb^{\leq n-1}$ is $\rho'(0)^n\id$.

Finally, we discuss holomorphic families of transformations. Let $X$ be a complex manifold and $\rho:X\rightarrow\mbb G,x\mapsto \rho_x$ a function. We say that $\rho$ is a \textbf{holomorphic family} of transformations if for any $x\in X$, there exists an open subset $V\subset X$ containing $x$ and an open $U\subset\Cbb$ containing $0$ such that  $(z,y)\in U\times V\mapsto \rho_y(z)$ is a holomorphic function on $U\times V$. Then it is clear that the coefficients $a_1,a_2,\dots$ in \eqref{eq48} depend holomorphically on the parameter $x\in X$. Hence the same true for $c_0,c_1,c_2,\dots$. Thus, by the formula \eqref{eq47}, for any $w\in\Wbb^{\leq n}$, $x\in X\mapsto \mc U(\rho_x)w$ is a $\Wbb^{\leq n}$-valued holomorphic function on $X$. Thus $\mc U(\rho)$ can be regarded as an isomorphism of  $\scr O_X$-modules 
\begin{align}
\mc U(\rho): \Wbb^{\leq n}\otimes_{\Cbb}\scr O_X\xrightarrow{\simeq}\Wbb^{\leq n}\otimes_{\Cbb}\scr O_X\label{eq64}
\end{align}
sending each constant function $w$ to the section $x\mapsto \mc U(\rho_x)w$. \index{U@$\mc U(\rho),\mc U(\eta_\blt)$} Its inverse is  $\mc U(\rho^{-1})$.

\begin{cv}
For any open subset $V\subset X$, any $v\in\Vbb^{\leq n}$  (resp. $w\in\Wbb^{\leq n}$) is also understood as the constant section $v\otimes 1$ (resp. $w\otimes 1$) in $(\Vbb^{\leq n}\otimes_{\Cbb}\scr O_X)(V)$ (resp. $(\Wbb^{\leq n}\otimes_{\Cbb}\scr O_X)(V)$).
\end{cv}

The following lemma will be used in Section \ref{lb50}. We let $\wtd L_n$ be $L_n$ if $n\neq 0$. \index{L0Ln@$\wtd L_0,\wtd L_n$}

\begin{lm}\label{lb29}
Let $T$ be an open subset of $\Cbb$ containing $0$. Let $\rho:T\rightarrow\Gbb,\zeta\mapsto\rho_\zeta$ be a holomorphic family of transformations satisfying $\rho_0(z)=z$. Then, for any $w\in\Wbb$, 
\begin{align}
\partial_\zeta\mc U(\rho_\zeta)w\Big|_{\zeta=0}=\sum_{n\geq 1}\frac 1{n!}\Big(\partial_\zeta\rho_\zeta^{(n)}(0)\Big|_{\zeta=0}\Big)\wtd L_{n-1}\cdot w\label{eq82}
\end{align}
where $\partial_\zeta\rho_\zeta^{(n)}(z)=\partial_z^n\partial_\zeta\rho_\zeta(z)=\partial_\zeta\partial_z^n\rho_\zeta(z)$.
\end{lm}

Note that when $w$ is a (non-necessarily constant) section of $\Wbb\otimes_\Cbb\scr O_T$, one needs to take $\partial_\zeta w$ into account when calculating the left hand side of \eqref{eq82}. Also, as the derivative of $w=\mc U(\rho_\zeta)\mc U(\rho_\zeta)^{-1}w$ is $0$, we obtain
\begin{align}
\partial_\zeta\big(\mc U(\rho_\zeta)^{-1}\big)w\Big|_{\zeta=0}=-\sum_{n\geq 1}\frac 1{n!}\Big(\partial_\zeta\rho_\zeta^{(n)}(0)\Big|_{\zeta=0}\Big)\wtd L_{n-1}\cdot w.\label{eq248}
\end{align}

\begin{proof}
Let $c_1,c_2,\dots\in\scr O_{\Cbb}(T)$ such that 
\begin{gather*}
\rho_\zeta(z)=\rho_\zeta'(0)\exp\Big(\sum_{n\geq 1}c_n(\zeta)z^{n+1}\partial_z\Big)(z).
\end{gather*}
Then $\rho_\zeta(z)$ equals
\begin{align*}
\rho_\zeta'(0)\Big(z+\sum_{n\geq 1}c_n(\zeta)z^{n+1}\Big)
\end{align*}
plus some polynomials of $z$ multiplied by at least two terms among $c_1(\zeta),c_2(\zeta),\dots$. Since $\rho_0(z)=z$, we have $\rho_0'(0)=1$ and $c_1(0)=c_2(0)=\cdots=0$. Therefore,
\begin{align*}
\partial_\zeta\rho_\zeta(z)\Big|_{\zeta=0}=\partial_\zeta\rho_\zeta'(0)z\Big|_{\zeta=0}+\Big(\sum_{n\geq 1}\partial_\zeta c_n(0)z^{n+1}\Big),
\end{align*}
which implies that when $n\geq 2$,
\begin{align*}
\frac 1{n!}\partial_\zeta\rho_\zeta^{(n)}(0)\Big|_{\zeta=0}=\partial_\zeta c_{n-1}(0).
\end{align*}
Thus, using \eqref{eq47} and again $\rho_0'(0)=1,c_1(0)=c_2(0)=\cdots=0$, we compute
\begin{align*}
\partial_\zeta\mc U(\rho_\zeta)w\Big|_{\zeta=0}=&\partial_\zeta\Big(\rho_\zeta'(0)^{\wtd L_0}\cdot \exp\Big(\sum_{n\geq1}c_n(\zeta) \wtd L_n\Big)\Big)w\Big|_{\zeta=0}\\
=&\partial_\zeta\rho_\zeta'(0)\wtd L_0w\Big|_{\zeta=0}+\Big(\sum_{n\geq 1}\partial_\zeta c_n(0)\wtd L_nw\Big).
\end{align*}
Now equation \eqref{eq82} follows from the last two equations.
\end{proof}

\begin{rem}\label{lb140}
Let $A=L_0-\wtd L_0$. Then using the fact that $\Res_{z=0}z^n Y_\Wbb(\cbf,z)=Y_\Wbb(\cbf)_n=L_{n-1}$, we can write \eqref{eq82} and \eqref{eq248} in the following form:
\begin{align}
&\partial_\zeta\mc U(\rho_\zeta)w\big|_{\zeta=0}=-\partial_\zeta\mc U(\rho_\zeta^{-1})w\big|_{\zeta=0}\nonumber\\
=&\Res_{z=0}~\partial_\zeta\rho_\zeta(z) Y_\Wbb(\cbf,z) w dz\big|_{\zeta=0}-\partial_\zeta\rho_\zeta'(0)Aw\big|_{\zeta=0}.
\end{align}
\end{rem}

\section{Sheaves of VOAs on complex curves}\label{lb23}

Let $C$ be a (non-necessarily compact) Riemann surface. Let $U,V$ be open subsets of $C$, and  $\eta:U\rightarrow\Cbb,\mu:V\rightarrow\Cbb$ univalent maps.  Define a holomorphic family $\varrho(\eta|\mu):U\cap V\rightarrow\Gbb$ \index{zz@$\varrho(\eta\lvert\mu)$} as follows. For any $p\in U\cap V$, $\eta-\eta(p)$ and $\mu-\mu(p)$ are local coordinates at $p$. We set $\varrho(\eta|\mu)_p\in\Gbb$ satisfying
\begin{align}
\eta-\eta(p)=\varrho(\eta|\mu)_p(\mu-\mu(p)).\label{eq49}
\end{align}
Let $z\in\scr O_{\Cbb,0}$ be the standard coordinate. Then, by composing both sides of \eqref{eq49} with $\mu^{-1}(z+\mu(p))$, we find the equivalent formula
\begin{align}
\varrho(\eta|\mu)_p(z)=\eta\circ\mu^{-1}(z+\mu(p))-\eta(p),
\end{align}
which justifies that $\varrho(\eta|\mu)$ is analytic. It is also clear that if $\eta_1,\eta_2,\eta_3$ are three local coordinates, then on their common domain the following cocycle condition holds:
\begin{align}
\varrho(\eta_3|\eta_1)=\varrho(\eta_3|\eta_2)\varrho(\eta_2|\eta_1).
\end{align}

Note that the linear map \eqref{eq51} can be extended to a homomorphism of $\Cbb((z))$-modules
\begin{gather}
\mbb V((z))\rightarrow(\text{End }\mbb W)[[z^{\pm1}]]\nonumber\\
u\otimes f\in \Vbb\otimes_{\Cbb}\Cbb((z))\mapsto Y_{\Wbb}(u\otimes f,z)=f(z)Y(u,z).\label{eq228}
\end{gather}
Moreover, for any $w\in\Wbb$ and $v\in\Vbb((z))$, it is clear that $Y(v,z)w\in\Wbb((z))$. The following theorem is one of the main results of \cite{Hua97}.

\begin{thm}\label{lb31}
Let $\Wbb$ be a $\Vbb$-module. Let $U\subset\Cbb$ be a neighborhood of $0$. Let $\alpha\in\scr O(U)$ be a local coordinate at $0$, let $\id_\Cbb\in\Gbb$ be the standard coordinate of $\Cbb$ (i.e.  the identity element of $\Gbb$), and let $z$ be the standard complex variable of $\Cbb$ (different from $\id_\Cbb$). Then for any $v\in\Vbb$ and $w\in\Wbb$, we have the following equation of elements in $\Wbb((z))$:
\begin{align}
\mc U(\alpha)Y_\Wbb(v,z)\mc U(\alpha)^{-1}\cdot w=Y_{\Wbb}\big(\mc U(\varrho(\alpha|\id_\Cbb))v,\alpha(z)\big)\cdot w.\label{eq52}
\end{align}
Note that $\mc U(\varrho(\alpha|\id_\Cbb))v$ is in $\Vbb\otimes_{\Cbb}\scr O(U)$ and hence can be regarded as an element of $\Vbb((z))$. Of course, \eqref{eq52} also holds in an obvious way for any $v\in\Vbb((z))$.
\end{thm}

For instance, take $\alpha(z)=\lambda z$ where $\lambda\in\Cbb^\times$. Then $\varrho(\alpha|\id_\Cbb)$ is constantly $\lambda$. It follows that
\begin{align}
\lambda^{\wtd L_0}Y_\Wbb(v,z)\lambda^{-\wtd L_0}=Y_\Wbb(\lambda^{L_0}v,\lambda z).\label{eq107}
\end{align}

The \textbf{sheaf of VOA} associated to $\Vbb$ is an $\scr O_C$-module $\scr V_C$ \index{VCVX@$\scr V_C,\scr V_C^{\leq n},\scr V_{\fk X},\scr V_{\fk X}^{\leq n}$} defined by 
\begin{align}
\scr V_C=\varinjlim_{n\in\Nbb}\scr V_C^{\leq n},\label{eq53}
\end{align}
where for each $n\in\Nbb$, $\scr V_C^{\leq n}$ is a locally free sheaf of rank $\dim\Vbb^{\leq n}$ described as follows. For any open subset $U\subset C$ and a univalent $\eta:U\rightarrow\Cbb$, we have an isomorphism of $\scr O_U$-modules \index{UV@$\mc U_\varrho(\eta),\mc U_\varrho(\varphi),\mc V_\varrho(\eta),\mc V_\varrho(\varphi)$}
\begin{align}\label{eq90}
\mc U_\varrho(\eta):\scr V_C^{\leq n}|_U\xrightarrow{\simeq} \Vbb^{\leq n}\otimes_{\Cbb}\scr O_U.
\end{align}
These isomorphisms are defined in such a way that if  $\mu:V\rightarrow\Cbb$ is also univalent, then on $U\cap V$ we have
\begin{align}
\boxed{~~\mc U_\varrho(\eta)\mc U_\varrho(\mu)^{-1}=\mc U(\varrho(\eta|\mu))~~}\quad\in\End_{\scr O_{U\cap V}}(\Vbb^{\leq n}\otimes_{\Cbb}\scr O_{U\cap V}).\label{eq55}
\end{align}
From \eqref{eq56}, we can compute that for any section $v$ of $\Vbb^{\leq n}\otimes_{\Cbb}\scr O_{U\cap V}$,
\begin{align}
\mc U_\varrho(\eta)\mc U_\varrho(\mu)^{-1}\cdot v=(\partial_\mu\eta)^n\cdot v~~\mod ~~\Vbb^{\leq n-1}\otimes_{\Cbb}\scr O_{U\cap V}.\label{eq58}
\end{align}

Now assume that $C$ is a (compact and possibly nodal) complex curve. We define for each $n\in\Nbb$ a locally free sheaf $\scr V_C^{\leq n}$ as follows. ($\scr V_C$ is defined again using \eqref{eq53}.) Let $\Sigma=\{x_1',x_2',\dots,x_M'\}$ be the set of nodes, and let $C_0=C-\Sigma$. Then $\scr V_{C_0}^{\leq n}$ is defined as above. Let $U$ be an open subset of $C$ containing only one of $\Sigma$, say $x_j'$. Let $\nu:\wtd C\rightarrow C$ be the normalization of $C$ as in Section \ref{lb20}, and let $\{y_j',y_j''\}=\nu^{-1}(x_j')$. Then $\nu^{-1}(U)$ is a disjoint union of two open subsets $V'\ni y_j',V''\ni y_j''$ of $\wtd C$. We assume that $U$ is small enough such that there exit local coordinates $\xi_j:V'\rightarrow\mbb C$ and $\varpi_j:V''\rightarrow\mbb C$ at $y_j'$ and $y_j''$ respectively. This means that $\xi_j,\varpi_j$ are univalent, and $\xi_j(y')=\varpi_j(y'')=0$. We also identify
\begin{align}
U-\{x_j'\} \simeq (V'-\{y_j'\})\sqcup (V''-\{y_j''\})
\end{align}
via $\nu$. Now, let $\scr V_C^{\leq n}(U)$ be the $\scr O_C(U)$-submodule of $\scr V_{C_0}^{\leq n}(U-\{x_j'\})$ generated by
\begin{align}
\boxed{~~\mc U_\varrho(\xi_j)^{-1}\big(\xi_j^{L_0}v\big)+\mc U_\varrho(\varpi_j)^{-1}\big(\varpi_j^{L_0}\mc U(\upgamma_1)v\big)~~}\qquad (\forall v\in\Vbb^{\leq n}),\label{eq54}
\end{align}
where we recall from example \ref{lb25} that $\mc U(\upgamma_1)=e^{L_1}(-1)^{L_0}$. To be more precise, \eqref{eq53} is a section on $(V'-\{y_j'\})\sqcup (V''-\{y_j''\})$ which equals $\mc U_\varrho(\xi_j)^{-1}\big(\xi_j^{L_0}v\big)$ on $(V'-\{y_j'\})$ and $\mc U_\varrho(\varpi_j)^{-1}\big(\varpi_j^{L_0}\mc U(\upgamma_1)v\big)$ on $V''$. Also, $\xi_j^{L_0}$ is an element of $\scr O_{C_0}(V'-\{y_j'\})$ acting on the constant section $v\in\Vbb^{\leq n}\otimes_{\Cbb}\scr O_{C_0}(V'-\{y_j'\})$, and $\varpi_j^{L_0}\mc U(\upgamma_1)v$ is understood in a similar way.  It is easy to see that $\scr V_C^{\leq n}(U)$ is generated freely by
\begin{align}
\mc U_\varrho(\xi_j)^{-1}\big(\xi_j^{\wt (v)}v\big)+\mc U_\varrho(\varpi_j)^{-1}\big((-\varpi_j)^{\wt(v)}\mc U(\upgamma_1)v\big)\label{eq225}
\end{align}
for all $v\in E$ where $E$ is any basis of $\Vbb^{\leq n}$ whose elements are homogenous. Since $\upgamma_1=\upgamma_1^{-1}$, $\scr V_C^{\leq n}(U)$ is also generated freely by
\begin{align*}
\mc U_\varrho(\xi_j)^{-1}\big(\xi_j^{L_0}\mc U(\upgamma_1)v\big)+\mc U_\varrho(\varpi_j)^{-1}\big(\varpi_j^{L_0}v\big)
\end{align*}
for all $v\in E$. By the gluing construction, we obtain the locally free $\scr O_C$-module $\scr V_C^{\leq n}$.

\begin{pp}\label{lb22}
Let $C$ be a complex curve and $n\in\Nbb$. Then we have the following isomorphism of  $\scr O_C$-modules:
\begin{align}
\scr V_C^{\leq n}/\scr V_C^{\leq n-1}\simeq\Vbb(n)\otimes_{\Cbb}\Theta_C^{\otimes n}.\label{eq267}
\end{align}
Under this isomorphism, if $U\subset C$ is open and smooth, and $\eta\in\scr O(U)$ is univalent, then for any $v\in\Vbb(n)$, $v\otimes \partial_\eta^n$ is identified with the equivalence class of $\mc U_\varrho(\eta)^{-1}v$.
\end{pp}

\begin{proof}
Recall that $\Sigma$ is the set of nodes. By the transition function\eqref{eq58}, we obtain a surjective $\scr O_{C-\Sigma}$-module morphism $\Psi:\scr V_{C-\Sigma}^{\leq n}\rightarrow\Vbb(n)\otimes\Theta_{C\setminus\{x\}}^{\otimes n}$ sending $\mc U_\varrho(\eta)^{-1}v$ to $v\otimes \partial_\eta^n$ if $v\in\Vbb(n)$, and to $0$ if $v\in\Vbb^{\leq n-1}$. $\Psi$ has kernel $\scr V_{C-\Sigma}^{\leq n-1}$. Now let $U$ be a neighborhood of $x_j'$ as in the setting of \eqref{eq54}. Then $\Psi$ sends \eqref{eq54} to $v\otimes \xi_j^n\partial_{\xi_j}^n|_{V'-\{y'\}}+v\otimes(-\varpi_j)^n\partial_{\varpi_j}^n|_{V''-\{y''\}}$ whenever $v\in\Vbb(n)\cap E$. (Recall that $E$ is a homogeneous basis of $\Vbb^{\leq n}$.) From this and \eqref{eq57} we see that $\Psi$ restricts to a surjective $\scr O_C$-module morphism $\Psi:\scr V_C^{\leq n}\rightarrow\Vbb(n)\otimes\Theta_C^{\otimes n}$ and that $\Ker\Psi(U)$ is $\scr O_C(U)$-generated by  \eqref{eq54} for all $v\in \Vbb^{\leq n-1}\cap E$. Thus $\Psi$ descends to an isomorphism \eqref{eq267}.
\end{proof}

As a consequence, we now prove a vanishing theorem for the sheaf of VOA.

\begin{thm}\label{lb21}
Let $\fk X=(C;x_1,\dots,x_N)$ be an $N$-pointed complex curve with $M$ nodes.  Let  $\wtd C$ be the normalization of $C$, and let $\wtd g$ be the largest genus of the connected components of $\wtd C$.  Then for any $n\in\Nbb$, there exists $k_0\in\Zbb$ depending only on $n,\wtd g,M$ such that 
\begin{align}
H^1\big(C,\scr V_C^{\leq n}\otimes\omega_C(k\SX)\big)=0\label{eq59}
\end{align}
for any $k>k_0$.
\end{thm}

Recall that the divisor $\SX$ is defined by $x_1+\cdots+x_N$. Also, by our definition of pointed complex curves, each connected component of $\wtd C$ contains at least one of (the pre-image of) $x_1,\dots,x_N$. Recall also that the dualizing sheaf $\omega_C$ is the inverse of $\Theta_C$.

We will see from the proof that $k_0$ can be chosen to be $n|2\wtd g-2|+2M$. Note that by Proposition \ref{lb121}, one can always find $k_0\in\Nbb$ such that \eqref{eq59} holds for any $k>k_0$. The importance of the present theorem is, however, that we may find $k_0$ which is independent of the complex structure of $C$.

\begin{proof}
$\scr V^{\leq -1}$ is trivial since $\Vbb^{\leq -1}$ is so. Since the vacuum vector $\id$ is killed by $L_0,L_1,L_2,\dots$, it is fixed by the action of $\Gbb$. From this and the fact that $\Vbb(0)$ is spanned by $\id$, it is clear that
\begin{align}
\scr V_C^{\leq 0}/\scr V_C^{\leq -1}\simeq \scr V_C^{\leq 0}\simeq \scr O_C.
\end{align}
Thus, by Theorem \ref{lb10},  the vanishing property \eqref{eq59} holds for any $k>2M$.

We now prove the theorem by induction. By Proposition \ref{lb22}, there is a  short exact sequence 
\begin{align*}
0\rightarrow \scr V_C^{\leq n-1}\otimes\omega_C(k\SX)\rightarrow \scr V_C^{\leq n}\otimes\omega_C(k\SX) \rightarrow \Vbb(n)\otimes_{\Cbb}\Theta_C^{\otimes(n-1)}(k\SX)\rightarrow 0,
\end{align*}
which induces a long exact sequence
\begin{align}
H^1(C,\scr V_C^{\leq n-1}\otimes\omega_C(k\SX))&\rightarrow H^1(C,\scr V_C^{\leq n}\otimes\omega_C(k\SX))\nonumber\\
&\rightarrow H^1\big(C,\Vbb(n)\otimes_{\Cbb}\Theta_C^{\otimes(n-1)}(k\SX)\big).\label{eq60}
\end{align}
Suppose the statement in our theorem is true for $n-1$ and any $k>(n-1)|2\wtd g-2|+2M$. Then, by  induction and Theorem \ref{lb10},  the first and the last terms of \eqref{eq60} equal $0$ for any $k>n|2\wtd g-2|+2M$. So the same is true for the middle term.
\end{proof}

\begin{rem}\label{lb88}
Let $D$ be an effective divisor on $C$. Using the same argument, one can show that Theorem \ref{lb21} holds verbatim if $\scr V_C^{\leq n}$ is replaced by $\scr V_C^{\leq n}(-D)$, except that $k_0$ should now also depend on $\deg D$.
\end{rem}

\section{Sheaves of VOAs on families of complex curves}

Let $\fk X=(\pi:\mc C\rightarrow\mc B)$ be a family of complex curves. Recall that $\Sigma$ is the critical locus. Let $U,V$ be open subsets of $\mc C-\Sigma$, and let $\eta:U\rightarrow\Cbb,\mu:V\rightarrow\Cbb$ be holomorphic functions such that $(\eta,\pi)$ and $(\mu,\pi)$ are biholomorphic maps from $U$ resp. $V$ to open subsets of $\Cbb\times\mc B$. This requirement is equivalent to that $\eta$ and $\mu$ are univalent on each fiber of $U$ and $V$ respectively. \index{zz@$\varrho(\eta\lvert\mu)$} For each $p\in U\cap V$, we define $\varrho(\eta|\mu)_p\in\scr O_{\Cbb,0}$ by
\begin{align}
\varrho(\eta|\mu)_p(z)=\eta\circ(\mu,\pi)^{-1}\big(z+\mu(p),\pi(p)\big)-\eta(p).\label{eq62}
\end{align}
Then $\varrho(\eta|\mu)_p$ is a holomorphic function of $z$ on  $\mu\big((U\cap V)_{\pi(p)}\big)$ where $(U\cap V)_{\pi(p)}$ is the fiber $U\cap V\cap \pi^{-1}(\pi(p))$. It is easy to check that for each $n\in\Nbb$,
\begin{align}
\partial_z^n\varrho(\eta|\mu)_p(0)=\partial_\mu^n \eta(p)\label{eq226},
\end{align}
where the partial derivative $\partial_\mu$ is defined to be vertical to $d\pi$. From this, we see that $\varrho(\eta|\mu)_p(0)=0$ and  $\partial_z\varrho(\eta|\mu)_p(0)\neq0$. So $\varrho(\eta|\mu)_p$ is an element of $\Gbb$. We thus obtain a family of transformations $\varrho(\eta|\mu):U\cap V\rightarrow\Gbb,p\mapsto \varrho(\eta|\mu)_p$, which  is clearly holomorphic.

As in Section \ref{lb23}, $\varrho(\eta|\mu)$ is also described by
\begin{align}
\eta-\eta(p)\big|_{(U\cap V)_{\pi(p)}}=\varrho(\eta|\mu)_p\big(\mu-\mu(p)\big|_{(U\cap V)_{\pi(p)}}\big).\label{eq63}
\end{align}
To see this, one composes both sides of \eqref{eq63} with $(\mu,\pi)^{-1}\big(z+\mu(p),\pi(p)\big)$. This relation shows that,   if for $j=1,2,3$ we have $\eta_j\in\scr O(U_j)$ which is univalent on each fiber of $U_j$, then  on $U_1\cap U_2\cap U_3$ we have
\begin{align}
\varrho(\eta_3|\eta_1)=\varrho(\eta_3|\eta_2)\varrho(\eta_2|\eta_1).
\end{align}



\subsection*{Definition of $\scr V_{\fk X}$}

We set \index{VCVX@$\scr V_C,\scr V_C^{\leq n},\scr V_{\fk X},\scr V_{\fk X}^{\leq n}$} 
\begin{align}
\scr V_{\fk X}=\varinjlim_{n\in\Nbb}\scr V_{\fk X}^{\leq n},
\end{align}
where for each $n\in\Nbb$, $\scr V_{\fk X}^{\leq n}$ is a locally free $\scr O_{\mc C}$-module  defined as follows. Suppose that $U$ is an open subset of $\mc C-\Sigma$ and $\eta\in\scr O(U)$ is univalent on each fiber of $U$. Then we have an isomorphism of $\scr O_U$-modules \index{UV@$\mc U_\varrho(\eta),\mc U_\varrho(\varphi),\mc V_\varrho(\eta),\mc V_\varrho(\varphi)$}
\begin{align}
\mc U_\varrho(\eta):\scr V_{\fk X}^{\leq n}|_U\xrightarrow{\simeq} \Vbb^{\leq n}\otimes_{\Cbb}\scr O_U.
\end{align}
These isomorphisms are defined in such a way that if $V$ is another open subset of $\mc C-\Sigma$ and $\mu\in\scr O_{\mc C}(V)$ is also univalent on each fiber, then on $U\cap V$ we have
\begin{align}
\boxed{~~\mc U_\varrho(\eta)\mc U_\varrho(\mu)^{-1}=\mc U(\varrho(\eta|\mu))~~}\quad\in\End_{\scr O_{U\cap V}}(\Vbb^{\leq n}\otimes_{\Cbb}\scr O_{U\cap V}).\label{eq76}
\end{align}
Recall that $\varrho(\mu|\eta)$ is a family of transformations on $U\cap V$, and $\mc U(\varrho(\mu|\eta))$ is defined as by \eqref{eq64}. Thus, we can defined $\scr V_{\fk X}^{\leq n}|_{\mc C-\Sigma}$ by gluing. Note that by \eqref{eq56} and \eqref{eq226}, we can compute that for any section $v$ of $\Vbb^{\leq n}\otimes_{\Cbb}\scr O_{U\cap V}$,
\begin{align}
\mc U_\varrho(\eta)\mc U_\varrho(\mu)^{-1}\cdot v=(\partial_\mu\eta)^n\cdot v~~\mod ~~\Vbb^{\leq n-1}\otimes_{\Cbb}\scr O_{U\cap V}.\label{eq79}
\end{align}

To define $\scr V_{\fk X}^{\leq n}$ near $\Sigma$, let $x'\in\Sigma$, and assume that near the fiber $\mc C_{\pi(x')}$, the family $\fk X$ is obtained via sewing  a family $\wtd{\fk X}$ of Riemann surfaces with local coordinates  as in Section \ref{lb24}. Then, by \eqref{eq19} and \eqref{eq227}, one can identify a neighborhood $W_j$ of $x'$ as $\mc D_{r_j}\times\mc D_{\rho_j}\times\mc D_{r_\bullet\rho_\bullet\setminus j}\times\wtd{\mc B}$ such that, by setting $\wtd {\mc B}_j=\mc D_{r_\bullet\rho_\bullet\setminus j}\times\wtd{\mc B}$, the projection $\pi|_{W_j}$ equals
\begin{align}
\pi:W_j=\mc D_{r_j}\times\mc D_{\rho_j}\times\wtd{\mc B}_j\xrightarrow{~\pi_{r_j,\rho_j}\times\id~} \mc D_{r_j\rho_j}\times\wtd{\mc B}_j.\label{eq69}
\end{align}
We thus have
\begin{gather}
W_j=\mc D_{r_j}\times\mc D_{\rho_j}\times\wtd{\mc B}_j,\label{eq115}\\
W_j\cap\Sigma=(0,0)\times\wtd{\mc B}_j.\nonumber
\end{gather}
As in \eqref{eq66} and \eqref{eq67}, we have open subsets of $W_j$:
\begin{gather*}
W_j'=\mc D_{r_j}^\times\times\mc D_{\rho_j}\times\wtd{\mc B}_j,\qquad W_j''=\mc D_{r_j}\times\mc D_{\rho_j}^\times\times\wtd{\mc B}_j.
\end{gather*}
Then it is clear that
\begin{align*}
W_j-\Sigma=W_j'\cup W_j''.
\end{align*}
Let $\xi_j,\varpi_j$ be the standard coordinates of $\mc D_{r_j},\mc D_{\rho_j}$, and extend them constantly to $W_j\rightarrow\mc D_{r_j},W_j\rightarrow\mc D_{\rho_j}$ respectively. Then $(\xi_j,\pi)$ and $(\varpi_j,\pi)$ are holomorphic open embeddings of $W_{j'},W_{j''}$ respectively; equivalently, $\xi_j,\varpi_j$ are univalent on each fiber of $W_j'$ and $W_j''$ respectively.

We shall define $\scr V^{\leq n}_{\fk X}|_{W_j}$ to be an $\scr O_{W_j}$-submodule of $\scr V_{\fk X}^{\leq n}|_{W_j-\Sigma}$ generated (freely) by some sections on $W_j$ whose restrictions to $W_j'$ and $W_j''$ are described under the trivilizations $\mc U_\varrho(\xi_j)$ and $\mc U_\varrho(\varpi_j)$ respectively. For that purpose, we need to first calculate the transition function \begin{align*}
\mc U_\varrho(\varpi_j)\mc U_\varrho(\xi_j)^{-1}=\mc U(\varrho(\varpi_j|\xi_j)):\Vbb^{\leq n}\otimes_{\Cbb}\scr O_{W_j'\cap W_j''}\xrightarrow{\simeq} \Vbb^{\leq n}\otimes_{\Cbb}\scr O_{W_j'\cap W_j''}.
\end{align*}
Set $q_j=\pi_{r_j,\rho_j}=\xi_j\varpi_j$. 

\begin{lm}\label{lb123}
Choose any $p\in W_j'\cap W_j''$. Then we have
\begin{gather*}
\varrho(\varpi_j|\xi_j)_p(z)=q_j(p)\upgamma_{\xi_j(p)}(z)
\end{gather*}
and hence
\begin{align*}
\mc U(\varrho(\varpi_j|\xi_j)_p)=q_j(p)^{L_0}\mc U(\upgamma_{\xi_j(p)}).
\end{align*}
\end{lm}
\begin{proof}
Choose any $x\in (W_j'\cap W_j'')_{\pi(p)}$. Then $\pi(x)=\pi(p)$ and hence $q_j(x)=q_j(p)$.  Since $\varpi_j=\xi_j^{-1}q_j$, we have
\begin{align*}
\varpi_j(x)-\varpi_j(p)=q_j(p)(\xi_j(x)^{-1}-\xi_j(p)^{-1}).
\end{align*}
By \eqref{eq63}, we have
\begin{align*}
\varpi_j(x)-\varpi_j(p)=\varrho(\varpi_j|\xi_j)_p(\xi_j(x)-\xi_j(p)).
\end{align*}
If we compare these two equations and set $z=\xi_j(x)-\xi_j(p)$, we obtain
\begin{align*}
&\varrho(\varpi_j|\xi_j)_p(z)=\varrho(\varpi_j|\xi_j)_p(\xi_j(x)-\xi_j(p))=q_j(p)(\xi_j(x)^{-1}-\xi_j(p)^{-1})\\
=&q_j(p)\big((\xi_j(p)+z)^{-1}-\xi_j(p)^{-1} \big)=q_j(p)\upgamma_{\xi_j(p)}(z).
\end{align*}
\end{proof}

We   define $\scr V_{\fk X}^{\leq n}|_{W_j}$ to be the  $\scr O_{W_j}$-submodule of $\scr V_{\fk X-\Sigma}^{\leq n}|_{W_j-\Sigma}$ generated by any section on $W_j-\Sigma$ whose restrictions to $W_j'$ and $W_j''$ are
\begin{align}
\boxed{~~\mc U_\varrho(\xi_j)^{-1}\big(\xi_j^{L_0}v\big)\qquad\text{resp.}\qquad \mc U_\varrho(\varpi_j)^{-1}\big(\varpi_j^{L_0}\mc U(\upgamma_1)v\big)  ~~}\label{eq75}
\end{align}
where $v\in\Vbb^{\leq n}$. Since $\upgamma_1=\upgamma_1^{-1}$ and hence $\mc U(\upgamma_1)=\mc U(\upgamma_1)^{-1}$, this definition is symmetric with respect to $\xi_j$ and $\varpi_j$. To check that \eqref{eq75} is well-defined, we need:

\begin{lm}
The two sections defined in \eqref{eq75} agree on $W_j'\cap W_j''$.
\end{lm}

\begin{proof}
Using   \eqref{eq78} and Lemma \ref{lb123}, we check that
\begin{align*}
&\mc U_\varrho(\varpi_j)\mc U_\varrho(\xi_j)^{-1}\xi_j^{L_0}v=\mc U(\varrho(\varpi_j|\xi_j))\xi_j^{L_0}v=q_j^{L_0}\mc U(\upgamma_{\xi_j})\xi_j^{L_0}v\\
=&q_j^{L_0}\xi_j^{-L_0}\mc U(\upgamma_1)v=\varpi_j^{L_0}\mc U(\upgamma_1)v.
\end{align*}
\end{proof}

It is easy to see that, if we take $v\in E$ where $E$ is a basis of $\Vbb^{\leq n}$ consisting of homogeneous vectors, then $\scr V_{\fk X}^{\leq n}|_{W_j}$ is generated freely by sections defined by \eqref{eq75} for all $v\in E$. Thus, by the gluing construction, we obtain a locally free $\scr O_{\mc C}$-module $\scr V_{\fk X}^{\leq n}$.

\begin{rem}\label{lb46}
Since the vacuum vector $\id$ is annihilated by $L_n$ ($n\geq 0$), we see that $\id$ is fixed by any transition function $\mc U(\varrho(\eta|\mu))$. Thus, we can define unambiguously an element $\id\in\scr V_{\fk X}(\mc C-\Sigma)$ (the \textbf{vacuum section}) such that for any open $U\subset\mc C-\Sigma$ and any $\eta\in\scr O(U)$ univalent on each fiber, $\mc U_\varrho(\eta)\id$ is the vaccum vector $\id$ (considered as a constant function). \index{1@$\id$} Also, by \eqref{eq75}, it is clear that
\begin{align*}
\id\in\scr V_{\fk X}(\mc C).
\end{align*}
\end{rem}

\subsection*{Restriction to fibers}

As in Section \ref{lb23}, we may use \eqref{eq79}, \eqref{eq75}, and \eqref{eq80} to show:

\begin{pp}\label{lb56}
For any $n\in\Nbb$, we have the following isomorphism of  $\scr O_C$-modules:
\begin{align}
\scr V_{\fk X}^{\leq n}/\scr V_{\fk X}^{\leq n-1}\simeq\Vbb(n)\otimes_{\Cbb}\Theta_{\mc C/\mc B}^{\otimes n}.\label{eq153}
\end{align}
Under this isomorphism, if $U\subset \mc C-\Sigma$ is open and smooth, and $\eta\in\scr O(U)$ is univalent on each fiber of $U$, then for any $v\in\Vbb(n)$, $v\otimes \partial_\eta^n$ is identified with the equivalence class of $\mc U_\varrho(\eta)^{-1}v$.
\end{pp}

By comparing the transition functions and looking at the generating sections near the nodes, it is easy to see:
\begin{pp}\label{lb26}
For any $n\in\Nbb$ and $b\in\mc B$, we have a natural isomorphism
\begin{align}
\scr V_{\fk X}^{\leq n}|\mc C_b\simeq\scr V_{\mc C_b}^{\leq n}.
\end{align}
\end{pp}

The following theorem is a generalization of Theorem \ref{lb9}. 

\begin{thm}\label{lb27}
Let $\fk X=(\pi:\mc C\rightarrow\mc B;\sgm_1,\dots,\sgm_N)$ be a family of $N$-pointed complex curves.  Let $n\in\Nbb$. Then  there exists $k_0\in\mbb N$ such that for any $k\geq k_0$, the $\scr O_{\mc B}$-module  $\pi_*\big(\scr V_{\fk X}^{\leq n}\otimes\omega_{\mc C/\mc B}(k\SX)\big)$ is locally free, and for any $b\in\mc B$ there is a natural isomorphism of vector spaces
\begin{align}
\frac{\pi_*\big(\scr V_{\fk X}^{\leq n}\otimes\omega_{\mc C/\mc B}(k\SX)\big)_b}{~\fk m_b\cdot\pi_*\big(\scr V_{\fk X}^{\leq n}\otimes\omega_{\mc C/\mc B}(k\SX)\big)_b~}\simeq H^0\big(\mc C_b,\scr V_{\mc C_b}^{\leq n}\otimes\omega_{\mc C_b}(k\SX(b))\big)
\end{align} 
defined by  restriction of sections. In particular,  $\dim  H^0\big(\mc C_b,\scr V_{\mc C_b}^{\leq n}\otimes\omega_{\mc C_b}(k\SX(b))\big)$ is locally constant over $b$.
\end{thm}

\begin{proof}
Recall that by Ehresmann's result, if $\fk X$ is smooth, then by our assumption in Section \ref{lb1}, $\mc B$  has finitely many connected components, and all the fibers over a connected component  are diffeomorphic.  Thus, by theorems \ref{lb7} and \ref{lb21} and Remark \ref{lb124}, for sufficiently large $k$, $H^r\big(\mc C_b,\scr V_{\mc C_b}^{\leq n}\otimes\omega_{\mc C_b}(k\SX(b))\big)$ vanishes for any $b\in\mc B$ and $r\geq 1$. Since the restriction of $\omega_{\mc C/\mc B}$ to $\mc C_b$ is $\omega_{\mc C_b}$, by Proposition \ref{lb26}, the restriction of $\scr V_{\fk X}^{\leq n}\otimes\omega_{\mc C/\mc B}(k\SX)$ to $\mc C_b$ is equivalent to $\scr V_{\mc C_b}^{\leq n}\otimes\omega_{\mc C_b}(k\SX(b))$.  Thus, our theorem follows easily from Grauert's Theorem \ref{lb11}.
\end{proof}

\begin{co}
Let $\fk X=(\pi:\mc C\rightarrow\mc B;\sgm_1,\dots,\sgm_N)$ be a family of $N$-pointed complex curves, and let $n\in\Nbb$. Then for any  Stein open subset $V$ of $\mc B$, there is $k_0\in\Nbb$ such that for  any integer $k\geq k_0$ and any $b\in V$,   the elements of $\pi_*\big(\scr V_{\fk X}^{\leq n}\otimes\omega_{\mc C/\mc B}(k\SX)\big)(V)$ (more precisely, their germs at $b$) generate the stalk $\pi_*\big(\scr V_{\fk X}^{\leq n}\otimes\omega_{\mc C/\mc B}(k\SX)\big)_b$, and their restrictions to $\mc C_b$ form the vector space $H^0\big(\mc C_b,\scr V_{\mc C_b}^{\leq n}\otimes\omega_{\mc C_b}(k\SX(b))\big)$.
\end{co}

\begin{proof}
Apply Theorem \ref{lb27}  and Cartan's theorem A (see Section \ref{lb20}). 
\end{proof}

As a variant (and easy consequence) of the above corollary, we have:

\begin{co}\label{lb33}
Let $\fk X=(\pi:\mc C\rightarrow\mc B;\sgm_1,\dots,\sgm_N)$ be a family of $N$-pointed complex curves. Then for any  Stein open subset $V$ of $\mc B$ and any $b\in V$, the elements of $\pi_*\big(\scr V_{\fk X}^{\leq n}\otimes\omega_{\mc C/\mc B}(\blt\SX)\big)(V)$ (more precisely, their germs at $b$) generate the stalk $\pi_*\big(\scr V_{\fk X}^{\leq n}\otimes\omega_{\mc C/\mc B}(\blt\SX)\big)_b$, and their restrictions to $\mc C_b$ form the vector space $H^0\big(\mc C_b,\scr V_{\mc C_b}^{\leq n}\otimes\omega_{\mc C_b}(\blt\SX(b))\big)$.
\end{co}

\subsection*{The subsheaf $\svir_c$}

We now define an important $\scr O_{\mc C}$-submodule $\svir_c$ \index{Virc@$\svir_c$} of $\scr V_{\fk X}^{\leq 2}$ related to the conformal vector $\cbf\in\Vbb(2)$.  If $U$ is an open subset of $\mc C-\Sigma$ equipped with a holomorphic $\eta:U\rightarrow\Cbb$ univalent on each fiber, then $\svir_c|_U$ is the  $\scr O_U$-submodule of $\scr V_{\fk X}|_U$ generated (freely) by $\mc U_\varrho(\eta)^{-1}\cbf$ and the vacuum section $\id$, which is locally free of rank $2$. This definition is independent of the choice of $\eta$. Indeed, if $\mu:U\rightarrow\Cbb$ is also univalent on each fiber, then $\mc U_\varrho(\mu)\mc U_\varrho(\eta)^{-1}\cbf=\mc U(\varrho(\mu|\eta))\cbf$, which can be calculated using the actions of $L_n$ ($n\geq 0$) on $\cbf$, is an $\scr O_U$-linear combination of $\cbf$ and $\id$ by Remark \ref{lb58}. Thus, by gluing all such $U$, we get $\svir_c|_{\mc C-\Sigma}$. Now assume that $U$ is a small neighborhood of a point of $\Sigma$. We let $\svir_c|_U$ be the submodule generated by the sections described in \eqref{eq75}, in which we set $v$ to be $\cbf$ and $\id$. This completes the definition of the $\scr O_{\mc C}$-submodule $\svir_c$. Note that the action of $\mc U(\varrho(\mu|\eta))$ and $\mc U(\upgamma_1)$ on $\cbf$ and $\id$ depends only on the central charge $c$ since this is true for $L_n$ ($n\geq 0$) by Remark \ref{lb58}. Thus, the $\scr O_{\mc C}$-module $\svir_c$ depends only on the number $c$ but not on $\Vbb$ or $\scr V_{\fk X}$.

By Proposition \ref{lb56}, we have a short exact sequence
\begin{align*}
0\rightarrow\scr V_{\fk X}^{\leq 1} \rightarrow\scr V_{\fk X}^{\leq 2}\xrightarrow{\uplambda}\Vbb(2)\otimes_\Cbb\Theta_{\mc C/\mc B}^{\otimes 2}\rightarrow 0
\end{align*}
where $\uplambda$ is described locally (outside $\Sigma$)  by sending  $\mc U_\varrho^{-1}(\eta)v$ (where $v\in\Vbb(2)$) to  $v\cdot \partial_\eta^2$ and sending the submodule $\scr V_{\fk X}^{\leq 1}$ to $0$. Using this description of $\uplambda$, it is easy to see that the restriction of $\uplambda$ to the subsheaf $\svir_c$ has image $\cbf\otimes_\Cbb \Theta_{\mc C/\mc B}^{\otimes 2}\simeq \Theta_{\mc C/\mc B}^{\otimes 2}$, and that its kernel is $\scr V_{\fk X}^{\leq 0}=\id\otimes_\Cbb\scr O_{\mc C}\simeq\scr O_{\mc C}$. Thus, we obtain an exact sequence
\begin{align}
0\rightarrow\scr O_{\mc C} \rightarrow\svir_c\xrightarrow{\uplambda}\Theta_{\mc C/\mc B}^{\otimes 2}\rightarrow 0.\label{eq168}
\end{align}
If we choose $U\subset\mc C-\Sigma$ and $\eta\in\scr O(U)$ holomorphic on each fiber,  then
\begin{gather*}
\uplambda:\quad \mc U_\varrho(\eta)^{-1}\cbf\mapsto \partial_\eta^2,\qquad \id\mapsto 0.
\end{gather*}
By tensoring with $\omega_{\mc C/\mc B}$, we get an exact sequence
\begin{gather}
0\rightarrow\omega_{\mc C/\mc B} \rightarrow\svir_c\otimes \omega_{\mc C/\mc B}\xrightarrow{\uplambda} \Theta_{\mc C/\mc B}\rightarrow 0\label{eq154}
\end{gather}
whose local expression outside $\Sigma$ is
\begin{gather}
\uplambda:\quad \mc U_\varrho(\eta)^{-1}\cbf ~d\eta\mapsto \partial_\eta,\qquad \id ~d\eta\mapsto 0.\label{eq155}
\end{gather}

\section{Lie derivatives}\label{lb50}

Fix a family of compact Riemann surfaces $\fk X=(\pi:\mc C\rightarrow\mc B)$. From Proposition \ref{lb56}, we see that the sheaf of VOA $\scr V_{\fk X}$ can be viewed as a twisted version of a direct sum of $\Theta_{\mc C/\mc B}^{\otimes n}$. It is well known that Lie derivatives can be defined for sections of $\Theta_{\mc C/\mc B}^{\otimes n}$ (whose restriction to $\Theta_{\mc C/\mc B}$ is given by the usual Lie bracket of vector fields). In this section, we define Lie derivatives for sections of $\scr V_{\fk X}$. The results of this section can be generalized easily to sections of $\scr V_{\fk X}|_{\mc C-\Sigma}$ when $\fk X$ is a family of complex curves.

Let $\varphi:U\rightarrow V$ be a biholomorphic map where $U,V$ are open subsets of $\mc C$ and $\varphi(U)=V$. We assume that $\varphi$ preserves fibers, i.e. $\varphi(U_{\pi(p)})=V_{\pi\circ\varphi(p)}$ for each $p\in U$. Then we have an \index{zz@$\varphi_*,(\eta_i)_*,(\eta_i,\pi)_*$} equivalence
\begin{gather*}
\varphi_*:\scr O_U\rightarrow\scr O_V,\qquad f\mapsto f\circ\varphi^{-1},
\end{gather*}
which makes each $\scr O_V$-module also an $\scr O_U$-module. $\varphi_*$ can be extended to
\begin{gather*}
\varphi_*:\Vbb^{\leq n}\otimes_{\Cbb}\scr O_U\xrightarrow{\simeq}\Vbb^{\leq n}\otimes_{\Cbb}\scr O_V,\qquad v\mapsto v\circ\varphi^{-1}.
\end{gather*}
Choose any $\eta\in\scr O(V)$ univalent on each fiber. Then we have a similar equivalence $(\eta,\pi)_*:\scr O_V\xrightarrow{\simeq}\scr O_{(\eta,\pi)(V)}$. Recall $\mc U_\varrho(\eta):\scr V_{\fk X}^{\leq n}|_V\xrightarrow{\simeq}\Vbb^{\leq n}\otimes_{\Cbb}\scr O_V$. Define an isomorphism \index{UV@$\mc U_\varrho(\eta),\mc U_\varrho(\varphi),\mc V_\varrho(\eta),\mc V_\varrho(\varphi)$}
\begin{gather}
\mc V_\varrho(\varphi):\scr V_{\fk X}^{\leq n}|_U\xrightarrow{\simeq} \scr V_{\fk X}^{\leq n}|_V,\nonumber\\
\mc U_\varrho(\eta)\mc V_\varrho(\varphi)=\varphi_*\cdot \mc U_\varrho(\eta\circ\varphi),\label{eq86}
\end{gather}
noting that $\mc U_\varrho(\eta\circ\varphi):\scr V_{\fk X}^{\leq n}|_U\rightarrow\Vbb^{\leq n}\otimes_{\Cbb}\scr O_U$. 

The definition of $\mc V_\varrho(\varphi)$ is independent of the choice of $\eta$. Indeed, if $\mu\in\scr O(V)$ is also univalent on each fiber, then, using \eqref{eq63}, it is not hard to show
\begin{align*}
\mc U(\varrho(\eta\circ\varphi|\mu\circ\varphi))=\varphi_*^{-1}\cdot\mc U(\varrho(\eta|\mu))\cdot\varphi_*,
\end{align*}
and hence equivalently that
\begin{align}
\mc U_\varrho(\eta\circ\varphi)\mc U_\varrho(\mu\circ\varphi)^{-1}=\varphi_*^{-1}\cdot\mc U_\varrho(\eta)\mc U_\varrho(\mu)^{-1}\cdot\varphi_*,
\end{align}
The independence follows easily from the above relation and \eqref{eq86}. Moreover, using the definition \eqref{eq86}, it is also not hard to show
\begin{align}
\mc V_\varrho(\psi\circ\varphi)=\mc V_\varrho(\psi)\mc V_\varrho(\varphi)
\end{align}
where $\psi:V\rightarrow W$ is another such fiber-preserving biholomorphic map. In particular, we have $\mc V_\varrho(\varphi^{-1})=\mc V_\varrho(\varphi)^{-1}$.

Let now $W$ be an open subset of $\mc C$, and let $V$ be a precompact open subset of $W$ whose closure is also in $W$. Note that since $\pi$ is (clearly) an open map, $\pi(W)$ is open. Recall the short exact sequence \eqref{eq81}:
\begin{align*}
0\rightarrow \Theta_{\mc C/\mc B}\rightarrow \Theta_{\mc C}\xrightarrow{d\pi}\pi^*\Theta_{\mc B}\rightarrow 0.
\end{align*}
Let $\fk x\in\Theta_{\mc C}(W)$ such that $d\pi(\fk x)$ equals $\pi^*(\fk y)$ for some $\fk y\in\Theta_{\mc B}(\pi(W))$. In other words, $\fk x$ is a vector field on $W$ whose projection to $\mc B$ depends only on the points of $\mc B$. Suppose that $W$ is small enough so that  we can choose $\eta\in\scr O(W)$ univalent on fibers, and choose coordinates  $\tau_\blt=(\tau_1,\dots,\tau_m)$ of $\pi(W)$. Denote $\tau_\blt\circ\pi$ also by $\tau_\blt$ for simplicity. Then $(\eta,\tau_\blt):V\rightarrow\Cbb\times\Cbb^m$ is a coordinate of $V$, and $\fk x$ takes the form
\begin{align}
\fk x=h(\eta,\tau_\blt)\partial_\eta+\sum_{j=1}^m g_j(\tau_\blt)\partial_{\tau_j}\label{eq85}
\end{align}
for some holomorphic functions $h$ on $(\eta,\tau_\blt)(W)$ and $g_1,\dots,g_m$ on $\tau_\blt(W)$.  Choose $\varphi^{\fk x}\in\scr O_{T\times V}(T\times V),(\zeta,p)\mapsto \varphi^{\fk x}_\zeta(p)$, where $T$ is an open subset of $\Cbb$ containing $0$, and  the following conditions are satisfied for any $p\in V$:
\begin{gather}
\varphi^{\fk x}_0(p)=p.\label{eq83}\\
\partial_\zeta\varphi^{\fk x}_\zeta(p)\big|_{\zeta=0}=\fk x(p).\label{eq84}
\end{gather}
The second condition is equivalent to that for any (local) section $f$ of $\scr O_V$,
\begin{gather}
\partial_\zeta (f\circ\varphi^{\fk x}_\zeta)\Big|_{\zeta=0}=\fk x f.\label{eq88}
\end{gather}
The first conditions implies that $\mc V_\varrho(\varphi^{\fk x}_0)$ is the identity map on $\scr V^{\leq n}_{\fk X}|_V$.

\begin{df}
For any $v\in\scr V_{\fk X}^{\leq n}(W)$ and $\fk x$ as above, we define $\mc L_{\fk x}v\in\scr V_{\fk X}^{\leq n}(W)$ as follows. Choose any $V\subset W$ whose closure is compact and contained in $W$,  and choose $\varphi^{\fk x}$ as above. Then
\begin{align}
\mc L_{\fk x}v\big|_V=\lim_{\zeta\rightarrow 0}\frac{\mc V_\varrho(\varphi^{\fk x}_\zeta)^{-1}\big(v\big|_{\varphi^{\fk x}_\zeta(V)}\big)-v\big|_V}{\zeta}.\label{eq89}
\end{align}
\end{df}

We now give an explicit formula of $\mc L_{\fk x}v$, which shows in particular that the above definition is independent of the choice of $\varphi^{\fk x}$ satisfying \eqref{eq83} and \eqref{eq84}. If $u$ is a section of $\Vbb^{\leq n}\otimes_{\Cbb}\scr O_{W}$, we say that a section of $\scr V_{\fk X}^{\leq n}|_W$ equals $u$ in the $\eta$-coordinate if this section is $\mc U_\varrho(\eta)^{-1}u$.

\begin{thm}\label{lb53}
Suppose that $\eta\in\scr O(W)$ is univalent on each fiber of $W$, $\fk x$ takes the form \eqref{eq85}, and $v\in\scr V_{\fk X}^{\leq n}(W)$ equals $u\in\Vbb^{\leq n}\otimes_{\Cbb}\scr O(W)$ in the $\eta$-coordinate. Then in the $\eta$-coordinate,  $\mc L_{\fk x}v$ equals
\begin{align}
h(\eta,\tau_\blt)\partial_\eta u+\sum_{j=1}^m g_j(\tau_\blt)\partial_{\tau_j}u-\sum_{k\geq 1}\frac 1{k!}\partial_\eta^k h(\eta,\tau_\blt)L_{k-1}u.\label{eq87}
\end{align}
\end{thm}

\begin{proof}
Choose $V,\varphi^{\fk x}$ as above. We have $v=\mc U_\varrho(\eta)^{-1}u$. Then, in the $\eta$-coordinate, $\mc V_\varrho(\varphi^{\fk x}_\zeta)^{-1}\big(v\big|_{\varphi^{\fk x}_\zeta(V)}\big)$ equals
\begin{align*}
\mc U_\varrho(\eta)\mc V_\varrho(\varphi^{\fk x}_\zeta)^{-1}\big(v\big|_{\varphi^{\fk x}_\zeta(V)}\big)=\mc U_\varrho(\eta)\mc V_\varrho(\varphi^{\fk x}_\zeta)^{-1}\mc U_\varrho(\eta)^{-1}\big(u\big|_{\varphi^{\fk x}_\zeta(V)}\big),
\end{align*}
which by \eqref{eq86}  equals
\begin{align*}
&\mc U_\varrho(\eta)\mc U_\varrho(\eta\circ\varphi^{\fk x}_\zeta)^{-1}(\varphi^{\fk x}_\zeta)_*^{-1}\big(u\big|_{\varphi^{\fk x}_\zeta(V)}\big)\\
=&\mc U(\varrho(\eta|\eta\circ\varphi^\xk_\zeta))\big(u\circ\varphi^{\fk x}_\zeta\big)\big|_V
\end{align*}
It is easy to see that the derivative over  $\zeta$ of the above expression at $\zeta=0$ equals \eqref{eq87}. Indeed, the first two terms of \eqref{eq87} come from the derivative of $u\circ\varphi^{\fk x}_\zeta$. The last term comes from the derivative of $\mc U(\varrho(\eta\circ\varphi_\zeta^{\fk x}|\eta)^{-1})$. Identify $V$ with $(\eta,\tau_\blt)(V)\subset\Cbb\times\Cbb^m$ via $(\eta,\tau_\blt)$. Then by \eqref{eq62},  
\begin{align*}
\varrho(\eta|\eta\circ\varphi^\xk_\zeta)_p(z)=\eta\circ\varphi_{-\zeta}^\xk(z+\eta\circ\varphi_\zeta^\xk(p),\tau_\blt(p))-\eta(p).
\end{align*}
So, as $\varphi^\xk_0=\id$, using \eqref{eq88}, we get
\begin{align*}
&\partial_\zeta\varrho(\eta|\eta\circ\varphi^\xk_\zeta)_p(z)\big|_{\zeta=0}=\partial_\zeta~\eta\circ\varphi_{-\zeta}^\xk(z+\eta(p),\tau_\blt(p))\big|_{\zeta=0}+\partial_\zeta~(z+\eta\circ\varphi_\zeta^\xk(p))\big|_{\zeta=0}\\
=&-h(z+\eta(p),\tau_\blt(p))+h(p),
\end{align*}
noting that $p=(\eta(p),\tau_\blt(p))$. So
\begin{align*}
\partial_\zeta\varrho(\eta|\eta\circ\varphi^\xk_\zeta)^{(k)}_p(0)\big|_{\zeta=0}=-\partial_\eta^k h(p).
\end{align*}
Thus, by Lemma \ref{lb29}, we have
\begin{align*}
\partial_\zeta\mc U(\varrho(\eta\circ\varphi_\zeta^{\fk x}|\eta)^{-1})\big|_{\zeta=0}=-\sum_{k\geq 1}\frac 1{k!}\partial_\eta^k h(\eta,\tau_\blt)L_{k-1}.
\end{align*}
\end{proof}

\begin{rem}\label{lb54}
Note that if $\varphi:U\rightarrow V$  is biholomorphic and fiber-preserving, then $\varphi$ maps any (complex and holomorphic) path in each $U_b$ to one in $V_{\varphi(b)}$. Thus $\varphi_*\equiv d\varphi$ maps tangent vectors of $U_b$ to those of $V_{\varphi(b)}$. Thus one can define an isomorphism  $d\varphi:\Theta_{\mc C/\mc B}|_U\xrightarrow{\simeq}\Theta_{\mc C/\mc B}|_V$, and hence 
\begin{align*}
\varphi_*\equiv d\varphi:\Theta_{\mc C/\mc B}^{\otimes n}|_U\xrightarrow{\simeq}\Theta_{\mc C/\mc B}^{\otimes n}|_V
\end{align*}
for each $n\in\Zbb$. One can thus use \eqref{eq89} (with $\mc V_\varrho(\varphi_\zeta^{\fk x})$ replaced by $d\varphi_\zeta^{\fk x}$)  to define the Lie derivatives  on $\Theta_{\mc C/\mc B}^{\otimes n}$. When $n\in\Nbb$, it is easy to see that the Lie derivatives on $\scr V_{\fk X}^{\leq n}/\scr V_{\fk X}^{\leq n-1}$ is the same as those on $\Vbb(n)\otimes_{\Cbb}\Theta_{\mc C/\mc B}^{\otimes n}$. One can also define Lie derivatives on $\scr V_{\fk X}^{\leq n}\otimes\omega_{\mc C/\mc B}$ (recall that $\omega_{\mc C/\mc B}=\Theta_{\mc C/\mc B}^{-1}$), and the formula of Lie derivatives is exactly the same as \eqref{eq87}, except that  $L_0$ should be replaced by $L_0-1$. (In other words, there is an extra term $\partial_\eta h(\eta,\tau_\blt)u$ contributed by the Lie derivatives on $\omega_{\mc C/\mc B}$.) In the next chapter, we will use the Lie derivatives on $\scr V_{\fk X}^{\leq n}\otimes\omega_{\mc C/\mc B}$ to define connections on sheaves of conformal blocks.
\end{rem}

\chapter{Sheaves of conformal blocks}\label{lb91}

\section{Spaces of conformal blocks}\label{lb34}

Let $\Vbb$ be always a (CFT-type) VOA. Let $\fk X=(C;x_1,\dots,x_N;\eta_1,\dots,\eta_N)$ be an $N$-pointed complex curve with local coordinates. Recall that if $\Wbb$ is a $\Vbb$-module, the vertex operator $Y_\Wbb$ can be regarded as a $\Cbb((z))$-module homomorphism $\Vbb((z))\otimes_{\Cbb}\Wbb\rightarrow\Wbb((z))$ sending each $v\otimes w$ to $Y_{\Wbb}(v,z)w$. (See \eqref{eq228}; here $v\in V$ is considered as the constant section in $\Vbb((z))$.) Let $\Wbb_1,\Wbb_2,\dots,\Wbb_N$ be $\Vbb$-modules. Set \index{Ww@$\Wbb_\blt,w_\blt$} $\Wbb_\blt=\Wbb_1\otimes \Wbb_2\otimes\cdots\otimes \Wbb_N$.

\begin{cv}
By $w\in\Wbb_\blt$, we mean a vector of $\Wbb_1\otimes\cdots\otimes \Wbb_N$. \index{w@$w_\blt,\Wbb_\blt$} By $w_\blt\in \Wbb_\blt$, we mean a vector of the form $w_1\otimes w_2\otimes\cdots\otimes w_N$, where $w_1\in \Wbb_1,\dots,w_N\in\Wbb_N$.
\end{cv}

Recall that $S_{\fk X}$ is the divisor $x_1+x_2+\cdots+x_N$. For each $1\leq i\leq N$, we choose a neighborhood $U_i$ of $x_i$ on which $\eta_i$ is defined. Then, by tensoring with the identity map of $\omega_{U_i}$,  the map \eqref{eq90} induces naturally an $\scr O_{U_i}$-module isomorphism 
\begin{align*}
\mc U_\varrho(\eta_i):\scr V_C|_{U_i}\otimes\omega_{U_i}(\blt S_{\fk X})\xrightarrow{\simeq}\Vbb\otimes_{\Cbb}\omega_{U_i}(\blt S_{\fk X}).
\end{align*}
Let $(\eta_i)_*:\omega_{U_i}\xrightarrow{\simeq}\omega_{\eta_i(U_i)}$ \index{zz@$\varphi_*,(\eta_i)_*,(\eta_i,\pi)_*$} be the pushforward of differentials, i.e. $(\eta_i)_*=(\eta_i^{-1})^*$. It can be extended by linearity to $(\eta_i)_*:\Vbb\otimes_{\Cbb}\omega_{U_i}(\blt S_{\fk X})\xrightarrow{\simeq}\Vbb\otimes_{\Cbb}\omega_{\eta_i(U_i)}(\blt 0)$. Let
\begin{gather}
\mc V_\varrho(\eta_i):\scr V_C|_{U_i}\otimes\omega_{U_i}(\blt S_{\fk X})\xrightarrow{\simeq}\Vbb\otimes_{\Cbb}\omega_{\eta_i(U_i)}(\blt 0)\nonumber\\
\mc V_{\varrho}(\eta_i)=(\eta_i)_*\mc U_{\varrho}(\eta_i).\label{eq233}
\end{gather}
In the case that $U_i$ and $\eta_i(U_i)$ are identified by $\eta_i$, we have $\mc U_\varrho(\eta_i)=\mc V_\varrho(\eta_i)$. Let $z$ be the standard coordinate of $\Cbb$. If $v$ is a section of $\scr V_C|_{U_i}\otimes\omega_{U_i}(\blt S_{\fk X})$ defined near $x_i$, we define a linear action of $v$ on $\Wbb_i$ such that if $w_i\in\Wbb_i$, \index{vw@$v\cdot w_i,v\cdot w_\blt$} then
\begin{align}
\boxed{~~v\cdot w_i=\Res_{z=0}Y_{\Wbb_i}(\mc V_\varrho(\eta_i)v,z)w_i~~}\label{eq230}
\end{align}

Define a linear action of $H^0(C,\scr V_C\otimes\omega_C(\blt S_{\fk X}))$ on $\Wbb_\blt$ as follows.  
If $v\in H^0(C,\scr V_C\otimes\omega_C(\blt S_{\fk X}))$, the action of $v$ on any $w_\blt$ is
\begin{align}
v\cdot w_\blt=\sum_{i=1}^N w_1\otimes w_2\otimes\cdots \otimes (v|_{U_i})\cdot w_i \otimes\cdots\otimes w_N.\label{eq231}
\end{align}
We now define a \textbf{space of covacua} 
\begin{gather}
\scr T_{\fk X}(\Wbb_\blt)=\frac{\Wbb_\blt}{H^0(C,\scr V_C\otimes\omega_C(\blt S_{\fk X}))\cdot \Wbb_\blt}
\end{gather}
whose dual vector space is denoted by $\scr T_{\fk X}^*(\Wbb_\blt)$ \index{T@$\scr T_{\fk X}(\Wbb_\blt),\scr T_{\fk X}^*(\Wbb_\blt)$} and called a \textbf{space of conformal blocks} or \textbf{space of vacua}.

A conformal block $\upphi\in \scr T_{\fk X}^*(\Wbb_\blt)$ is understood as a chiral correlation function in physics. According to the above definition, $\upphi$ as a linear functional on $\Wbb_\blt$ should vanish on the subspace $H^0(C,\scr V_C\otimes\omega_C(\blt S_{\fk X}))\cdot \Wbb_\blt$. Such condition is similar to the Jacobi identity for VOAs. We now  interpret this condition in a similar fashion as Theorem \ref{lb19}. 

For each $\upphi\in \Wbb_\blt^*$ and $x_i$,  if $w_\blt\in \Wbb_\blt$, we define
\begin{align*}
\wr\upphi_{x_i}(w_\blt)\in \Vbb^*[[z^{\pm1}]]
\end{align*}
whose evaluation on each $v\in \Vbb$, written as $\wr\upphi_{x_i}(v,w_\blt)$, equals
\begin{align}
\wr\upphi_{x_i}(v,w_\blt)=\upphi(w_1\otimes w_2\otimes\cdots\otimes Y_{\Wbb_i}(v,z)w_i\otimes\cdots\otimes w_N).\label{eq125}
\end{align}
By the lower truncation property, the above expression is an element of $\Cbb((z))$. Also, the above expression makes sense when $v$ is a section of $\Vbb\otimes_{\Cbb}\scr O_{\eta_i(U_i)}$ defined near $z=0$. By linearity, we can define $\wr\upphi_{x_i}(v,w)$ for any $w\in\Wbb_\blt$. The following theorem is also true when $C$ is nodal; however, we will only be interested in the smooth case. We \index{CS@$C-\SX,\mc C-\SX,E-\SX$} understand
\begin{align*}
C-\SX=C-\{x_1,\dots,x_N\}.
\end{align*}

\begin{thm}\label{lb32}
Assume that $C$ is smooth. Let $\upphi\in \Wbb_\blt^*$ . Then the following are equivalent.
\begin{enumerate}[label=(\alph*)]
\item $\upphi$ is an element of $\scr T_{\fk X}^*(\Wbb_\blt)$. 
\item For each $w\in\Wbb_\blt$, there exists a (necessarily unique) element \index{zz@$\wr\upphi$}
\begin{align}
\wr\upphi(w)\in H^0(C-S_{\fk X},\scr V_C^*)\label{eq119}
\end{align}
such that for each $1\leq i\leq N$, if we identify $U_i\simeq \eta_i(U_i)$ via $\eta_i$ and identify $\scr V_C|_{U_i}\simeq\Vbb\otimes_\Cbb\scr O_{U_i}$ via $\mc U_\varrho(\eta_i)$, then the evaluation of $\wr\upphi(w)$ with any $v\in\scr V_C(U_i)$ (restricted to $U_i-x_i$), written as $\wr\upphi(v,w)$, is
\begin{align}
\wr\upphi(v,w)=\wr\upphi_{x_i}(v,w).\label{eq229}
\end{align}
\end{enumerate}	
\end{thm}

\begin{proof}
Choose any $v\in H^0(C,\scr V_C\otimes\omega_C(\blt S_{\fk X}))$. Assume $w=w_\blt=w_1\otimes\cdots\otimes w_N$. Then 
\begin{align}
\sum_{i=1}^N\Res_{z=0}~ {\wr\upphi}_{x_i}\big(v,w_\blt\big)=\upphi(v\cdot w_\blt).\label{eq91}
\end{align}
Suppose that (b) is true. Then  $\wr\upphi(v,w_\blt)$ is an element of $H^0(C,\omega_C(\blt S_{\fk X}))$, and the left hand side of the above expression equals $\sum_{i=1}^N \Res_{x_i}\wr\upphi(v,w_\blt)$, which, by residue theorem, equals $0$. Thus $\upphi$ vanishes on $v\cdot w_\blt$. This proves (a).

We now assume (a). Choose any $n\in\Nbb$, and restrict $\wr\upphi_{x_i}(w_\blt)$ to $\Vbb^{\leq n}$ (or $\Vbb^{\leq n}\otimes_{\Cbb}\scr O_{U_i}$ when considering non-constant sections), which gives
\begin{align*}
\wr\upphi_{x_i}^{\leq n}(w_\blt)\in (\Vbb^{\leq n})^*((z)).
\end{align*}
Then for any $w_\blt\in\Wbb_\blt$ and $v\in H^0(C,\scr V_C^{\leq n}\otimes\omega_C(\blt S_{\fk X}))$, since $\upphi$ vanishes on $H^0(C,\scr V_C^{\leq n}\otimes\omega_C(\blt S_{\fk X}))\cdot \Wbb_\blt$, we have $\sum_{i=1}^N \Res_{z=0}\wr\upphi_{x_i}(v,w_\blt)$ equals $0$ by \eqref{eq91}. Thus, by strong residue theorem, there exists
\begin{align*}
\wr\upphi^{\leq n}(w_\blt)\in H^0(C,(\scr V_C^{\leq n})^*\otimes\omega_C(\blt S_{\fk X}))
\end{align*}
whose series expansion near each $x_i$ is $\wr\upphi_{x_i}^{\leq n}$. Equivalently, \eqref{eq229} holds for any $i$ and any $v\in\scr V_C^{\leq n}(U_i)$. It is clear that $\wr\upphi^{\leq n'}(w_\blt)$ restricts to $\wr\upphi^{\leq n}(w_\blt)$ when $n'\geq n$. One can thus define $\wr\upphi(w_\blt)$ to be the projective limit of $\wr\upphi^{\leq n}(w_\blt)$ over $n$.
\end{proof}

Our next goal is to give a coordinate-free definition of conformal blocks. Let $\fk X=(C;x_1,\dots,x_N)$ be an $N$-pointed complex curve, and let $\Wbb_1,\dots,\Wbb_N$ be $\Vbb$-modules. Define a vector space $\scr W_{\fk X}(\Wbb_\blt)$ isomorphic to $\Wbb_\blt$ as follows. $\scr W_{\fk X}(\Wbb_\blt)$ is a (infinite rank) vector bundle on the $0$-dimensional manifold $\{C\}$ (consider as the base manifold of the family $C\rightarrow \{C\}$). For any choice of local coordinates $\eta_\blt=(\eta_1,\dots,\eta_N)$ of $x_1,\dots,x_N$ respectively, we have a trivialization \index{U@$\mc U(\rho),\mc U(\eta_\blt)$}
\begin{align}
\mc U(\eta_\blt):\scr W_{\fk X}(\Wbb_\blt)\xrightarrow{\simeq} \Wbb_\blt  \label{eq121}
\end{align}
such that if $\mu_\blt$ is another set of local coordinates, then
\begin{align}
\mc U(\eta_\blt)\mc U(\mu_\blt)^{-1}=&\mc U(\eta_\blt\circ\mu_\blt^{-1})\nonumber\\
:=&\mc U(\eta_1\circ\mu_1^{-1})\otimes \mc U(\eta_2\circ\mu_2^{-1})\otimes\cdots\otimes \mc U(\eta_N\circ\mu_N^{-1}).\label{eq123}
\end{align}
If $v\in H^0(C,\scr V_C\otimes\omega_C(\blt S_{\fk X}))$ and $w\in\scr W_{\fk X}(\Wbb_\blt)$, we set
\begin{align}
v\cdot w=\mc U(\eta_\blt)^{-1}\cdot v\cdot\mc U(\eta_\blt)\cdot w,\label{eq232}
\end{align}
where the action of $v$ on $\mc U(\eta_\blt)w$ (which depends on $\eta_\blt$) is defined by \eqref{eq230} and \eqref{eq231}.

\begin{lm}\label{lb125}
The definition of $v\cdot w$ in \eqref{eq232} is independent of the choice of $\eta_\blt$.
\end{lm}

\begin{proof}
We prove this lemma for the case $N=1$. The general cases can be proved in a similar way. Choose local coordinates $\eta,\mu$ at $x=x_1$ defined on a neighborhood $U$. We identify $U$ with $\mu(U)$ via $\mu$. So $\mu$  is identified with the standard coordinate $\id_\Cbb$ of $\Cbb$, and $\eta\in\Gbb$. (We will denote by $z$ the standard complex variable of $\Cbb$.) Also, identify $\scr W_{\fk X}(\Wbb)$ (where $\Wbb=\Wbb_1=\Wbb_\blt$) with $\Wbb$ via $\mc U(\mu)$. So $\mc U(\mu)=\mc U(\id_\Cbb)=\id$.  Choose any $w\in\Wbb$ , and assume that $v$ is of the form
\begin{gather*}
\mc V_\varrho(\mu)v=\mc U_\varrho(\mu)v=u(z)dz
\end{gather*}
where $u=u(z)\in\Vbb\otimes_\Cbb\scr O_\Cbb(\blt 0)(U)$. Then, by \eqref{eq55}, we have
\begin{gather*}
\mc U_\varrho(\eta)v=\mc U(\varrho(\eta|\id_\Cbb))u(z)dz=\mc U(\varrho(\eta|\id_\Cbb)_z)u(z)dz.
\end{gather*}
By \eqref{eq233}, we have
\begin{gather*}
\qquad \mc V_\varrho(\eta)v=\eta_*\big(\mc U(\varrho(\eta|\id_\Cbb)_z)u(z)dz \big)=\mc U(\varrho(\eta|\id_\Cbb)_{\eta^{-1}(z)})\cdot u(\eta^{-1}(z))\cdot d(\eta^{-1}(z)).
\end{gather*}
We calculate that
\begin{align*}
&\mc U(\eta)^{-1}\cdot v\cdot\mc U(\eta)\cdot w=\Res_{z=0}~\mc U(\eta)^{-1}Y_\Wbb\big(\mc V_\varrho(\eta)v,z\big)\mc U(\eta)w\\
=&\Res_{z=0}~\mc U(\eta)^{-1}Y_\Wbb\big(\mc U(\varrho(\eta|\id_\Cbb)_{\eta^{-1}(z)})\cdot u(\eta^{-1}(z)),z\big)\mc U(\eta)w\cdot d(\eta^{-1}(z)),
\end{align*}
which by Theorem \ref{lb31} equals
\begin{align*}
&\Res_{z=0}~Y_\Wbb\big( u(\eta^{-1}(z)),\eta^{-1}(z)\big)w\cdot d(\eta^{-1}(z))\\
=&\Res_{z=0}~Y_\Wbb\big( u(z),z\big)w\cdot dz=\mc U(\mu)^{-1}\cdot v\cdot\mc U(\mu)\cdot w.
\end{align*}
The proof is complete.
\end{proof}

Thus, we have a coordinate-independent linear action of $H^0(C,\scr V_C\otimes\omega_C(\blt S_{\fk X}))$ on $\scr W_{\fk X}(\Wbb_\blt)$. Then the space of conformal blocks $\scr T_{\fk X}^*(\Wbb_\blt)$ \index{T@$\scr T_{\fk X}(\Wbb_\blt),\scr T_{\fk X}^*(\Wbb_\blt)$} is the dual space of the space of covacua
\begin{gather}
\scr T_{\fk X}(\Wbb_\blt)=\frac{\scr W_{\fk X}(\Wbb_\blt)}{H^0(C,\scr V_C\otimes\omega_C(\blt S_{\fk X}))\cdot \scr W_{\fk X}(\Wbb_\blt)}.
\end{gather}

\begin{rem}\label{lb40}
We remark that Theorem \ref{lb32} still holds in this general setting. Indeed, by Theorem \ref{lb31}, for each $w\in\scr W_{\fk X}(\Wbb_\blt)$ one can define $\wr\upphi_{x_i}(w)$ whose expression is covariant under the change of local coordinates. Then $\upphi$ is a conformal block if and only if for each $w$, all these $\wr\upphi_{x_i}(w)$ can be extended to  a (necessarily unique) $\wr\upphi(w)$ which is independent of the choice of local coordinates.
\end{rem}

\begin{eg}
Let $C=\Pbb^1$. Let $\zeta$ be the standard coordinate of $\Cbb=\Pbb^1-\{\infty\}$, and let $\varpi=\zeta^{-1}$ be a coordinate of $\infty$ defined on $\Pbb^1-\{0\}$. Let $\fk X=(\Pbb^1;0,\infty;\zeta,\varpi)$. So the divisor $S_{\fk X}$ is $0+\infty$, and hence $\Pbb^1-S_{\fk X}=\Cbb-\{0\}=\Cbb^\times$. Choose $V$-modules $\Wbb_1,\Wbb_2$ associated to $0,\infty$ respectively. Let $\Wbb_\blt=\Wbb_1\otimes\Wbb_2$. We shall show that there is an isomorphism $\Hom_V(\Wbb_1,\Wbb_2')\simeq\scr T_{\fk X}^*(\Wbb_\blt)$. 

Define a linear map $\Hom_V(\Wbb_1,\Wbb_2')\rightarrow\scr T_{\fk X}^*(\Wbb_\blt)$ as follows. If $T\in \Hom_V(\Wbb_1,\Wbb_2')$, then the corresponding conformal block $\upphi_T$, as a linear functional on $\Wbb_\blt$, is defined by
\begin{align}
\upphi_T(w_\blt)=\bk{Tw_1,w_2}
\end{align}
for each $w_\blt:=w_1\otimes w_2\in \Wbb_1\otimes \Wbb_2$. We now verify that $\upphi_T$ is a conformal block by verifying (b) of Theorem \ref{lb32}. $\scr V_{\Pbb^1}|_{\Cbb^\times}$ is  generated by all $\mc U_\varrho(\zeta)^{-1}v$ where $v\in\Vbb$. Moreover, it is easy to define $\wr\upphi_T(w_\blt)\in H^0(\Cbb^\times,\scr V_{\Pbb^1}^*)$ such that for any $v$,
\begin{align*}
\wr\upphi_T\big(\mc U_\varrho(\zeta)^{-1}v,w_\blt\big)=\bk{TY_{\Wbb_1}(v,\zeta)w_1,w_2},
\end{align*}
considering $\mc U_\varrho(\zeta)^{-1}v$ as a section on $\Cbb^\times$. The series expansion of $\wr\upphi_T(w_\blt)$ near $0$ is clearly $\wr\upphi_{T,0}(w_\blt)$. Near $\infty$, we have
\begin{align*}
&\wr\upphi_{T,\infty}\big(\mc U_\varrho(\zeta)^{-1}v,w_\blt\big)=\upphi_{T,\infty}\big(w_1\otimes Y_{\Wbb_2}(\mc U_\varrho(\varpi)\mc U_\varrho(\zeta)^{-1}v,\varpi) w_2\big)\\
=&\bk{Tw_1,Y_{\Wbb_2}(\mc U_\varrho(\varpi)\mc U_\varrho(\zeta)^{-1}v,\varpi) w_2}=\bk{Tw_1,Y_{\Wbb_2}(\mc U(\varrho(\varpi|\zeta))v,\zeta^{-1}) w_2}.
\end{align*}
According to (the proof of) Lemma \ref{lb123} (note that we have $\xi=\zeta$ and $q=\xi\varpi=1$), we have
\begin{align*}
\mc U(\varrho(\varpi|\zeta))=\mc U(\upgamma_\zeta),
\end{align*}
where we recall by \eqref{eq73} that $\mc U(\upgamma_\zeta)=e^{\zeta L_1}(-\zeta^{-2})^{L_0}$ when acting on $V$. Thus, by \eqref{eq92} and the above calculation, we have
\begin{align*}
&\wr\upphi_{T,\infty}\big(\mc U_\varrho(\zeta)^{-1}v,w_\blt\big)=\bk{Tw_1,Y_{\Wbb_2}(\mc U(\upgamma_\zeta)v,\zeta^{-1}) w_2}\\
=&\bk{Tw_1,Y_{\Wbb_2}(v,\zeta)^\tr w_2}=\bk{Y_{\Wbb_2}(v,\zeta)Tw_1, w_2}=\bk{TY_{\Wbb_1}(v,\zeta)w_1,w_2}.
\end{align*}
Thus, the series expansion of $\wr\upphi_T(w_\blt)$ near $\infty$ is $\wr\upphi_{T,\infty}(w_\blt)$. So $\upphi_T$ is a conformal block.

It is obvious that the map $T\mapsto\upphi_T$ is injective. To show the surjectivity, we choose any conformal block $\upphi$. Define a linear map $T:W_1\rightarrow W_2^*$ satisfying  $\upphi(w_\blt)=\bk{Tw_1,w_2}$ for any $w_1,w_2$. Then for any $v\in \Vbb$, by the above calculation, we see that   $\bk{TY_{\Wbb_1}(v,\zeta)w_1,w_2}$ and $\bk{Tw_1,Y_{\Wbb_2}(v,\zeta)^\tr w_2}$ (which are elements of $\Cbb((\zeta))$ and $\Cbb((\zeta^{-1}))$ respectively) are expansions near $\zeta=0$ and $\zeta^{-1}=0$ of the same function in $H^0(\Pbb^1,\scr O_{\Pbb^1}(\blt 0+\blt\infty))$. Thus, they are equal to a polynomial of $\zeta$. Take $v$ to be the conformal vector $\cbf$. Then one sees that $T$ intertwines the actions of $L_0$. This shows that $T$ has image in $\Wbb_2'$, since, in particular, $T$ maps $L_0$-homogeneous vectors of $\Wbb_1$ to those of $\Wbb_2'$ with the same weights. It is now obvious that $T$ intertwines the actions of $\Vbb$, and that $\upphi=\upphi_T$.
\end{eg}

\begin{eg}\label{lb44}
Choose any $z\in\Cbb^\times$.  Let again $\zeta$ be the standard coordinate of $\Cbb$. Then $\zeta,\zeta-z,\zeta^{-1}$ are local coordinates of $0,z,\infty$ respectively. Set  $\fk X=(\Pbb^1;0,z,\infty;\zeta,\zeta-z,\zeta^{-1})$. Choose any $\Vbb$-module $\Wbb$, and set  $\Wbb_\blt=\Wbb\otimes\Vbb\otimes\Wbb'$. Then the vertex operator $Y_\Wbb$ can be viewd as an element $\upphi$ of $\scr T^*{\fk X}(\Wbb_\blt)$ by sending each $w_\blt=w\otimes v\otimes w'\in\Wbb\otimes\Vbb\otimes\Wbb'$ to the scalar
\begin{align*}
\upphi(w_\blt)=\bk{Y_\Wbb(v,z)w,w'}.
\end{align*}
One can define $\wr\upphi(w_\blt)\in H^0(\Cbb-\{z\},\scr V_{\Pbb^1}^*)$ such that for any $u\in\Vbb$, the section $\mc U_\rho(\zeta)^{-1}u$ on $\Cbb-\{z\}$ is sent to the function $f$ in Theorem \ref{lb19}. Using the argument in the previous example, it is not hard to check that $\wr\upphi_0(w_\blt),\wr\upphi_z(w_\blt),\wr\upphi_\infty(w_\blt)$ are the series expansions of $\wr\upphi(w_\blt)$ near $0,z,\infty$ respectively. This proves that $\upphi$ is a conformal block.
\end{eg}

\section{Sheaves of conformal blocks}

Let $\fk X=(\pi:\mc C\rightarrow\mc B;\sgm_1,\dots,\sgm_N;\eta_1,\dots,\eta_N)$ be a family of $N$-pointed complex curves with local coordinates. Let $z$ be the standard coordinate of $\Cbb$. Let $\scr O_{\mc B}((z))$ be the $\scr O_{\mc B}$-module associating to each open $V\subset\mc B$ the algebra $\scr O(V)((z))$. If $\Wbb$ is a $\Vbb$-module, we have a homomorphism of $\scr O_{\mc B}((z))$-modules \index{Y@$Y,Y_\Wbb$}
\begin{gather}
Y_\Wbb:\big(\Vbb\otimes_\Cbb\scr O_{\mc B}((z))\big)\otimes_{\scr O_{\mc B}} (\Wbb\otimes_\Cbb\scr O_{\mc B})\rightarrow \Wbb\otimes_\Cbb\scr O_{\mc B}((z)),\nonumber\\
v\otimes w\mapsto Y_{\Wbb}(v,z)w 
\end{gather}
where $v=v\otimes 1,w=w\otimes 1$ are constant sections. Note a section of $\Vbb\otimes_\Cbb\scr O_{\Cbb\times\mc B}$ on a neighborhood of $\{0\}\times V$ (where $V$ is an open subset of $\mc B$) can be regarded as an element of $\Vbb\otimes_\Cbb\scr O_{\mc B}(V)((z))$ by taking series expansion.

If $\Wbb_1,\dots,\Wbb_N$ are $\Vbb$-modules, we set $\Wbb_\blt=\Wbb_1\otimes \cdots\otimes \Wbb_N$ as usual, and define
\begin{align}
\scr W_{\fk X}(\Wbb_\blt)=\Wbb_\blt\otimes_{\Cbb}\scr O_{\mc B}.
\end{align}
Choose mutually disjoint neighborhoods $U_1,\dots,U_N$ of $\sgm_1(\mc B),\dots,\sgm_N(\mc B)$ on which $\eta_1,\dots,\eta_N$ are defined respectively. For each $i$ we have
\begin{gather*}
\mc U_\varrho(\eta_i):\scr V_{\fk X}\otimes\omega_{\mc C/\mc B}(\blt S_{\fk X})\big|_{U_i}\xrightarrow{\simeq}\Vbb\otimes_{\Cbb}\omega_{\mc C/\mc B}(\blt S_{\fk X})\big|_{U_i},
\end{gather*}
and the \index{zz@$\varphi_*,(\eta_i)_*,(\eta_i,\pi)_*$} pushforward
\begin{gather*}
(\eta_i,\pi)_*:\omega_{\mc C/\mc B}\big|_{U_i}\xrightarrow{\simeq}\omega_{(\eta_i,\pi)(U_i)/\mc B}
\end{gather*}
equaling $((\eta_i,\pi)^{-1})^*$, which, by tensoring with $\id_\Vbb$, gives rise to (identifying $\mc B$ with $\{0\}\times\mc B$)
\begin{gather*}
(\eta_i,\pi)_*:\Vbb\otimes_{\Cbb}\omega_{\mc C/\mc B}(\blt S_{\fk X})\big|_{U_i}\xrightarrow{\simeq}\Vbb\otimes_{\Cbb}\omega_{(\eta_i,\pi)(U_i)/\mc B}(\blt \mc B).
\end{gather*}
Define \index{UV@$\mc U_\varrho(\eta),\mc U_\varrho(\varphi),\mc V_\varrho(\eta),\mc V_\varrho(\varphi)$}
\begin{gather}
\mc V_\varrho(\eta_i):\scr V_{\fk X}\otimes\omega_{\mc C/\mc B}(\blt S_{\fk X})\big|_{U_i}\xrightarrow{\simeq}\Vbb\otimes_{\Cbb}\omega_{(\eta_i,\pi)(U_i)/\mc B}(\blt \mc B),\nonumber\\
\mc V_{\varrho}(\eta_i)=(\eta_i,\pi)_*\cdot \mc U_{\varrho}(\eta_i).
\end{gather}
If $V$ is an open subset of $\mc B$, and $v$ is a section of $\scr V_{\fk X}\otimes\omega_{\mc C/\mc B}(\blt S_{\fk X})\big|_{U_i}$ defined near $\sgm_i(V)$, then $\mc V_\varrho(\eta_i)v$, which is a section defined near $\{0\}\times V$, can be viewed (by taking series expansion) as an element of $\Vbb\otimes_\Cbb\scr O_{\mc B}(V)((z))dz$. If $w_i\in\Wbb_i\otimes_\Cbb\scr O_{\mc B}(V)$, we set \index{vw@$v\cdot w_i,v\cdot w_\blt$}
\begin{align}
v\cdot w_i=\Res_{z=0}Y_{\Wbb_i}(\mc V_\varrho(\eta_i)v,z)w_i.\label{eq116}
\end{align}

One can now define an $\scr O_{\mc B}$-linear action of $\pi_*\big(\scr V_{\fk X}\otimes\omega_{\mc C/\mc B}(\blt S_{\fk X})\big)$ on $\scr W_{\fk X}(\Wbb_\blt)$ as follows. If $V\subset\mc B$ is open, for any $v$ in $\pi_*\big(\scr V_{\fk X}\otimes\omega_{\mc C/\mc B}(\blt S_{\fk X})\big)(V)=\big(\scr V_{\fk X}\otimes\omega_{\mc C/\mc B}(\blt S_{\fk X})\big)(\pi^{-1}(V))$ and any $w\in\scr W_{\fk X}(\Wbb_\blt)(V)$, we set
\begin{align}
v\cdot w_\blt=\sum_{i=1}^N w_1\otimes w_2\otimes\cdots \otimes (v|_{U_i\cap\pi^{-1}(V)})\cdot w_i \otimes\cdots\otimes w_N.\label{eq152}
\end{align}
Then for each $b\in V$, the value $(v\cdot w_\blt)(b)$ (which is a vector inside the fiber $\scr W_{\fk X}(\Wbb_\blt)|b\simeq \Wbb_\blt$) equals
\begin{align}
(v\cdot w_\blt)(b)=(v|{\mc C_b})\cdot w_\blt(b)\label{eq96}
\end{align}
where $w_\blt(b)=w_1(b)\otimes w_2(b)\otimes\cdots\otimes w_N(b)$, and $v|{\mc C_b}$, the restriction of $v$ to the fiber $\mc C_b$, is in $\scr V_{\mc C_b}\otimes\omega_{\mc C_b}(\blt S_{\fk X}(b))$. The action on the right hand side of \eqref{eq96} is defined by \eqref{eq230} and \eqref{eq231}.

Define a \textbf{sheaf of covacua} \index{T@$\scr T_{\fk X}(\Wbb_\blt),\scr T_{\fk X}^*(\Wbb_\blt)$}
\begin{gather}
\scr T_{\fk X}(\Wbb_\blt)=\frac{\scr W_{\fk X}(\Wbb_\blt)}{\pi_*\big(\scr V_{\fk X}\otimes\omega_{\mc C/\mc B}(\blt S_{\fk X})\big)\cdot \scr W_{\fk X}(\Wbb_\blt)}\label{eq99}
\end{gather}
whose dual sheaf is denoted by $\scr T_{\fk X}^*(\Wbb_\blt)$ and called a \textbf{sheaf of conformal blocks} or \textbf{sheaf of vacua}. $\pi_*\big(\scr V_{\fk X}\otimes\omega_{\mc C/\mc B}(\blt S_{\fk X})\big)\cdot \scr W_{\fk X}(\Wbb_\blt)$ is the sheaf of $\scr O_{\mc B}$-modules associated to the presheaf whose sections on any open $V\subset\mc B$ are (linear combinations of) those in $\pi_*\big(\scr V_{\fk X}\otimes\omega_{\mc C/\mc B}(\blt S_{\fk X})\big)(V)\cdot \scr W_{\fk X}(\Wbb_\blt)(V)$. 

For each $b\in\mc B$, we let $\fk X_b$ be the restriction of $\fk X$ to $\mc C_b=\pi^{-1}(b)$, i.e. \index{XX@$\wtd {\fk X}$ and $\fk X$!$\fk X_b$}
\begin{align*}
\fk X_b=(\mc C_b;\sgm_1(b),\dots,\sgm_N(b);\eta_1|_{U_1\cap\mc C_b},\dots,\eta_N|_{U_N\cap\mc C_b}).
\end{align*}
We show that each fiber of the sheaf $\scr T_{\fk X}(\Wbb_\blt)$ is isomorphic to the space of covacua $\scr T_{\fk X_b}(\Wbb_\blt)$. Note the obvious $\scr O_{\mc B,b}$-module isomorphism of  stalks
\begin{gather}
\scr T_{\fk X}(\Wbb_\blt)_b\simeq \frac{\scr W_{\fk X}(\Wbb_\blt)_b}{\pi_*\big(\scr V_{\fk X}\otimes\omega_{\mc C/\mc B}(\blt S_{\fk X})\big)_b\cdot \scr W_{\fk X}(\Wbb_\blt)_b}.\label{eq95}
\end{gather}
Recall that $\fk m_b$ is the maximal ideal of $\scr O_{\mc B,b}$.

\begin{thm}\label{lb69}
For any $b\in\mc B$, the evaluation map at $b$:
\begin{gather}
\scr W_{\fk X}(\Wbb_\blt)_b=\Wbb_\blt\otimes_\Cbb\scr O_{\mc B,b}\rightarrow \Wbb_\blt,\qquad w\mapsto w(b)\label{eq94}
\end{gather}
descends to an isomorphism
\begin{gather}
\scr T_{\fk X}(\Wbb_\blt)|b=\frac{\scr T_{\fk X}(\Wbb_\blt)_b}{\fk m_b\cdot \scr T_{\fk X}(\Wbb_\blt)_b}\quad\xlongrightarrow{\simeq} \quad\scr T_{\fk X_b}(\Wbb_\blt).\label{eq97}
\end{gather}
\end{thm}

\begin{proof}
By \eqref{eq95}, the fiber $\scr T_{\fk X}(\Wbb_\blt)|b$ equals $\scr W_{\fk X}(\Wbb_\blt)_b$ modulo the subspace spanned by $\fk m_b\cdot\scr W_{\fk X}(\Wbb_\blt)_b$ and $\pi_*\big(\scr V_{\fk X}\otimes\omega_{\mc C/\mc B}(\blt S_{\fk X})\big)_b\cdot \scr W_{\fk X}(\Wbb_\blt)_b$. The first one is sent by the map \eqref{eq94} to $0$, and the second one into $H^0\big(\mc C_b,\scr V_{\mc C_b}\otimes\omega_{\mc C_b}(\blt S_{\fk X}(b))\big) \Wbb_\blt$ according to the relation \eqref{eq96}. Thus the linear map \eqref{eq97} is well-defined. It is clearly surjective. To show that \eqref{eq97} is injective, it suffices to show that the map  \eqref{eq94} sends $\pi_*\big(\scr V_{\fk X}\otimes\omega_{\mc C/\mc B}(\blt S_{\fk X})\big)_b\cdot \scr W_{\fk X}(\Wbb_\blt)_b$ onto $H^0\big(\mc C_b,\scr V_{\mc C_b}\otimes\omega_{\mc C_b}(\blt S_{\fk X}(b))\big) \Wbb_\blt$. This follows from Corollary \ref{lb33}.
\end{proof}

\begin{rem}
For any open subset $V\subset\mc B$, an element $\upphi\in\scr T_{\fk X}^*(\Wbb_\blt)(V)$ is an homomorphism $\scr W_{\fk X}(\Wbb_\blt)|_V\rightarrow\scr O_V$ vanishing on $\pi_*\big(\scr V_{\fk X}\otimes\omega_{\mc C/\mc B}(\blt S_{\fk X})\big)(W)\cdot \scr W_{\fk X}(\Wbb_\blt)(W)$ for any open subset $W\subset V$. By \eqref{eq96}, it is clear that this vanishing requirement $\Leftarrow$ for any $b\in V$, the fiber
\begin{align*}
\upphi(b):\scr W_{\fk X}(\Wbb_\blt)|b\simeq\scr W_{\fk X_b}(\Wbb_\blt)\rightarrow \scr O_{\mc B}|b\simeq \Cbb
\end{align*}
vanishes on $H^0(C,\scr V_{\mc C_b}\otimes\omega_{\mc C_b}(\blt S_{\fk X_b}))\cdot\Wbb_\blt$. By Corollary \ref{lb33}, we also have $\Rightarrow$. This proves:
\end{rem}

\begin{pp}\label{lb35}
Let $V\subset\mc B$ be open, and let $\upphi:\scr W_{\fk X}(\Wbb_\blt)|_V\rightarrow \scr O_V$ be a homomorphism of $\scr O_V$-modules. Then $\upphi$ is a conformal block if and only if its restriction to each fiber is a conformal block. More precisely, $\upphi\in\scr T_{\fk X}^*(\Wbb_\blt)(V)$ if and only if $\upphi(b)\in\scr T_{\fk X_b}^*(\Wbb_\blt)$ for any $b\in V$.
\end{pp}

We now define the sheaves of covacua and conformal blocks for any family $\fk X=(\pi:\mc C\rightarrow\mc B;\sgm_1,\dots,\sgm_N)$ of $N$-pointed complex curves whose local coordinates are not specified. If $V\subset\mc B$ is open, then $\fk X_V$ denotes the subfamily \index{XX@$\wtd {\fk X}$ and $\fk X$!$\fk X_V,\mc C_V=\pi^{-1}(V)$} \index{CV@$\mc C_V=\pi^{-1}(V)$}
\begin{align}
\fk X_V=(\pi:\mc C_V=\pi^{-1}(V)\rightarrow V;\sgm_1|_V,\dots,\sgm_N|_V).\label{eq140}
\end{align}
Define $\scr W_{\fk X}(\Wbb_\blt)$ to be an infinite rank locally free sheaf on $\mc B$ as follows. For any connected open subset $V\subset \mc B$ together with local coordinates $\eta_1,\dots,\eta_N$ of the family $\pi:\mc C_V\rightarrow V$ defined near $\sgm_1(V),\dots,\sgm_N(V)$ respectively, we have a trivialization \index{U@$\mc U(\rho),\mc U(\eta_\blt)$}
\begin{align}
\mc U(\eta_\blt)\equiv \mc U(\eta_1)\otimes\cdots\otimes\mc U(\eta_N):\scr W_{\fk X}(\Wbb_\blt)|_V\xrightarrow{\simeq}\Wbb_\blt\otimes_\Cbb\scr O_V\label{eq120}
\end{align}
such that if $\mu_\blt$ is another set of local coordinates, then
\begin{align*}
\mc U(\eta_\blt)\mc U(\mu_\blt)^{-1}:\Wbb_\blt\otimes_\Cbb\scr O_V\xrightarrow{\simeq}\Wbb_\blt\otimes_\Cbb\scr O_V
\end{align*}
is defined such that for any constant section $w_\blt=w_1\otimes\cdots\otimes w_N\in\Wbb_\blt$, $\mc U(\eta_\blt)\mc U(\mu_\blt)^{-1}w_\blt$, as a $\Wbb_\blt$-valued holomorphic function, satisfies
\begin{align}
&\Big(\mc U\big(\eta_\blt\big)\mc U\big(\mu_\blt\big)^{-1}w_\blt\Big )(b)\equiv \Big(\mc U\big(\eta_\blt\big|\mu_\blt\big)w_\blt\Big )(b)\nonumber\\
=&\mc U\big((\eta_1|\mu_1)_b \big)w_1\otimes\mc U\big((\eta_2|\mu_2)_b \big)w_2\otimes\cdots\otimes \mc U\big((\eta_N|\mu_N)_b \big)w_N\label{eq124}
\end{align}
for any $b\in V$. Here, for each $1\leq i\leq N$, $(\eta_i|\mu_i)_b$ is the element in $\Gbb$ satisfying
\begin{align}
(\eta_i|\mu_i)_b(z)=\eta_i\circ(\mu_i,\pi)^{-1}(z,b).\label{eq244}
\end{align}
If we compare the transition functions \eqref{eq123} and \eqref{eq124}, we see that there is a natural and coordinate-independent isomorphism of vector spaces
\begin{align*}
\scr W_{\fk X}(\Wbb_\blt)|b\simeq \scr W_{\fk X_b}(\Wbb_\blt)
\end{align*}
where $\scr W_{\fk X_b}(\Wbb_\blt)$ is defined near \eqref{eq121}. We shall identify these two spaces in the following.

The action of $\pi_*\big(\scr V_{\fk X}\otimes\omega_{\mc C/\mc B}(\blt S_{\fk X})\big)$ on $\scr W_{\fk X}(\Wbb_\blt)$ is defined fiberwisely by the action of $\pi_*\big(\scr V_{\fk X}\otimes\omega_{\mc C/\mc B}(\blt S_{\fk X})\big)(b)=H^0\big(\mc C_b,\scr V_{\mc C_b}\otimes\omega_{\mc C_b}(\blt S_{\fk X}(b))\big)$ (recall again Corollary \ref{lb33}) on $\scr W_{\fk X_b}(\Wbb_\blt)$. By \eqref{eq96}, it is clear that if we choose a set of local coordinates $\eta_\blt$ near $\sgm_1(V),\dots,\sgm_N(V)$ as above, and if we  identify $\scr W_{\fk X}(\Wbb_\blt)|_V\simeq\Wbb_\blt\otimes_\Cbb\scr O_V$ via $\mc U(\eta_\blt)$, then this action is described by \eqref{eq116} and \eqref{eq152}.

The \textbf{sheaf of covacua} $\scr T_{\fk X}(\Wbb_\blt)$ is defined still by \eqref{eq99}, and the \textbf{sheaf of conformal blocks} $\scr T_{\fk X}^*(\Wbb_\blt)$ is its dual sheaf.

\begin{pp}\label{lb41}
Assume that $\mc B$ is connected, and let $\upphi:\scr W_{\fk X}(\Wbb_\blt)\rightarrow \scr O_{\mc B}$ be a homomorphism of $\scr O_{\mc B}$-modules. Suppose that $V$ is a non-empty open subset of $\mc B$, and the restriction $\upphi|_V:\scr W_{\fk X}(\Wbb_\blt)|_V\rightarrow \scr O_V$ is  an element in $\scr T_{\fk X}^*(\Wbb_\blt)(V)$. Then $\upphi\in\scr T_{\fk X}^*(\Wbb_\blt)(\mc B)$.
\end{pp}

\begin{proof}
We first assume that $\mc B$ is small enough such that $\fk X$ can be equipped with $N$ local coordinates,  and that the connect manifold $\mc B$ is biholomorphic to a polydisc which in particular is Stein. Then for each $b\in\mc B$, the restriction of $\pi_*\big(\scr V_{\fk X}\otimes\omega_{\mc C/\mc B}(\blt S_{\fk X})\big)(\mc B)\cdot \scr W_{\fk X}(\Wbb_\blt)(\mc B)$ to the fiber $\mc C_b$ is $H^0(\mc C_b,\scr V_{\fk X_b}\otimes\omega_{\mc C_b}(\blt S_{\fk X}))\cdot \scr W_{\fk X_b}(\Wbb_\blt)$ by Corollary \ref{lb33}. Thus, $\upphi$ is a conformal block if and only if its evaluation with any element of $\pi_*\big(\scr V_{\fk X}\otimes\omega_{\mc C/\mc B}(\blt S_{\fk X})\big)(\mc B)\cdot \scr W_{\fk X}(\Wbb_\blt)(\mc B)$ (which is a holomorphic function on $\mc B$) is zero. By our assumption, such holomorphic function vanishes on $V$. Thus it is zero on $\mc B$.

In general, we let $A$ be the set of all $b\in\mc B$ such that $b$ has a neighborhood $W$ such that the restriction $\upphi|_W$ is a conformal block. Then $A$ is open and (by assumption) non-empty. For any $b\in\mc B-A$, let $W$ be a connected neighborhood of $b$ small enough as in the first paragraph. Then $\upphi|_W$ is a conformal block if $W$ has a non-zero open subset $V$ such that $\upphi|_V$ is a conformal block. Therefore $W$ must be disjoint from $A$. This shows that $\mc B-A$ is open. Thus $\mc B=A$.
\end{proof}

Next, we generalize Theorem \ref{lb32} to sheaves of conformal blocks. Assume that $\fk X$ has local coordinates $\eta_1,\dots,\eta_N$.  Let $\upphi:\scr W_{\fk X}(\Wbb_\blt)\rightarrow \scr O_{\mc B}$ be a homomorphism of $\scr O_{\mc B}$-modules. Then for each $1\leq i\leq N$, we have an $\scr O_{\mc B}((z))$-module homomorphism
\begin{gather*}
\wr\upphi_{\sgm_i(\mc B)}:\big(\Vbb\otimes_\Cbb\scr O_{\mc B}((z))\big)\otimes_{\scr O_{\mc B}} (\Wbb_\blt\otimes_\Cbb\scr O_{\mc B})\rightarrow \scr O_{\mc B}((z))
\end{gather*}
such that for each $v\in\Vbb,w_\blt\in\Wbb_\blt$, considered as constant sections on $\mc B$ of  $\Vbb\otimes_\Cbb\scr O_{\mc B}((z))$ and $\Wbb_\blt\otimes_\Cbb\scr O_{\mc B}$ respectively, $v\otimes w_\blt$ is sent to
\begin{align}
\wr\upphi_{\sgm_i(\mc B)}(v\otimes w_\blt)\equiv\wr\upphi_{\sgm_i(\mc B)}(v,w_\blt)=\upphi(w_1\otimes w_2\otimes\cdots\otimes Y_{\Wbb_i}(v,z)w_i\otimes\cdots\otimes w_N).
\end{align}
For any $b\in\mc B$, for the fiber map $\upphi(b):\scr W_{\fk X_b}(\Wbb_\blt)\rightarrow\Cbb$ we can define $\wr\upphi(b)_{\sgm_i(b)}$ as in Section \ref{lb34}. Then it is clear that the following elements in $\Cbb((z))$ are equal:
\begin{align}
\wr\upphi(b)_{\sgm_i(b)}(v,w_\blt)=\big(\wr\upphi_{\sgm_i(\mc B)}(v,w_\blt)\big)(b).\label{eq100}
\end{align}
Shortly speaking, the restriction of $\wr\upphi_{\sgm_i(\mc B)}$ to each fiber $\mc C_b$ equals $\wr\upphi(b)_{\sgm_i(b)}$.

For any (non-necessarily open) subset $E\subset\mc C$, \index{CS@$C-\SX,\mc C-\SX,E-\SX$} we set \index{SX@$\SX$, $\SX(b)$!$E-S_{\fk X}$}
\begin{align}
E-S_{\fk X}=E-\bigcup_{i=1}^N\sgm_i(\mc B).
\end{align}
Choose mutually disjoint neighborhoods $U_1,\dots,U_N$ of $\sgm_1(\mc B),\dots,\sgm_N(\mc B)$ on which $\eta_1,\dots,\eta_N$ are defined respectively. Note that $\scr V_{\fk X}$ and $\scr V_{\fk X}|_{\mc C-S_{\fk X}}$ are $\scr O_{\mc B}$ modules realized by pulling back $\scr O_{\mc B}$ to $\mc C$. The following theorem is clearly true if $\mc B$ is replaced by an open subset $V$ and $\fk X$ by the the subfamily $\fk X_V$.

\begin{thm}\label{lb45}
Assume that $\fk X$ is a smooth family. Let $\upphi:\scr W_{\fk X}(\Wbb_\blt)\rightarrow \scr O_{\mc B}$ be a homomorphism of $\scr O_{\mc B}$-modules. Then the following are equivalent.
\begin{enumerate}[label=(\alph*)]
	\item $\upphi$ is an element of $\scr T_{\fk X}^*(\Wbb_\blt)(\mc B)$. 
	\item  For each open $V\subset \mc B$ and $w\in\scr W_{\fk X}(\Wbb_\blt)(V)\simeq \Wbb_\blt\otimes_\Cbb\scr O_{\mc B}(V)$, there is a (necessarily unique) element \index{zz@$\wr\upphi$}
	\begin{align*}
	\wr\upphi(w)\in H^0(\mc C_V-S_{\fk X},\scr V_{\fk X}^*)
	\end{align*}
satisfying that for each $1\leq i\leq N$, if we identify $U_i\simeq (\eta_i,\pi)(U_i)$ via $(\eta_i,\pi)$ and identify $\scr V_{\fk X}|_{U_i}\simeq\Vbb\otimes_\Cbb\scr O_{U_i}$ via $\mc U_\varrho(\eta_i)$, then the evaluation of $\wr\upphi(w)$ with any $v\in\scr V_{\fk X}(U_i\cap\mc C_V)$ (restricted to $U_i\cap\mc C_V-S_{\fk X}$) is
\begin{align}
\wr\upphi(v,w)=\wr\upphi_{\sgm_i(\mc B)}(v,w).\label{eq101}
\end{align}
\end{enumerate}	
\end{thm}
Note that  $v\in\scr V_{\fk X}(U_i\cap\mc C_V)$ can be regarded as an element of $\Vbb\otimes_\Cbb\scr O_{\mc B}(V)((z))$.

\begin{proof}
Suppose that (b) holds. Then by \eqref{eq100}, for any $b\in\mc B$ and any section $w$ of $\scr W_{\fk X}(\Wbb_\blt)$ defined near $b$, the restriction $\wr\upphi(w)|_{\mc C_b-S_{\fk X}}=\wr\upphi(w)|_{\mc C_b-S_{\fk X_b}}$, which is an element in $H^0(\mc C_b-S_{\fk X_b},\scr V_{\fk X_b}^*)$, has series expansion 
$\wr\upphi(b)_{\sgm_i(b)}$ near each $\sgm_i(b)$. Thus, by Theorem \ref{lb32}, $\upphi(b)$ is a conformal block on the fiber $\mc C_b$. Since this is true for each $b\in\mc B$, by Proposition \ref{lb35}, $\upphi$ is a conformal block on the family $\fk X$. This proves (a).

Now assume (a). Choose  open $V\subset \mc B$ and $w\in\scr W_{\fk X}(\Wbb_\blt)(V)$. For any $n\in\Nbb$, the restriction of $\wr\upphi_{\sgm_i(\mc B)}(w)$ to $\Vbb^{\leq n}\otimes_\Cbb\scr O_V$ gives a homomorphism of $\scr O_V$-modules
\begin{align*}
s_i:=\wr\upphi^{\leq n}_{\sgm_i(\mc B)}(w):\Vbb^{\leq n}\otimes_\Cbb\scr O_V\rightarrow\scr O_V((z)).
\end{align*}
which can also be considered as
\begin{align*}
s_i\in\big((\Vbb^{\leq n})^*\otimes_\Cbb\scr O(V)\big)((z)).
\end{align*}
Let $\scr E=(\scr V_{\fk X}^{\leq n})^*|_{\mc C_V}$. By the fact due to Proposition \ref{lb35} that $\upphi(b)$ is a conformal block on $\mc C_b$ for each $b\in V$, it is easy to see that  $s_1,\dots,s_N$ satisfy (b) of Theorem \ref{lb18}. Thus they also satisfy (a) of that theorem, namely, that $s_1,\dots,s_N$ can be extended to an element
\begin{align*}
s:=\wr\upphi^{\leq n}(w)\in H^0(\mc C_V,\scr V_{\fk X}^*(\blt S_{\fk X})|_{\mc C_V}),
\end{align*}
which can also be regarded as in $H^0(\mc C_V-S_{\fk X},\scr V_{\fk X}^*|_{\mc C_V})$. It is clear that $\wr\upphi^{\leq n'}(w)$ restricts to $\wr\upphi^{\leq n}(w)$ when $n'\geq n$. Thus $\wr\upphi(w)$, the projective limit over $n\in\Nbb$ of $\wr\upphi^{\leq n}(w)$, satisfies \eqref{eq101}. This proves (b).
\end{proof}

\begin{rem}
In the case that local coordinates are not assigned to the smooth family $\fk X$, one can still define $\wr\upphi$ for each $\upphi\in\scr T_{\fk X}^*(\Wbb_\blt)(\mc B)$. Indeed, for each $1\leq i\leq N$,  the restriction of $\wr\upphi_{\sgm_i(\mc B)}$ to each fiber $\mc C_b$ is independent of the choice of local coordinates by (the proof of) Lemma \ref{lb125}. So is $\wr\upphi$. Thus, one can define $\wr\upphi$ locally, and glue them together to obtain the global section.
\end{rem}

\section{Sewing  conformal blocks}\label{lb76}

\subsection*{Formal conformal blocks}

Let $\fk X=(\pi:\mc C\rightarrow\mc B;\sgm_1,\dots,\sgm_N;\eta_1,\dots,\eta_N)$ be a family of $N$-pointed complex curves with local coordinates obtained via sewing the following smooth family
\begin{align*}
\wtd{\fk X}=(\wtd\pi:\wtd{\mc C}\rightarrow\wtd{\mc B};\sgm_1,\dots,\sgm_N;\sgm_1',\dots,\sgm_M';\sgm_1'',\dots,\sgm_M'';\eta_1,\dots,\eta_N;\xi_1,\dots,\xi_M;\varpi_1,\dots,\varpi_M).
\end{align*}
(See Section \ref{lb24}.) As suggested by the notations, we require as in Remark \ref{lb39} that the $N$-points  $\sgm_1,\dots,\sgm_N$ and the local coordinates $\eta_1,\dots,\eta_N$ of $\fk X$ are constant with respect to sewing. In this section, \emph{we only assume that each connected component of each fiber $\wtd{\mc C}_b$ contains at least one of the $N+2M$ marked points of $\wtd{\fk X}_b$}. This is slightly weaker than the assumption in Rem. \ref{lb39}. Choose $\Vbb$-modules $\Wbb_1,\dots,\Wbb_N$, and $\Mbb_1,\dots,\Mbb_M$ whose contragredient modules are $\Mbb_1',\dots,\Mbb_M'$, which are associated to $\sgm_1,\dots,\sgm_N$, $\sgm_1',\dots,\sgm_M'$, $\sgm_1'',\dots,\sgm_M''$ respectively.

To simplify discussions, we assume throughout this section that $\wtd{\mc B}$ is Stein. This assumption allows us to work with modules instead of sheaves of modules. Recall that we have
\begin{align*}
\mc D_{r_\bullet\rho_\bullet}=\mc D_{r_1\rho_1}\times\cdots\times\mc D_{r_M\rho_M},\qquad\mc B=\wtd{\mc B}\times \mc D_{r_\bullet\rho_\bullet}.
\end{align*}
Then $\mc B$ is also Stein.  We have the identification $\scr W_{\fk X}(\Wbb_\blt)\simeq \Wbb_\blt\otimes_{\Cbb}\scr O_{\mc B}$ realized by $\mc U(\eta_\blt)$. By taking series expansions,  $\scr O(\mc B)$ can be regarded as a subalgebra of
\begin{align*}
\scr O(\wtd{\mc B})[[q_\blt]]=\scr O(\wtd{\mc B})[[q_1,q_2,\dots,q_M]],
\end{align*}
whose elements are formal power series of $q_1,\dots,q_M$ whose coefficients are in $\scr O(\wtd{\mc B})$. In particular, $\scr O(\wtd{\mc B})[[q_\blt]]$ is an $\scr O(\mc B)$-module. Choose a homomorphism of $\scr O(\mc B)$-modules
\begin{align}
\upphi:\scr W_{\fk X}(\Wbb_\blt)(\mc B)=\Wbb_\blt\otimes_\Cbb\scr O(\mc B)\rightarrow\scr O(\wtd{\mc B})[[q_\blt]].\label{eq102}
\end{align}
We say that $\upphi$ is a \textbf{formal conformal block} if $\upphi$ vanishes on $\pi_*\big(\scr V_{\fk X}\otimes\omega_{\mc C/\mc B}(\blt S_{\fk X})\big)(\mc B)\cdot \scr W_{\fk X}(\Wbb_\blt)(\mc B)$. We say that $\upphi$ \textbf{converges absolutely and locally uniformly (a.l.u.)} if the image of $\upphi$ is in $\scr O(\mc B)$. This name is explained below.

\begin{rem}\label{lb77}
Write
\begin{align*}
q_\blt^{n_\blt}=q_1^{n_1}q_2^{n_2}\cdots q_M^{n_M}
\end{align*}
for any $n_\blt=(n_1,n_2,\dots,n_M)\in\Nbb^M$. For each $w\in \scr W_{\fk X}(\Wbb_\blt)(\mc B)$, we have the series expansion
\begin{align*}
\upphi(w)=\sum_{n_\blt\in\Nbb^M}\upphi(w)_{n_\blt}\cdot q_\blt^{n_\blt}
\end{align*}
where each $\upphi(w)_{n_\blt}$ is a holomorphic function on $\wtd{\mc B}$. Then it is clear that  $\upphi$ converges a.l.u. on $\mc B$ if and only if for any $w$ (which is sufficient to be constant) and any  compact subsets $K\subset\wtd{\mc B}$ and $Q\subset\mc D_{r_\blt\rho_\blt}$,  there exists $C>0$ such that
\begin{align}
\sum_{n_\blt\in\Nbb^M}\big|\upphi(w)_{n_\blt}(b)\big|\cdot |q_\blt^{n_\blt}|\leq C
\end{align}
for any $b\in K$ and $q_\blt=(q_1,\dots,q_M)$ in $Q$. 
\end{rem}

If $\upphi$ converges a.l.u. on $\mc B$, one can regard $\upphi$ as a homomorphism of $\scr O_{\mc B}$-modules
\begin{align*}
\upphi:\scr W_{\fk X}(\Wbb_\blt)=\Wbb_\blt\otimes_\Cbb\scr O_{\mc B}\rightarrow\scr O_{\mc B}
\end{align*}
whose values at the global sections of $\scr W_{\fk X}(\Wbb_\blt)$  are given by \eqref{eq102}.

\begin{pp}\label{lb43}
Let $\upphi$ in \eqref{eq102} be a formal conformal block, and assume that $\upphi$ converges a.l.u..
\begin{enumerate}
\item If for every $b\in\wtd{\mc B}$, each connected component of  $\wtd{\mc C}_b$ contains at least one $\sgm_1(b),\dots,\sgm_N(b)$ (cf. Rem. \ref{lb39}), then  $\upphi$ is a conformal block, i.e. $\upphi\in\scr T_{\fk X}^*(\Wbb_\blt)(\mc B)$.
\item If  for every $b\in\mc B-\Delta=\wtd{\mc B}\times\mc D_{r_\blt\rho_\blt}^\times$, each connected component of $\mc C_b$ contains at least one of $\sgm_1(b),\dots,\sgm_N(b)$, then $\upphi$ is a conformal block on $\mc B-\Delta$.
\end{enumerate}

\end{pp}

\begin{proof}
Case 1 follows easily from Corollary \ref{lb33}, Proposition \ref{lb35}, and our assumption that $\wtd{\mc B}$ (and hence $\mc B$) is Stein.

In case 2, note that Proposition \ref{lb35} and Theorem \ref{lb27} apply to the restriction $\fk X_{\mc B-\Delta}$. So it suffices to prove for all $n\in\Nbb$ that for sufficiently large $k\in\Nbb$, the elements of $\pi_*\big(\scr V_{\fk X}\otimes\omega_{\mc C/\mc B}(k S_{\fk X})\big)(\mc B)$  generate the stalk $\pi_*\big(\scr V_{\fk X}^{\leq n}\otimes\omega_{\mc C/\mc B}(k\SX)\big)_b$. This follows from Cartan's theorem A and that $\pi_*\big(\scr V_{\fk X}\otimes\omega_{\mc C/\mc B}(k S_{\fk X})\big)$ is coherent (by Grauert's direct image theorem).
\end{proof}

\subsection*{Sewing conformal blocks}

Let $\Wbb_\blt\otimes \Mbb_\blt\otimes\Mbb_\blt'$ be
\begin{align}
\Wbb_1\otimes\cdots\otimes \Wbb_N\otimes\Mbb_1\otimes\Mbb_1'\otimes\cdots\otimes\Mbb_M\otimes\Mbb_M'.
\end{align}
We have switched the orders and put each $\Mbb_j$ and its contragredient module $\Mbb_j'$ together, which are associated to $\sgm_j'(\mc B)$ and $\sgm_j''(\mc B)$ respectively. Our goal is to define a formal conformal block from each element of $\scr T_{\wtd{\fk X}}^*(\Wbb_\blt\otimes \Mbb_\blt\otimes\Mbb_\blt')$. Notice that for each  $j=1,\dots,M$, $(\Mbb_j'\otimes\Mbb_j)^*$ can be regarded as the algebraic completion of $\Mbb_j\otimes\Mbb_j'$. Define \index{zzz@$\btr\otimes_j\btl$}
\begin{align*}
\btr\otimes_j\btl\in (\Mbb_j'\otimes\Mbb_j)^*
\end{align*}
such that for any $m'\in\Mbb_j',m\in\Mbb_j$,
\begin{align*}
\bk{\btr\otimes_j\btl,m'\otimes m}=\bk{m',m}.
\end{align*}
Let $A\in\End(\Mbb_j)$ whose transpose $A^\tr\in\End(\Mbb_j')$ exists, i.e.,
\begin{align}
\bk{Am,m'}=\bk{m,A^\tr m'}\label{eq103}
\end{align}
for any $m'\in\Mbb_j',m\in\Mbb_j$. Then we have an element
\begin{align}
A\btr\otimes_j\btl\equiv \btr\otimes_j A^\tr\btl\quad\in (\Mbb_j'\otimes\Mbb_j)^*\label{eq108}
\end{align}
whose value at each $m'\otimes m$ is \eqref{eq103}.

More explicitly, for each $n\in\Nbb$ we choose a basis $\{m(n,a)\}_a$ of the finite dimensional vector space $\Mbb_j(n)$ (which we recall is the $n$-eigenspace of $\wtd L_0$; recall also Convention \ref{lb36}). Its dual basis $\{\wch m(n,a)\}_a$ is a basis of $\Mbb_j'(n)=\Mbb_j(n)^*$ satisfying $\bk{m(n,a),\wch m(n,b)}=\delta_{a,b}$. Then we have
\begin{gather*}
\btr\otimes_j\btl=\sum_{n\in\Nbb}\sum_am(n,a)\otimes\wch m(n,a),
\end{gather*}
and
\begin{align*}
&A\btr\otimes_j\btl=\sum_{n\in\Nbb}\sum_a A\cdot m(n,a)\otimes\wch m(n,a)\\
=& \btr\otimes_j A^\tr\btl=\sum_{n\in\Nbb}\sum_am(n,a)\otimes A^\tr\cdot \wch m(n,a).
\end{align*}

Let $P(n)$ be the projection of $\Mbb_j$ onto $\Mbb_j(n)$. Its transpose, which is the projection of $\Mbb_j'$ onto $\Wbb_j'(n)$, is also denoted by $P(n)$. Then we clearly have
\begin{align*}
P(n)\btr\otimes_j\btl=\btr\otimes_j P(n)\btl=\sum_am(n,a)\otimes\wch m(n,a)\qquad\in\Mbb_j\otimes\Mbb_j'.
\end{align*}
Recall $\wtd L_0^\tr=\wtd L_0$. Define
\begin{align}
q_j^{\wtd L_0}=\sum_{k\in\Nbb}P(n)q_j^n\qquad\in\End(\Mbb_j)[[q_j]].
\end{align} 
Then we have
\begin{align}
q_j^{\wtd L_0}\btr\otimes_j\btl=\btr\otimes_j q_j^{\wtd L_0} \btl\qquad \in (\Mbb_j\otimes\Mbb_j')[[q_j]].\label{eq109}
\end{align}

For any $\uppsi\in\scr T_{\wtd{\fk X}}^*(\Wbb_\blt\otimes\Mbb_\blt\otimes\Mbb_\blt')(\wtd{\mc B})$, we define its (normalized) \textbf{sewing} $\wtd{\mc S}\uppsi$ which is an $\scr O(\mc B)$-module homomorphism \index{S@$\wtd{\mc S}\uppsi,\mc S\uppsi$}
\begin{align*}
\wtd{\mc S}\uppsi:\scr W_{\fk X}(\Wbb_\blt)(\mc B)=\Wbb_\blt\otimes_\Cbb\scr O(\mc B)\rightarrow\scr O(\wtd{\mc B})[[q_\blt]],
\end{align*}
and the (standard) \textbf{sewing}
\begin{align*}
\mc S\uppsi:\scr W_{\fk X}(\Wbb_\blt)(\mc B)=\Wbb_\blt\otimes_\Cbb\scr O(\mc B)\rightarrow\scr O(\wtd{\mc B})\{q_\blt\},
\end{align*}
as follows. Regard $\uppsi$ as an $\scr O(\wtd{\mc B})$-module homomorphism $\Wbb_\blt\otimes\Mbb_\blt\otimes\Mbb_\blt'\otimes_\Cbb\scr O(\wtd{\mc B})\rightarrow \scr O(\wtd{\mc B})$. $\wtd{\mc S}\uppsi$ is defined such that for any constant section $w\in\Wbb_\blt$,
\begin{align}
\wtd{\mc S}\uppsi(w)=\uppsi\Big(w\otimes (q_1^{\wtd L_0}\btr\otimes_1\btl)\otimes\cdots\otimes (q_M^{\wtd L_0}\btr\otimes_M\btl)\Big).
\end{align}
$\mc S\uppsi$ is defined similarly, except that the normalized energy operator $\wtd L_0$ is replaced by the standard one $L_0$. When $\Mbb_1,\dots,\Mbb_M$ are irreducible, $\wtd{\mc S}\upphi$ differs from $\mc S\upphi$ by a factor $q_1^{\lambda_1}\cdots q_M^{\lambda_M}$ for some $\lambda_1,\dots,\lambda_M\in\Cbb$. Our goal is to show that $\wtd{\mc S}\uppsi$ is a formal conformal block. In the case that $\Mbb_1,\dots \Mbb_M$ are semisimple (which is sufficient for our purpose), this will show that $\mc S\upphi$ is also a formal conformal block in a suitable sense. We begin our proof with the following observation, in which we have omitted the subscript $j$ of $\xi,\varpi,q$ for simplicity.

\begin{lm}\label{lb37}
Let $R$ be any unital commutative $\Cbb$-algebra. (For instance, $R=\scr O(X)$ where $X$ is a complex manifold.) For any $u\in\Vbb$ and $f\in R[[\xi,\varpi]]$, the following two elements in $(\Mbb_j\otimes\Mbb_j'\otimes R)[[q]]$ (where the tensor products are over $\Cbb$) are equal:
\begin{align}
&\Res_{\xi=0}~Y_{\Mbb_j}\big(\xi^{L_0}u,\xi\big)q^{\wtd L_0}\btr\otimes_j\btl\cdot f(\xi,q/\xi)\frac{d\xi}{\xi}\nonumber\\
=&\Res_{\varpi=0}~q^{\wtd L_0}\btr\otimes_j Y_{\Mbb_j'}\big(\varpi^{L_0}\mc U(\upgamma_1)u,\varpi\big)\btl\cdot f(q/\varpi,\varpi)\frac{d\varpi}{\varpi}.\label{eq110}
\end{align}
\end{lm}

\begin{rem}\label{lb38}
We  explain the meaning of the left hand side; the other side can be understood in a similar way. As   $q^{\wtd L_0}\btr\otimes_j\btl$ is an element of $(\Mbb_j\otimes\Mbb_j')[[q]]$, $Y_{\Mbb_j}\big(\xi^{L_0}u,\xi\big)q^{\wtd L_0}\btr\otimes_j\btl$ is an element of $(\Mbb_j\otimes\Mbb_j')((\xi))[[q]]$, i.e. it is a formal \emph{power} series of $q$ whose coefficients are in $(\Mbb_j\otimes\Mbb_j')((\xi))$. (Note that one cannot switch the order of $((\xi))$ and $[[q]]$.) Identify $(\Mbb_j\otimes\Mbb_j')((\xi))[[q]]\simeq (\Mbb_j\otimes\Mbb_j'\otimes 1)((\xi))[[q]]$, which is a subspace of the $R((\xi))[[q]]$-module $(\Mbb_j\otimes\Mbb_j'\otimes R)((\xi))[[q]]$.  On the other hand,  write $f(\xi,\varpi)=\sum_{m,n\in\Nbb}f_{m,n}\xi^m\varpi^n$ where each $f_{m,n}$ is in $R$. Then
\begin{align*}
f(\xi,q/\xi)=\sum_{n\geq 0}\sum_{k\geq -n}f_{n+k,n}\xi^k q^n,
\end{align*}
which shows $f(\xi,q/\xi)\in R((\xi))[[q]]$. Thus, the term in the residue on the left hand side is an element in
\begin{align*}
(\Mbb_j\otimes\Mbb_j'\otimes R)((\xi))[[q]]d\xi,
\end{align*}
whose residue is in $(\Mbb_j\otimes\Mbb_j'\otimes R)[[q]]$.
\end{rem}

\begin{proof}[Proof of Lemma \ref{lb37}]
Consider $Y_{\Mbb_j}\big(\xi^{L_0}u,\xi\big)q^{\wtd L_0}$ as an element of $\End(\Mbb_j)[[\xi^{\pm1},q]]$. Since $\wtd L_0^\tr=\wtd L_0$,  we have the following relations of elements of $\End(\Mbb_j')[[\xi^{\pm1},q^{\pm 1}]]$:
\begin{align*}
&\big(Y_{\Mbb_j}\big(\xi^{L_0}u,\xi\big)q^{\wtd L_0}\big)^\tr=q^{\wtd L_0}\big(Y_{\Mbb_j}\big(\xi^{L_0}u,\xi\big)\big)^\tr\xlongequal{\eqref{eq92}}q^{\wtd L_0}Y_{\Mbb_j'}\big(\mc U(\upgamma_\xi)\xi^{L_0}u,\xi^{-1}\big)\\
\xlongequal{\eqref{eq78}}&q^{\wtd L_0}Y_{\Mbb_j'}\big(\xi^{-L_0}\mc U(\upgamma_1)u,\xi^{-1}\big)\xlongequal{\eqref{eq107}}Y_{\Mbb_j'}\big((q/\xi)^{L_0}\mc U(\upgamma_1)u,q/\xi\big)q^{\wtd L_0}.
\end{align*}
(Note that due to the appearance of $q/\xi$ in the vertex operator,  it was not known a priori that the right hand side contains no negative powers of $q$.) Thus, by \eqref{eq108}, we have the following equations of elements in $(\Mbb_j'\otimes\Mbb_j)^*[[\xi^{\pm1},q^{\pm 1}]]$:
\begin{align}
&Y_{\Mbb_j}\big(\xi^{L_0}u,\xi\big)q^{\wtd L_0}\btr\otimes_j\btl=\btr\otimes_j \big(Y_{\Mbb_j}\big(\xi^{L_0}u,\xi\big)q^{\wtd L_0}\big)^\tr\btl\nonumber\\
=&\btr\otimes_j~Y_{\Mbb_j'}\big((q/\xi)^{L_0}\mc U(\upgamma_1)u,q/\xi\big)q^{\wtd L_0}\btl=q^{\wtd L_0}\btr\otimes_j~Y_{\Mbb_j'}\big((q/\xi)^{L_0}\mc U(\upgamma_1)u,q/\xi\big)\btl.\label{eq111}
\end{align}
Since for each $n$, $P(n)\btr\otimes_j\btl$ is in $\Mbb_j\otimes\Mbb_j'$, \eqref{eq111} is actually an element in $(\Mbb_j\otimes\Mbb_j')[[\xi^{\pm1},q^{\pm 1}]]$.

Let 
\begin{gather*}
A(\xi,q)=Y_{\Mbb_j}\big(\xi^{L_0}u,\xi\big)q^{\wtd L_0}\btr\otimes_j\btl,\\
B(\varpi,q)=q^{\wtd L_0}\btr\otimes_j~Y_{\Mbb_j'}\big(\varpi^{L_0}\mc U(\upgamma_1)u,\varpi\big)\btl,
\end{gather*}
considered as elements of $(\Mbb_j\otimes\Mbb_j')[[\xi^{\pm1},q^{\pm 1}]]$ and $(\Mbb_j\otimes\Mbb_j')[[\varpi^{\pm1},q^{\pm 1}]]$ respectively. Then \eqref{eq111} says $A(\xi,q)=B(q/\xi,q)$. Let $C(\xi,\varpi)\in(\Mbb_j\otimes\Mbb_j')[[\xi^{\pm 1},\varpi^{\pm1}]]$  be $A(\xi,\xi\varpi)$, which also equals $B(\varpi,\xi\varpi)$. Since $A(\xi,q)$ contains only non-negative  powers of $q$, so does $A(\xi,\xi\varpi)$ for $\varpi$.   Similarly, since $B(\varpi,q)$ contains only non-negative powers of $q$, so does $B(\varpi,\xi\varpi)$ for $\xi$. Therefore $C(\xi,\varpi)$ is an element in  $(\Mbb_j\otimes\Mbb_j')[[\xi,\varpi]]$, where the latter  can be identified with the subspace $(\Mbb_j\otimes\Mbb_j'\otimes 1)[[\xi,\varpi]]$ of the $R[[\xi,\varpi]]$-module $(\Mbb_j\otimes\Mbb_j'\otimes R)[[\xi,\varpi]]$. Thus $D(\xi,\varpi):=f(\xi,\varpi)C(\xi,\varpi)$ is well-defined as an element in $(\Mbb_j\otimes\Mbb_j'\otimes R)[[\xi,\varpi]]$. It is easy to check that
\begin{align*}
\Res_{\xi=0}\bigg(D(\xi,q/\xi)\frac{d\xi}{\xi}\bigg)=\Res_{\varpi=0}\bigg(D(q/\varpi,\varpi)\frac{d\varpi}{\varpi}\bigg).
\end{align*}
(Indeed, they both equal $\sum_{n\in\Nbb}D_{n,n}q^n$ if we write $D(\xi,\varpi)=\sum_{m,n\in\Nbb}D_{m,n}\xi^m\varpi^n$.) This proves \eqref{eq110}.
\end{proof}

\begin{thm}\label{lb42}
Let $\uppsi\in\scr T_{\wtd{\fk X}}^*(\Wbb_\blt\otimes\Mbb_\blt\otimes\Mbb_\blt')(\wtd{\mc B})$. Then $\wtd{\mc S}\uppsi$ is a formal conformal block.
\end{thm}

\begin{proof}
Step 1. Note that we have $S_{\fk X}=\sum_{i=1}^N\sgm_i(\mc B)$ and $S_{\wtd{\fk X}}=\sum_{i=1}^N\sgm_i(\wtd{\mc B})+\sum_{j=1}^M(\sgm_j'(\wtd{\mc B})+\sgm_j''(\wtd{\mc B}))$. Choose any $v$ in $\pi_*\big(\scr V_{\fk X}\otimes\omega_{\mc C/\mc B}(\blt S_{\fk X})\big)(\mc B)=\big(\scr V_{\fk X}\otimes\omega_{\mc C/\mc B}(\blt S_{\fk X})\big)(\mc C)$.  In this first step, we would like to construct a formal power series expansion
\begin{align}
v=\sum_{n_\blt\in\Nbb^M}v_{n_\blt}q_\blt^{n_\blt}\label{eq114}
\end{align}
where each $v_{n_\blt}$ is in $\big(\scr V_{\wtd{\fk X}}\otimes\omega_{\wtd{\mc C}/\wtd{\mc B}}(\blt S_{\wtd{\fk X}})\big)(\wtd{\mc C})$.  

First, choose any precompact open subset $\wtd U$ of $\wtd{\mc C}$ disjoint from the double points $\sgm_j'(\wtd{\mc B})$ and $\sgm_j''(\wtd{\mc B})$ for all $1\leq j\leq M$. Then for each $j$ one can find small enough positive numbers $\epsilon_j<r_j,\lambda_j<\rho_j$  such that $\wtd U\times\mc D_{\epsilon_\blt\lambda_\blt}$ is an open subset of  $\wtd {\mc C}\times\mc D_{r_\blt\rho_\blt}-\bigcup_j F_j'-\bigcup_j F_j''$ in \eqref{eq112}, and hence an open subset of $\mc C$.   Moreover, by \eqref{eq234}, the projection $\pi:\mc C\rightarrow\mc B$ equals $\wtd \pi\times \id:\wtd {\mc C}\times\mc D_{r_\blt\rho_\blt}\rightarrow \wtd{\mc B}\times \mc D_{r_\blt\rho_\blt}$ when restricted to  $\wtd U\times \mc D_{\epsilon_\blt\lambda_\blt}$. It follows that the section $v|_{\wtd U\times \mc D_{\epsilon_\blt\lambda_\blt}}$  of $\scr V_{\fk X}\otimes\omega_{\mc C/\mc B}(\blt S_{\fk X})$ can be regarded as a section of $\scr V_{\wtd{\fk X}\times\mc D_{r_\blt\rho_\blt}}\otimes\omega_{\wtd{\mc C}\times\mc D_{r_\blt\rho_\blt}/\wtd{\mc B}\times\mc D_{r_\blt\rho_\blt}}(\blt S_{\fk X})$, which, by taking power series expansions at $q_\blt=0$, is in turn an element of $\scr V_{\wtd{\fk X}}\otimes\omega_{\wtd{\mc C}/\wtd{\mc B}}(\blt S_{\wtd{\fk X}})(\wtd U)[[q_\blt]]$. The coefficient before $q_\blt^{n_\blt}$ defines $v_{n_\blt}|_{\wtd U}$. This defines the section $v_{n_\blt}$ of $\scr V_{\wtd{\fk X}}\otimes\omega_{\wtd{\mc C}/\wtd{\mc B}}(\blt S_{\wtd{\fk X}})$ on $\wtd{\mc C}-\bigcup_{j=1}^M(\sgm_j'(\wtd{\mc B})\cup \sgm_j''(\wtd{\mc B}))$ satisfying \eqref{eq114}.

We now show that $v_{n_\blt}$ has poles of orders at most $n_j+1$ at $\sgm_j'(\wtd{\mc B})$ and $\sgm_j''(\wtd{\mc B})$. This will imply that $v_{n_\blt}$ extends to a section of $\scr V_{\wtd{\fk X}}\otimes\omega_{\wtd{\mc C}/\wtd{\mc B}}(\blt S_{\wtd{\fk X}})$ on $\wtd{\mc C}$. Let  $\wtd{\mc B}_j=\mc D_{r_\blt\rho_\blt\setminus j}\times\wtd{\mc B}$, and choose open sets $W_j,W_j',W_j''$  as in the paragraph containing equation \eqref{eq115}. Define coordinates $\xi_j,\varpi_i,q_j$ as in the beginning of Section \ref{lb24}.  Then, by \eqref{eq75} and \eqref{eq80}, $v|_{W_j-\Sigma}$ is a sum of those whose restrictions to $W_j',W_j''$  under the trivializations $\mc U_\varrho(\xi_j),\mc U_\varrho(\varpi_j)$ are
\begin{align}
f(\xi_j,q_j/\xi_j,\cdot)\xi_j^{L_0}u\cdot \frac{d\xi_j}{\xi_j}  \qquad\text{resp.}\qquad -f(q_j/\varpi_j,\varpi,\cdot)\varpi_j^{L_0}\mc U(\upgamma_1)u\cdot \frac{d\varpi_j}{\varpi_j}\label{eq118}
\end{align}
where $u\in\Vbb$ and $f=f(\xi_j,\varpi_j,\cdot)\in\scr O(W_j)$, and the coordinates of $\wtd {\mc B}_j$ are suppressed as the dot. In the above two terms, if we take power series expansions of $q_\blt$, then it is obvious that the coefficients before  $q_\blt^{n_\blt}$ have poles of orders at most $n_j+1$ at $\xi_j=0$ and $\varpi_j=0$ respectively. This proves the claim.

Step 2. By \eqref{eq115}, we can regard $f(\xi_j,\varpi_j,\cdot)$ as an element of $\scr O(\wtd{\mc B}_j)[[\xi_j,\varpi_j]]$, which in turn is an element of $\scr O(\wtd{\mc B})[[q_{\blt\setminus j},\xi_j,\varpi_j]]$. Thus, by Lemma \ref{lb37} (applied to $R=\scr O(\wtd{\mc B})[[q_{\blt\setminus j}]]$) and the fact that $v|_{W_j-\Sigma}$ is a (finite) sum of those of the form \eqref{eq118}, we have the following equation of elements in $(\Mbb_j\otimes\Mbb_j'\otimes \scr O(\wtd{\mc B}))[[q_\blt]]$:
\begin{align}
\sum_{n_\blt\in\Nbb^M}\big(v_{n_\blt}\cdot q_j^{\wtd L_0}\btr\otimes_j\btl+ q_j^{\wtd L_0}\btr\otimes_j~ v_{n_\blt}\cdot\btl\big)q_\blt^{n_\blt} =0\label{eq117}
\end{align}
where the actions of $v_{n_\blt}$ on $\Mbb_j$ and $\Mbb_j'$ are defined by \eqref{eq116} using the  local coordinates $\xi_j,\varpi_j$ of $\wtd{\fk X}$. On the other hand, since $\uppsi$ is  conformal block, for each $n_\blt$, the element $A_{n_\blt}\in\scr O(\wtd{\mc B})$ defined by
\begin{align*}
A_{n_\blt}:=&\uppsi\Big(v_{n_\blt}\cdot w\otimes (q_1^{\wtd L_0}\btr\otimes_1\btl)\otimes\cdots\otimes (q_M^{\wtd L_0}\btr\otimes_M\btl)\Big)\nonumber\\
&+\sum_{j=1}^M\uppsi\Big(w\otimes (q_1^{\wtd L_0}\btr\otimes_1\btl)\otimes\cdots\otimes v_{n_\blt}\cdot q_j^{\wtd L_0}\btr\otimes_j\btl \otimes\cdots\otimes (q_M^{\wtd L_0}\btr\otimes_M\btl)\Big)\nonumber\\
&+\sum_{j=1}^M\uppsi\Big(w\otimes (q_1^{\wtd L_0}\btr\otimes_1\btl)\otimes\cdots\otimes  q_j^{\wtd L_0}\btr\otimes_j~v_{n_\blt}\cdot\btl \otimes\cdots\otimes (q_M^{\wtd L_0}\btr\otimes_M\btl)\Big)
\end{align*}
equals $0$. Here, similarly, the action of $v_{n_\blt}$ on $w$ is defined by \eqref{eq116} and \eqref{eq152} using the local coordinates $\eta_\blt$.  By \eqref{eq117}, we have
\begin{align*}
0=\sum_{n_\blt\in\Nbb^M}A_{n_\blt}q_\blt^{n_\blt}=\sum_{n_\blt\in\Nbb^M}\uppsi\Big(v_{n_\blt}\cdot w\otimes (q_1^{\wtd L_0}\btr\otimes_1\btl)\otimes\cdots\otimes (q_M^{\wtd L_0}\btr\otimes_M\btl)\Big)q_\blt^{n_\blt},
\end{align*}
which is exactly $\wtd{\mc S}\uppsi(v\cdot w)$. This finishes the proof that $\wtd{\mc S}\uppsi$ is a formal conformal block.
\end{proof}

\section{Propagation of conformal blocks}\label{lb149}

\subsection*{$\wr\upphi$ as a conformal block}

Let $\fk X=(C;x_1,\dots,x_N)$ be an $N$-pointed compact Riemann surface. Recall the divisor $S_{\fk X}=x_1+\cdots+x_N$. As in previous sections, we write $C-\{x_1,\dots,x_N\}$ as $C-S_{\fk X}$ for brevity. Then the projection onto the second component $C\times(C-S_{\fk X})\rightarrow (C-S_{\fk X})$ is the family over $C-S_{\fk X}$ with constant  fiber $C$. We understand each $x_i$ as the constant section $x_i:C-S_{\fk X}\rightarrow C\times (C-S_{\fk X})$, i.e., its value at each $y\in C-S_{\fk X}$ is $(x_i,y)$. Let
\begin{align*}
\iota:C-S_{\fk X}\rightarrow C\times (C-S_{\fk X})
\end{align*}
be the diagonal map, i.e., sending each $y$ to $(y,y)$. We thus have a family of $(N+1)$-pointed curve \index{X@$\wr\fk X,\wr\wr\fk X$}
\begin{gather}
\wr\fk X=(C\times(C-S_{\fk X})\rightarrow (C-S_{\fk X});\iota,x_1,\dots,x_N).
\end{gather}
Let $\Wbb_1,\dots,\Wbb_N$ be $\Vbb$-modules associated to the $N$ points of $\fk X$. Choose $\upphi\in\scr T_{\fk X}^*(\Wbb_\blt)$. Our first goal in this section is to show that the $\wr\upphi$  in \eqref{eq119} can be identified naturally with an element of $\scr T_{\wr\fk X}^*(\Vbb\otimes\Wbb_\blt)(C-S_{\fk X})$. 

We make the identification of $\scr O_{C-S_{\fk X}}$-modules
\begin{align}
\scr W_{\wr\fk X}(\Vbb\otimes\Wbb_\blt)\simeq \scr V_{C-S_{\fk X}}\otimes_\Cbb\scr W_{\fk X}(\Wbb_\blt)\label{eq122}
\end{align}
as follows. Choose any open subset $V\subset C-S_{\fk X}$ and a univalent map $\mu:V\rightarrow \Cbb$, and choose local coordinates $\eta_1,\dots,\eta_N$ of $x_1,\dots,x_N$ respectively. Then we have a local coordinate $\varrho(\mu)$ of $\iota(V)$ defined on $V\times V$ to be
\begin{align}
\varrho(\mu)_y(x)\equiv\varrho(\mu)(x,y)=\mu(x)-\mu(y)\label{eq128}
\end{align}
for any $x,y\in V$. One can also regard each $\eta_i$ as the local coordinate of $x_i(C-S_{\fk X})$ constant over the base $C-S_{\fk X}$. Write $(\varrho(\mu),\eta_\blt)=(\varrho(\mu),\eta_1,\dots,\eta_N)$. We then have trivialization
\begin{align*}
\mc U(\varrho(\mu),\eta_\blt):\scr W_{\wr\fk X}(\Vbb\otimes\Wbb_\blt)|_V\xrightarrow{\simeq} \Vbb\otimes_\Cbb\Wbb_\blt\otimes_\Cbb\scr O_V
\end{align*}
as in \eqref{eq120}. On the other hand, we have $\mc U_\varrho(\mu):\scr V_C|_V\xrightarrow{\simeq}\Vbb\otimes_\Cbb\scr O_V$ and $\mc U(\eta_\blt):\scr W_{\fk X}(\Wbb_\blt)\xrightarrow{\simeq}\Wbb_\blt$ as in \eqref{eq90} and \eqref{eq121}, which give
\begin{align*}
\mc U_\varrho(\mu)\otimes \mc U(\eta_\blt):\scr V_C|_V\otimes_\Cbb\scr W_{\fk X}(\Wbb_\blt)\xrightarrow{\simeq} \Vbb\otimes_\Cbb\Wbb_\blt\otimes_\Cbb\scr O_V.
\end{align*}
Then the isomorphism \eqref{eq122} on $V$ is given by
\begin{align}
\big(\mc U_\varrho(\mu)\otimes \mc U(\eta_\blt)\big)^{-1}\mc U(\varrho(\mu),\eta_\blt):\scr W_{\wr\fk X}(\Vbb\otimes\Wbb_\blt)|_V\xrightarrow{\simeq} \scr V_C|_V\otimes_\Cbb\scr W_{\fk X}(\Wbb_\blt).\label{eq129}
\end{align}
We identify the above two sheaves of modules via the above map, so that we also have identification
\begin{align}
\mc U(\varrho(\mu),\eta_\blt)=\mc U_\varrho(\mu)\otimes \mc U(\eta_\blt).\label{eq126}
\end{align}
Using \eqref{eq55}, \eqref{eq123}, and \eqref{eq124}, it is not hard to see that this isomorphism is independent of the choice of $\mu$ and $\eta_\blt$. Thus \eqref{eq122} can be defined globally.

Recall from \eqref{eq119} and Remark \ref{lb40} that for each $w\in\scr W_{\fk X}(\Wbb_\blt)$, $\wr\upphi(w)$ is an element of $H^0(C-S_{\fk X},\scr V_C^*)=\scr V_C^*(C-S_{\fk X})=\Hom_{\scr O_{C-S_{\fk X}}}(\scr V_{C-S_{\fk X}},\scr O_{C-S_{\fk X}})$ whose evaluation $\wr\upphi(v,w_\blt)$ with any $v\in\scr V_C(V)$ (where $V\subset C-S_{\fk X}$ is open) is also written as $\wr\upphi(v\otimes w_\blt)$. This notation suggests that we regard $\wr\upphi$ as an element of $\Hom_{\scr O_{C-S_{\fk X}}}(\scr V_{C-S_{\fk X}}\otimes_\Cbb\scr W_{\fk X}(\Wbb_\blt),\scr O_{C-S_{\fk X}})$, which, through the isomorphism \eqref{eq122}, becomes a homomorphism of $\scr O_{C-S_{\fk X}}$-modules \index{zz@$\wr\upphi$}
\begin{align*}
\wr\upphi:\scr W_{\wr\fk X}(\Vbb\otimes\Wbb_\blt)\rightarrow\scr O_{C-S_{\fk X}}.
\end{align*}
We would like to show that $\wr\upphi$ is a conformal block.

\begin{thm}[Propagation of conformal blocks]\label{lb126}
For any $\upphi\in\scr T_{\fk X}^*(\Wbb_\blt)$, $\wr\upphi$ is an element of $\scr T_{\wr\fk X}^*(\Vbb\otimes\Wbb_\blt)$. 
\end{thm}

\begin{proof}
Assume without loss of generality that $C$ is connected. Then so is the base manifold $C-S_{\fk X}$ of the family $\wr\fk X$. Choose any of $x_1,\dots,x_N$, say $x_1$. Choose local coordinates $\eta_\blt$ of $x_\blt$ such that  $\eta_1$ is defined on a neighborhood $V\ni x_1$. Assume that under the coordinate $\eta_1$, $V$ is an open disc centered at $x_1$ with radius $r>0$. Identify $\scr W_{\fk X}(\Wbb_\blt)$ with $\Wbb_\blt$ via the triviliazation $\mc U(\eta_\blt)$ (which means we set $\mc U(\eta_\blt)=\id$). By Proposition \ref{lb41}, it suffices to prove that the restriction of $\wr\upphi$ to $V-\{x_1\}$ is a conformal block. Identify $V$ with an open subset of $\Cbb$ via $\eta_1$, which makes  $\eta_1$ equal to the standard coordinate $z$ of $\Cbb$. We also have $\mc U_\varrho(\eta_1)=\mc V_\varrho(\eta_1)$. Choose any $w_\blt\in W_\blt$ and $v\in\Vbb$. By the fact that $\wr\upphi$ equals $\wr\upphi_{x_1}$ (defined by \eqref{eq125}) near $x_1$, we have
\begin{align}
&~~~~\wr\upphi\big(\mc U(\varrho(\eta_1),\eta_\blt)^{-1}(v\otimes w_\blt)\big)\xlongequal{\eqref{eq126}}\wr\upphi\big(\mc U_\varrho(\eta_1)^{-1}v, w_\blt\big)\nonumber\\
&=\upphi\big(Y_{\Wbb_1}(v,z)w_1\otimes w_2\otimes \dots\otimes w_N\big)\nonumber\\
&\xlongequal{\eqref{eq107}}\upphi\big(z^{\wtd L_0}Y_{\Wbb_1}(z^{-L_0}v,1)z^{-\wtd L_0}w_1\otimes w_2\otimes \dots\otimes w_N\big).\label{eq127}
\end{align}

Define
\begin{align*}
\wtd{\fk Y}=(\Pbb^1\sqcup C;0,1,x_2,\dots,x_N;\infty;x_1;z,z-1,\eta_2,\dots,\eta_N;z^{-1};\eta_1).
\end{align*}
Its sewing (near $\infty\in\Pbb^1,x_1\in C$ controlled by $|z^{-1}|<1,|\eta_1|<r$) gives a family $\fk Y$ with base manifold the open disc $\mc D_r$ with radius $r$ (which is equivalent to $V$). Associate $\Wbb_1,\Vbb,\Wbb_2,\dots,\Wbb_N$ to $0,1,x_2,\dots,x_N$ and $\Wbb_1',\Wbb_1$ to $\infty,x_1$ respectively. Let $\uppsi$ be a linear functional on $\Wbb_1\otimes\Vbb\otimes\Wbb_2\otimes\dots\otimes\Wbb_N\otimes \Wbb_1'\otimes \Wbb_1$, for which we switch the order to $\Wbb_1\otimes\Wbb_2\otimes\cdots\otimes\Wbb_N\otimes \Wbb_1\otimes\Vbb\otimes\Wbb_1'$ (where the order of two $\Wbb_1$ are also switched), is defined by $\upphi\in\Wbb_\blt^*$ (a linear functional on the first $N$ components) tensor the conformal block on $(\Pbb;0,1,\infty)$ (which is a linear functional on the last three components) defined as in example \ref{lb44}. Then the sewing $\wtd{\mc S}\uppsi$ (which is a formal conformal block by Theorem \ref{lb42}), evaluated with the constant section $v\otimes w_\blt$ (under the local coordinates $z,z-1,\eta_2,\dots,\eta_N$ constant with respect to sewing), is
\begin{align*}
\upphi\big(q^{\wtd L_0}Y_{\Wbb_1}(v,1)w_1\otimes w_2\otimes \dots\otimes w_N\big)\qquad (\in\Cbb[[q]]).
\end{align*}
This expression is an element of $\scr O(\mc D_r^\times)$ (and hence of $\scr O(\mc D_r)$) since \eqref{eq127} is an element of $\scr O(\mc D_r^\times)$. If we scale the local coordinates of $0,1\in\Pbb^1$ by $q^{-1}$, then the above expression becomes \eqref{eq127} (with $z$ replaced by $q$), and the restriction $\fk Y_{\mc D_r^\times}$ of $\fk Y$ to the punctured disc $\mc D_r^\times$, including the sections and the local coordinates, is exactly  $\wr\fk X|_{V-\{x_1\}}$. This shows that, after scaling the coordinates, $\wtd {\mc S}\uppsi$ becomes exactly $\wr\upphi$. So $\wtd {\mc S}\uppsi$ converges a.l.u..    By Proposition \ref{lb43}-2, $\wtd{\mc S}\uppsi$ is an actual conformal block associated to the restricted family $\fk Y_{\mc D_r^\times}$. So $\wr\upphi|_{V-\{x_1\}}$ is a conformal block associated to $\wr\fk X|_{V-\{x_1\}}$. This finishes the proof. 
\end{proof}

\begin{co}\label{lb127}
For each $y\in C-\SX$, if we set $\wr\fk X_y=(C;y,x_1,\dots,x_N)$, then the value of $\wr\upphi$ at $y$, written as $\wr\upphi|_y$, is an element of $\scr T_{\wr\fk X_y}^*(\Vbb\otimes\Wbb_\blt)$.
\end{co}

\begin{rem}
We describe the explicit form of $\wr\upphi|_y$ as an element of $\scr T_{\wr\fk X_y}^*(\Vbb\otimes\Wbb_\blt)$. Choose  local coordinates  $\eta_\blt$ at $x_\blt$. Choose  a univalent map $\mu$ defined on a neighborhood $V$ of $y$. Then, $\varrho(\mu)_y$ (defined in \eqref{eq128}) is a local coordinate at $y$.  Identify $\mc U(\eta_\blt):\scr W_{\fk X}(\Wbb_\blt)\xrightarrow{\simeq}\Wbb_\blt$.  Then for any $v\in\Vbb$ and $w_\blt\in\Wbb_\blt$,
\begin{align}
\wr\upphi\big|_y\big(\mc U(\varrho(\mu)_y)^{-1}v\otimes w_\blt\big)=\wr\upphi\big(\mc U_\varrho(\mu)^{-1}v,w_\blt\big)\big|_y,
\end{align}
where we recall that $\mc U_\varrho(\mu)^{-1}v$ is in $\scr V_C(V)$.

Note that Theorem \ref{lb126} is stronger than Corollary \ref{lb127}, in that it says also that $\wr\upphi|_y$ varies holomorphically over $y$.
\end{rem}

\subsection*{Double propagation}

We now want to propagate the conformal block $\wr\upphi$. For any two sets $A,B$, we set \index{Conf@$\Conf$}
\begin{gather}
\Conf(A,B)=\{(a,b)\in A\times B:a\neq b\},\qquad \Conf_2(A)=\Conf(A,A).\label{eq135}
\end{gather}
Fix trivialization $\scr W_{\fk X}(\Wbb_\blt)\simeq\Wbb_\blt$ via $\mc U(\eta_\blt)$. Let $V\subset C-S_{\fk X}$ be open. For each $v\in\scr V_C(V)$ and $w_\blt\in\Wbb_\blt$, $v\otimes w_\blt$ can be regarded as an element of $\scr W_{\wr\fk X}(\Vbb\otimes\Wbb_\blt)(V)$ by \eqref{eq129}. Therefore, by Theorem \ref{lb45}, we have $\wr\wr\upphi(v\otimes w_\blt)$, written as $\wr\wr\upphi(v,w_\blt)$ in the following, \index{zz@$\wr\upphi$!$\wr\wr\upphi,\wr^n\upphi$} sending each $u\in\scr V_C(U)$ (where $U$ is an open subset of $C-S_{\fk X}$ disjoint from $V$) to an element $\wr\wr\upphi(u,v,w_\blt)\in\scr O(U\times V)$. (To apply Theorem \ref{lb45}, we extend $u$ to an element of $\scr V_{\wr\fk X}(U\times V)$ constantly over $V$.) This map is compatible with the restrictions to subsets of $U$ and of $V$. Using this compatibility, for any (non-necessarily disjoint) open subsets $U,V\subset C$, $\wr\wr\upphi(w_\blt)$ can be extended to a homomorphism of $\scr O(U)$-$\scr O(V)$ bimodules
\begin{align}
\wr\wr\upphi(w_\blt):\scr V_C(U)\otimes_\Cbb \scr V_C(V)\rightarrow\scr O_{\Conf_2(C-S_{\fk X})}(\Conf(U-S_{\fk X},V-S_{\fk X}))
\end{align}
compatible with the restrictions to subsets of $U$ and $V$, such that for any $u\in\scr V_C(U),v\in\scr V_C(V)$ and open subsets $U_0\subset U,V_0\subset V$ satisfying $U_0\cap V_0=\emptyset$, $\wr\wr\upphi(u,v,w_\blt)|_{U_0\times V_0}$ is the element $\wr\wr\upphi(u|_{U_0-S_{\fk X}},v|_{V_0-S_{\fk X}},w_\blt)$ describe above (which is in $\scr O_{\Conf_2(C-S_{\fk X})}((U_0-S_{\fk X})\times (V_0-S_{\fk X}))$). For brevity, such compatibility is summarized by saying that $\wr\wr\upphi(w_\blt)$ is a homomorphism of $\scr O_C\boxtimes\scr O_C$-modules
\begin{align}
\boxed{~~\wr\wr\upphi(w_\blt):\scr V_C\boxtimes\scr V_C\rightarrow\scr O_{\Conf_2(C-S_{\fk X})}~~}
\end{align}
Similar to this description, we can regard $\wr\upphi(w_\blt)$ as a homomorphism of $\scr O_C$-modules
\begin{align}
\boxed{~~\wr\upphi(w_\blt):\scr V_C\rightarrow\scr O_{C-S_{\fk X}}~~}
\end{align}
whose value at each $v\in\scr O(V)$ is equal to the one at $v|_{V-S_{\fk X}}$. 

The next theorem is just the restatement of the description of $\wr\upphi$ in Theorem \ref{lb32}. We assume that for each $i=1,\dots,N$, $\eta_i$ is defined on $W_i$, and that $W_i\cap W_j=\emptyset$ if $i\neq j$. Let $z$ be the standard coordinate of $\Cbb$.
\begin{thm}\label{lb47}
Choose any  $w_\blt\in\Wbb_\blt$.  Choose $V\subset W_i$ an open disc centered at $x_i$ (under the coordinate $\eta_i$), identify $V$ with a neighborhood of $0\in\Cbb$ via $\eta_i$, and identify $\scr V_V\simeq \Vbb\otimes_\Cbb\scr O_V$ via the trivialization $\mc U_\varrho(\eta_i)=\mc V_\varrho(\eta_i)$. Choose $v\in\scr V_C(V)$, and choose  $y\in V-S_{\fk X}=V-\{x_i\}$. Then
\begin{align}
\wr\upphi(v,w_\blt)\big|_y=\upphi\big(w_1\otimes\cdots\otimes Y_{\Wbb_i}(v,z)w_j\otimes\cdots\otimes w_N\big)\big|_{z=\eta_i(y)}.
\end{align}
Moreover, we have
\begin{align}
\wr\upphi(\id,w_\blt)|_y=\upphi(w_\blt).\label{eq130}
\end{align}
\end{thm}

In this theorem, $\id\in \scr V_C(C)$ is the vacuum section defined in Remark \ref{lb46}. Since $Y_{\Wbb_i}(\id,z)=\id_{\Wbb_i}$, \eqref{eq130} is clearly true when $y$ is in a neighborhood of $x_1,\dots,x_N$. Thus \eqref{eq130} is true for any $y\in C$.

We now generalize this theorem to $\wr\wr\upphi$. 

\begin{thm}\label{lb48}
Choose  any $w_\blt\in\Wbb_\blt$. Choose any $U,V$ open subsets of $C$ with $\eta:U\rightarrow\Cbb,\mu:V\rightarrow\Cbb$ univalent maps,  identify $\scr V_U\simeq \Vbb\otimes_\Cbb\scr O_U,\scr V_V\simeq \Vbb\otimes_\Cbb\scr O_V$ via trivializations $\mc U_\varrho(\eta)=\mc V_\varrho(\eta),\mc U_\varrho(\mu)=\mc V_\varrho(\mu)$ respectively. Choose $u,v$   in $\scr V_C(U),\scr V_C(V)$ respectively, and choose $x\in U-S_{\fk X},y\in V-S_{\fk X}$ satisfying $x\neq y$. Then the following are true.
\begin{enumerate}[label=(\arabic*)]
\item If $U$ is an open disc in $W_i$ centered at $x_i$ (under the coordinate $\eta_i$) and does not contain $y$, and if $\eta=\eta_i$. Then
\begin{align}
\wr\wr\upphi(u,v,w_\blt)\big|_{x,y}=\wr\upphi\big(v,w_1\otimes\cdots\otimes Y_{\Wbb_i}(u,z)w_i\otimes\cdots\otimes w_N\big)\big|_y~\big|_{z=\eta_i(x)}
\end{align}
where the series of $z$ on the right hand side converges absolutely, and we regard $u\in\Vbb\otimes\Cbb[[z]]$ by taking Taylor series expansion of the variable $\eta_i$ at $x_i$.
\item If $U=V$ and do not contain $x_1,\dots,x_N$, if $\eta=\mu$, and if $U$ contains the closed disc with center $y$ and radius $|\eta(x)-\eta(y)|$ (under the coordinate $\eta$), then
\begin{align}
\wr\wr\upphi(u,v,w_\blt)\big|_{x,y}=\wr\upphi\big(Y(u,z)v,w_1\otimes\cdots\otimes w_N\big)\big |_y~\big |_{z=\eta(x)-\eta(y)}
\end{align}
where the series of $z$ on the right hand side converges absolutely, and we regard $u\in\Vbb\otimes\Cbb[[z]]$ by taking Taylor series expansion of the variable $\eta-\eta(y)$ at $y$.
\item  We have
\begin{align}
\wr\wr\upphi(\id,v,w_\blt)=\wr\upphi(v,w_\blt).
\end{align}
\item We have
\begin{align}
\wr\wr\upphi(u,v,w_\blt)\big|_{x,y}=\wr\wr\upphi(v,u,w_\blt)\big|_{y,x}.
\end{align}
\end{enumerate}
\end{thm}

\begin{proof}
It is easy to see that
\begin{align}
\wr\wr\upphi_{x,y}=\wr(\wr\upphi|_y)|_x.
\end{align}
Thus (1) (2) (3) follow directly from Theorem \ref{lb47} and relation \eqref{eq130} (with $\upphi$ replaced by $\wr\upphi|_y$). We now prove (4). It suffices to assume that $C$ is connected.

Assume first of all that $N>1$. Let $U$ and $V$ be open discs in $W_1,W_2$ centered at $x_1,x_2$ and identified with open subsets of $\Cbb$ via $\eta_1,\eta_2$ respectively. (Note that under this identification, we have $x_1=0$ and $x_2=0$.) Let $\zeta$ be also the standard coordinate of $\Cbb$. Then from (1) and Theorem \ref{lb47}, $\wr\wr\upphi(u,v,w_\blt)\big|_{x,y}$ equals the evaluation of 
\begin{align}
g(z,\zeta):=\upphi\big(Y_{\Wbb_1}(u,z)w_1\otimes Y_{\Wbb_2}(v,\zeta)w_2\otimes w_3\otimes\cdots\otimes w_N\big)\label{eq131}
\end{align}
(which is an element of $\Cbb((z,\zeta))$) first at $\zeta=\eta_2(y)$ and then at $z=\eta_1(x)$. By varying $x$ and $y$, $\wr\wr\upphi(u,v,w_\blt)$ is clearly a two-variable holomorphic function $f(z,\zeta)$ on $(U-\{0\})\times (V-\{0\})$. Thus, we have for any $z_0\in U-\{0\},\zeta_0\in V-\{0\}$ that $f(z,\zeta)|_{\zeta=\zeta_0}|_{z=z_0}=g(z,\zeta)|_{\zeta=\zeta_0}|_{z=z_0}$.

By taking Laurant series expansions, we may regard $f(z,\zeta)$ as an element of $\Cbb[[z^{\pm 1},\zeta^{\pm 1}]]$.  By applying $\Res_{\zeta=0}\Res_{z=0}(\cdots)z^m\zeta^ndzd\zeta$ to $f$ and $g$ for any $m,n\in\Zbb$ (note the order of the two residues), we see that $f(z,\zeta)$ and $g(z,\zeta)$ can be regarded as identical elements of $\Cbb[[z^{\pm 1},\zeta^{\pm1}]]$. Since $g$ is in $\Cbb((z,\zeta))$, so is $f$. Since the double series $f(z,\zeta)$  converges absolutely when $z\in U-\{0\}$ and $\zeta\in V-\{0\}$, so does $g(z,\zeta)$. Therefore, the evaluations of \eqref{eq131} at $\zeta=\eta_2(y)$ and at $z=\eta_1(x)$ commute. This proves (4) for the above chosen $U,V$. By analytic continuation, one may use the argument in the proof of Proposition \ref{lb41} to show (4) when $U$ is as above and $V$ is any open subset of $C$. Another application of this argument proves (4) for any open $U,V\subset C$.

Finally, we assume that $N=1$. Then, by \eqref{eq130}, $\upphi$ is the restriction of a conformal block $\uppsi$ on $C$ with two marked points. Since (4) is true for $\uppsi$, it holds also for $\upphi$.
\end{proof}

\subsection*{Multiple propagations}

Although single and double propagations are sufficient for proving the main results of this monograph, it would be interesting to generalize them to multiple propagation.

Choose $n\in\Nbb$ and define \index{Conf@$\Conf$}
\begin{gather*}
\Conf(A_1,\dots,A_n)=\{(a_1,\dots,a_n)\in A_1\times\cdots A_n:a_i\neq a_j\text{ for any }1\leq i<j\leq n\},\\
\Conf_n(A)=\Conf(A,A,\dots,A).
\end{gather*}
One can apply the propagation $n$-times  to obtain, for each $w_\blt\in\Wbb_\blt$, a homomorphism of $\scr O_C^{\boxtimes n}$-modules
\begin{align}
\wr^n\upphi(w_\blt):\scr V_C^{\boxtimes n}\rightarrow\scr O_{\Conf_n(C-S_{\fk X})},
\end{align}
which means that for any open subsets $U_1,\dots,U_n$ of $C$, we have  an $\scr O(U_1)\otimes_\Cbb\cdots\otimes_\Cbb(U_n)$-module \index{zz@$\wr\upphi$!$\wr\wr\upphi,\wr^n\upphi$} homomorphism
\begin{align*}
\wr^n\upphi(w_\blt):\scr V_C(U_1)\otimes_\Cbb\cdots\otimes_\Cbb\scr V_C(U_n)\rightarrow\scr O_{\Conf_n(C-S_{\fk X})}(\Conf(U_1-S_{\fk X},\dots,U_n-S_{\fk X}))
\end{align*}
compatible with respect to restrictions. We have the following generalization of Theorem \ref{lb48}. Recall again that $W_1,\dots,W_N$ are mutually disjoint neighborhoods of $x_1,\dots,x_N$ on which $\eta_1,\dots,\eta_N$ are defined respectively.

\begin{thm}\label{lb92}
Choose any $w_\blt\in\Wbb_\blt$. For each $1\leq k\leq n$, choose an open subset  $U_k$  of $C$ equipped with a univalent map $\mu_k:U_k\rightarrow\Cbb$, identify $\scr V_{U_k}\simeq \Vbb\otimes_\Cbb\scr O_{U_k}$ via trivialization $\mc U_\varrho(\mu_k)=\mc V_\varrho(\mu_k)$,  choose $v_k\in\scr V_C(U_k)=\Vbb\otimes_\Cbb\scr O(U_k)$, and choose $y_k\in U_k-S_{\fk X}$ satisfying $y_j\neq y_k$ for any $1\leq j<k\leq n$. Then the following are true.
\begin{enumerate}[label=(\arabic*)]
\item If $U_1$ is an open disc of $W_i$ centered at $x_i$ (under the coordinate $\eta_i$) and does not contain $y_2,\dots,y_n$, and if $\mu_1=\eta_i$, then
	\begin{align}
	&\wr^n\upphi(v_1,v_2,\dots,v_n,w_\blt)\big|_{y_1,y_2,\dots,y_n}\nonumber\\
	=&\wr^{n-1}\upphi\big(v_2,\dots,v_n,w_1\otimes\cdots\otimes Y_{\Wbb_i}(v_1,z)w_i\otimes\cdots\otimes w_N\big)\big|_{y_2,\dots,y_n}~\big|_{z=\eta_i(y_1)}
	\end{align}
where the series of $z$ on the right hand side converges absolutely, and $v_1$ is considered as an element of $\Vbb\otimes\Cbb((z))$ by taking Taylor series expansion with respect to the variable $\eta_j$ at $x_j$.
\item If $U_1=U_2$ and do not contain $x_1,\dots,x_N,y_3,\dots,y_n$, if $\mu_1=\mu_2$, and if $U_2$ contains the closed disc with center $y_2$ and radius $|\mu_2(y_1)-\mu_2(y_2)|$ (under the coordinate $\mu_2$), then
	\begin{align}
	&\wr^n\upphi(v_1,v_2,v_3,\dots,v_n,w_\blt)\big|_{y_1,y_2,\dots,y_n}\nonumber\\
	=&\wr^{n-1}\upphi\big(Y(v_1,z)v_2,v_3,\dots,v_n,w_1\otimes\cdots\otimes w_N\big)\big |_{y_2,\dots,y_n}~\big |_{z=\mu_2(y_1)-\mu_2(y_2)}
	\end{align}
where the series of $z$ on the right hand side converges absolutely, and $v_1$ is considered as an element of $\Vbb\otimes\Cbb((z))$ by taking Taylor series expansion with respect to the variable $\mu_2-\mu_2(y_2)$ at $y_2$.
	\item  We have
	\begin{align}
\wr^n\upphi(\id,v_2,v_3,\dots,v_n,w_\blt)=\wr^{n-1}\upphi(v_2,\dots,v_n,w_\blt).
	\end{align}
	\item For any permutation $\upsigma$ of the set $\{1,2,\dots,n\}$, we have
	\begin{align}
	\wr^n\upphi(v_{\upsigma(1)},\dots,v_{\upsigma(n)},w_\blt)\big|_{y_{\upsigma(1)},\dots,y_{\upsigma(n)}}=\wr^n\upphi(v_1,\dots,v_n,w_\blt)\big|_{y_1,\dots,y_n}.
	\end{align}
\end{enumerate}
\end{thm}

\begin{proof}
(1) (2) (3) follow from
\begin{align}
\wr^n\upphi|_{y_1,y_2,\dots,y_n}=\wr(\wr^{n-1}\upphi|_{y_2,\dots,y_n})|_{y_1}.
\end{align}
If $\upsigma$ fixes $3,\dots,n$ and exchanges $1,2$, then (4) follows from Theorem \ref{lb48}-(4). If $1\leq k<n$ and $\upsigma$ exchanges $k,k+1$ and fixes the others, then the first two components of $\wr^{n-k+1}\upphi$ are exchangeable. Thus, by propagating $\wr^{n-k+1}\upphi$ for $k-1$ times, we see that the $k$-th and the $(k+1)$-th components of $\wr^n\upphi$ are exchangeable. This proves (4) in general.
\end{proof}

\section{A commutator formula}

The results of this section will be used in Section \ref{lb84} to define a logarithmic connection on sheaves of conformal blocks.

Let $\fk X=(\pi:\mc C\rightarrow\mc B;\sgm_1,\dots,\sgm_N)$ be a family of $N$-pointed compact Riemann surfaces. We are going to define a sheaf action $\scr V_{\fk X}\otimes\omega_{\mc C/\mc B}(\blt S_{\fk X})\curvearrowright \scr V_{\fk X}(\blt S_{\fk X})$ which is $\scr O_{\mc B}$-linear on $\scr V_{\fk X}\otimes\omega_{\mc C/\mc B}(\blt S_{\fk X})$ and $\scr O_{\mc C}$-linear on $\scr V_{\fk X}(\blt S_{\fk X})$. In other words, we shall define a homomorphism of $\scr O_{\mc B}$-modules
\begin{align*}
\mbf L:\scr V_{\fk X}\otimes\omega_{\mc C/\mc B}(\blt S_{\fk X})\rightarrow\underline{\End}_{\scr O_{\mc C}}(\scr V_{\fk X}(\blt S_{\fk X})).
\end{align*}
Thus, for any open subset $W\subset\mc C$, any element of $\scr V_{\fk X}\otimes\omega_{\mc C/\mc B}(\blt S_{\fk X})(W)$ gives an element of $\End_{\scr O_W}(\scr V_{\fk X}(\blt S_{\fk X})|_W)$.

Choose open $W\subset\mc C$  together with $\eta\in\scr O(W)$ univalent on each fiber. We assume that $W$ is small enough such that $\pi(W)$ has coordinates $\tau_\blt=(\tau_1,\tau_2,\dots)$. Write also $\tau_\blt\circ\pi$ as $\tau_\blt$ for simplicity. Identify $W$ with an open subset of $\Cbb\times\mc B$ via $(\eta,\tau_\blt)$, and identify 
\begin{gather*}
\scr V_{\fk X}|_W\simeq \Vbb\otimes_\Cbb\scr O_W,\qquad \scr V_{\fk X}\otimes\omega_{\mc C/\mc B}|_W\simeq \Vbb\otimes_\Cbb\omega_{\mc C/\mc B}|_W
\end{gather*}
via the trivialization $\mc U_\varrho(\eta)=\mc V_\varrho(\eta)$. Let $z$ be the standard coordinates of $\Cbb$ (which is identified with $\eta$). Then for any $udz=u(z,\tau_\blt)dz$ in $\Vbb\otimes_\Cbb\omega_{\mc C/\mc B}(\blt S_{\fk X})(W)$, open subset $U\subset W$, and $v=v(z,\tau_\blt)$ in $\Vbb\otimes_\Cbb\scr O_W(\blt S_{\fk X})(U)$, we define the action of $udz$ on $v$,  written as  $\mbf L_{udz} v=(\mbf L_{udz} v)(z,\tau_\blt)$, to be
\begin{align}
\mbf L_{udz}~v=\Res_{\zeta-z=0}~Y\big(u(\zeta,\tau_\blt),\zeta-z\big)v(z,\tau_\blt) d\zeta \label{eq148}
\end{align}
where $\zeta$ is another distinct standard complex variable of $\Cbb$. Note that $u,v$ are $\Vbb$-valued meromorphic functions on $W$ with possible poles at $\sgm_1(\mc B),\dots,\sgm_N(\mc B)$. That $\mbf L_{udz}v$ has finite poles at $\SX$ follows from the easy calculation
\begin{align}
\mbf L_{udz}~v=\sum_{n\geq 0}\frac 1{n!}Y\big(\partial_z^n u(z,\tau_\blt) \big)_n v(z,\tau_\blt)\label{eq235}
\end{align}
where the sum is finite by the lower truncation property.

Similar to \eqref{eq116}, the definition of this action is independent of the choice of $\eta$ thanks to Theorem \ref{lb31}. Thus,  it can be extended to $\scr V_{\fk X}\otimes\omega_{\mc C/\mc B}(\blt S_{\fk X})\curvearrowright \scr V_{\fk X}(\blt S_{\fk X})$. (We will not this fact; in Section \ref{lb84}, we shall only use the local expression of $\mbf L$ as in \eqref{eq235}.) By tensoring with the identity map of $\omega_{\mc C/\mc B}$, we get a homomorphism of $\scr O_{\mc B}$-modules
\begin{align*}
\mbf L:\scr V_{\fk X}\otimes\omega_{\mc C/\mc B}(\blt S_{\fk X})\rightarrow\underline{\End}_{\scr O_{\mc C}}(\scr V_{\fk X}\otimes\omega_{\mc C/\mc B}(\blt S_{\fk X}))
\end{align*}
whose local expression under $\eta$ is
\begin{align}
\mbf L_{udz}~vdz=&\Big(\Res_{\zeta-z=0}~Y\big(u(\zeta,\tau_\blt),\zeta-z\big)v(z,\tau_\blt) d\zeta\Big)dz\nonumber\\
=&\sum_{n\geq 0}\frac 1{n!}Y\big(\partial_z^n u(z,\tau_\blt) \big)_n v(z,\tau_\blt)dz.\label{eq149}
\end{align}

Now, we assume that $\mc B$ is small enough such that we have a family of $N$-pointed compact Riemann surfaces with local coordinates $\fk X=(\pi:\mc C\rightarrow\mc B;\sgm_1,\dots,\sgm_N;\eta_1,\dots,\eta_N)$. Let $\Wbb_1,\dots,\Wbb_N$ be $\Vbb$-modules associated to the $N$-points. Let now $W$ be a neighborhood of $\sgm_i(\mc B)$ on which $\eta_i$ is defined,  set $\eta=\eta_i$, and take the identifications mentioned above. Recall that by \eqref{eq116}, for each $w_i\in\Wbb_i\otimes_\Cbb\scr O(\mc B)$ and $vdz\in\Vbb\otimes_\Cbb\omega_{\mc C/\mc B}(\blt S_{\fk X})(W)$,
\begin{align*}
vdz\cdot w_i=\Res_{z=0}Y_{\Wbb_i}(v,z)w_idz.
\end{align*}

The following is the main result of this section.

\begin{pp}\label{lb51}
For any $udz,vdz\in\Vbb\otimes_\Cbb\omega_{\mc C/\mc B}(\blt S_{\fk X})(W)$ and $w_i\in\Wbb_i\otimes_\Cbb\scr O(\mc B)$,
\begin{align*}
udz\cdot vdz\cdot w_i-vdz\cdot udz\cdot w_i=(\mbf L_{udz}vdz)\cdot w_i.
\end{align*}
\end{pp}
Note that the same identity holds when $w_i$ is replaced by any $w\in\Wbb_\blt\otimes_\Cbb\scr O(\mc B)$.
\begin{proof}
We write $w_i$ as $w$ for simplicity. Since the actions of $udz,vdz$ and the definition of $\mbf L$ can be defined fiberwisely, it suffices to assume that $\mc B$ is a single point. We thus suppress the symbol $\tau_\blt$. Note that $W$ is identified with $\eta(W)$ under $\eta$. We assume that $W=\eta(W)$ is an open disc. Choose any $w'\in\Wbb_i'$. Then there exists a holomorphic function $f=f(z,\zeta)$ on $\Conf_2(W-\{0\})$ (recall \eqref{eq135}) such that for each fixed $z$, the series expansion (with respect to the variable $\zeta$) near $0,z$ are respectively
\begin{gather*}
\alpha_z(\zeta)=\bk{w',Y_{\Wbb_i}(v,z)Y_{\Wbb_i}(u,\zeta)w}\qquad\in\Cbb((\zeta)),\\
\gamma_z(\zeta-z)=\bk{w',Y_{\Wbb_i}(Y(u,\zeta-z)v,z)w}\qquad\in\Cbb((\zeta-z)),
\end{gather*}
and that for each fixed $\zeta$, the series expansion with respect to $z$ at $0$ is
\begin{align*}
\beta_\zeta(z)=\bk{w',Y_{\Wbb_i}(u,\zeta)Y_{\Wbb_i}(v,z)w}\qquad\in\Cbb((z)).
\end{align*}
Indeed, let $\upphi$ be the conformal block on $(\Pbb^1;0,\infty)$ defined in example \ref{lb44}. Then $f=\wr\wr\upphi(u,v,w\otimes w')$. That its series expansions are described as above follows from theorems \ref{lb47} and \ref{lb48}.

Choose circles $C_1,C_2,C_3\subset W$ centered at $0$ with radii $r_1,r_2,r_3$ respectively satisfying $r_1<r_2<r_3$. For each $z\in C_2$, choose a circle $C(z)$ centered at $z$ whose radius is less than $r_2-r_1$ and $r_3-r_2$. Then 
\begin{align}
&\bk{w',vdz\cdot udz\cdot w}=\Res_{z=0}(\Res_{\zeta=0}\alpha_z(\zeta)d\zeta)dz=\Res_{z=0}(\Res_{\zeta=0}f(z,\zeta)d\zeta)dz\nonumber\\
=&\oint_{C_2}\Bigg(\oint_{C_1}f(z,\zeta)d\zeta\Bigg)dz.\label{eq136}
\end{align}
Similarly,
\begin{align}
\bk{w',udz\cdot vdz\cdot w}=\oint_{C_3}\Bigg(\oint_{C_2}f(z,\zeta)dz\Bigg)d\zeta=\oint_{C_2}\Bigg(\oint_{C_3}f(z,\zeta)d\zeta\Bigg)dz.\label{eq137}
\end{align}
Note that in the last equation, the two contour integrals are interchangeable since $f$ is holomorphic and in particular continuous on $\Conf_2(W-\{0\})$. Also, 
\begin{align*}
&\bk{w',(\mbf L_{udz}vdz)\cdot w_i}=\Res_{z=0}(\Res_{\zeta-z=0}\gamma_z(\zeta-z)d\zeta)dz\nonumber\\
=&\oint_{C_2}\Bigg(\oint_{C(z)}f(z,\zeta)d\zeta \Bigg)dz,
\end{align*}
which, by Cauchy integral theorem (applied to the function $\zeta\mapsto f(z,\zeta)$ for each fixed $z$), equals the difference of \eqref{eq137} and \eqref{eq136}. This finishes the proof.
\end{proof}

The following useful observation will be used in constructing connections on sheaves of conformal blocks.

\begin{lm}\label{lb55}
For any $vdz\in\Vbb\otimes_\Cbb\omega_{\mc C/\mc B}(\blt S_{\fk X})(W)$ and $w_i\in\Wbb_i\otimes_\Cbb\scr O(\mc B)$,
\begin{align*}
\big((\partial_zv)dz\big)\cdot w_i+\big((L_{-1}v)dz\big)\cdot w_i=0.
\end{align*}
\end{lm}

\begin{proof}
Recall that when $v$ is a constant section, we have the $L_{-1}$-derivative property \eqref{eq134}. Since we do not assume here that $v$ is constant, we have
\begin{align*}
\partial_z Y_{\Wbb_i}(v,z)=Y_{\Wbb_i}(\partial_zv,z)+Y_{\Wbb_i}(L_{-1}v,z).
\end{align*}
Thus, the left hand side of the equation we want to prove equals
\begin{align*}
\Res_{z=0}~\partial_z\big(Y_{\Wbb_i}(v,z)w_i\big)dz.
\end{align*}
This residue must be $0$ since the series expansion of $\partial_z(\cdots)$ with respect to $z$ does not contain $z^{-1}$.
\end{proof}

\section{The logarithmic connections}\label{lb84}

Let $\fk X=(\pi:\mc C\rightarrow\mc B)$ be a family of complex curves. If $\scr E$ is an $\scr O_{\mc B}$-module, then a \textbf{logarithmic connection} $\nabla$ on $\scr E$ associates to each open subset $U\subset \mc B$ a bilinear map
\begin{gather*}
\nabla:\Theta_{\mc B}(-\log\Delta)(U)\times\scr E(U)\rightarrow\scr E(U),\qquad (\yk,s)\mapsto \nabla_{\yk} s
\end{gather*}
satisfying conditions (a) and (b) of Definition \ref{lb49}, namely,
\begin{enumerate}[label=(\alph*)]
\item If $V$ is an open subset of $U$, then $\nabla_{\yk|_V}s|_V=(\nabla_{\yk} s)|_V$.
\item If $f\in\scr O_{\mc B}(U)$, then
\begin{gather*}
\nabla_{f\yk}s=f\nabla_{\yk}s,\\
\nabla_{\yk}(fs)=\yk(f)s+f\nabla_{\yk}s.
\end{gather*}
\end{enumerate}
Thus, $\nabla$ is a connection if $\fk X$ is a family of compact Riemann surfaces (equivalently, $\Delta=\emptyset$).

Given a logarithmic connection $\nabla$ on $\scr E$, one can define  on the dual sheaf $\scr E^*$ the dual connection (also denoted by $\nabla$) as follows. Choose any open $U\subset\mc B$, $\yk\in\Theta_{\mc B}(-\log\Delta)(U)$, and  $\upsigma\in\scr E^*(U)=\Hom_{\scr O_U}(\scr E|_U,\scr O_U)$, then $\nabla_{\yk}\upsigma$, which is an element of $\Hom_{\scr O_U}(\scr E|_U,\scr O_U)$, is defined such that for any open subset $V\subset U$ and $s\in\scr E(U)$,
\begin{align}
\bk{\nabla_{\yk}\upsigma,s}=\yk\bk{\upsigma,s}-\bk{\upsigma,\nabla_{\yk}s}.
\end{align}

We now assume that $\fk X$ is a family of $N$-pointed complex curves $\fk X=(\pi:\mc C\rightarrow\mc B;\sgm_1,\dots,\sgm_N)$ equipped with $\Vbb$-modules $\Wbb_1,\dots,\Wbb_N$. Our goal of this section is to define locally a logarithmic connection $\nabla$ on $\scr T_{\fk X}^*(\Wbb_\blt)$ near each point of $\mc B$. Since our task is local, we assume that $\mc B$ is small enough such that the following hold. 
\begin{enumerate}[label=(\roman*)]
\item $\fk X$ is either smooth or is obtained by sewing a smooth family.
\item If $\fk X$ is smooth,  then $\mc B$ is biholomorphic to a  Stein open subset of a complex coordinate space $\Cbb^m$  ($m\in\Nbb$); if $\fk X$ is obtained by sewing a smooth family $\wtd{\fk X}$, then  $\wtd{\mc B}$ is biholomorphic to a  Stein open subset of $\Cbb^m$.
\item We can equip $\sgm_1,\dots,\sgm_N$ with local coordinates:
\begin{align*}
\fk X=(\pi:\mc C\rightarrow\mc B;\sgm_1,\dots,\sgm_N;\eta_1,\dots,\eta_N).
\end{align*}
\end{enumerate}
Note that condition (ii) implies that $\mc B$ is Stein. (If $\fk X$ is obtained by sewing $\wtd{\fk X}$, then $\mc B$ is a product of the Stein manifold $\mc B$ and some open discs, which is therefore Stein.) Recall  the description of $\Theta_{\mc B}(-\log\Delta)$ near \eqref{eq133}, which shows that, due to condition (ii), $\Theta_{\mc B}(-\log\Delta)$ is a free $\scr O_{\mc B}$-module, i.e., it is generated freely by finitely many global sections $\fk y_1,\fk y_2,\dots\in \Theta_{\mc B}(-\log\Delta)(\mc B)$.

We shall construct a (global) logarithmic connection $\nabla$ over $\mc B$ whenever the above three conditions are satisfied.  For that purpose, we shall define, for each  $\fk y\in\{\fk y_1,\fk y_2,\dots \}$, a sheaf map (not a homomorphism of $\scr O_{\mc B}$-modules !)
\begin{align*}
\nabla_{\fk y}:\scr T_{\fk X}(\Wbb_\blt)\rightarrow\scr T_{\fk X}(\Wbb_\blt)
\end{align*}
satisfying that for any open $U\subset \mc B$, $w\in\scr T_{\fk X}(\Wbb_\blt)(U)$, and $f\in\scr O(U)$,
\begin{align*}
\nabla_{\fk y}(fw)=\fk y(f)w+f\nabla_{\fk y}w.
\end{align*}
(Indeed, we will do this for any $\fk y\in \Theta_{\mc B}(-\log\Delta)(\mc B)$.) Then the  \textbf{differential operators} $\nabla_{\fk y_1},\nabla_{\fk y_2},\dots$ extend to a logarithmic connection $\nabla$ on $\scr T_{\fk X}(\Wbb_\blt)$, whose dual  connection  is the one $\nabla$ on $\scr T_{\fk X}^*(\Wbb_\blt)$.

\subsection*{Defining $\nabla$ on $\scr W_{\fk X}(\Wbb_\blt)$}

We now fix $\fk y\in \Theta_{\mc B}(-\log\Delta)(\mc B)$. Recall that $\scr T_{\fk X}(\Wbb_\blt)$ is the quotient of $\scr W_{\fk X}(\Wbb_\blt)=\Wbb_\blt\otimes_\Cbb \scr O_{\mc B}$ (identified via $\mc U(\eta_\blt)$) by the $\scr O_{\mc B}$-submodule $\pi_*\big(\scr V_{\fk X}\otimes\omega_{\mc C/\mc B}(\blt S_{\fk X})\big)\cdot \scr W_{\fk X}(\Wbb_\blt)$. Our plan is to  define the action of $\nabla_{\fk y}$ on $\scr W_{\fk X}(\Wbb_\blt)$, and then to show that $\nabla_{\fk y}$ preserves that submodule.

Choose any $k\in\Nbb$. Then, due to \eqref{eq138}, we have a short exact sequence of $\scr O_{\mc C}$-modules
\begin{align*}
0\rightarrow \Theta_{\mc C/\mc B}(kS_{\fk X})\rightarrow \Theta_{\mc C}(-\log \mc C_\Delta+kS_{\fk X})\xrightarrow{d\pi}\big(\pi^*\Theta_{\mc B}(-\log \Delta)\big)(kS_{\fk X})\rightarrow 0.
\end{align*}
where $\Theta_{\mc C}(-\log \mc C_\Delta+kS_{\fk X})$ is short for $\Theta_{\mc C}(-\log \mc C_\Delta)(kS_{\fk X})$. By Theorem \ref{lb10}, there exists $k_0\in\Nbb$ such that for any $k\geq k_0$ and $b\in \mc B$, we have $H^1(\mc C_b,\Theta_{\mc C_b}(kS_{\fk X}))=0$. Thus, when $k\geq k_0$, we have $R^1\pi_*\Theta_{\mc C/\mc B}(kS_{\fk X})=0$ due to Grauert's
Theorem \ref{lb11}. Therefore, if we apply \eqref{eq4} to the above short exact sequence, we get an exact sequence of $\scr O_{\mc B}$-modules
\begin{align*}
0\rightarrow \pi_*\Theta_{\mc C/\mc B}(k S_{\fk X})\rightarrow \pi_*\Theta_{\mc C}(-\log \mc C_\Delta+k S_{\fk X})\xrightarrow{d\pi}\pi_*\big(\pi^*\Theta_{\mc B}(-\log \Delta)\big)(k S_{\fk X})
\rightarrow 0.
\end{align*}
Consider this as a short exact sequence of $\scr O_{\mc B}$-modules. Since $\pi_*\Theta_{\mc C/\mc B}(k S_{\fk X})$ is coherent by Grauert direct image theorem (indeed it is locally free by Theorem \ref{lb11}), and since $\mc B$ is assumed to be Stein, we have $H^1(\mc B,\pi_*\Theta_{\mc C/\mc B}(k S_{\fk X}))=0$ by Cartan's theorem B. Therefore, we have an exact sequence of vector spaces $0\rightarrow H^0\big(\mc B,\pi_*\Theta_{\mc C/\mc B}(k S_{\fk X})\big)\rightarrow H^0\big(\mc B,\pi_*\Theta_{\mc C}(-\log \mc C_\Delta+k S_{\fk X})\big)\xrightarrow{d\pi}H^0\big(\mc B,\pi_*\big(\pi^*\Theta_{\mc B}(-\log \Delta)\big)(k S_{\fk X})\big)
\rightarrow 0$. To simplify notations, we take the direct limit over all $k\geq k_0$ to obtain an exact sequence 
\begin{align}
0&\rightarrow H^0\big(\mc B,\pi_*\Theta_{\mc C/\mc B}(\blt S_{\fk X})\big)\rightarrow H^0\big(\mc B,\pi_*\Theta_{\mc C}(-\log \mc C_\Delta+\blt S_{\fk X})\big)\nonumber\\
&\xrightarrow{d\pi}H^0\big(\mc B,\pi_*\big(\pi^*\Theta_{\mc B}(-\log \Delta)\big)(\blt S_{\fk X})\big)
\rightarrow 0.\label{eq151}
\end{align}

Choose any  $\fk y\in\Theta_{\mc B}(-\log\Delta)(\mc B)$. Then, its pull back $\pi^*\fk y$ is in $\pi^*\Theta_{\mc B}(-\log\Delta)(\mc C)$ (recall \eqref{eq139}), which can be viewed as an element of $\pi_*\big(\pi^*\Theta_{\mc B}(-\log \Delta)\big)(\mc B)$ and hence of $\pi_*\big(\pi^*\Theta_{\mc B}(-\log \Delta)\big)(\blt S_{\fk X})(\mc B)$. Since the $d\pi$ in the above exact sequence is a surjective, there exists  a \textbf{lift} $\wtd{\fk y}$ of $\fk y$, i.e., an element $\wtd{\fk y}$ satisfying
\begin{gather}
\wtd{\fk y}\in \Theta_{\mc C}(-\log \mc C_\Delta+\blt S_{\fk X})(\mc C),\nonumber\\
d\pi(\wtd{\fk y})=\pi^*\fk y.\label{eq141}
\end{gather}
(Recall that  $\Theta_{\mc C}(-\log \mc C_\Delta+\blt S_{\fk X})(\mc C)$ equals $\pi_*\Theta_{\mc C/\mc B}(-\log \mc C_\Delta+\blt S_{\fk X})(\mc B)$.)

We are going to use $\wtd{\fk y}$ to define $\nabla_{\fk y}$. Let $\tau_\blt=(\tau_1,\tau_2,\dots)$ be coordinates of $\pi(\mc B)$. For each $1\leq i\leq N$, choose a neighborhood $W_i$ of $\sgm_i(\mc B)$ on which $\eta_i$ is defined, such that $W_i\cap W_j=\emptyset$ when $i\neq j$.  Write $\tau_\blt\circ\pi$ also as $\tau_\blt$ for simplicity, so that $(\eta_i,\tau_\blt)$ is a coordinate of $W_i$. Identify $W_i$ with its image via $(\eta_i,\tau_\blt)$, so that $\eta_i$ is identified with the standard coordinate $z$ of $\Cbb$.  Write
\begin{align}
\fk y=\sum_jg_j(\tau_\blt)\partial_{\tau_j}.
\end{align}
Then by \eqref{eq21} and \eqref{eq141}, when restricted to $W_i$, $\wtd{\fk y}$ can be written as
\begin{align}
\wtd{\fk y}|_{W_i}=h_i(z,\tau_\blt)\partial_z+\sum_jg_j(\tau_\blt)\partial_{\tau_j}\label{eq143}
\end{align}
where, due to the divisor $\blt S_{\fk X}$, $z^kh_i(z,\tau_\blt)$ is a holomorphic function on $W_i$ for some $k\in\Nbb$. Recall that $\cbf$ is the conformal vector of $\Vbb$. We set
\begin{gather}
\upnu(\wtd\yk)\in \scr V_{\fk X}\otimes\omega_{\mc C/\mc B}(\blt S_{\fk X})(W_1\cup\cdots\cup W_N)\nonumber\\
\mc U_\varrho(\eta_i)\upnu(\wtd\yk)|_{W_i}=h_i(z,\tau_\blt)\cbf ~dz.\label{eq147}
\end{gather}
Now, for any open subset $U\subset\mc B$ and any $w_i\in\Wbb_i\otimes_\Cbb\scr O_{\mc B}(U)$, we define
\begin{align}
\boxed{~~\nabla_{\fk y}w_i=\sum_jg_j(\tau_\blt)\partial_{\tau_j} w_i-\upnu(\wtd\yk)\cdot w_i~~}\label{eq142}
\end{align}
For any $w_\blt=w_1\otimes\cdots\otimes w_N\in\Wbb_\blt\otimes_\Cbb\scr O_{\mc B}(U)$, we set
\begin{align}
\nabla_\yk w_\blt=\sum_{i=1}^Nw_1\otimes w_2\otimes\cdots\otimes \nabla_\yk w_i\otimes\cdots\otimes w_N.\label{eq185}
\end{align}
This finishes the definition of $\nabla_\yk$. 

\begin{rem}
It is easy to check that the definition of $\nabla_\yk$ is independent of the choice of the coordinate $\tau_\blt$ of $\mc B$. Thus, if we assume in (ii) just that $\mc B$ is a Stein manifold, then one still has a well defined  $\nabla_\yk$ whose local expression is given by \eqref{eq142} and \eqref{eq185} when choosing a coordinate $\tau_\blt$ for a small enough open subset of $\mc B$. The requirement that $\mc B$ is also biholomorphic to an open subset of $\Cbb^m$ is used to define $\nabla$, for which the freeness of $\Theta_{\mc B}(-\log\Delta)$ is needed.
\end{rem}

\begin{rem}
As mentioned before, we define $\nabla$ by first defining $\nabla_{\yk_1},\nabla_{\yk_2},\dots$ where $\yk_1,\yk_2,\dots$ generate freely $\Theta_{\mc B}(-\log\Delta)$, and then extending it $\scr O_{\mc B}$-linearly to a logarithmic connection over $\mc B$. $\nabla_{\yk_1},\nabla_{\yk_2},\dots$ are defined by formula \eqref{eq142} using the lifts $\wtd\yk_1,\wtd\yk_2,\dots$. Then for any section $\yk$ of $\Theta_{\mc B}(-\log\Delta)$, the operator $\nabla_\yk$ can also be defined by \eqref{eq142} by choosing the lift $\wtd\yk$ in the following way: choose unique sections $f_1,f_2,\dots$ of $\scr O_{\mc B}$ such that $\yk=f_1\yk_1+f_2\yk_2+\cdots$. Then $\wtd\yk=(f_1\circ\pi)\wtd\yk_1+(f_2\circ\pi)\wtd\yk_2+\cdots$.
\end{rem}

\begin{rem}
One can write down the explicit formula of $\upnu(\wtd\yk)\cdot w_i$. Write
\begin{align}
h_i(z,\tau_\blt)=\sum_{k\in\Zbb}\wht h_i(k,\tau_\blt)z^k,\label{eq246}
\end{align}
noting that $\wht h_i(k,\tau_\blt)$ vanishes when $k$ is sufficiently small. Then, using $Y_{\Wbb_i}(\cbf)_k=L_{k-1}$, we compute
\begin{align}
\upnu(\wtd\yk)\cdot w_i=&\Res_{z=0}h_i(z,\tau_\blt)Y_{\Wbb_i}(\cbf,z)w_idz=\sum_{k\in\Zbb}\wht h_i(k,\tau_\blt)\Res_{z=0}z^kY_{\Wbb_i}(\cbf,z)w_idz\nonumber\\
=&\sum_{k\in\Zbb}\wht h_i(k,\tau_\blt) Y_{\Wbb_i}(\cbf)_kw_i=\sum_{k\in\Zbb}\wht h_i(k,\tau_\blt)L_{k-1}w_i.\label{eq247}
\end{align}
We conclude:
\begin{align}
\nabla_{\fk y}w_i=\sum_jg_j(\tau_\blt)\partial_{\tau_j} w_i-\sum_{k\in\Zbb}\wht h_i(k,\tau_\blt)L_{k-1}w_i.\label{eq144}
\end{align}
\end{rem}

\begin{rem}
We give a heuristic explanation of our definition of $\nabla_\yk$. Assume that $\fk X$ is a smooth family. For brevity, we also assume $N=i=1$, and write $\Wbb_i=\Wbb,w_i=w,\sgm_i=\sgm,\eta_i=\eta$. Let $\zeta\mapsto\varphi_\zeta^{\yk}$ and  $\zeta\mapsto\varphi_\zeta^{\wtd\yk}$  be the (complex) one-parameter flows (in $\mc C$ and in $\mc B$) integrated from the vector fields $\yk$ and $\wtd\yk$ respectively. Fix $b\in\mc B$ and set $b(\zeta)=\varphi_\zeta^{\yk}(b)$. Choose a closed disc $D(b)$ in the fiber $\mc C_b$ centered at $\sgm(b)$. Then we have an equivalence of open Riemann surfaces
\begin{align*}
\varphi_\zeta^{\wtd\yk}:\mc C_b-D(b)\xrightarrow{\simeq}\mc C_{b(\zeta)}-D(b(\zeta))
\end{align*}
where $D(b(\zeta))$ is a closed disc inside $\mc C_b$. Now, $\eta$ gives local coordinates of $\mc C_b-D(b)$ and $\mc C_{b(\zeta)}-D(b(\zeta))$ near the circles $\partial D(b)$ and $\partial(D(b(\zeta)))$ respectively. Pull back the coordinate $\eta$ near $\partial(D(b(\zeta)))$ to one near $\partial D(b)$ through the bihomolorphic map $\varphi_\zeta^{\wtd\yk}$, and call this new coordinate $\eta_\zeta$. Then, we have two local coordinates of $\mc C_b-D(b)$: they are $\eta,\eta_\zeta$, both defined near $\partial D(b)$.

Now, we shall find the condition of $w$ to be a parallel section in the direction of $\yk$. Then $w(b)$ in the $\eta$ coordinate should be equal to $w(b(\zeta))$ in the $\eta_\zeta$ coordinate. So one should expect the following identity of elements in $\scr W_{\fk X}(\Wbb)$:
\begin{align*}
\mc U(\eta)^{-1}w(b)=\mc U(\eta_\zeta)^{-1}w(b(\zeta)).
\end{align*}
Here, $\mc U(\eta)$ and $\mc U(\eta_\zeta)$ are the trivilizations of the vector space $\scr W_{\fk X}(\Wbb)$ induced by the local coordinates $\eta,\eta_\zeta$ near the circles. (We have defined such trivializations in the paragraphs near \eqref{eq121} when the local coordinates are defined near (neighborhoods of) points. It is reasonable to expect that they can be generalized, at least formally, to those defined near circles.) Thus, formally, we have
\begin{align*}
w(b(\zeta))=\mc U(\eta_\zeta\circ\eta^{-1})w(b).
\end{align*}
We now want to take the derivative of this equation. Set $h=h_i,\wht h=\wht h_i$. Then $\partial_\zeta(\eta_\zeta\circ\eta^{-1})|_{\zeta=0}$, the derivative at $\zeta=0$ of the transformation $\eta_\zeta\circ\eta^{-1}$, should equal the vector field $h(z,\tau_\blt)\partial_z$ at $\tau_\blt=b$ due to \eqref{eq143}.  Thus, the derivative of $\mc U(\eta_\zeta\circ\eta^{-1})$ at $\tau_\blt=b$ should equal the action of the vector field $h(z,b)\partial_z=\sum_{k\in\Zbb}\wht h(k,b)z^k\partial_z$ on $\Wbb$, which, by the correspondence $z^k\partial_z\leftrightarrow L_{k-1}$ (see \eqref{eq145})\footnote{Here, unlike in Section \ref{lb52}, the action of $z\partial_z$ is $L_0$ but not $\wtd L_0$.}, should be 
$\sum_{k\in\Zbb}\wht h(k,b)L_{k-1}$. So
\begin{align*}
\partial_\zeta w(b(\zeta))\Big|_{\zeta=0}=\sum_{k\in\Zbb}\wht h(k,b)L_{k-1}w(b).
\end{align*}
By \eqref{eq144}, this is equivalent to $\nabla_{\yk}w_i=0$ at $b$.
\end{rem}

\subsection*{Defining $\nabla$ on $\scr T_{\fk X}(\Wbb_\blt)$}

Assume as before that $\mc B$ satisfies conditions (i)-(iii).
\begin{thm}
The logarithmic connection $\nabla$ on $\scr W_{\fk X}(\Wbb_\blt)$ descends to one on $\scr T_{\fk X}(\Wbb_\blt)$.
\end{thm}

To prove this theorem, we need to show that $\nabla_{\fk y}$ descends to one on $\scr T_{\fk X}(\Wbb_\blt)$. This means that we need to check that $\nabla_{\fk y}$ preserves the $\scr O_{\mc B}$-submodule $\pi_*\big(\scr V_{\fk X}\otimes\omega_{\mc C/\mc B}(\blt S_{\fk X})\big)\cdot \scr W_{\fk X}(\Wbb_\blt)$.

\begin{proof}
Choose any open subset $U\subset\mc B$. Choose any $w_i\in\Wbb_i\otimes_\Cbb\scr O_{\mc B}(U)$ for each $1\leq i\leq N$ and set $w_\blt=w_1\otimes\cdots\otimes w_N$ which is in $\scr W_{\fk X}(\Wbb_\blt)(U)=\Wbb_\blt\otimes_\Cbb\scr O_{\mc B}(U)$. Choose $v\in\pi_*\big(\scr V_{\fk X}\otimes\omega_{\mc C/\mc B}(\blt S_{\fk X})\big)(U)$. So $v\in\scr V_{\fk X}\otimes\omega_{\mc C/\mc B}(\blt S_{\fk X})(\mc C_U)$. (Recall that by our notation, $\mc C_U=\pi^{-1}(U)$.) We shall show that $[\nabla_\yk,v]=\mc L_{\wtd\yk}v$,  i.e.,
\begin{align}
\nabla_\yk (v\cdot w_i)=v\cdot \nabla_\yk w_i+\mc L_{\wtd\yk}v\cdot w_i.\label{eq146}
\end{align}
Then the same relation holds with $w_i$ replaced by $w_\blt$, which will finish the proof of the theorem. (Note that $\mc L_{\wtd\yk}$ is defined as in Section \ref{lb50}, and that $\mc L_{\wtd\yk}v$ is also in $\scr V_{\fk X}\otimes\omega_{\mc C/\mc B}(\blt S_{\fk X})(\mc C_U)$.) Since both sides of \eqref{eq146} are $\Wbb_i$-valued holomorphic functions on $U$, to prove this relation on $U$, it suffices to verify it on $U-\Delta$. Therefore, we  assume in the following that $\Delta=\emptyset$ and $\mc B=U$. So $\fk X$ is a smooth family.

Recall that $\eta_i$ is defined on $W_i$, and that we have identified $W_i$ with $(\eta_i,\tau_\blt)(W_i)$ and hence $\eta_i$ with the standard coordinate $z$.   Since, in \eqref{eq146}, $v$ is acting on $\Wbb_i$, we use $\mc U_\varrho(\eta_i)=\mc V_\varrho(\eta_i)$ to identify $\scr V_{\fk X}|_{W_i}\simeq\Vbb\otimes_\Cbb\scr O_{W_i}$ and hence $\scr V_{\fk X}\otimes\omega_{\mc C/\mc B}|_{W_i}\simeq \Vbb\otimes_\Cbb \scr O_{W_i}dz$. Then the section $v|_{W_i}\in\scr V_{\fk X}\otimes\omega_{\mc C/\mc B}(\blt\SX)(W_i)$ can be written as $v=udz$ where $u=u(z,\tau_\blt)\in \Vbb\otimes_\Cbb\scr O_{W_i}(\blt\SX)(W_i)$.

We now use \eqref{eq142} to calculate that
\begin{align*}
\nabla_\yk (v\cdot w_i)=\nabla_\yk (udz\cdot w_i)=\sum_jg_j(\tau_\blt)\partial_{\tau_j} (udz\cdot w_i)-\upnu(\wtd\yk)\cdot udz\cdot w_i,
\end{align*}
which, by Proposition \ref{lb51}, equals
\begin{align*}
&\sum_jg_j(\tau_\blt)udz\cdot\partial_{\tau_j}  w_i-udz\cdot\upnu(\wtd\yk)\cdot  w_i+\sum_jg_j(\tau_\blt) ((\partial_{\tau_j}u)dz\cdot w_i)-(\mbf L_{\upnu(\wtd\yk)}udz)\cdot w_i\\
=&v\cdot \nabla_\yk w_i+\sum_jg_j(\tau_\blt) ((\partial_{\tau_j}u)dz\cdot w_i)-(\mbf L_{\upnu(\wtd\yk)}udz)\cdot w_i.
\end{align*}
Thus, one can finish proving \eqref{eq146} by showing
\begin{align*}
\mc L_{\wtd\yk}udz\cdot w_i=\sum_jg_j(\tau_\blt) ((\partial_{\tau_j}u)dz\cdot w_i)-(\mbf L_{\upnu(\wtd\yk)}udz)\cdot w_i.
\end{align*}
By Lemma \ref{lb55}, the above equation follows if we can show
\begin{align}
\mc L_{\wtd\yk}udz=\sum_jg_j(\tau_\blt) (\partial_{\tau_j}u)dz-\mbf L_{\upnu(\wtd\yk)}udz+(\partial_z+L_{-1})(\cdots)dz.\label{eq150}
\end{align}

By \eqref{eq147} and \eqref{eq149}, we have
\begin{align*}
&\mbf L_{\upnu(\wtd\yk)}udz=\sum_{k\in\Nbb}\frac 1{k!}\partial_z^kh_i(z,\tau_\blt)Y(\cbf)_ku(z,\tau_\blt)dz\\
=&\sum_{k\in\Nbb}\frac 1{k!}\partial_z^kh_i(z,\tau_\blt)L_{k-1}u(z,\tau_\blt).
\end{align*}
By Theorem \ref{lb53} and also Remark \ref{lb54} (which explains the appearance of $\partial_zh_i(z,\tau_\blt) udz$ below), we have
\begin{align*}
&\mc L_{\wtd\yk}udz=h_i(z,\tau_\blt)\partial_z udz+\sum_{j=1}^m g_j(\tau_\blt)\partial_{\tau_j}udz-\sum_{k\geq 1}\frac 1{k!}\partial_z^k h_i(z,\tau_\blt)L_{k-1}udz+\partial_zh_i(z,\tau_\blt) udz\\
=&\partial_z(h_i(z,\tau_\blt) u)dz+\sum_{j=1}^m g_j(\tau_\blt)\partial_{\tau_j}udz-\sum_{k\geq 0}\frac 1{k!}\partial_z^k h_i(z,\tau_\blt)L_{k-1}udz+h_i(z,\tau_\blt)L_{-1}udz\\
=&\sum_{j=1}^m g_j(\tau_\blt)\partial_{\tau_j}udz-\sum_{k\geq 0}\frac 1{k!}\partial_z^k h_i(z,\tau_\blt)L_{k-1}udz+(\partial_z+L_{-1})(h_i(z,\tau_\blt) u)dz.
\end{align*}
This proves \eqref{eq150}.
\end{proof}

\subsection*{Projective uniqueness of $\nabla$}

We shall show that the definition of $\nabla_\yk$ on $\scr T_{\fk X}(\Wbb_\blt)$ is independent, up to $\scr O_{\mc B}$-scalar multiples, of the choice of the lift $\wtd\yk$.

\begin{pp}\label{lb73}
Suppose that  $\fk y\in\Theta_{\mc B}(-\log\Delta)(\mc B)$ has two lifts $\wtd\yk,\wtd\yk'\in\Theta_{\mc C}(-\log\mc C_\Delta+\blt S_{\fk X})(\mc C)$ which together with $\eta_\blt$ define $\nabla_\yk,\nabla_\yk'$  respectively. Then there exists  $f\in\scr O(\mc B)$ such that 
\begin{align*}
\nabla_\yk s-\nabla_\yk's =-fs
\end{align*}
for any section $s$ of $\scr T_{\fk X}(\Wbb_\blt)$. Consequently, for  any section $\upphi$ of $\scr T_{\fk X}^*(\Wbb_\blt)$, we have
\begin{align}
\nabla_\yk \upphi-\nabla_\yk'\upphi =f\upphi.
\end{align}
\end{pp}

\begin{rem}\label{lb143}
A different proof of this proposition will be given in Section \ref{lb139}, and we will see that  $f=\#(\wtd\yk-\wtd\yk')$ which is calculated by \eqref{eq179}. 
\end{rem}

\begin{proof}
Set $\wtd 0=\wtd\yk-\wtd\yk'$ and $\upnu(\wtd 0)=\upnu(\wtd\yk)-\upnu(\wtd\yk')$. Then $d\pi(\wtd 0)=0$, i.e., $\wtd 0$ is a lift of the zero tangent field of $\mc B$. Hence, by the exact sequence \eqref{eq151}, we have
\begin{gather*}
\wtd 0\in\pi_*\Theta_{\mc C/\mc B}(\blt S_{\fk X})(\mc B).
\end{gather*}
Let $W=W_1\cup\cdots\cup W_N$. Then $\upnu(\wtd 0)$ is an element of $ \scr V_{\fk X}\otimes\omega_{\mc C/\mc B}(\blt S_{\fk X})(W)$. We shall show that the action of $\upnu(\wtd 0)$ on $\scr W_{\fk X}(\Wbb_\blt)=\Wbb_\blt\otimes_\Cbb\scr O_{\mc B}$ (defined by acting on each component as in \eqref{eq152}) descends to  an $\scr O_{\mc B}$-scalar multiplication on $\scr T_{\fk X}(\Wbb_\blt)$. 

On each $W_i$ and under the previously mentioned trivializations,  $\wtd 0$ and $\upnu(\wtd 0)$ can be written as
\begin{gather*}
\wtd 0|_{W_i}=a_i(z,\tau_\blt)\partial_z,\qquad \mc U_\varrho(\eta_i)\upnu(\wtd 0)|_{W_i}=a_i(z,\tau_\blt)\cbf~dz.
\end{gather*}
So
\begin{align*}
\upnu(\wtd 0)\in \svir_c\otimes\omega_{\mc C/\mc B}(\blt S_{\fk X})(W).
\end{align*}
By Theorem \ref{lb10}, there exists $k_0\in\Nbb$ such that for each $b\in\mc B$ and $k\geq k_0$, we have $H^1(\mc C_b,\omega_{\mc C_b}(kS_{\fk X}))=0$. Thus, as argued for \eqref{eq151}, by constructing long exact sequences twice and using the fact that $\mc B$ is Stein, we see that the short exact sequence \eqref{eq154} gives rise to an exact sequence of vector spaces
\begin{align}
0\rightarrow H^0\big(\mc B,\pi_*\omega_{\mc C/\mc B}(\blt S_{\fk X})\big) \rightarrow  H^0\big(\mc B,\pi_*\big(\svir_c\otimes \omega_{\mc C/\mc B}(\blt S_{\fk X})\big)\big)\xrightarrow{\uplambda}  H^0\big(\mc B,\pi_*\Theta_{\mc C/\mc B}(\blt S_{\fk X})\big)\rightarrow 0.\label{eq176}
\end{align}
Since $\uplambda$ is  surjective,  there is an element $u\in \pi_*\big(\svir_c\otimes \omega_{\mc C/\mc B}(\blt S_{\fk X})\big)(\mc B)=\svir_c\otimes \omega_{\mc C/\mc B}(\blt S_{\fk X})(\mc C)$ such that $\uplambda(u)=\wtd 0$.  By \eqref{eq155} (note the identification that $\eta=z$), on  $W_i$ we can express $u$ by
\begin{align*}
\mc U_\varrho(\eta_i)u|_{W_i}=a_i(z,\tau_\blt)\cbf~dz+b_i(z,\tau_\blt)\id~dz.
\end{align*}
Hence
\begin{align*}
\upnu(\wtd 0)|_{W_i}-u|_{W_i}=-b_i(z,\tau_\blt)\id~dz.
\end{align*}
Set $f(\tau_\blt)=-\sum_i \Res_{z=0}~b_i(z,\tau_\blt)dz$. Then $\upnu(\wtd 0)-u|_W$ acts as $f\cdot\id$ on $\scr W_{\fk X}(\Wbb_\blt)$ and hence on  $\scr T_{\fk X}(\Wbb_\blt)$. By the definition of $\scr T_{\fk X}(\Wbb_\blt)$ as a quotient of $\scr W_{\fk X}(\Wbb_\blt)$, $u$ acts trivially on $\scr T_{\fk X}(\Wbb_\blt)$. So $\upnu(\wtd 0)$ acts as $f\cdot\id$ on $\scr T_{\fk X}(\Wbb_\blt)$.
\end{proof}

\section{Local freeness}

We say that the VOA $\Vbb$ is \textbf{$C_2$-cofinite} if the subspace of $\Vbb$ spanned by $C_2(\Vbb):=\{Y(u)_{-2}v:u,v\in\Vbb \}$ \index{C2V@$C_2(\Vbb)$} has finite codimension. The following important result is due to Buhl.

\begin{thm}[Cf.  \cite{Buhl02} Thm. 1]\label{lb122}
Assume $\Vbb$ is $C_2$-cofinite. Then there exist $Q\in\Nbb$ and a finite set $\Ebb$ of homogeneous vectors of $\Vbb$   satisfying the following condition: For any weak $\Vbb$-module $\Wbb$ generated by a vector $w_0$, there exists $L\in\Nbb$ such that $\Wbb$ is spanned by elements of the form 
\begin{align}
Y_\Wbb(v_k)_{-n_k}Y_\Wbb(v_{k-1})_{-n_{k-1}}\cdots Y_\Wbb(v_1)_{-n_1} w_0\label{eq156}
\end{align}
where $n_k\geq n_{k-1}\geq\cdots\geq n_1>-L$ and $v_1,v_2,\dots,v_k\in\Ebb$. In addition, for any $1\leq j\leq k$, if $n_j>0$ then $n_j>n_{j-1}$; if $n_j\leq 0$ then $n_j=n_i$ for at most $Q$ different $i$.
\end{thm}

We will fix this $\Ebb$ in this section.

\begin{co}\label{lb65}
Assume that $\Vbb$ is $C_2$-cofinite. Let $\Wbb$ be a finitely generated admissible $\Vbb$-module. Then for any $n\in\Nbb$, there exists $\nu(n)\in\Nbb$ such that any $\wtd L_0$-homogeneous vector $w\in\Wbb$ whose weight $\wtd\wt(w)>\nu(n)$ is a sum of vectors in   $Y_\Wbb(v)_{-l}\Wbb(\wtd\wt(w)-\wt(v)-l+1)$ where $v\in\Ebb$ and $l>n$.
\end{co}

\begin{proof}
Assume without loss of generality that $\Wbb$ is generated by a single $\wtd L_0$-homogeneous vector $w_0$. Let $T$ be the set of all vectors of the form \eqref{eq156} where $n_k\leq n$. Then, by the above theorem, $T$ is a finite subset of $\Wbb$. Set $\nu(n)=\max\{\wtd\wt(w_1):w_1\in T\}$. If $w\in\Wbb$ is $\wtd L_0$-homogeneous with weight $\wtd\wt(w)>\nu(n)$, then we can also write $w$ as a sum of $\wtd L_0$-homogeneous vectors of the form \eqref{eq156}, but now the $n_k$ must be greater than $\nu$ since such vector is not in $T$. This proves that $w$ is a sum of $\wtd L_0$-homogeneous vectors of the form $Y_\Wbb(v)_{-l}w_2$ where $v\in\Ebb$, $l>n$, and $w_2\in\Wbb$ is $\wtd L_0$-homogeneous. The same is true if $\Wbb$ is assumed finitely generated. By \eqref{eq203}, we have $\wtd\wt(w_2)=\wtd\wt(w)-\wt(v)-l+1$.
\end{proof}

Recall the definition of sheaf of conformal blocks $\scr T_{\fk X}(\Wbb_\blt)=\scr W_{\fk X}(\Wbb_\blt)/\scr N$ where $\scr N=\pi_*\big(\scr V_{\fk X}\otimes\omega_{\mc C/\mc B}(\blt S_{\fk X})\big)\cdot \scr W_{\fk X}(\Wbb_\blt)$. For any section $\sigma$ of $\scr W_{\fk X}(\Wbb_\blt)$, its equivalence class (considered as a section of $\scr T_{\fk X}(\Wbb_\blt)$) is written as $[\sigma]$.

\begin{thm}\label{lb67}
Let $\Vbb$ be $C_2$-cofinite, let $\Wbb_1,\dots,\Wbb_N$ be finitely generated $\Vbb$-modules, and let $\fk X=(\pi:\mc C\rightarrow\mc B;\sgm_1,\dots,\sgm_N;\eta_1,\dots,\eta_N)$ be a family of $N$-pointed complex curves with local coordinates. Assume that $\mc B$ is a  Stein manifold. Then there exist finitely many elements $s_1,s_2,\dots$ of $\scr W_{\fk X}(\Wbb_\blt)(\mc B)$ such that  for any element $\sigma\in\scr W_{\fk X}(\Wbb_\blt)(\mc B)$, its equivalence class $[\sigma]$ in $\scr T_{\fk X}(\Wbb_\blt)(\mc B)$ is an $\scr O(\mc B)$-linear combination of $[s_1],[s_2],\dots$.
\end{thm}

\begin{proof}
Since local coordinates are chosen, we identify $\scr W_{\fk X}(\Wbb_\blt)$ with $\Wbb_\blt\otimes_\Cbb\scr O_{\mc B}$.  Let $E=\max\{\wt(v):v\in\Ebb\}$. By our assumption on $\fk X$ and $\mc B$, and by Theorem \ref{lb21}, there exists $k_0\in\Nbb$ such that
\begin{align}
H^1(\mc C_b,\scr V_{\mc C_b}^{\leq E}\otimes\omega_{\mc C_b}(kS_{\fk X}))=0\label{eq163}
\end{align}
for any $b\in\mc B$ and $k\geq k_0$. We fix an arbitrary $k\in\Nbb$ satisfying $k\geq E+k_0$.  

Introduce a weight $\wtd\wt$ on $\Wbb_\blt$ such that $\wtd\wt(w_\blt)=\wtd\wt(w_1)+\wtd\wt(w_2)+\cdots+\wtd\wt(w_N)$ when $w_1,\dots,w_N$ are $\wtd L_0$-homogeneous. For each $n\in\Nbb$, $\Wbb_\blt^{\leq n}$ (resp. $\Wbb_\blt(n)$) denotes the subspace spanned by all $\wtd L_0$-homogeneous homogeneous vectors $w\in\Wbb_\blt$ satisfying $\wtd\wt(w)\leq n$ (resp. $\wtd\wt(w)=n$). Note that $\Wbb_\blt^{\leq N\nu(k)}$ is finite dimensional. We shall prove  the claim of our theorem by choosing  $s_1,s_2,\dots$ to be a basis of $\Wbb_\blt^{\leq N\nu(k)}$. By induction, it suffices  to show that for any $n>N\nu(k)$,  any  vector of $\Wbb_\blt(n)$ (considered as constant sections of $\Wbb_\blt\otimes_\Cbb\scr O(\mc B)$) is a (finite) sum of elements of $\Wbb_\blt^{\leq n-1}\otimes_\Cbb\scr O(\mc B)$ mod  $\scr N(\mc B)$.

Choose any $w_\blt=w_1\otimes\cdots\otimes w_N\in \Wbb_\blt(n)$ such that $w_1,\dots,w_N$ are $\wtd L_0$-homogeneous. Then one of $w_1,\dots,w_N$ must have $\wtd L_0$-weight greater than $\nu(k)$. Assume, without loss of generality, that $\wtd\wt(w_1)>\nu(k)$. Then, by Corollary \ref{lb65}, $w_1$ is a sum of non-zero $\wtd L_0$-homogeneous vectors of the form $Y_{\Wbb_1}(u)_{-l}w_1^\circ$ where $u\in\Ebb$, $l>k$, $w_1^\circ\in\Wbb_1$, and $\wtd\wt(w_1^\circ)=\wtd\wt(w_1)-\wt(u)-l+1$. Thus $\wtd\wt(w_1)-\wtd\wt(w_1^\circ)\geq l-1\geq k\geq E+k_0.$

It suffices to show that each $Y_{\Wbb_1}(u)_{-l}w_1^\circ\otimes w_2\otimes\cdots\otimes w_N$ is a sum of elements of $\Wbb_\blt^{\leq n-1}\otimes_\Cbb\scr O(\mc B)$ mod  $\scr N(\mc B)$. Thus, we assume for simplicity that $w_1=Y_{\Wbb_1}(u)_{-l}w_1^\circ$.  Then
\begin{align*}
w_\blt=Y_{\Wbb_1}(u)_{-l}w_1^\circ\otimes w_2\otimes\cdots\otimes w_N.
\end{align*}
Set also  
\begin{align*}
w_\blt^\circ=w_1^\circ\otimes w_2\otimes\cdots\otimes w_N.
\end{align*}
Then $n-\wtd\wt(w_\blt^\circ)= \wtd\wt(w_\blt)-\wtd\wt(w_\blt^\circ)\geq E+k_0.$ Thus
\begin{align}
\wtd\wt(w_\blt^\circ)\leq n-E-k_0.\label{eq164}
\end{align}

Consider the short exact sequence of $\scr O_{\mc C}$-modules
\begin{align*}
0\rightarrow \scr V_{\mc C}^{\leq E}\otimes\omega_{\mc C/\mc B}(k_0S_{\fk X})\rightarrow \scr V_{\mc C}^{\leq E}\otimes\omega_{\mc C/\mc B}(lS_{\fk X})\rightarrow\scr G\rightarrow 0
\end{align*}
where $\scr G$ is the quotient of the previous two sheaves. By \eqref{eq163} and Grauert's Theorem \ref{lb11}, we see that $R^1\pi_*(\scr V_{\mc C}^{\leq E}\otimes\omega_{\mc C/\mc B}(k_0S_{\fk X}))=0$, and $\pi_*(\scr V_{\mc C}^{\leq E}\otimes\omega_{\mc C/\mc B}(k_0S_{\fk X}))$ is locally free. Thus, we may apply  \eqref{eq4} to obtain an exact sequence of $\scr O_{\mc B}$-modules
\begin{align*}
0\rightarrow \pi_*\big(\scr V_{\mc C}^{\leq E}\otimes\omega_{\mc C/\mc B}(k_0S_{\fk X})\big)\rightarrow \pi_*\big(\scr V_{\mc C}^{\leq E}\otimes\omega_{\mc C/\mc B}(lS_{\fk X})\big)\rightarrow \pi_*\scr G\rightarrow 0.
\end{align*}
Since $\mc B$ is assumed to be a Stein manifold, by Cartan's theorem B, we know $H^1(\mc B,\pi_*(\scr V_{\mc C}^{\leq E}\otimes\omega_{\mc C/\mc B}(k_0S_{\fk X})))=0$. Thus, there is an exact sequence
\begin{align}
&0\rightarrow H^0\big(\mc B,\pi_*\big(\scr V_{\mc C}^{\leq E}\otimes\omega_{\mc C/\mc B}(k_0S_{\fk X})\big)\big)\rightarrow H^0\big(\mc B,\pi_*\big(\scr V_{\mc C}^{\leq E}\otimes\omega_{\mc C/\mc B}(lS_{\fk X})\big)\big)\nonumber\\
&\rightarrow H^0\big(\mc B,\pi_*\scr G\big)\rightarrow 0. \label{eq170}
\end{align}

Note that $H^0\big(\mc B,\pi_*\scr G\big)$ is exactly $\scr G(\mc C)$. Choose mutually disjoint neighborhoods $W_1,\dots,W_N$ of $\sgm_1(\mc B),\dots,\sgm_N(\mc B)$ respectively. For each $1\leq i\leq N$, identify $\scr V_{\mc C}^{\leq E}\otimes\omega_{\mc C/\mc B}|_{W_i}$ with $\Vbb^{\leq E}\otimes_\Cbb\omega_{\mc C/\mc B}|_{W_i}$ via $\mc U_\varrho(\eta_i)$, and identify $\eta_i$ with the standard coordinate $z$ by identifying $W_i$ with $(\eta_i,\pi)(W_i)$. Define an element $\upupsilon\in\scr G(\mc C)$ as follows. $\upupsilon|_{W_1}$ is the equivalence class represented by $uz^{-l}dz$, and $\upupsilon|_{\mc C-\sgm_1(\mc B)}=0$. Since the second map in the above exact sequence is surjective, $\upupsilon$ lifts to an element $\wht\upupsilon$ of $H^0\big(\mc B,\pi_*\big(\scr V_{\mc C}^{\leq E}\otimes\omega_{\mc C/\mc B}(lS_{\fk X})\big)\big)$, i.e., of $\big(\scr V_{\mc C}^{\leq E}\otimes\omega_{\mc C/\mc B}(lS_{\fk X})\big)(\mc C)$. Moreover, by the definition of $\scr G$ as a quotient, for each $1\leq i\leq N$ we have an element $v_i$ of $\Vbb^{\leq E}\otimes_\Cbb\scr O_{\mc C}(k_0 S_{\fk X})(W_i)$ (and hence of $\Vbb^{\leq E}\otimes_\Cbb\scr O_{W_i}(k_0 \sgm_i(\mc B))(W_i)$) such that
\begin{gather*}
\wht\upupsilon|_{W_1}=uz^{-l}dz+v_1dz,\\
\wht\upupsilon|_{W_i}=v_idz\qquad(2\leq i\leq N).
\end{gather*}
Notice that $\Res_{z=0}Y(\cdot,z)z^ndz=Y(\cdot)_n$. It follows that the element $\wht\upupsilon\cdot w_\blt^\circ$, which is clearly in $\scr N(\mc B)$,  equals $w_\blt+w_\triangle$ where
\begin{align*}
w_\triangle=(v_1dz)\cdot w_1^\circ\otimes w_2\otimes\cdots \otimes  w_i\otimes\cdots\otimes w_N+\sum_{i=2}^Nw_1^\circ\otimes w_2\otimes\cdots \otimes (v_idz)\cdot w_i\otimes\cdots\otimes w_N.
\end{align*}
Thus $[w_\blt]=-[w_\triangle]$. For each $1\leq i\leq N$, $v_i$ has pole at $z=0$ with order at most $k_0$. Thus, by \eqref{eq203}, the action of $v_idz$ on $\Wbb_i$ increases the $\wtd L_0$-weight by at most $E+k_0-1$. It follows from \eqref{eq164} that $w_\triangle\in\Wbb_\blt^{\leq {n-1}}\otimes_\Cbb\scr O(\mc B)$. The proof is complete.
\end{proof}

\begin{thm}\label{lb68}
Let $\Vbb$ be $C_2$-cofinite, let $\Wbb_1,\dots,\Wbb_N$ be finitely generated $\Vbb$-modules, and let $\fk X=(\pi:\mc C\rightarrow\mc B;\sgm_1,\dots,\sgm_N)$ be a family of $N$-pointed compact Riemann surfaces. Then $\scr T_{\fk X}(\Wbb_\blt)$ and (hence) $\scr T_{\fk X}^*(\Wbb_\blt)$ are locally free. 
\end{thm}

\begin{proof}
We have seen that $\scr T_{\fk X}(\Wbb_\blt)$ admits a connection near any point of $\mc B$. Thus, to prove the local freeness, it suffices to verify that $\scr T_{\fk X}(\Wbb_\blt)$ satisfies the two conditions of Theorem \ref{lb66}. Assume that $\mc B$ is a Stein manifold, and $\fk X$ is equipped with local coordinates.  Consider  $A=\pi_*(\scr V_{\fk X}\otimes\omega_{\mc C/\mc B}(\blt S_{\fk X}))(\mc B)$ as a vector space. Then $\scr W_{\fk X}(\Wbb_\blt)=\Wbb_\blt\otimes_\Cbb\scr O_{\mc B}$ and $\scr F=A\otimes_\Cbb\Wbb_\blt\otimes_\Cbb\scr O_{\mc B}$ are clearly locally free with infinite rank. Define a homomorphism of $\scr O_{\mc B}$-modules $\varphi:\scr F\rightarrow\scr W_{\fk X}(\Wbb_\blt)$ such that for each $v\in A$, $w\in\Wbb_\blt$, $U\subset \mc B$ being open, and $f\in\scr O_{\mc B}(U)$,
\begin{align*}
\varphi: (v,w,f)\mapsto f\cdot(v\cdot w)|_U.
\end{align*}
Then both $\varphi(\scr F)$ and $\scr N$ are $\scr O_{\mc B}$-submodules of $\scr W_{\fk X}(\Wbb_\blt)$. It is clear that $\varphi(\scr F(U))\subset\scr N(U)$ and hence $\varphi(\scr F)(U)\subset\scr N(U)$. By Corollary \ref{lb33}, we have $\varphi(\scr F)_b=\scr N_b$ for each $b\in\mc B$. Thus $\varphi(\scr F)=\scr N(U)$. Thus $\scr T_{\fk X}(\Wbb_\blt)=\coker(\varphi)$, which verifies condition (a) of Theorem \ref{lb66}. By Theorem \ref{lb67}, $\scr T_{\fk X}(\Wbb_\blt)$ satisfies condition (b). 
\end{proof}

\begin{co}\label{lb119}
In the setting of Theorem \ref{lb68}, the dimension function
\begin{align*}
b\in\mc B\mapsto \dim_\Cbb\scr T_{\fk X_b}(\Wbb_\blt)
\end{align*}
is finite and  constant on each connected component of $\mc B$.
\end{co}

\begin{proof}
This follows from theorems \ref{lb68} and \ref{lb69}.
\end{proof}

\begin{rem}
For each $g\in\Nbb$ and $N\in\Zbb_+$, one has a universal family $\fk X$ of $N$-pointed connected   compact Riemann surfaces with genus $g$, where $\mc B$ is the Teichm\"uller space $\mc T_{g,N}$ of all such curves, and any such curve is biholomorphic to one of the fibers of $\fk X$. Since $\mc T_{g,N}$ is well known to be connected, it follows that the dimensions of the spaces of conformal blocks are finite and depend only on $\Vbb,\Wbb_\blt,g,N$ but not on the complex structures of the Riemann surfaces.
\end{rem}

\chapter{Sewing and factorization}

\section{Projective structures}

\subsection*{Schwarzian derivatives and projective structures}
Let $\fk X=(\pi:\mc C\rightarrow\mc B)$ be a family of compact Riemann surfaces. Let us calculate the transition functions of $\svir_c$ (which is a subsheaf of $\scr V_{\fk X}$). Choose an open subset $U\subset \mc C$ and holomorphic functions $\eta,\mu:U\rightarrow\Cbb$ univalent on each fiber. If $f\in\scr O(U)$ and $\partial_\eta f$ is nowhere zero, we define the \textbf{Schwarzian derivative} of $f$ over $\eta$ to be \index{S@$\Sbf_\eta f,\Sbf_\eta\fk P$}
\begin{align}
\Sbf_\eta f=\frac{\partial_\eta^3f}{\partial_\eta f}-\frac 32 \Big(\frac{\partial_\eta^2f}{\partial_\eta f} \Big)^2
\end{align}
where the partial derivative $\partial_\eta$ is defined with respect to $(\eta,\pi)$, i.e., it is annihilated by $d\pi$ and restricts to $d/d\eta$ on each fiber. Similarly, one can define $\Sbf_\mu f$. The change of variable formula is easy to calculate:
\begin{align}\label{eq166}
\Sbf_\mu f =(\partial_\mu\eta)^2 \Sbf_\eta f+\Sbf_\mu\eta.
\end{align}
Take $f=\mu$ and notice $\Sbf_\mu\mu=0$, we have
\begin{align}
\Sbf_\mu\eta=-(\partial_\mu\eta)^2\Sbf_\eta\mu.\label{eq167}
\end{align}
Assuming $f$ is also univalent on each fiber, we obtain the cocycle relation.
\begin{align}
\Sbf_\mu\eta\cdot d\mu^2=-\Sbf_\eta\mu\cdot d\eta^2,\qquad \Sbf_\mu f\cdot d\mu^2 + \Sbf_f \eta\cdot df^2+\Sbf_\eta\mu\cdot d\eta^2=0.\label{eq260}
\end{align}

By \eqref{eq76}, the transition function from $\mu$ to $\eta$ is given by $\mc U_\varrho(\eta)\mc U_\varrho(\mu)^{-1}=\mc U(\varrho(\eta|\mu))$. The vacuum vector $\id$ is clearly fixed by the transition function. So we just need to calculate $\mc U(\varrho(\eta|\mu))\cbf$. By Remark \ref{lb58} and formula \eqref{eq47}, if $\rho=\rho(z)\in\Gbb$, then $\mc U(\rho)\cbf=\rho'(0)^{L_0}e^{c_2L_2}\cbf=\rho'(0)^{L_0}(\cbf+\frac c2c_2\id)=\rho'(0)^2\cbf+\frac c2c_2\id$ where $c$ is the central charge, and $c_2$, which is given by \eqref{eq165}, is $\frac 16\Sbf_z\rho(0)$. Replace $\rho$ by $\varrho(\eta|\mu):U\rightarrow\Gbb$. Then $\rho^{(n)}(0)$ should be replaced by $\partial_\mu^n\eta$. Thus the transition function $\mc U(\varrho(\eta|\mu))$  is described by
\begin{gather}
\mc U(\varrho(\eta|\mu))\id=\id,\qquad  \mc U(\varrho(\eta|\mu))\cbf=(\partial_\mu\eta)^2\cbf+\frac {c}{12}\Sbf_\mu\eta \cdot \id.\label{eq169}
\end{gather}

We collect some useful properties of the Schwarzian derivatives.  

\begin{pp}\label{lb70}
The following are true.
\begin{enumerate}[label=(\arabic*)]
\item If the restriction of $\eta$ to each fiber $U_b=U\cap\pi^{-1}(b)$ (where $b\in\mc B$) is a M\"obius transformation of $\mu$, i.e., of the form $\frac{a\mu+b}{c\mu+d}$ where $ad-bc\neq 0$, then $\Sbf_\mu \eta=0$.
\item  Let $Q\in\scr O(U)$. Then, for each $x\in U$, one can find a neighborhood $V\subset U$ of $x$ and a function $f\in\scr O(V)$ univalent on each fiber $V_b=V\cap\pi^{-1}(b)$, such that $\Sbf_\eta f=Q$.
\item If $f,g\in\scr O(U)$ are univalent on each fiber, then $\Sbf_\eta f=\Sbf_\eta g$ if and only if $\Sbf_f g=0$.

\item If $f,Q\in\scr O(U)$, $f$ is univalent on each fiber, and $\Sbf_\mu f=Q$, then 
\begin{align}
\Sbf_\eta f=(\partial_\mu\eta)^{-2}(Q-\Sbf_\mu \eta).\label{eq172}
\end{align}
\end{enumerate}
\end{pp}

We remark that the converse of (1) is also true: If $f$ is univalent on each fiber, and if $S_\eta f=0$, then the restriction of $f$ to each fiber is a M\"obius transformation of $\eta$.

\begin{proof}
(1) can be verified directly. To prove (2), we identify $U$ with an open subset of $\Cbb\times\mc B$ via $(\eta,\pi)$. So $\eta$ is identified with the standard coordinate $z$. We choose a neighborhood $V\subset U$ of $x$ of the form $\mc D\times W$ where $W\subset\mc B$ is open, and $\mc D$ is an open disc centered at point $p=\eta(x)\in\Cbb$. Consider the differential equation
\begin{align*}
\partial_z^2h+Qh/2=0
\end{align*}
which can be transformed to an $\Cbb^2$-valued $1$-st order differential equation
\begin{align*}
\partial_z
\begin{pmatrix}
\alpha\\
\beta
\end{pmatrix}
=
\begin{pmatrix}
0&1\\
-Q/2&0
\end{pmatrix}
\begin{pmatrix}
\alpha\\
\beta
\end{pmatrix}
\end{align*}
where $\alpha,\beta$ and $h$ are related by $\alpha=h,\beta=\partial_zh$. By Theorem \ref{lb12}, there exist solutions $h_1,h_2\in\scr O(V)$ satisfying  the initial conditions $h_1(\cdot,p)=1,\partial_z h_1(\cdot,p)=0$ and $h_2(\cdot,p)=0,\partial_z h_2(\cdot,p)=1$. It is easy to check that $f:=h_2/h_1$ satisfies $\Sbf_z f=Q$, and is defined and satisfies $\partial_zf\neq 0$ near $\{p\}\times W$.

(3) follows from \eqref{eq166}, which says $\Sbf_\eta g=(\partial_\eta f)^2\Sbf_fg+\Sbf_\eta f$. (4) follows from \eqref{eq166} and \eqref{eq167}.
\end{proof}

\begin{df}
An open cover $(U_\alpha,\eta_\alpha)_{\alpha\in\fk A}$ of $\mc C$, where each open set $U_\alpha$ is equipped with a function $\eta_\alpha\in\scr O(U_\alpha)$ holomorphic on each fiber, is called a (family of) \textbf{projective chart} of $\fk X$, if for any $\alpha,\beta\in\fk A$, we have $\Sbf_{\eta_\beta}\eta_\alpha=0$ on $U_\alpha\cap U_\beta$. Two projective charts are called equivalent if their union is a projective chart. An equivalence class of projective charts is called a \textbf{projective structure}. Equivalently, a projective structure is a maximal projective chart.

Suppose that $\fk X=(\pi:\mc C\rightarrow\mc B;\sgm_1,\dots,\sgm_N;\eta_1,\dots,\eta_N)$ is a family of $N$-pointed compact Riemann surfaces with local coordinates. We say that the local coordinates $\eta_1,\dots,\eta_N$ \textbf{admit a projective structure} if, by choosing mutually disjoint neighborhoods $W_1,\dots,W_N$ of $\sgm_1(\mc B),\dots,\sgm_N(\mc B)$ on which $\eta_1,\dots,\eta_N$ are defined respectively, there is a projective chart of $\fk X$ which contains $(W_1,\eta_1),\dots,(W_N,\eta_N)$.
\end{df}

\begin{rem}\label{lb75}
	Let $\fk P$ be a projective chart on $\fk X$. Choose an open subset $W\subset\mc C$ and a fiberwisely univalent $\eta\in\scr O(W)$. One can \index{S@$\Sbf_\eta f,\Sbf_\eta\fk P$} define an element 
	\begin{align*}
	\Sbf_\eta\fk P\in \scr O(W)
	\end{align*}
	as follows. Choose any $(U,\mu)\in\fk P$. Then $\Sbf_\eta\fk P=\Sbf_\eta\mu$ on $W\cap U$. To check that $\Sbf_\eta\fk P$ is well defined, suppose there is another $(V,\zeta)\in\fk P$. Then $\Sbf_\mu\zeta=0$ on $U\cap V$. Thus $\Sbf_\eta\mu=\Sbf_\eta\zeta$ on $U\cap V\cap W$ by Proposition \ref{lb70}.
\end{rem}

\subsection*{Existence and classification of projective structures}

In what follows we assume for simplicity that each fiber $\mc C_b$ of the family $\fk X$ is connected. We also assume that $c$ is a non-zero central charge. Recall the exact sequence \eqref{eq168}, and tensor it by the identity map of  $\omega_{\mc C/\mc B}^{\otimes 2}$. We obtain a short  exact sequence
\begin{align}
0\rightarrow \omega_{\mc C/\mc B}^{\otimes 2}\rightarrow\svir_c\otimes \omega_{\mc C/\mc B}^{\otimes 2}\xrightarrow{\uplambda}\scr O_{\mc C}\rightarrow 0,\label{eq174}
\end{align}
which induces a long one
\begin{align}
0\rightarrow\pi_*(\omega_{\mc C/\mc B}^{\otimes 2})\rightarrow \pi_*(\svir_c\otimes \omega_{\mc C/\mc B}^{\otimes 2})\xrightarrow{\uplambda} \scr O_{\mc B}\xrightarrow{\delta} R^1\pi_*(\omega_{\mc C/\mc B}^{\otimes 2}).\label{eq236}
\end{align}
Here, we have used the obvious equivalence $\scr O_{\mc B}\xrightarrow{\simeq}\pi_*\scr O_{\mc C},f\mapsto f\circ\pi$. We thus obtain a linear map
\begin{align}
\uplambda:H^0\big(\mc B,\pi_*(\svir_c\otimes \omega_{\mc C/\mc B}^{\otimes 2})\big)\rightarrow  H^0(\mc B,\scr O_{\mc B}).
\end{align}
Consider $1\in H^0(\mc B,\scr O_{\mc B})$, i.e. the constant function on $\mc B$ with value $1$.

\begin{thm}\label{lb72}
There is a one to one correspondence between the subset $\uplambda^{-1}(1)$ of $H^0\big(\mc B,\pi_*(\svir_c\otimes \omega_{\mc C/\mc B}^{\otimes 2})\big)$ and  the projective structures of $\fk X$.
\end{thm}

\begin{proof}
First of all, observe that by \eqref{eq155}, for any open $U\subset \mc C$ and $\mu\in\scr O(U)$ univalent on each fiber, we have
\begin{gather*}
\uplambda:\quad \mc U_\varrho(\mu)^{-1}\cbf ~d\mu^2\mapsto 1,\qquad \id ~d\mu^2\mapsto 0
\end{gather*}
where $\id$ is the vacuum section. If $\eta\in\scr O(U)$ is also univalent on each fiber, then by \eqref{eq169} and $d\mu=(\partial_\mu\eta)^{-1} d\eta$, we have
\begin{align}
\mc U_\varrho(\mu)^{-1}\cbf~d\mu^2=\mc U_\varrho(\eta)^{-1}\cbf~d\eta^2+\frac {c}{12}(\partial_\mu\eta)^{-2}\Sbf_\mu\eta \cdot \id~d\eta^2.\label{eq171}
\end{align}

Choose any $\nu\in\uplambda^{-1}(1)$. The corresponding projective chart $\fk P$ is defined as follows. For any $x\in\mc C$, we choose a neighborhood $U$ of $x$ and  $\mu\in\scr O(U)$ univalent on each fiber $U\cap\pi^{-1}(b)$. Since $\nu$ is an element of $\svir_c\otimes \omega_{\mc C/\mc B}^{\otimes 2}(\mc C)$ sent by $\uplambda$ to $1$, we have
\begin{gather}
\nu|_U=\mc U_\varrho(\mu)^{-1}\cbf ~d\mu^2-\frac{c}{12}Q\cdot \id~d\mu^2 \label{eq173}
\end{gather} 
for some $Q\in\scr O(U)$. By Proposition \ref{lb70}, we may find an open subset $U_x\subset U$ containing $x$ so that there exists $f_x\in\scr O(U)$ which is univalent on each fiber and satisfies $\Sbf_\mu f_x=Q$. We claim that $\fk P:=\{(U_x,f_x)\}_{x\in\mc C}$ is a projective chart. Let $y\in \mc C$, and let $(U_y,f_y)$ be obtained in the same way through a function $\eta\in\scr O(U_y)$ univalent on each fiber. Let $V=U_x\cap U_y$. Then, by \eqref{eq171}, we have
\begin{gather*}
\nu|_V=\mc U_\varrho(\eta)^{-1}\cbf ~d\eta^2-\frac{c}{12}(\partial_\mu\eta)^{-2}\big(Q-\Sbf_\mu\eta\big)\cdot\id~d\eta^2.
\end{gather*}
Thus,  we have $\Sbf_\eta f_y|_V=(\partial_\mu\eta)^{-2}(Q-\Sbf_\mu\eta)$. On the other hand, by \eqref{eq172}, we also have $\Sbf_\eta f_x|_V=(\partial_\mu\eta)^{-2}(Q-\Sbf_\mu\eta)$. So $\Sbf_\eta f_x|_V=\Sbf_\eta f_y|_V$. This implies, by Proposition \ref{lb70}, that $\Sbf_{f_x}f_y=0$ on $V$. Thus $\fk P$ is projective. By maximizing $\fk P$, we obtain the projective structure.

If $\nu'\neq \nu$ defines another projective chart $\fk P'=\{U'_x,f'_x\}$, then for the corresponding $Q'$ defined similarly by \eqref{eq173} on $U'_x$, we have $\Sbf_\mu f_x'=Q'$. Since $\nu'\neq \nu$, we may find $x$ such that $Q\neq Q'$. Thus, on $U_x\cap U_x'$ we have $\Sbf_\mu f_x\neq \Sbf_\mu f_x'$. Hence $\Sbf_{f_x} f_x'\neq 0$ by Proposition \ref{lb70}. So $\{U'_x,f'_x\}$ is not equivalent to $\{U_x,f_x\}$. Thus, the map $\nu\mapsto\fk P$ is injective. 

Finally, we show that this map is also surjective. Choose projective structure  $\fk P$. If $(U,\mu)$ and $(U,\eta)$ belongs to this projective structure, then the transition function \eqref{eq171} becomes $\mc U_\varrho(\mu)^{-1}\cbf~d\mu^2=\mc U_\varrho(\eta)^{-1}\cbf~d\eta^2$. By this formula, it is clear that one can find $\nu$ such that $\uplambda(\nu)=1$, and that on each $(U,\mu)$ belonging to the projective structure, we have $\nu|_U=\mc U_\varrho(\mu)^{-1}\cbf ~d\mu^2$. Namely, the $Q$ for $\nu|_U$ is zero. It is obvious that the $\nu$ corresponds to $\fk P$.
\end{proof}

\begin{pp}
The sheaf map $\delta:\scr O_{\mc B}\rightarrow R^1\pi_*(\omega_{\mc C/\mc B}^{\otimes 2})$ in \eqref{eq236} is zero. Consequently, we have a long exact sequence
\begin{align}
0\rightarrow H^0\big(\mc B,\pi_*(\omega_{\mc C/\mc B}^{\otimes 2})\big)\rightarrow H^0\big(\mc B,\pi_*(\svir_c\otimes \omega_{\mc C/\mc B}^{\otimes 2})\big)\xrightarrow{\uplambda} H^0\big(\mc B,\scr O_{\mc B}\big)\rightarrow H^1\big(\mc B,\pi_*(\omega_{\mc C/\mc B}^{\otimes 2})\big).\label{eq175}
\end{align}
\end{pp}

As we will see in the proof, $\pi_*(\omega_{\mc C/\mc B}^{\otimes 2})$ is a locally free $\scr O_{\mc B}$-module.

\begin{proof}
Recall that by Ehresmann's theorem, the (assumed connected) fibers of $\fk X$ are diffeomorphic. We let $g$ be the genus. If $g>1$, then for each $b\in\mc B$, we have $H^1(\mc C_b,\omega_{\mc C_b}^{\otimes 2})=0$ by Corollary \ref{lb5}. Thus, by Grauert's Theorem \ref{lb2}, $\pi_*(\omega_{\mc C/\mc B}^{\otimes 2})$ is locally free, and  $R^1\pi_*(\omega_{\mc C/\mc B}^{\otimes 2})=0$, which shows $\delta=0$.

We now treat the case $g=0,1$. Assume first of all that $\mc B$ is a single point. Then $\fk X$ admits a projective structure.  Indeed, when $g=0$, $\fk X$ is equivalent to $\Pbb^1$, which obviously admits a projective structure (e.g. $\{(\Pbb^1-\{\infty\},z),(\Pbb^1-\{0\},1/z)\}$). If $g=1$, then it is well known that $\fk X$ is equivalent to $\Cbb/\Lambda$ where $\Lambda$ is a real rank $2$ lattice in $\Cbb$ generated by $1$ and $\tau$ in the upper half plane. (See e.g. \cite{Hain08}.) A projective structure on $\Cbb$ preserved by $\Lambda$ clearly exists, which descends to one of $\fk X$. Thus, by Theorem \ref{lb72}, $\uplambda^{-1}(1)$ is nonempty. Therefore, by the exactness of \eqref{eq236}, $\scr O_{\mc B}\simeq \Cbb$ is in the kernel of $\delta$, which shows $\delta=0$.

We now treat the general case where $\mc B$ a complex manifold. Since $\scr O_{\mc B}$ (as an $\scr O_{\mc B}$-module) is generated by the global section $1$, it suffices to prove $\delta(1)=0$.  We first claim that $R^1\pi_*(\omega_{\mc C/\mc B}^{\otimes 2})$ is locally free. Indeed,  when $g=1$, we have $\Theta_{\mc C_b}\simeq\scr O_{\mc C_b}$ by the lattice realization. Thus, by Serre duality, we have $H^1(\mc C_b,\omega_{\mc C_b}^{\otimes 2})\simeq H^0(\mc C_b,\Theta_{\mc C_b})\simeq H^0(\mc C_b,\scr O_{\mc C_b})$, which has constant (over $b\in\mc B$) dimension $1$. Thus $R^1\pi_*(\omega_{\mc C/\mc B}^{\otimes 2})$ (and also $\pi_*(\omega_{\mc C/\mc B}^{\otimes 2})$) is locally free (of rank $1$) by Grauert's Theorem \ref{lb2}. When $g=0$, the same is true since all fibers are equivalent to $\Pbb^1$. Moreover, Grauert's theorem tells us that the fiber of $R^1\pi_*(\omega_{\mc C/\mc B}^{\otimes 2})$ at $b$ is naturally equivalent to $H^1(\mc C_b,\omega_{\mc C_b}^{\otimes 2})$. Thus, it suffices to show that for any $b\in\mc B$, the restriction of $\delta(1)$ to the fiber $\mc C_b$ is the zero element of $H^1(\mc C_b,\omega_{\mc C_b}^{\otimes 2})$. This follows from the previous paragraph.
\end{proof}

\begin{thm}\label{lb71}
Let $\fk X=(\pi:\mc C\rightarrow\mc B)$ be a family of compact Riemann surfaces. Suppose that $\mc B$ is a Stein manifold. Then there is a projective chart on $\fk X$.
\end{thm}

\begin{proof}
Since $\mc B$ is Stein and $\pi_*(\omega_{\mc C/\mc B}^{\otimes 2})$ is locally free, we have $H^1(\mc B,\pi_*(\omega_{\mc C/\mc B}^{\otimes 2}))=0$ by Cartan's theorem B. Therefore, the map $\uplambda$ in \eqref{eq175} is surjective. Thus, a projective structure exists by Theorem \ref{lb72}.
\end{proof}

\begin{co}
Let $\fk X=(\pi:\mc C\rightarrow\mc B;\sgm_1,\dots,\sgm_N)$ be a family of $N$-pointed compact Riemann surfaces. Then for any $b\in\mc B$, there is a neighborhood $V\ni b$ such that the restricted family $\fk X_V$ can be equipped with local coordinates $\eta_1,\dots,\eta_N$ which admit a projective structure.
\end{co}

\begin{proof}
Assume without loss of generality that $\mc B$ is Stein. Thus we can choose a projective chart $\{(U_\alpha,\mu_\alpha)\}_{\alpha\in\fk A}$ of $\fk X$. Choose $b\in \mc B$. Then one may shrink $\mc B$ such that $b$ is still in $\mc B$, and that for any $1\leq i\leq N$, we can find $\alpha\in\fk A$ such that $U_\alpha$ contains $\sgm_i(\mc B)$. Set $\eta_i=\mu_\alpha-\mu_\alpha\circ\sgm_i\circ\pi$, defined near $\sgm_i(\mc B)$. It is clear that $\eta_1,\dots,\eta_N$ are compatible with the chosen projective chart.
\end{proof}

\section{Actions of $\pi_*\Theta_{\mc C/\mc B}(\blt S_{\fk X})$}\label{lb139}

In this section, we fix  $\fk X=(\pi:\mc C\rightarrow\mc B;\sgm_1,\dots,\sgm_N;\eta_1,\dots,\eta_N)$ to be a family of $N$-pointed compact Riemann surfaces with local coordinates. We assume for simplicity that $\mc B$ is a Stein manifold with coordinates $\tau_\blt=(\tau_1,\dots,\tau_N)$. Let $\Vbb$ be a VOA with central charge $c$, and let $\Wbb_1,\dots,\Wbb_N$ be $\Vbb$-modules. By definition, $\pi_*(\scr V_{\fk X}\otimes\omega_{\mc C/\mc B}(\blt S_{\fk X}))$ and hence its subsheaf $\pi_*(\svir_c\otimes\omega_{\mc C/\mc B}(\blt S_{\fk X}))$ act trivially on the sheaf of covacua $\scr T_{\fk X}(\Wbb_\blt)$.

Recall the exact sequence \eqref{eq176}
\begin{align*}
0\rightarrow H^0\big(\mc B,\pi_*\omega_{\mc C/\mc B}(\blt S_{\fk X})\big) \rightarrow  H^0\big(\mc B,\pi_*\big(\svir_c\otimes \omega_{\mc C/\mc B}(\blt S_{\fk X})\big)\big)\xrightarrow{\uplambda}  H^0\big(\mc B,\pi_*\Theta_{\mc C/\mc B}(\blt S_{\fk X})\big)\rightarrow 0.
\end{align*}
We shall define an action of $H^0\big(\mc B,\pi_*\Theta_{\mc C/\mc B}(\blt S_{\fk X})\big)$ on $\scr T_{\fk X}(\Wbb_\blt)$, which turns out to be an $\scr O(\mc B)$-scalar multiplication. Such definition already appears in the proof of Theorem \ref{lb73}, where the action of $\wtd 0\in H^0\big(\mc B,\pi_*\Theta_{\mc C/\mc B}(\blt S_{\fk X})\big)$ is $f$. Our goal in this section is to express $f$ in terms of a projective structure $\fk P$.

Choose mutually disjoint neighborhoods $W_1,\dots,W_N$ of $\sgm_1(\mc B),\dots,\sgm_N(\mc B)$ on which $\eta_1,\dots,\eta_N$ are defined respectively. Write each $\tau_j\circ\pi$ as $\tau_j$ for short, so that $(\eta_i,\tau_\blt)$ is a coordinate of $W_i$. Set $W=W_1\cup\cdots\cup W_N$. Choose any $\theta\in H^0\big(\mc B,\pi_*\Theta_{\mc C/\mc B}(\blt S_{\fk X})\big)$, which, in each $W_i$, is expressed as
\begin{align}
\theta|_{W_i}=a_i(\eta_i,\tau_\blt)\partial_{\eta_i}.\label{eq177}
\end{align}
As in \eqref{eq147}, we define $\upnu(\theta)\in \svir_c\otimes \omega_{\mc C/\mc B}(\blt S_{\fk X})(W)$ such that
\begin{align}
\mc U_\varrho(\eta_i)\upnu(\theta)|_{W_i}=a_i(\eta_i,\tau_\blt)\cbf~d{\eta_i}.\label{eq178}
\end{align}
The action of $\theta$ on $\scr T_{\fk X}(\Wbb_\blt)$ is defined to be the action of $\upnu(\theta)$, namely, is determined by
\begin{align}
\upnu(\theta)\cdot w_\blt=\sum_{i=1}^N w_1\otimes\cdots\otimes \upnu(\theta)\cdot w_i\otimes\cdots \otimes w_N
\end{align}
for any $w_\blt=w_1\otimes\cdots\otimes w_N\in\Wbb_\blt$. (Recall \eqref{eq116}.) Such definition depends on the choice of local coordinates $\eta_1,\dots,\eta_N$.

\begin{lm}\label{lb74}
Assume $\eta_1,\dots,\eta_N$ admit a projective structure. Then the action of $\theta$ on $\scr T_{\fk X}(\Wbb_\blt)$ is zero.
\end{lm}

\begin{proof}
Assume that $(W_1,\eta_1),\dots,(W_N,\eta_N)$ belong to a projective structure $\fk P$. Then, by \eqref{eq169}, the transition function for $\cbf \otimes \omega_{\mc C/\mc B}$  between two projective coordinates is the same as that for $\Theta_{\mc C/\mc B}$, namely, when $\Sbf_\mu\eta=0$, $\partial_\mu$ changes to $\partial_\mu\eta\cdot \partial_\eta$, and $\cbf d\mu$ changes to $\partial_\mu\eta\cdot \cbf d\eta$, sharing the same transition function $\partial_\mu\eta$. Thus, as $\theta$ is over $\mc C$, $\upnu(\theta)$ can be extended to a section of $\svir_c\otimes \omega_{\mc C/\mc B}(\blt S_{\fk X})$ on $\mc C$. In particular, $\upnu(\theta)$ is in $\big(\scr V_{\fk X}\otimes \omega_{\mc C/\mc B}(\blt \SX)\big)(\mc C)=\pi_*(\scr V_{\fk X}\otimes \omega_{\mc C/\mc B}(\blt \SX))(\mc B)$. Thus, $\upnu(\theta)$ acts trivially on $\scr T_{\fk X}(\Wbb_\blt)$.
\end{proof}

\begin{co}
In Proposition \ref{lb73},  if $\eta_1,\dots,\eta_N$ admit a projective structure outside  $\mc C_{\Delta}$ (the union of all the nodal fibers), then the definition of the logarithmic connection $\nabla$ is independent of the lifts, i.e., we have $\nabla_\yk=\nabla'_\yk$.
\end{co}

Thus, the projectiveness of $\nabla$ is controlled by the projective structures of $\fk X$.

\begin{proof}
In the proof of Proposition \ref{lb73}, we need to show that $f=0$. It suffices to prove this outside the discriminant locus $\Delta$. Thus, we may assume $\Delta=\emptyset$ and hence $\fk X$ is a smooth family. We know that the action of $\upnu(\wtd 0)$ on $\scr T_{\fk X}(\Wbb_\blt)$ is the multiplication by $f$. By Lemma \ref{lb74}, the action is trivial. Thus $f=0$.
\end{proof}

\begin{pp}\label{lb80}
Suppose that $\fk X$ has a projective structure $\fk P$. Choose $\theta\in H^0\big(\mc B,\pi_*\Theta_{\mc C/\mc B}(\blt S_{\fk X})\big)$ whose local expression is given by \eqref{eq177}. Then the action of $\upnu(\theta)$ on $\scr T_{\fk X}(\Wbb_\blt)$ (defined by the local coordinates $\eta_\blt$) is the $\scr O(\mc B)$-scalar multiplication  by
\begin{align}
\#(\theta):=\frac{c}{12}\sum_{i=1}^N \Res_{\eta_i=0}~ \Sbf_{\eta_i}\fk P\cdot a_i(\eta_i,\tau_\blt)~d\eta_i.\label{eq179}
\end{align}
\end{pp}
Note that each $\Sbf_{\eta_i}\fk P$ (defined in Remark \ref{lb75}) is an element of $\scr O(W_i)$. Also, the residue $\Res_{\eta_i=0}$ is taken with respect to the coordinate $(\eta_i,\tau_\blt)$.

\begin{proof}
It suffices to prove that the claim is locally true. Thus, we may shrinking $\mc B$ and $W_1,\dots,W_N$ such  that for each $1\leq i\leq N$, there exists a coordinate  $\mu_i\in\scr O(W_i)$ at $\sgm_i(\mc B)$ such that $(W_i,\mu_i)\in\fk P$. Then
\begin{align*}
\theta|_{W_i}=a_i(\eta_i,\tau_\blt)\cdot(\partial_{\mu_i}\eta_i)^{-1}\partial_{\mu_i}.
\end{align*}
Our strategy is to compare the action of $\upnu(\theta)$ defined by the coordinates $\mu_\blt$ (which is trivial by Lemma \ref{lb74}) with $\wtd \upnu(\theta)$ defined by $\eta_\blt$. Set $\wtd\upnu(\theta)\in \svir_c\otimes \omega_{\mc C/\mc B}(\blt S_{\fk X})(W)$ such that $\mc U_\varrho(\mu_i)\wtd\upnu(\theta)|_{W_i}=a_i(\eta_i,\tau_\blt)\cdot(\partial_{\mu_i}\eta_i)^{-1}\cbf~d{\mu_i}$. Then
\begin{align*}
\mc U_\varrho(\mu_i)\wtd\upnu(\theta)|_{W_i}=a_i(\eta_i,\tau_\blt)\cdot(\partial_{\mu_i}\eta_i)^{-2}\cbf~d{\eta_i}.
\end{align*}
By Lemma \ref{lb74}, the action of $\wtd\upnu(\theta)$ on $\scr T_{\fk X}(\Wbb_\blt)$ is zero. Notice that the action of $\wtd\upnu(\theta)$ is independent of the choice of local coordinates. (See Lemma \ref{lb125} and the paragraphs before Proposition \ref{lb41}.) By \eqref{eq169}, we have
\begin{align*}
&\mc U_\varrho(\eta_i)\wtd\upnu(\theta)|_{W_i}=\mc U(\varrho(\eta_i|\mu_i))\mc U_\varrho(\mu_i)\wtd\upnu(\theta)|_{W_i}\\
=&a_i(\eta_i,\tau_\blt)\cbf~d{\eta_i}+\frac{c}{12} a_i(\eta_i,\tau_\blt)\cdot(\partial_{\mu_i}\eta_i)^{-2}\Sbf_{\mu_i}\eta_i\cdot \id~d{\eta_i}
\end{align*}
By \eqref{eq178} and \eqref{eq167}, we have
\begin{align*}
\mc U_\varrho(\eta_i)\wtd\upnu(\theta)|_{W_i}=&\mc U_\varrho(\eta_i)\upnu(\theta)|_{W_i}-\frac{c}{12} a_i(\eta_i,\tau_\blt)\cdot\Sbf_{\eta_i}\mu_i\cdot \id~d{\eta_i}\\
=&\mc U_\varrho(\eta_i)\upnu(\theta)|_{W_i}-\frac{c}{12} a_i(\eta_i,\tau_\blt)\cdot\Sbf_{\eta_i}\fk P\cdot \id~d{\eta_i}.
\end{align*}
Since the action of $\wtd\upnu(\theta)$ is zero, the action of $\upnu(\theta)$ equals the sum over $i$ of the actions of $\frac{c}{12} a_i(\eta_i,\tau_\blt)\cdot\Sbf_{\eta_i}\fk P\cdot \id~d{\eta_i}$, which is exactly the scalar multiplication by \eqref{eq179}.
\end{proof}

\section{Convergence of sewing}\label{lb147}

In this section, we assume the setting of Section \ref{lb76}. In particular, $\fk X$ is a family of $N$-pointed complex curves with local coordinates obtained by sewing the smooth family
\begin{align*}
\wtd{\fk X}=(\wtd\pi:\wtd{\mc C}\rightarrow\wtd{\mc B};\sgm_1,\dots,\sgm_N;\sgm_1',\dots,\sgm_M';\sgm_1'',\dots,\sgm_M'';\eta_1,\dots,\eta_N;\xi_1,\dots,\xi_M;\varpi_1,\dots,\varpi_M),
\end{align*}
to which the $\Vbb$-modules $\Wbb_1,\dots,\Wbb_N$, $\Mbb_1,\dots,\Mbb_M'$, and their contragredient modules $\Mbb_1',\dots,\Mbb_M'$ are associated.  Also, we assume throughout this section that \emph{$\Vbb$ is $C_2$-cofinite,  $\Wbb_1,\dots,\Wbb_N$ are finitely generated $\Vbb$-modules, and $\Mbb_1,\dots,\Mbb_M$ (and hence their contragredient modules) are semi-simple}, i.e., they are finite direct sums of irreducible modules.

For each $n\in\Cbb$, let $P_n$ be the projection of each $\Vbb$-module onto its $L_0$-weight $n$ subspace. Recall the notation $q_\blt^{n_\blt}=q_1^{n_1}\cdots q_M^{n_M}$. Given $\uppsi\in\scr T_{\wtd{\fk X}}^*(\Wbb_\blt\otimes\Mbb_\blt\otimes\Mbb_\blt')(\wtd{\mc B})$, we say that $\mc S\uppsi$ \textbf{converges absolutely and locally uniformly (a.l.u.)} if $\mc S\uppsi$ converges in the sense of Remark \ref{lb77}, i.e., if for any $w_\blt\in\Wbb_\blt$ and any compact subsets $K\subset\wtd{\mc B}$ and $Q\subset\mc D_{r_\blt\rho_\blt}^\times$, there exists $C>0$ such that
\begin{align*}
\sum_{n_\blt\in\Cbb^M}\Big|\uppsi\Big(w_\blt\otimes (P_{n_1}\btr\otimes_1\btl)\otimes\cdots\otimes (P_{n_M}\btr\otimes_M\btl)\Big)(b)   \Big|\cdot |q_\blt^{n_\blt}|\leq C
\end{align*}
for any $b\in K$ and $q_\blt=(q_1,\dots,q_M)\in Q$. When $\Mbb_1,\dots,\Mbb_M$ are irreducible (i.e. simple), since $L_0$ and $\wtd L_0$ differ by a scalar multiplication, we have $\wtd{\mc S}\uppsi=q_\blt^{\lambda_\blt}\mc S\uppsi$  where $\lambda_1,\dots,\lambda_M$ are constants. Thus $\mc S\uppsi$ converges a.l.u. if and only if $\wtd{\mc S}\uppsi$ does. The same is true when $\Mbb_1,\dots,\Mbb_M$ are semisimple. Recall also that $\mc S\uppsi$ is a formal conformal block.

As in the proof of Theorem \ref{lb67}, for each $k\in\Nbb$, $\Wbb_\blt^{\leq k}$ (resp. $\Wbb_\blt(k)$) denotes the (finite dimensional) subspace spanned by all $\wtd L_0$-homogeneous homogeneous vectors $w\in\Wbb_\blt$ satisfying $\wtd\wt(w)\leq k$ (resp. $\wtd\wt(w)=k$). This gives a filtration (resp. grading) of $\Wbb_\blt$. We define
\begin{align*}
\mc S\uppsi^{\leq k}\in (\Wbb_\blt^{\leq k})^*\otimes_\Cbb\scr O(\wtd{\mc B})\{q_\blt\}
\end{align*}
whose evaluation with each $w\in\Wbb_\blt^{\leq k}$ is $\mc S\uppsi(w)$. Clearly, the a.l.u. convergence of $\mc S\uppsi$ holds if and only if $\mc S\uppsi^{\leq k}\in (\Wbb_\blt^{\leq k})^*\otimes_\Cbb\scr O(\mc B-\Delta)$ for any large enough $k$. Recall $\mc B=\mc D_{r_\blt\rho_\blt}\times \wtd{\mc B}$.

\begin{thm}\label{lb79}
Assume that $\wtd{\mc B}$ (and hence $\mc B$) is a Stein manifold.  Then there exists $k_0\in\Zbb_+$ such that for any $k\geq k_0$, there exist
\begin{align*}
A_1,\dots,A_M\in \End_{\Cbb}\big((\Wbb_\blt^{\leq k})^* \big)\otimes_\Cbb\scr O(\mc B)
\end{align*}
not depending on $\Mbb_1,\dots,\Mbb_M$, such that for any $1\leq j\leq M$,
\begin{align}
q_j\partial_{q_j} (\mc S\uppsi^{\leq k})=A_j\cdot\mc S\uppsi^{\leq k}.\label{eq180}
\end{align}
\end{thm}

For simplicity, we shall prove this theorem for $M=1$. For general $M$ the idea of the proof is the same.  We set $q=q_1=q_\blt$, $\Mbb=\Mbb_1=\Mbb_\blt$, etc.. In this case, we have $\mc B=\mc D_{r\rho}\times\wtd{\mc B}$ and $\Delta=\{0\}\times\wtd{\mc B}$. We assume that $\Mbb$ and hence $\Mbb'$ are irreducible, so that $\wtd L_0$ and $L_0$ are equal up to a constant. Recall the assumptions of $U',U''$  in Section \ref{lb24}:    $U'$ (resp. $U''$) is a neighborhood of $\sgm'(\wtd {\mc B})$ (resp. $\sgm''(\wtd {\mc B})$) such that
\begin{gather*}
(\xi,\wtd\pi):U'\rightarrow \mc D_{r}\times\wtd{\mc B}\qquad\text{resp.}\qquad (\varpi,\wtd\pi):U''\rightarrow \mc D_{\rho}\times\wtd{\mc B}
\end{gather*}
is a biholomorphic map.  Recall $\SX=\sum_i\sgm_i(\mc B)$. We set  $S_{\wtd{\fk X}}=\sum_{i=1}^N\sgm_i(\wtd{\mc B})+\sgm'(\wtd{\mc B})+\sgm''(\wtd{\mc B})$ to be a divisor of $\wtd{\mc C}$. Let
\begin{align*}
\Gamma=\sgm'(\wtd{\mc B})\cup\sgm''(\wtd{\mc B}).
\end{align*}

Our first step is to show that $\uppsi$ is a formal parallel section in the direction of $q$. Define $\fk y\in\Theta_{\mc B}(-\log\Delta)(\mc B)$ to be $\fk y=q\partial_q$, regarded as constant over $\wtd{\mc B}$. Choose $\wtd{\fk y}\in\Theta_{\mc C}(-\log\mc C_\Delta+\blt S_{\fk X})(\mc C)$ satisfying $d\pi(\wtd\yk)=q\partial_q$ as in \eqref{eq141}. We shall take the series expansion of the vertical part of $\wtd\yk$.


For any open precompact subset $\wtd V\subset \wtd{\mc C}-\Gamma$ and an $\eta\in\scr O(\wtd V)$ univalent on each fiber, choose an open subdisc $\mc D\subset\mc D_{r\rho}$ centered at $0$ with standard coordinate $q$, and assume that $\mc D$ is small enough such that $\mc D\times\wtd V\simeq\wtd V\times\mc D$  can be regarded as an open subset of $\wtd {\mc C}\times\mc D_{r\rho}-F_1'-F_1''$ (recall \eqref{eq112}) and hence of $\mc C$.  Consider $\eta$ also as an element of $\scr O(\mc D\times\wtd V)$ which is constant over $\mc D$. Thus $\partial_q\eta=0$. Then there exists $h\in\scr O(\blt \SX)(\mc D\times\wtd V)$ such that 
\begin{align}
\wtd\yk|_{\mc D\times\wtd V}=h\partial_\eta+q\partial_q.\label{eq239}
\end{align}
Write $h=\sum_{n\in\Nbb}h_nq^n$ where $h_n\in\scr O(\blt S_{\wtd{\fk X}})(\wtd V)$. For each $n\in\Nbb$, set an element $\wtd\yk^\perp_n\in\Theta_{\wtd{\mc C}/\wtd{\mc B}}(\blt S_{\wtd{\fk X}})(\wtd V)$ by
\begin{align}
\wtd\yk^\perp_n|_{\wtd V}=h_n\partial_\eta.\label{eq237}
\end{align}

\begin{lm}
The locally defined $\wtd\yk^\perp_n$ is independent of the choice of $\eta$, and hence can be extended to an element of $\Theta_{\wtd{\mc C}/\wtd{\mc B}}(\blt S_{\wtd{\fk X}})(\wtd {\mc C}-\Gamma)$
\end{lm}

\begin{proof}
Suppose we have another $\mu\in\scr O(\wtd V)$ univalent on each fiber, which is extended constantly to $\mc D\times\wtd V$. Then $\partial_q\mu=0$ and hence $\wtd\yk|_{\mc D\times\wtd V}=h\cdot\partial_\eta\mu\cdot\partial_\mu+q\partial_q.$ Note that $\partial_\eta\mu$ is constant over $q$. Thus, if we define $\wtd\yk^\perp_n|_{\wtd V}$ using $\mu$, then $\wtd\yk^\perp_n|_{\wtd V}=h_n\cdot\partial_\eta\mu\cdot\partial_\mu$, which agrees with \eqref{eq237}.
\end{proof}

We shall show that $\wtd\yk^\perp_n$ has poles of finite orders at $\Gamma$. For that purpose, we need to describe explicitly $\wtd\yk$ near the critical locus $\Sigma$. Let us first recall the geometry of $\fk X$ near  $\Sigma$. By the paragraph containing \eqref{eq115}, any $x'\in\Sigma$ is contained in a neighborhood $W$ of the form
\begin{gather*}
W= \mc D_r\times\mc D_\rho\times\wtd{\mc B},\\
W\cap \Sigma\simeq (0,0)\times\wtd{\mc B},\\
\pi:W=\mc D_r\times\mc D_\rho\times\wtd{\mc B}\xrightarrow{\pi_{r,\rho}\times\id}\mc D_{r\rho}\times\wtd{\mc B}=\mc B,
\end{gather*}
where $\pi_{r,\rho}:\mc D_r\times\mc D_\rho\rightarrow\mc D_{r\rho}$ is the multiplication map. As usual, we let $\xi,\varpi$ be respectively the standard coordinates of $\mc D_r,\mc D_\rho$. Then $(\xi,\varpi)$ is a coordinate of $\mc D_r\times\mc D_\rho$. Set $q=\pi_{r,\rho}=\xi\varpi$. 

In the following, we let $\tau_\blt$ be any biholomorphic map from $\wtd{\mc B}$ to an open subset of a complex manifold. If $\wtd{\mc B}$ is small enough, then $\tau_\blt$ can be a set of coordinates of $\wtd{\mc B}$. The purpose of introducing $\tau_\blt$ is only to indicate the dependence of certain functions on the points of $\wtd{\mc B}$.  Thus, $(\xi,q,\tau_\blt)$ and $(\varpi,q,\tau_\blt)$ are respectively biholomorphic maps of 
\begin{align*}
W'=\mc D_r^\times\times\mc D_\rho\times\wtd{\mc B},\qquad W''=\mc D_r\times\mc D_\rho^\times\times\wtd{\mc B}
\end{align*}
to complex manifolds. By \eqref{eq22}, we can find $a,b\in\scr O((\xi,\varpi,\tau_\blt)(W))$ such that
\begin{align*}
\wtd\yk|_W=a(\xi,\varpi,\tau_\blt)\xi\partial_\xi+b(\xi,\varpi,\tau_\blt)\varpi\partial_\varpi.
\end{align*}
Since $d\pi(\xi\partial_\xi)=d\pi(\varpi\partial_\varpi)=q\partial_q$ by \eqref{eq182}, we must have
\begin{align}
a+b=1.
\end{align}
This relation, together with \eqref{eq74}, shows that under the coordinates $(\xi,q,\tau_\blt)$ and $(\varpi,q,\tau_\blt)$ respectively,
\begin{gather}
\wtd\yk|_{W'}=a(\xi,q/\xi,\tau_\blt)\xi\partial_\xi+q\partial_q,\qquad \wtd\yk|_{W''}=b(q/\varpi,\varpi,\tau_\blt)\varpi\partial_\varpi+q\partial_q.\label{eq238}
\end{gather}

\begin{lm}
For each $n\in\Nbb$, $\wtd\yk_n^\perp$ has poles of orders at most $n-1$ at $\sgm'(\wtd{\mc B})$ and $\sgm''(\wtd{\mc B})$. Consequently, $\wtd\yk^\perp_n$ is an element of $\Theta_{\wtd {\mc C}/\wtd {\mc B}}(\blt S_{\wtd{\fk X}})(\wtd{\mc C})$.
\end{lm}

\begin{proof}
Let us write
\begin{gather*}
a(\xi,\varpi,\tau_\blt)=\sum_{m,n\in\Nbb}a_{m,n}(\tau_\blt)\xi^m\varpi^n,\qquad b(\xi,\varpi,\tau_\blt)=\sum_{m,n\in\Nbb}b_{m,n}(\tau_\blt)\xi^m\varpi^n
\end{gather*}
where $a_{m,n},b_{m,n}\in\scr O(\tau_\blt(\wtd{\mc B}))$. Then
\begin{gather}
a(\xi,q/\xi,\tau_\blt)=\sum_{n\geq 0,l\geq -n}a_{l+n,n}(\tau_\blt)\xi^lq^n,\qquad b(q/\varpi,\varpi,\tau_\blt)=\sum_{m\geq 0,l\geq -m}b_{m,l+m}(\tau_\blt)\varpi^lq^m.\label{eq188}
\end{gather}
Combine these two relations with \eqref{eq237} and \eqref{eq238}, and take the coefficients before $q^n$. We obtain
\begin{align}
\wtd\yk_n^\perp\Big|_{U'-\sgm'(\wtd{\mc B})}=\sum_{l\geq -n}a_{l+n,n}(\tau_\blt)\xi^{l+1}\partial_\xi,\qquad \wtd\yk_n^\perp\Big|_{U''-\sgm''(\wtd{\mc B})}=\sum_{l\geq -n}b_{n,l+n}(\tau_\blt)\varpi^{l+1}\partial_\varpi,\label{eq187}
\end{align}
which finishes the proof.
\end{proof}

One can then let $\upnu(\wtd\yk_n^\perp)$  be a section of $\svir_c\otimes\omega_{\wtd{\mc C}/\wtd{\mc B}}(\blt S_{\wtd{\fk X}})$ defined on $U'\cup U''$ (near $\sgm'(\wtd{\mc B}),\sgm''(\wtd{\mc B})$) and near $\sgm_1(\wtd{\mc B}),\dots,\sgm_N(\wtd{\mc B})$ as in  \eqref{eq178}. Note that $\svir_c$ is defined over $\wtd{\mc C}$. Also, $\upnu(\wtd\yk_n^\perp)$ depends on the local coordinates $\eta_1,\dots,\eta_N,\xi,\varpi$. Recall the correspondence $\partial_\xi\mapsto \cbf d\xi,\partial_\varpi\mapsto \cbf d\varpi$. We calculate the actions of $\upnu(\wtd\yk_n^\perp)$ on $\Mbb$ and on $\Mbb'$ to be respectively
\begin{align}
\Res_{\xi=0}\sum_{l\geq -n}a_{l+n,n}Y_\Mbb(\cbf,\xi)\xi^{l+1}d\xi,\qquad \Res_{\varpi=0}\sum_{l\geq -n}b_{n,l+n}Y_{\Mbb'}(\cbf,\varpi)\varpi^{l+1}d\varpi.\label{eq183}
\end{align}
In the following proofs, we will suppress the symbol $\tau_\blt$ when necessary.

\begin{lm}\label{lb78}
The following equation of elements of $(\Mbb\otimes\Mbb')[[q]]$ is true.
\begin{align}
L_0q^{\wtd L_0}\btr\otimes~\btl=\sum_{n\in\Nbb}\upnu(\wtd\yk_n^\perp)q^{n+\wtd L_0}\btr\otimes~\btl+\sum_{n\in\Nbb}q^{n+\wtd L_0}\btr\otimes ~\upnu(\wtd\yk_n^\perp)\btl\label{eq184}
\end{align}
\end{lm}
Note that as $\Mbb$ is assumed to be irreducible, the equation still holds if $\wtd L_0$ is replaced by $L_0$.

\begin{proof}
It is obvious that $\mc U(\upgamma_1)\cbf=\cbf$, $\xi^{L_0}\cbf=\xi^2\cbf$, $\varpi^{L_0}\cbf=\varpi^2\cbf$. Notice Remark \ref{lb38}. We have  
\begin{align*}
&Y_\Mbb(\xi^{L_0}\cbf,\xi) q^{\wtd L_0}\btr\otimes~\btl\cdot a(\xi,q/\xi)\frac{d\xi}{\xi}\\
=&\sum_{n\geq 0}\sum_{l\geq -n}Y_\Mbb(\cbf,\xi) q^{n+\wtd L_0}\btr\otimes~\btl\cdot a_{l+n,n}\xi^{l+1}d\xi
\end{align*}	
as elements of $(\Mbb\otimes\Mbb'\otimes\scr O(\wtd{\mc B}))((\xi))[[q]]d\xi$.  Take $\Res_{\xi=0}$ and notice \eqref{eq183}. Then, the above expression becomes the first summand on the right hand side of \eqref{eq184}. A similar thing could be said about the second summand. Thus, the right hand side of \eqref{eq184} equals
\begin{align*}
&\Res_{\xi=0}Y_\Mbb(\xi^{L_0}\cbf,\xi) q^{\wtd L_0}\btr\otimes~\btl\cdot a(\xi,q/\xi)\frac{d\xi}{\xi}\\
+&\Res_{\varpi=0} q^{\wtd L_0}\btr\otimes~Y_{\Mbb'}(\varpi^{L_0}\mc U(\upgamma_1)\cbf,\varpi)\btl\cdot b(q/\varpi,\varpi)\frac{d\varpi}{\varpi}.
\end{align*}
By Lemma \ref{lb37} and that $a+b=1$, it equals
\begin{align*}
&\Res_{\xi=0}Y_\Mbb(\xi^{L_0}\cbf,\xi) q^{\wtd L_0}\btr\otimes~\btl\cdot \frac{d\xi}{\xi}=\Res_{\xi=0}Y_\Mbb(\cbf,\xi) q^{\wtd L_0}\btr\otimes~\btl\cdot \xi d\xi\\
=&Y_\Mbb(\cbf)_1 q^{\wtd L_0}\btr\otimes~\btl=L_0 q^{\wtd L_0}\btr\otimes~\btl.
\end{align*}
\end{proof}

\begin{lm}\label{lb81}
For any $w_\blt\in\Wbb_\blt$, we have the following relation of elements of $\scr O(\wtd{\mc B})\{q\}$.
\begin{align*}
q\partial_q\mc S\uppsi(w_\blt)=\sum_{n\in\Nbb}\uppsi(w_\blt\otimes \upnu(\wtd\yk_n^\perp)q^{n+L_0}\btr\otimes~\btl)+\sum_{n\in\Nbb}\uppsi(w_\blt\otimes q^{n+L_0}\btr\otimes ~\upnu(\wtd\yk_n^\perp)\btl).
\end{align*}
\end{lm}

\begin{proof}
We have
\begin{align*}
q\partial_q\mc S\uppsi(w_\blt)=q\partial_q\uppsi(w_\blt\otimes q^{L_0}\btr\otimes~\btl)=\uppsi(w_\blt\otimes L_0 q^{L_0}\btr\otimes~\btl).
\end{align*}
By the Lemma \ref{lb78} (with $\wtd L_0$ replaced by $L_0$), the desired equation is proved.
\end{proof}

As usual, we let $\upnu(\wtd\yk_n^\perp)w_\blt$ denote $\sum_i w_1\otimes\cdots\otimes \upnu(\wtd\yk_n^\perp) w_i\otimes\cdots\otimes w_N$.  Recall $\mc B=\mc D_{r\rho}\times\wtd{\mc B}$. For any $w_\blt\in\Wbb_\blt$, one can define $\nabla_{q\partial_q}w_\blt\in\Wbb_\blt\otimes_\Cbb\scr O(\mc D_{r\rho}\times\wtd {\mc B})$ using \eqref{eq142} and \eqref{eq185}, which equals $\nabla_{q\partial_q}w_\blt=-\upnu(\wtd\yk)w_\blt$. The action of $\upnu(\wtd\yk)$ clearly equals that of $\sum_{n\in\Nbb}q^n \upnu(\wtd\yk^\perp_n)$. Thus, we obtain
\begin{align}
\nabla_{q\partial_q}w_\blt=-\sum_{n\in\Nbb}q^n \upnu(\wtd\yk^\perp_n)w_\blt.\label{eq186}
\end{align}
In particular, the series on the right hand side converges absolutely.

The  following lemma claims that up to a formal projective term, $\mc S\uppsi$ is parallel in the direction of $q\partial_q$, where the connection is defined by the chosen lift $\wtd\yk$. Recall that $\wtd{\mc B}$ is Stein. Thus, we can choose a projective structure $\fk P$ on $\wtd{\fk X}$, which exists due to Theorem \ref{lb71}.

\begin{pp}\label{lb82}
There exists $\#(\wtd\yk_n^\perp)\in\scr O(\wtd{\mc B})$ for each $n\in\Nbb$, such that for any $w_\blt\in\Wbb_\blt$, we have the following equation of elements of $\scr O(\wtd{\mc B})\{q\}$:
\begin{align*}
q\partial_q\mc S\uppsi(w_\blt)=\mc S\uppsi(\nabla_{q\partial_q}w_\blt)+\sum_{n\in\Nbb}\#(\wtd\yk_n^\perp)q^n\cdot\mc S\uppsi(w_\blt).
\end{align*}
\end{pp}

\begin{proof}
$\#(\wtd\yk_n^\perp)$ is defined by Proposition \ref{lb80}. Moreover, by that proposition, we have
\begin{align*}
&w_\blt\otimes \upnu(\wtd\yk_n^\perp)q^{L_0}\btr\otimes~\btl+w_\blt\otimes q^{L_0}\btr\otimes ~\upnu(\wtd\yk_n^\perp)\btl+\upnu(\wtd\yk_n^\perp)w_\blt\otimes q^{L_0}\btr\otimes ~\btl\\
=&\#(\wtd\yk_n^\perp)\cdot w_\blt\otimes q^{L_0}\btr\otimes ~\btl.
\end{align*}
By Lemma \ref{lb81} and relation \eqref{eq186}, it is easy to prove the desired equation.
\end{proof}

We still need one more result before we can prove Theorem \ref{lb79}: the projective term $\sum_{n\in\Nbb}\#(\wtd\yk_n^\perp)q^n$ converges absolutely.

\begin{pp}\label{lb83}
$\sum_{n\in\Nbb}\#(\wtd\yk_n^\perp)q^n$ is an element of $\scr O(\mc D_{r\rho}\times\wtd{\mc B})=\scr O(\mc B)$.
\end{pp}

\begin{proof}
Let $\wtd V_1,\dots,\wtd V_N$ be mutually disjoint neighborhoods of $\sgm_1(\wtd{\mc B}),\dots,\sgm_N(\wtd{\mc B})$ on which $\eta_1,\dots,\eta_N$ are defined. Assume that they are disjoint from $U',U''$. Then $\mc D_{r\rho}\times\wtd V_1,\dots,\mc D_{r\rho}\times\wtd V_N$ are neighborhoods of $\sgm_1(\mc B),\dots,\sgm_N(\mc B)$. Write $\tau_\blt\circ\pi$ also as $\tau_\blt$ for simplicity. Recall \eqref{eq239}. We may write
\begin{align*}
\wtd\yk|_{\mc D_{r\rho}\times\wtd V_i}=h_i(q,\eta_i,\tau_\blt)\partial_{\eta_i}+q\partial_q
\end{align*}
where $h_i(q,\eta_i,\tau_\blt)\in\scr O_{\mc C}(\blt\SX)(\mc D_{r\rho}\times\wtd V_i)$. Write $h_i=\sum_{n}h_{i,n}q^n$. Then by \eqref{eq237},
\begin{align}
\wtd\yk_n^\perp|_{\wtd V_i}=h_{i,n}(\eta_i,\tau_\blt)\partial_{\eta_i}.\label{eq240}
\end{align}

Combine \eqref{eq187} and \eqref{eq240}, and apply Proposition \ref{lb80} to the family $\wtd{\fk X}$. We obtain
\begin{align*}
\#(\wtd\yk_n^\perp)=\frac{c}{12}\big(A_n+B_n+\sum_{i=1}^NC_{i,n}\big)
\end{align*}
where
\begin{gather*}
A_n=\sum_{l\geq -n}\Res_{\xi=0}~\Sbf_\xi\fk P\cdot a_{l+n,n}(\tau_\blt)\xi^{l+1}d\xi,\\
B_n=\sum_{l\geq -n}\Res_{\varpi=0}~\Sbf_\varpi\fk P\cdot b_{n,l+n}(\tau_\blt)\varpi^{l+1}d\varpi,\\
C_{i,n}=\Res_{\eta_i=0}~\Sbf_{\eta_i}\fk P\cdot h_{i,n}(\eta_i,\tau_\blt)d\eta_i.
\end{gather*}
Notice that $\Sbf_{\eta_i}\fk P=\Sbf_{\eta_i}\fk P(\eta_i,\tau_\blt)$, $\Sbf_\varpi\fk P=\Sbf_\varpi\fk P(\varpi,\tau_\blt)$,  $\Sbf_\xi\fk P=\Sbf_\xi\fk P(\xi,\tau_\blt)$ are  holomorphic functions on $\wtd V_i,U',U''$ which are identified with their images under $(\eta_i,\tau_\blt),(\xi,\tau_\blt),(\varpi,\tau_\blt)$   respectively.

We have
\begin{align}
\sum_{n\geq 0} A_nq^n=\sum_{n\geq 0}\sum_{l\geq -n}\Res_{\xi=0}~\Sbf_\xi\fk P\cdot a_{l+n,n}(\tau_\blt)\xi^{l+1}q^nd\xi.\label{eq189}
\end{align}
We claim that \eqref{eq189} is an element of $\scr O(\mc D_{r\rho}\times\wtd{\mc B})$.  Note that $a(\xi,q/\xi,\tau_\blt)$ is defined when $|q|/\rho<|\xi|<r$. Choose any $\epsilon\in (0,r\rho)$. Choose a circle $\gamma'$ surrounding $\mc D_{\epsilon/\rho}$ and inside $\mc D_r$. Then, when $\xi$ is on $\gamma$, $a(\xi,q/\xi,\tau_\blt)$ can be defined whenever $|q|<\epsilon$. Thus, 
\begin{align*}
A:=\frac{1}{2\im\pi}\oint_{\gamma'}\Sbf_\xi\fk P(\xi,\tau_\blt)\cdot a(\xi,q/\xi,\tau_\blt)\xi d\xi
\end{align*}
is a holomorphic function defined whenever  $|q|<\epsilon$. Recall the first equation of \eqref{eq188}, and note that the series converges absolutely and uniformly when $\xi\in\gamma'$ and $|q|\leq \epsilon$, by the double Laurent series expansion of $a(\xi,q/\xi,\tau_\blt)$. Therefore, the above contour integral equals
\begin{align*}
\sum_{n\geq 0}\sum_{l\geq -n}\oint_{\gamma'}\frac{1}{2\im\pi}\Sbf_\xi\fk P(\xi,\tau_\blt)\cdot a_{l+n,n}(\tau_\blt)\xi^{l+1}q^nd\xi,
\end{align*}
which clearly equals \eqref{eq189} as an element of $\scr O(\wtd{\mc B})[[q]]$. Thus \eqref{eq189} is an element of $\scr O(\mc D_\epsilon\times\wtd{\mc B})$ whenever $\epsilon<r\rho$, and hence when $\epsilon=r\rho$.

A similar argument shows  $\sum B_nq^n$ converges a.l.u. to
\begin{align*}
B:=\frac{1}{2\im\pi}\oint_{\gamma''}\Sbf_\varpi\fk P(\varpi,\tau_\blt)\cdot b(q/\varpi,\varpi,\tau_\blt)\varpi d\varpi
\end{align*}
where $\gamma''$ is any circle in $\mc D_\rho$ surrounding $0$. Finally, we compute
\begin{align*}
C_i:=&\sum_{n\geq 0}C_{i,n}q^n=\sum_{n\geq 0}\Res_{\eta_i=0}~\Sbf_{\eta_i}\fk P(\eta_i,\tau_\blt)\cdot h_{i,n}(\eta_i,\tau_\blt)q^nd\eta_i\\
=&\Res_{\eta_i=0}~\Sbf_{\eta_i}\fk P(\eta_i,\tau_\blt)\cdot h_i(q,\eta_i,\tau_\blt)d\eta_i
\end{align*}
which is clearly inside $\scr O(\mc D_{r\rho}\times\wtd{\mc B})$. The proof is now complete. We summarize that the projective term equals
\begin{align}
\sum_{n\in\Nbb}\#(\wtd\yk_n^\perp)q^n=\frac c{12}\big(A+B+\sum_{i=1}^N C_i\big).\label{eq241}
\end{align}
\end{proof}

\begin{proof}[Proof of Theorem \ref{lb79}]
We choose $k_0$ such that the $s_1,s_2,\dots$ found in Theorem \ref{lb67} are in $\Wbb_\blt^{\leq k_0}\otimes_\Cbb\scr O(\mc B)$. We may assume $s_1,s_2,\dots$ form a basis of $\Wbb_\blt^{\leq k_0}$, regarded as constant sections of $\Wbb_\blt^{\leq k_0}\otimes_\Cbb\scr O(\mc B)$. Fix any $k\geq k_0$, and extend $s_1,s_2,\dots$ to a basis of $\Wbb_\blt^{\leq k}$.

By propositions \ref{lb82} and \ref{lb83},  for each $s_i$ of $s_1,s_2,\dots$,  we have the following equation of elements of $\scr O(\wtd{\mc B})\{q\}$:
\begin{align*}
q\partial_q\mc S\uppsi(s_i)=\mc S\uppsi(\nabla_{q\partial_q}s_i)+g\mc S\uppsi(s_i)
\end{align*}
where $g\in\scr O(\mc D_{r\rho}\times\wtd{\mc B})=\scr O(\mc B)$ equals \eqref{eq241} and is hence independent of $s_1,s_2,\dots$. By Theorem \ref{lb67}, we can find $f_{i,j}\in\scr O(\mc B)$ such that $\nabla_{q\partial_q}s_i$ equals $\sum_j f_{i,j}s_j$ mod  sections of $\pi_*\big(\scr V_{\fk X}\otimes\omega_{\mc C/\mc B}(\blt S_{\fk X})\big)\cdot \scr W_{\fk X}(\Wbb_\blt)$. (Indeed, the proof of that theorem shows that the relation holds mod a sum of elements of the form $\wht\upupsilon\cdot w_\blt^\circ\in\pi_*\big(\scr V_{\fk X}\otimes\omega_{\mc C/\mc B}(\blt S_{\fk X})\big)(\mc B)\cdot \scr W_{\fk X}(\Wbb_\blt)(\mc B)$). Since, by Theorem \ref{lb42}, $\mc S\uppsi$ is a formal conformal block, we must have
\begin{align*}
q\partial_q\mc S\uppsi(s_i)=\sum_j f_{i,j}\mc S\uppsi(s_j)+g\mc S\uppsi(s_i).
\end{align*}
The proof is  completed by setting the matrix-valued function $A_1=A$ to be $(f_{i,j}+g\delta_{i,j})_{i,j}$.
\end{proof}

\begin{thm}
Choose $\uppsi\in\scr T_{\wtd{\fk X}}^*(\Wbb_\blt\otimes\Mbb_\blt\otimes\Mbb_\blt')(\wtd{\mc B})$. Then $\mc S\uppsi$ and $\wtd{\mc S}\uppsi$ converge a.l.u..
\end{thm}

We are not assuming $\wtd{\mc B}$ to be Stein.
\begin{proof}
Assume without loss of generality that $\Mbb_1,\dots,\Mbb_M$ are irreducible. Then $\mc S\uppsi$ and $\wtd{\mc S}\uppsi$ differ by $q_\blt^{\lambda_\blt}$ for some $\lambda_1,\dots,\lambda_M\in\Cbb$.

When $M=1$, the a.l.u. convergence follows directly from theorems \ref{lb79} and \ref{lb14}. The general case can be proved by induction. For simplicity, we assume $M=2$ and explain the idea. According to the base case, for any $w_\blt\in\Wbb_\blt$, we know that $\wtd{\mc S}\uppsi(w_\blt)(q_1,q_2,\tau_\blt)$ is an element of $\scr O(\mc D_{r_1\rho_1}\times\wtd{\mc B})[[q_2]]$. Also, by the base case (applied to the smooth family over $\mc D_{r_1\rho_1}^\times\times\wtd{\mc B}$), we know that $\wtd{\mc S}\uppsi(w_\blt)(q_1,q_2,\tau_\blt)$ converges a.l.u. to an element of $\scr O(\mc D_{r_1\rho_1}^\times\times\mc D_{r_2\rho_2}\times\wtd{\mc B})$. This finishes the proof.
\end{proof}

Recall that the logarithmic connection on $\scr T_{\fk X}^*(\Wbb_\blt)$ is dual to the one on $\scr T_{\fk X}(\Wbb_\blt)$. From propositions \ref{lb82} and \ref{lb83} and their proofs, we have:

\begin{thm}\label{lb85}
Assume $\wtd{\mc B}$ is Stein. For each $1\leq k\leq M$, define $\nabla_{q_k\partial_{q_k}}$ on $\scr T_{\fk X}(\Wbb_\blt)$ and hence on $\scr T_{\fk X}^*(\Wbb_\blt)$ using $\eta_\blt$ and a lift $\wtd\yk$ of $q_k\partial_{q_k}$ as in Section \ref{lb84}. Then there exists a projective term $f_k\in\scr O(\mc B)$ such that
\begin{align*}
\nabla_{q_k\partial_{q_k}}\mc S\uppsi=f_k\cdot  \mc S\uppsi.
\end{align*}
\end{thm}

Recall $\mc B=\mc D_{r_\blt\rho_\blt}\times\wtd{\mc B}$. When $M=1$, the projective term $f_1$ is \eqref{eq241}. In the following remark, we give the formula of $f_k$ for a general $M$.

\begin{rem}\label{lb138}
Let $\wtd\yk\in\Theta_{\mc C}(-\log\mc C_\Delta+\blt\SX)(\mc C)$ be a lift of $\yk=q_k\partial_{q_k}$. Choose neighborhoods $U'_1,\dots,U'_M,U''_1,\dots,U''_M$ of $\sgm_1'(\wtd{\mc B}),\dots,\sgm_M'(\wtd{\mc B}),\sgm_1''(\wtd{\mc B}),\dots,\sgm_M''(\wtd{\mc B})$ as in \eqref{eq113}. Choose $\wtd V_1,\dots,\wtd V_N$ to be neighborhoods of $\sgm_1(\wtd{\mc B}),\dots,\sgm_N(\wtd{\mc B})$ as in the proof of Proposition \ref{lb83}. Assume that they are disjoint from $U'_1,\dots,U'_M,U''_1,\dots,U''_M$. Then $\mc D_{r_\blt\rho_\blt}\times\wtd V_1,\dots,\mc D_{r_\blt\rho_\blt}\times\wtd V_N$ are neighborhoods of $\sgm_1(\mc B),\dots,\sgm_N(\mc B)$ in $\mc C$. For each $1\leq i\leq N,1\leq j\leq M$, write
\begin{gather*}
\wtd\yk|_{\mc D_{r_\blt\rho_\blt}\times\wtd V_i}=h_i(q_\blt,\eta_i,\tau_\blt)\partial_{\eta_i}+q_k\partial_{q_k},\\
\wtd\yk|_{\mc D_{r_j}\times\mc D_{\rho_j}\times\mc D_{r_\blt\rho_\blt\setminus j}\times\wtd{\mc B}}=a_j(\xi_j,\varpi_j,q_{\blt\setminus j},\tau_\blt)\xi_j\partial_{\xi_j}+b_j(\xi_j,\varpi_j,q_{\blt\setminus j},\tau_\blt)\varpi_j\partial_{\varpi_j}+(1-\delta_{j,k})q_k\partial_{q_k}.
\end{gather*}
where $h_i,a_j,b_j$ are holomorphic functions on suitable domains, and
\begin{align*}
a_j+b_j=\delta_{j,k}.
\end{align*}
Choose a projective structure $\fk P$ of $\fk X$. Then $\Sbf_{\eta_i}\fk P=\Sbf_{\eta_i}\fk P(\eta_i,\tau_\blt)$, $\Sbf_{\xi_j}\fk P=\Sbf_{\xi_j}\fk P(\xi_j,\tau_\blt)$, $\Sbf_{\varpi_j}\fk P=\Sbf_{\varpi_j}\fk P(\varpi_j,\tau_\blt)$  are  holomorphic functions on $\wtd V_i,U_j',U''_j$ which are identified with their images under $(\eta_i,\tau_\blt),(\xi_j,\tau_\blt),(\varpi_j,\tau_\blt)$   respectively. Choose circles $\gamma_j'\subset\mc D_{r_j}$ and $\gamma_j''\subset\mc D_{\rho_j}$ surrounding $0$. Then the projective term $f_k$ in Theorem \ref{lb85} equals
\begin{align*}
f_k=\frac{c}{12}\cdot\bigg(\sum_{1\leq j\leq M}A_j+\sum_{1\leq j\leq M}B_j+\sum_{1\leq i\leq N}C_i\bigg),
\end{align*}
where
\begin{align*}
A_j=\oint_{\gamma_j'}\Sbf_{\xi_j}\fk P(\xi_j,\tau_\blt)\cdot a_j(\xi_j,q_j/\xi_j,q_{\blt\setminus j},\tau_\blt)\cdot\xi_jd\xi_j,\\
B_j=\oint_{\gamma_j''}\Sbf_{\varpi_j}\fk P(\varpi_j,\tau_\blt)\cdot b_j(q_j/\varpi_j,\varpi_j,q_{\blt\setminus j},\tau_\blt)\cdot \varpi_jd\varpi_j,\\
C_i=\Res_{\eta_i=0}~\Sbf_{\eta_i}\fk P(\eta_i,\tau_\blt)\cdot h_i(q_\blt,\eta_i,\tau_\blt)d\eta_i.
\end{align*}
\end{rem}

\section{Linear independence of sewing}\label{lb104}

We continue the study of sewing, but assume that $\wtd{\fk X}$ is a single compact Riemann surface, and $M=1$, i.e.,
\begin{align*}
\wtd{\fk X}=(\wtd C;x_1,\dots,x_N;x';x'';\eta_1,\dots,\eta_N;\xi;\varpi).
\end{align*}
The main result of this section can be generalized, by induction, to any $M\in\Zbb_+$. As usual, we assume that each connected component of $\wtd C$ contains at least one of $x_1,\dots,x_N$. Let $\mc E$ be a finite set of mutually inequivalent irreducible $\Vbb$-modules. Choose open discs $W'\simeq \mc D_r,W''\simeq\mc D_\rho$ (with respect to the local coordinates $\xi,\varpi$) around $x',x''$ which do not contain $x_1,\dots,x_N$. Assume $\Vbb$ is $C_2$-cofinite and $\Wbb_1,\dots,\Wbb_N$ are finitely-generated $\Vbb$-modules associated to $x_1,\dots,x_N$. 

Let
\begin{align*}
\fk X=(\pi:\mc C\rightarrow\mc D_{r\rho};x_1,\dots,x_N;\eta_1,\dots,\eta_N)
\end{align*}
be the family of complex curves obtained by sewing $\wtd{\fk X}$, where $x_1,\dots,x_N,\eta_1,\dots,\eta_N$ are extended from those of $\wtd C$ and are constant over $\mc D_{r\rho}$. For each $q\in\mc D_{r\rho}^\times$, we define a linear map
\begin{gather}
\fk S_q:\bigoplus_{\Mbb\in\mc E}\scr T_{\wtd{\fk X}}^*(\Wbb_\blt\otimes\Mbb\otimes\Mbb')\rightarrow\scr T_{\fk X_q}^*(\Wbb_\blt),\label{eq191}\\
\bigoplus_\Mbb\uppsi_\Mbb\mapsto \sum_\Mbb\mc S\uppsi_\Mbb(q).\nonumber
\end{gather}
Similarly, one can define $\wtd{\fk S}_q$ by replacing $\mc S$ with $\wtd{\mc S}$. Notice that $\sum_\Mbb\wtd{\mc S}\uppsi_\Mbb(q)=\sum_\Mbb q^{\lambda_\Mbb}\mc S\uppsi_\Mbb(q)$ for some constants $\lambda_\Mbb$ depending only on $\Mbb$. Thus $\fk S_q$ is injective if and only if $\wtd{\fk S}_q$ is.

\begin{thm}\label{lb114}
Choose any $q\in\mc D_{r\rho}^\times$. Then $\fk S_q$ and $\wtd{\fk S}_q$ are injective linear maps.
\end{thm}

\begin{proof}
Let us fix $q_0\in\mc D_{r\rho}^\times$ and let $q$ denote a complex variable. Let us prove that $\fk S_{q_0}$ is injective. Suppose that the finite sum $\sum_\Mbb\mc S\uppsi_\Mbb(q_0)$ equals $0$. We shall prove by contradiction that $\uppsi_\Mbb=0$ for any $\Mbb\in\mc E$. 

Suppose this is not true. Let $\mc F$ be the (finite) subset of all $\Mbb\in\mc E$ satisfying $\uppsi_\Mbb\neq0$. Then $\mc F$ is not an empty set. We first show that $\upphi:=\sum_\Mbb\mc S\uppsi_\Mbb$ (which is a multivalued holomorphic function on $\mc D_{r\rho}^\times$) satisfies $\upphi(q)=0$ for each $q\in\mc D_{r\rho}^\times$. Choose any large enough  $k\in\Nbb$. Then, by Theorem \ref{lb79},  $\upphi^{\leq k}$ satisfies a linear differential equation on $\mc D_{r\rho}^\times$ of the form $\partial_q\upphi^{\leq k}=q^{-1}A\cdot \upphi^{\leq k}$. Moreover, it satisfies the initial condition $\upphi^{\leq k}(q_0)=0$. Thus, by Theorem \ref{lb12}, $\upphi^{\leq k}$ is constantly $0$. So is $\upphi$.

Consider the $\Vbb\times\Vbb$-module $\mbb X:=\bigoplus_{\Mbb\in\mc F}\Mbb\otimes\Mbb'$. Define a linear map $\kappa:\mbb X\rightarrow \Wbb_\blt^*$ as follows. If $m\otimes m'\in\Mbb\otimes\Mbb'$, then the evaluation of $\kappa(m\otimes m')$ with any $w_\blt\in\Wbb_\blt$ is
\begin{align*}
\bk{\kappa(m\otimes m'),w_\blt}=\uppsi_\Mbb(w_\blt\otimes m\otimes m').
\end{align*}
We claim that $\Ker(\kappa)$ is a non-zero subspace of $\mbb X$ invariant under the action of $\Vbb\times \Vbb$. If this can be proved, then, by Lemma \ref{lb87}, we have $\Mbb\otimes\Mbb'\subset \ker(\kappa)$ for some $\Mbb\in\mc F$. Therefore, $\uppsi_\Mbb(w_\blt\otimes m\otimes m')=0$ for any $w_\blt\in\Wbb_\blt$ and $m\otimes m'\in\Mbb\otimes\Mbb'$. Namely, $\uppsi_\Mbb=0$. So $\Mbb\notin\mc F$, which gives a contradiction. 

For any $n\in\Cbb$, let $P_n$ be the projection of $\Mbb$ onto its $L_0$-weight $n$ subspace. Then 
\begin{align*}
\upphi(w_\blt)(q)=\sum_{\Mbb\in\mc F}\sum_{n\in\Cbb} \uppsi_\Mbb(w_\blt\otimes P_n\btr\otimes~\btl)q^n.
\end{align*}
Since this multivalued function is always $0$, by Lemma \ref{lb86}, any coefficient before $q^n$ is $0$. Thus $P_n\btr\otimes~\btl$ (which is an element of $\Mbb_{(n)}\otimes\Mbb_{(n)}'$) is in $\ker(\kappa)$ for any $n$. Thus $\ker(\kappa)$ must be non-empty.

Suppose now that $\sum_j m_j\otimes m'_j\in\Ker(\kappa)$ where each $m_j\otimes m'_j$ belongs to some $\Mbb\otimes\Mbb'$. We set $\uppsi_{\wtd\Mbb}(w_\blt\otimes m_j\otimes m_j')=0$ if $\Mbb,\wtd\Mbb\in\mc F$ and $\Mbb\neq \wtd\Mbb$. Choose any $n\in\Nbb$ and $l\in\Zbb$.  We shall show that $\sum_j Y(u)_l m_j\otimes m_j'\in\Ker(\kappa)$ for any $u\in\Vbb^{\leq n}$. (Here $Y$ denotes $Y_\Mbb$ for a suitable $\Mbb$.) Thus $\Ker(\kappa)$  is $\Vbb\times\id$-invariant. A similar argument will show that $\Ker(\kappa)$ is $\id\times\Vbb$-invariant, and hence $\Vbb\times\Vbb$-invariant.

Set divisors $D_1=x_1+\cdots+x_N$ and $D_2=x'+x''$. Choose a natural number $k_2\geq l$ such that $Y(u)_km_j=Y(u)_km_j'=0$ for any $j$, any $k\geq k_2$,  and any $u\in\Vbb^{\leq n}$. This is possible by the lower truncation property. By Corollary \ref{lb6}, we can choose $k_1\in\Nbb$ such that $H^1(\wtd C,\scr V_{\wtd C}^{\leq n}\otimes\omega_{\wtd C}(k_1D_1-k_2D_2))=0$. Thus, the short exact sequence
\begin{align*}
0\rightarrow \scr V_{\wtd C}^{\leq n}\otimes\omega_{\wtd C}(k_1D_1-k_2D_2) \rightarrow \scr V_{\wtd C}^{\leq n}\otimes\omega_{\wtd C}(k_1D_1-lD_2)\rightarrow\scr G\rightarrow 0
\end{align*}
(where $\scr G$ is the quotient of the previous two sheaves) induces another one
\begin{align*}
&0\rightarrow H^0\big(\wtd C,\scr V_{\wtd C}^{\leq n}\otimes\omega_{\wtd C}(k_1D_1-k_2D_2)\big) \rightarrow H^0\big(\wtd C,\scr V_{\wtd C}^{\leq n}\otimes\omega_{\wtd C}(k_1D_1-lD_2)\big)\\
&\rightarrow H^0(\wtd C,\scr G)\rightarrow 0.
\end{align*}
Choose any $u\in\Vbb^{\leq n}$. Choose $v\in H^0(\wtd C,\scr G)$ to be $\mc U_\varrho(\xi)^{-1}u\xi^ld\xi$ in $W'$ and $0$ in $\wtd C-\{x'\}$. Then $v$ has a lift $\nu$ in $H^0\big(\wtd C,\scr V_{\wtd C}^{\leq n}\otimes\omega_{\wtd C}(k_1D_1-lD_2)\big)$, which must be of the form
\begin{gather*}
\mc U_\varrho(\xi)\nu|_{W'}=u\xi^ld\xi+\xi^{k_2}(\mathrm{elements~of~}\Vbb^{\leq n}\otimes_\Cbb\scr O(W'))d\xi,\\
\mc U_\varrho(\varpi)\nu|_{W''}=\varpi^{k_2}(\mathrm{elements~of~}\Vbb^{\leq n}\otimes_\Cbb\scr O(W''))d\varpi.
\end{gather*}
It is clear that $\nu\cdot m_j=Y(u)_lm_j$ and $\nu\cdot m_j'=0$. Thus, as each $\uppsi_\Mbb$ vanishes on $\nu\cdot(\Wbb_\blt\otimes\Mbb\otimes\Mbb')$, we have
\begin{align*}
&\sum_{\Mbb\in\mc F}\sum_j\uppsi_\Mbb(w_\blt\otimes Y(u)_lm_j\otimes m_j')=-\sum_{\Mbb\in\mc F}\sum_j\uppsi_\Mbb((\nu\cdot w_\blt)\otimes m_j\otimes m_j')\\
=&-\sum_j\bk{\kappa(m_j\otimes m_j'),\nu\cdot w_\blt}=0.
\end{align*}
So $\sum_j Y(u)_l m_j\otimes m_j'\in\Ker(\kappa)$.
\end{proof}

\begin{lm}\label{lb87}
If $\Ker(\kappa)$ is non-zero subspace of $\mbb X$ invariant under the action of $\Vbb\times\Vbb$, then it contains an irreducible summand $\Mbb\otimes\Mbb'$.
\end{lm}

\begin{proof}
Set $\mbb A=\Ker(\kappa)$. By basic representation theory, $\mbb A$ contains an irreducible $\Vbb\times\Vbb$-submodule $\mbb A_0$.  (See for instance \cite{Lang} section XVII.2.) Let $\iota:\mbb A_0\rightarrow\mbb X$ be the inclusion, and let $p_\Mbb:\mbb X\rightarrow\Mbb\otimes\Mbb'$ be the standard projection. Then $p_\Mbb\circ\iota:\mbb A_0\rightarrow\Mbb\otimes\Mbb'$ is non-zero for some $\Mbb\in\mc F$. Since $\Mbb\otimes\Mbb'$ is an irreducible $\Vbb\times\Vbb$-module (see Proposition \ref{lb113} and Theorem \ref{lb108}), $p_\Mbb\circ\iota$ is surjective. Since $\mbb A_0$ is irreducible and $p_\Mbb\circ\iota$ is non-zero, the kernel of $p_\Mbb\circ\iota$ must be empty. So $p_\Mbb\circ\iota$ is an isomorphism. Hence $\mbb A_0$ is isomorphic to $\Mbb\otimes\Mbb'$. If $p_{\wtd\Mbb}\circ\iota$ is non-zero for another $\wtd\Mbb\in\mc F$, then the same argument shows $\mbb A_0\simeq \wtd\Mbb\otimes\wtd\Mbb'$, which is impossible since $\Mbb\otimes\Mbb'$ is not isomorphic to $\wtd\Mbb\otimes\wtd\Mbb'$. So $\mbb A_0$ must be exactly $\Mbb\otimes\Mbb'$.
\end{proof}

\begin{lm}\label{lb86}
Let $E$ be a finite subset of $\Cbb$. Choose an element
\begin{align*}
f(z)=\sum_{\alpha\in E+\Nbb}c_\alpha z^\alpha
\end{align*}
of $\Cbb\{z\}$. Let $\epsilon>0$. Assume that $f(z)$ converges absolutely to $0$ on $\mc D_\epsilon^\times$. Namely, for any  $z\in\mc D_r^\times$, there is $C>0$ such that
\begin{align}
\sum_{\alpha\in E+\Nbb}|c_\alpha z^\alpha|\leq C,\label{eq190}
\end{align} 
and the infinite sum $\sum_{\alpha\in E+\Nbb}c_\alpha z^\alpha$ converges to $0$. Then $c_\alpha=0$ for each $\alpha$.
\end{lm}

\begin{rem}\label{lb89}
Note that $f(z)$ can be written as $z^{\alpha_1}f_1(z)+\cdots+z^{\alpha_n}f_n(z)$ where $f_1,\dots,f_n\in\Cbb[[z]]$, and any two of $\alpha_1,\dots,\alpha_n$ do not differ by an integer. It is easy to see that $f(z)$ converges absolutely on $\mc D_\epsilon^\times$ if and only if $f_1,\dots,f_n\in\scr O(\mc D_\epsilon)$. Indeed, the if part is obvious; the only if part follows from the root test. Moreover, it is clear that $f_1,\dots,f_n\in\scr O(\mc D_\epsilon)$ implies  $f(z)$ converges a.l.u. on $\mc D_\epsilon^\times$, i.e.,  for compact subset $K\subset\mc D_\epsilon^\times$, there is $C>0$ such that \eqref{eq190} holds for any $z\in K$.
\end{rem}

\begin{proof}[Proof of Lemma \ref{lb86}]
Assume that the coefficients of $f$ are not all $0$. We can let $r\in\Rbb$ be the smallest number such that $c_\alpha\neq 0$ for some $\alpha\in E+\Nbb$ satisfying $\mathrm{Re}(\alpha)=r$. Let $\beta_1=r+\im s_1,\dots,\beta_k=r+\im s_k$ be all the elements of $E+\Nbb$ whose real parts are $r$. (So $s_1,\dots,s_k\in\Rbb$.) Notice $s_i\neq s_j$ when $i\neq j$. Then 
\begin{align*}
z^{-r}f(z)=c_{\beta_1}z^{\im s_1}+\cdots +c_{\beta_k}z^{\im s_k}+g(z)
\end{align*}
where $g(z)\in\Cbb\{z\}$ can be written as
\begin{align*}
g(z)=z^{\gamma_1}h_1(z)+\cdots+z^{\gamma_m}h_m(z)
\end{align*}
for some $h_1,\dots,h_m\in\Cbb[[z]]$,  the real parts of $\gamma_1,\dots,\gamma_m\in\Cbb$ are all $>0$, and any two of $\gamma_1,\dots,\gamma_m$ do not differ by an integer. Since $z^{-r}f(z)$ converges absolutely on $\mc D_\epsilon^\times$, we have $h_1,\dots,h_m\in\scr O(\mc D_\epsilon)$ by Remark \ref{lb89}. Let $t$ be a real variable.  Then it is easy to see by induction  that for any $j\in\Nbb$,
\begin{align*}
\lim_{t\rightarrow-\infty}\partial_t^jg(e^t)=0.
\end{align*}
Since the $j$-th derivative of $e^{-rt}f(e^t)$ over $t$ is constantly $0$, we have
\begin{align*}
\lim_{t\rightarrow-\infty}(c_{\beta_1}s_1^{j-1}e^{\im s_1 t}+\cdots +c_{\beta_k}s_k^{j-1}e^{\im s_k t})=0.
\end{align*}
Let $A=(s_i^{j-1})_{1\leq i,j\leq k}\in M_{k\times k}(\Cbb)$.  Then
\begin{align*}
\lim_{t\rightarrow-\infty} (c_{\beta_1} e^{\im s_1t},\dots,c_{\beta_k} e^{\im s_k t})\cdot A=0.
\end{align*} 
Since $A$ is a Vandermonde matrix which is invertible, we conclude that $c_{\beta_1} e^{\im s_1t},\dots,c_{\beta_k} e^{\im s_k t}$ all converge to $0$, which forces $c_{\beta_1},\dots,c_{\beta_k}$ to be $0$. This gives a contradiction.
\end{proof}

\begin{co}\label{lb118}
There are finitely many equivalence classes of irreducible $\Vbb$-modules.
\end{co}

Recall that $\Vbb$ is assumed to be $C_2$-cofinite.

\begin{proof}
We let $\wtd{\fk X}=(\Pbb^1;1;0;\infty;z-1;z;1/z)$. Set $r=\rho=1$. Let $\fk X$ be obtained by sewing $\wtd{\fk X}$. Then, for each $q\in\mc D_1^\times$, $\mc C_q$ is of genus $1$. Choose $\mc E$ to be an arbitrary finite set of mutually inequivalent irreducible $\Vbb$-modules. We claim that the cardinality of $\mc E$ is bounded by dimension of $\scr T_{\fk X_q}^*(\Vbb)$ which is finite by Theorem \ref{lb67}. Indeed, consider the map $\fk S_q$ in \eqref{eq191}, where we set $\Wbb_\blt$ to be $\Vbb$. For each $\Mbb\in\mc E$, the vector space $\scr T_{\wtd {\fk X}}^*(\Vbb\otimes\Mbb\otimes\Mbb')$ is nontrivial by the construction in example \ref{lb44}. Thus, the dimension of the domain of $\fk S_q$ must be no less than the cardinality of $\mc E$, which is bounded by the dimension of $\scr T_{\fk X_q}^*(\Vbb)$ since $\fk S_q$ is injective.
\end{proof}

\section{More on VOA modules}

We fix a  VOA $\Vbb$. In this section, we do not assume $\Vbb$ is $C_2$-cofinite.

\subsection*{A criterion on weak $\Vbb$-modules}

Let $\Wbb$ be a vector space. Let $Y_\Wbb$ associates to each $v\in\Vbb$ and $n\in\Zbb$ an operator $Y_\Wbb(v)_n\in\End(\Wbb)$, and assume that the map $v\in\Vbb\mapsto Y_\Wbb(v)_n$ is linear. Set $Y_\Wbb(v,z)=\sum_{n\in\Zbb}Y_\Wbb(v)_nz^{-n-1}\in \End(\Wbb)[[z,z^{-1}]]$. We assume that for each $v\in\Vbb$ and $w\in\Wbb$, we have $Y_\Wbb(v)_nw=0$ for any $n$ small enough, i.e.,
\begin{align*}
Y_\Wbb(v,z)w\in\Wbb((z)).
\end{align*}
We say that $(\Wbb,Y_\Wbb)$ (or simply $\Wbb$) is a (lower-truncated) \textbf{linear representation} of $\Vbb$. 

Let $\Wbb^\circ$ be a subspace of the dual space $\Wbb^*$ of $\Wbb$. We say that $\Wbb^\circ$ is \textbf{dense}, if for any $w\in\Wbb$, we have $w=0$ iff $\bk{w',w}=0$ for any $w'\in\Wbb^\circ$. 

\begin{pp}\label{lb94}
Let $\Wbb$ be a linear representation of $\Vbb$. Assume that  $Y_\Wbb(\id,z)=\id_\Wbb$. Assume  that $\Wbb^*$ has a dense subspace $\Wbb^\circ$ satisfying the following condition: For each $w'\in\Wbb^\circ$, there exist $\epsilon>0$ such that
\begin{enumerate}[label=(\arabic*)]
	\item for each $v\in\Vbb,w\in\Wbb$, $\bk{Y_\Wbb(v,z)w,w'}\in\Cbb((z))$ is the laurent series expansion of an element of $\scr O(\mc D_\epsilon^\times)$;
	\item for each $u,v\in\Vbb,w\in\Wbb$, there exists $f=f(\zeta,z)\in\scr O(\Conf_2(\mc D_\epsilon^\times))$ such that for any $n\in\Zbb$ and $z\in\mc D_\epsilon^\times$,
	\begin{gather}
		\bk{Y_\Wbb(v,z)Y_\Wbb(u)_nw,w'}=\Res_{\zeta=0}~f(\zeta,z)\zeta^nd\zeta,\label{eq193}\\
		\bk{Y_\Wbb(Y(u)_nv,z)w,w'}=\Res_{\zeta-z=0}~f(\zeta,z)(\zeta-z)^nd\zeta,\label{eq194}
	\end{gather}
	and for any $n\in\Zbb$ and $\zeta\in\mc D_\epsilon^\times$,
	\begin{gather}
		\bk{Y_\Wbb(u,\zeta)Y_\Wbb(v)_nw,w'}=\Res_{z=0}~f(\zeta,z)z^ndz.\label{eq195}
	\end{gather}
\end{enumerate}
Then $\Wbb$ is a weak $\Vbb$-module.
\end{pp}

(2) says that for each $z\in\mc D_\epsilon^\times$, $\bk{w',Y_\Wbb(v,z)Y_\Wbb(u,\zeta)w}\in\Cbb((\zeta))$ and $\bk{w',Y_\Wbb(Y(u,\zeta-z)v,z)w}\in\Cbb((\zeta-z))$ are respectively the Laurent series expansions of $f(\zeta,z)$ at $\zeta=0$ and $\zeta=z$; for each $\zeta\in\mc D_\epsilon^\times$
, $\bk{Y_\Wbb(u,\zeta)Y_\Wbb(v,z)w,w'}\in\Cbb((z))$ is the Laurent series expansion of $f(\zeta,z)$ at $z=0$.

We remark that in practice, we can often choose $\epsilon$ independent of $w'$.

\begin{proof}
Choose circles $C_1,C_2,C_3$ in $\mc D_\epsilon^\times$ surrounding $0$ with radii $r_1<r_2<r_3$.  For each $z\in C_2$, choose a circle $C(z)$ centered at $z$ whose radius is less than $r_2-r_1$ and $r_3-r_2$. Choose any $m,n\in\Zbb$. Then $P(z)=Q(z)-R(z)$ where
\begin{gather*}
P(z)=\oint_{C(z)}f(\zeta,z)\zeta^m(\zeta-z)^nd\zeta,\\
Q(z)=\oint_{C_3}f(\zeta,z)\zeta^m(\zeta-z)^nd\zeta,\\
R(z)=\oint_{C_1}f(\zeta,z)\zeta^m(\zeta-z)^nd\zeta.
\end{gather*}
As in the proof of Theorem \ref{lb19}, we may use \eqref{eq193} and \eqref{eq194} to calculate that
\begin{gather*}
P(z)=\sum_{l\in\Nbb}{m\choose l}\bk{Y_\Wbb(Y(u)_{n+l}v,z)w,w'}z^{m-l},\\
R(z)=\sum_{l\in\Nbb}{n\choose l}(-1)^{n-l}\bk{Y_\Wbb(v,z)Y_\Wbb(u)_{m+l}w,w'}z^{n-l}.
\end{gather*}
Thus, for any $h\in\Zbb$,
\begin{gather*}
\oint_{C_2} P(z)z^hdz=\sum_{l\in\Nbb}{m\choose l}\bk{Y_\Wbb(Y(u)_{n+l}v)_{m+h-l}w,w'}z^{m-l},\\
\oint_{C_2}R(z)z^hdz=\sum_{l\in\Nbb}{n\choose l}(-1)^{n-l}\bk{Y_\Wbb(v)_{n+h-l}Y_\Wbb(u)_{m+l}w,w'}z^{n-l}.
\end{gather*}
Also, by \eqref{eq195}, it is not hard to see that
\begin{align*}
&\oint_{C_2} Q(z)z^hdz=\oint_{C_2}\oint_{C_3}f(\zeta,z)z^h\zeta^m(\zeta-z)^nd\zeta dz\\
=&\oint_{C_3}\oint_{C_2}f(\zeta,z)z^h\zeta^m(\zeta-z)^ndz d\zeta\\
=&\oint_{C_3}\sum_{l\in\Nbb}{n\choose l}(-1)^{l}\bk{Y_\Wbb(u,\zeta)Y_\Wbb(v)_{h+l}w,w'}\zeta^{m+n-l}d\zeta\\
=&\sum_{l\in\Nbb}{n\choose l}(-1)^{l}\bk{Y_\Wbb(u)_{m+n-l}Y_\Wbb(v)_{h+l}w,w'}.
\end{align*}
This proves the Jacobi identity \eqref{eq192} since we have $\oint_{C_2}P(z)z^hdz=\oint_{C_2}Q(z)z^hdz-\oint_{C_2}R(z)z^hdz$.
\end{proof}

\subsection*{Lowest weight weak $\Vbb$-modules}

Most results in this section is well-known.

Let $(\Wbb,Y_\Wbb)$ be a weak $\Vbb$-module. We define its \textbf{lowest weight subspace} to be \index{zz@$\Omega(\Wbb),\Omega_{+-}(\Wbb)$}
\begin{align*}
\Omega(\Wbb)=\{w\in\Wbb:Y_\Wbb(v)_nw=0\text{ for any homogeneous }v\in\Vbb,n\geq \wt(v)\}.
\end{align*}
Using \eqref{eq242}, it is easy to see that $Y(v)_n\Omega(\Wbb)\subset\Omega(\Wbb)$ when $\wt(v)=n+1$.  If $\Wbb$ is admissible, then the lowest $\wtd L_0$-weight space of $\Wbb$ is a subspace of $\Omega(\Wbb)$. In particular, $\Omega(\Wbb)$ is non-trivial if $\Wbb$ is so. Also, if $(\Wbb_i)_i$ is a collection of weak $\Vbb$-modules, then
\begin{align}
\Omega\Big(\bigoplus_i \Wbb_i\Big)=\bigoplus_i\Omega(\Wbb_i).\label{eq205}
\end{align}
Indeed, it is clear that a vector $w$ in the direct sum is annihilated by any $Y_\Wbb(v)_n$ (where $n\geq\wt(v)$) if and only if each component of $w$ is so. We say that $\Wbb$ is a \textbf{lowest weight weak $\Vbb$-module}, if $\Wbb$ is generated by the lowest weight vectors, i.e., vectors in $\Omega(\Wbb)$.

It will be interesting to know if a weak module has non-trivial lowest weight subspace. The following lemma provides a criterion.

\begin{lm}\label{lb96}
Assume that $\Wbb$ is admissible, and let $\Mbb$ be a non-trivial weak $\Vbb$-submodule of $\Wbb$. Then $\Omega(\Mbb)$ is non-trivial.
\end{lm}

Note that in general, we always have the obvious relation $\Omega(\Mbb)=\Omega(\Wbb)\cap \Mbb$.

\begin{proof}
Choose any $w\in\Mbb$. For each $k\in\Nbb$, let $\Wbb^{\leq k}$ be the subspace spanned by the $\wtd L_0$-homogeneous vectors with weights $\leq k$. Since $\Mbb=\bigcup_{k\in\Nbb}(\Mbb\cap\Wbb^{\leq k})$, we can find the smallest $k$ such that $\Mbb\cap\Wbb^{\leq k}$ is non-trivial. If $v\in\Vbb$ is homogeneous, $n\in\Zbb$, and $n\geq\wt(v)$. Then $Y(v)_n\Wbb^{\leq k}\subset\Wbb^{\leq k-1}$ by \eqref{eq203}. Thus $Y(v)_n(\Mbb\cap\Wbb^{\leq k})\subset\Mbb\cap\Wbb^{\leq k-1}=\{0\}$. So $\Mbb\cap\Wbb^{\leq k}$ is a non-trivial subspace of $\Omega(\Mbb)$.
\end{proof}

Let $\fk U(\Vbb)$ be the universal unital associative algebra generated freely by the elements $(v,n)$ where $v\in\Vbb$ is homogeneous and $n\in\Zbb$. \index{UV@$\fk U(\Vbb),\fk U_0(\Vbb)$} Then, $\Wbb$ is a representation of $\fk U(\Vbb)$ defined by $(v,n)\cdot w=Y_\Wbb(v)_nw$. We say that $(v,n)$ is \textbf{raising} (resp. \textbf{lowering}, \textbf{strictly raising}, \textbf{strictily lowering}) if $\wt(v)-n-1$ is $\geq 0$ (resp. $\leq 0$, $>0$, $<0$). Then $\Omega(\Wbb)$ is precisely the subspace of vectors annihilated by all strictly lowering elements.

We let $\fk U_0(\Vbb)$ be the subalgebra of $\fk U(\Vbb)$ of generated by $1$ and $(v_1,n_1)\cdots (v_k,n_k)$ where $k\in\Zbb_+$ and 
\begin{align*}
\sum_{i=1}^k (\wt(v_i)-n_i-1)=0.
\end{align*}
Then $\fk U_0(\Vbb)\Omega(\Wbb)\subset\Omega(\Wbb)$. Thus, $\Omega(\Wbb)$ is a representation of $\fk U_0(\Vbb)$. Using the commutator formula \eqref{eq242} to move all the strictly lowering elements to the right, we see that the action of $\fk U_0(\Vbb)$ on $\Omega(\Wbb)$ is determined by $(v,n)$ where $v$ is homogeneous and $n=\wt(v)-1$.

\begin{pp}\label{lb98}
Assume that $\Wbb$ is an irreducible admissible $\Vbb$-module. Then $\Omega(\Wbb)$ equals the lowest nontrivial $\wtd L_0$-weight space, and  is an irreducible $\fk U_0(\Vbb)$-module.
\end{pp}
\begin{proof}
Assume without loss of generality that $\Wbb(0)$ is the lowest nontrivial $\wtd L_0$-weight space. Clearly $\Wbb(0)\subset \Omega(\Wbb)$. If $w\in\Omega(\Wbb)$, we let $w(k)$ be the (non-zero) component of $w$ in $\bigoplus_{n\in\Nbb}\Wbb(n)$ with the largest weight $k$. We claim $k=0$, which shows $w\in\Wbb(0)$, and hence $\Omega(\Wbb)\subset\Wbb(0)$. By the irreducibility of $\Wbb$, there exists  $x\in\fk U(\Vbb)$ lowering the $\wtd L_0$-weights by $k$, such that $xw(k)$ is a non-zero vector of $\Wbb(0)$. We must have $xw=xw(k)$, which shows $xw$ is  non-zero.  Suppose $k>0$. Then $x$ must contain strictly lowering components. By using the commutator formula \eqref{eq242} to move all the strictly lowering components of $x$ to the right, we can find a strictly lowering element whose action on $w$ is non-zero. This contradicts $w\in\Omega(\Wbb)$.

Let $\mbf U$ be any non-trivial $\fk U_0(\Vbb)$-invariant subspace of $\Wbb(0)$. By \eqref{eq203}, it is easy to see that $\fk U(\Vbb)\mbf U\cap \Wbb(0)\subset\fk U_0(\Vbb)\mbf U$. So $\fk U(\Vbb)\mbf U\cap \Wbb(0)\subset\mbf U$. $\fk U(\Vbb)\mbf U$ is clearly a non-trivial weak $\Vbb$-submodule of $\Wbb$. Thus, by the irreducibility of $\Wbb$, we have $\fk U(\Vbb)\mbf U=\Wbb$. Hence $\Wbb(0)=\mbf U$. 
\end{proof}

\begin{co}\label{lb130}
Let $\Wbb$ be an irreducible (ordinary)  $\Vbb$-module. Then $\Omega(\Wbb)$ equals the lowest non-trivial eigenspace of $L_0$. In particular, $\Omega(\Wbb)$ is non-trivial and finite-dimensional.
\end{co}

\begin{proof}
Set $\wtd L_0=L_0$ and apply Proposition \ref{lb98}.
\end{proof}

\begin{rem}\label{lb111}
Let $\Wbb$ be an irreducible (ordinary) $\Vbb$-module. Choose $\wtd L_0$ whose lowest weight is $0$. We know that $\Omega(\Wbb)=\Wbb(0)$ and $\Omega(\Wbb')=\Wbb'(0)$. Note that $\Wbb'(0)$ is the dual space of $\Wbb(0)$. Moreover, if $v\in\Vbb$ is homogeneous,  by \eqref{eq211}, we know that for any $w\in\Wbb(0),w'\in\Wbb'(0)=\Wbb(0)^*$,
\begin{align}
\bk{Y_{\Wbb'}(v)_{\wt(v)-1}w',w}=\sum_{l\in\Nbb}\frac{(-1)^{\wt(v)}}{l!}\bk{w',Y_\Wbb(L_1^lv)_{\wt(v)-l-1}w}.\label{eq215}
\end{align}
Thus, the action of $\fk U_0(\Vbb)$ on $\Wbb'(0)$ is completely determined by that of $\fk U_0(\Vbb)$ on $\Wbb(0)$. We will see a stronger result in Section \ref{lb110}, that the irreducible $\Vbb$-module $\Wbb$ is completely determined by the $\fk U_0(\Vbb)$-module $\Wbb(0)$.
\end{rem}

\begin{pp}\label{lb97}
Assume that $(\Wbb,Y_\Wbb)$ is a lowest weight weak $\Vbb$-module with finite dimensional $\Omega(\Wbb)$. Then $\Wbb$ is an admissible $\Vbb$-module.
\end{pp}

\begin{proof}
For any $\lambda\in\Cbb$, we let $\Wbb_{[\lambda]}$ be the subspace of all $w\in \Wbb$ satisfying $(L_0-\lambda)^Nw=0$ for some $N\in\Zbb_+$. For any homogeneous $v\in\Vbb$ and $m\in\Zbb$, we have by \eqref{eq104} that
\begin{align*}
Y_\Wbb(v)_m(L_0-\lambda)=(L_0-(\wt(v)+\lambda-m-1))Y_\Wbb(v)_m.
\end{align*}
Thus
\begin{align}
Y_\Wbb(v)_m(L_0-\lambda)^N=(L_0-(\wt(v)+\lambda-m-1))^NY_\Wbb(v)_m,
\end{align}
which shows 
\begin{align}
Y_\Wbb(v)_m\Wbb_{[\lambda]}\subset\Wbb_{[\lambda+\wt(v)-m-1]}.\label{eq204}
\end{align}
That $\bigvee_{\lambda\in\Cbb}\Wbb_{[\lambda]}=\bigoplus_{\lambda\in\Cbb}\Wbb_{[\lambda]}$ follows as in the finite dimensional case: Suppose  $v_1+\cdots+v_n=0$ where $(L_0-\lambda_i)^Nv_i=0$ for each $1\leq i\leq n$ and $\lambda_i\neq \lambda_j$ when $i\neq j$. Set polynomials $p(x)=(x-\lambda_1)^N,q(x)=(x-\lambda_2)^N\cdots(x-\lambda_n)^N$. Then $p(L_0)v_1=0$ and $q(L_0)v_1=q(L_0)(v_1+v_2+\cdots+v_N)=0$. Since $p(x),q(x)$ have no common divisor other than $1$, there exist polynomials $a(x),b(x)$ such that $a(x)p(x)+b(x)q(x)=1$. So $v_1=a(L_0)p(L_0)v_1+b(L_0)q(L_0)v_1=0$. Similar argument shows $v_2=\cdots=v_n=0$.

We say that a vector $w\in\Wbb$ is a generalized eigenvector of $L_0$ if $w\in\Wbb_{[\lambda]}$ for some $\lambda\in\Cbb$. Since $L_0\Omega(\Wbb)\subset\Omega(\Wbb)$ and $\Omega(\Wbb)$ is finite dimension, by the Jordan canonical form for $L_0|_{\Omega(\Wbb)}$,  $\Omega(\Wbb)$ must be spanned by generalized eigenvectors of $L_0$. By \eqref{eq204}, the same is true for $\Wbb$. We thus have grading $\Wbb=\bigoplus_{\lambda\in\Cbb}W_{[\lambda]}$. Moreover, by \eqref{eq204} and that $\Omega(\Wbb)$ is finite-dimensional, we may find a finite subset $E\subset\Cbb$ such that $\Wbb=\bigoplus_{\lambda\in E+\Nbb}W_{[\lambda]}$, and that any two elements of $E$ do not differ by an integer.  Set $\Wbb(n)$ to be $\Wbb_{[\lambda]}$ if $\lambda-n\in E$ for some $\lambda\in E+\Nbb$ (such $\lambda$ must be unique); otherwise, set $\Wbb(n)=0$. Then \eqref{eq203} is satisfied, and $\Wbb=\bigoplus_{n\in\Nbb}\Wbb(n)$. Thus, $\Wbb$ is admissible.
\end{proof}

\begin{thm}\label{lb102}
Assume that $(\Wbb,Y_\Wbb)$ is a lowest weight weak $\Vbb$-module with finite dimensional $\Omega(\Wbb)$.
\begin{enumerate}[label=(\arabic*)] 
\item $\Mbb$ is an irreducible weak $\Vbb$-submodule of $\Wbb$ if and only if $\Mbb$ is generated by $\Omega(\Mbb)$, and $\Omega(\Mbb)$ is an irreducible $\fk U_0(\Vbb)$-module.  In that case, $\Mbb$ is an admissible $\Vbb$-module.

\item There is a 1-1 correspondence between irreducible weak $\Vbb$-submodules $\Mbb$ of $\Wbb$ and irreducible $\fk U_0(\Vbb)$-submodules $\mbf U$ of $\Omega(\Wbb)$. The relation is given by $\mbf U=\Omega(\Mbb)$ and $\Mbb=\fk U(\Vbb)\mbf U$.

\item $\Wbb$ is completely reducible if and only if $\Omega(\Wbb)$ is so. In that case, there are finitely many irreducible weak $\Vbb$-submodules  $\Mbb_1,\dots,\Mbb_n$ (which are admissible) such that
\begin{gather*}
\Mbb=\bigoplus_{i=1}^n\Mbb_i,\qquad\Omega(\Mbb)=\bigoplus_{i=1}^n\Omega(\Mbb_i).
\end{gather*}
\end{enumerate}	
\end{thm}

\begin{proof}
By Proposition \ref{lb97}, $\Wbb$ is admissible.

(1): Let $\Mbb$ be a weak $\Vbb$-submodule of $\Wbb$. Note that by Lemma \ref{lb96}, $\Omega(\Mbb)$ is nontrivial. Assume first of all that $\Mbb$ is an irreducible (non-trivial) weak $\Vbb$-module. Then $\Mbb$ is generated by any non-trivial subspace, and hence by $\Omega(\Mbb)$. Since $\Omega(\Mbb)\subset\Omega(\Wbb)$, $\Omega(\Mbb)$ is finite-dimensional. By Proposition \ref{lb97}, $\Mbb$ is admissible. Thus, by Proposition \ref{lb98}, $\Omega(\Mbb)$ is irreducible. Conversely, assume $\Omega(\Mbb)$ is irreducible. Let $\Mbb_1$ be a non-trivial weak $\Vbb$-submodule of $\Mbb$. By Lemma \ref{lb96}, $\Omega(\Mbb_1)$ is a non-trivial $\fk U_0(\Vbb)$-submodule of $\Omega(\Mbb)$. Thus $\Omega(\Mbb_1)=\Omega(\Mbb)$. Suppose $\Mbb$ is generated by $\Omega(\Mbb)$. Then $\Mbb$ is generated by $\Omega(\Mbb_1)$. So $\Mbb\subset\Mbb_1$. Hence $\Mbb$ is irreducible.

(2): Note that by part (1), $\fk U(\Vbb)\mbf U$ is irreducible. We shall show $\fk U(\Vbb)\Omega(\Mbb)=\Mbb$ and $\Omega(\fk U(\Vbb)\mbf U)=\mbf U$.  Since $\Omega(\Mbb)$ is a nontrivial subspace of $\Mbb$, we must have $\fk U(\Vbb)\Omega(\Mbb)\subset\Mbb$ and hence, by the irreducibility of $\Mbb$, that $\fk U(\Vbb)\Omega(\Mbb)=\Mbb$. Since $\fk U(\Vbb)\mbf U$ is clearly nontrivial, $\Omega(\fk U(\Vbb)\mbf U)$ is a nontrivial $\fk U_0(\Vbb)$-submodule of $\mbf U$. So $\Omega(\fk U(\Vbb)\mbf U)=\mbf U$.

(3) If $\Wbb$ is completely reducible, i.e., equivalent to $\bigoplus_i\Wbb_i$ where each $\Wbb_i$ is an irreducible weak $\Vbb$-module, then by \eqref{eq205}, $\Omega(\Wbb)$ is equivalent to $\bigoplus_i\Omega(\Wbb_i)$ where each $\Omega(\Wbb_i)$ is irreducible. So $\Omega(\Wbb)$ is completely reducible.

Now, assume that $\Omega(\Wbb)$ is completely reducible. Since  $\Omega(\Wbb)$ is finite dimensional, it is semisimple, i.e., $\Omega(\Wbb)=\bigoplus_{i=1}^N\mbf U_i$ where $N\in\Zbb_+$, and each $\mbf U_i$ is an irreducible $\fk U_0(\Vbb)$-submodule of $\Omega(\Wbb)$. Since $\Wbb$ is generated by $\Omega(\Wbb)$, it is clear that $\Wbb=\bigvee_{i=1}^n \fk U(\Vbb)(\mbf U_i)$. Thus, $\Wbb$ is a sum of irreducible modules. By basic representation theory (see for instance \cite{Lang} section XVII.2.), $\Wbb$ is completely reducible.
\end{proof}

\begin{df}
Let $\mbf U$ be a finite dimensional representation of $\fk U_0(\Vbb)$. We say that $\mbf U$ is \textbf{$\Vbb$-admissible} if there exists a weak $\Vbb$-module such that the $\fk U_0(\Vbb)$-module $\Omega(\Wbb)$ is equivalent to $\mbf U$. By restricting $\Wbb$ to $\fk U(\Vbb)\Omega(\Wbb)$, we may assume that $\Wbb$ is a lowest weight weak $\Vbb$-module.
\end{df}

By the previous results, we clearly have:

\begin{pp}\label{lb132}
The following are equivalent.
\begin{enumerate}[label=(\alph*)]
\item Every lowest weight admissible $\Vbb$-module $\Wbb$  with finite dimensional $\Omega(\Wbb)$ is completely reducible.
\item Every lowest weight weak  $\Vbb$-module $\Wbb$  with finite dimensional $\Omega(\Wbb)$ is  a finite direct sum of irreducible admissible $\Vbb$-modules.
\item Every finite-dimensional $\Vbb$-admissible $\fk U_0(\Vbb)$-module is semi-simple.
\end{enumerate}
\end{pp}

\begin{df}\label{lb133}
If $\Vbb$ satisfies one of the three conditions of Proposition \ref{lb132}, we say that $\Vbb$ is \textbf{rational}.
\end{df}

Note that by Corollary \ref{lb99}, if $\Vbb$ is $C_2$-cofinite and rational, each lowest weight weak  $\Vbb$-module $\Wbb$  with finite dimensional $\Omega(\Wbb)$ is a finite direct sum of irreducible (ordinary) $\Vbb$-modules.

\begin{rem}
Our definition of rationality is weaker than the usual one, which says any admissible $\Vbb$-module is completely reducible. Assuming $\Vbb$ is $C_2$-cofinite, then the two notations are equivalent. Indeed, our rationality is equivalent to the semisimplicity of the Zhu's algebra $A(\Vbb)$ of $\Vbb$. The latter is equivalent to the usual rationality due to \cite{McR21}.
\end{rem}

Suppose now that $(\Wbb,Y_+,Y_-)$ is a weak $\Vbb\times\Vbb$-module. We let $\Omega_+(\Mbb)$ (resp. $\Omega_-(\Mbb)$) be the lowest weight subspace of $(\Wbb,Y_+)$ (resp. $(\Wbb,Y_-)$). Set \index{zz@$\Omega(\Wbb),\Omega_{+-}(\Wbb)$}
\begin{align}
\Omega_{+-}(\Wbb)=\Omega_+(\Wbb)\cap \Omega_-(\Wbb).\label{eq207}
\end{align}
Then, it is clear that $\Omega_{+-}(\Wbb)$ is $\fk U_0(\Vbb)\times \fk U_0(\Vbb)$-invariant.

\subsection*{Some results for associative algebras}

Let $A$ be an associative algebra and $\mbf U$ be a representation of $A$. It is clear that if $\mbf U_1,\mbf U_2$ are inequivalent irreducible representations of $A$, then $\Hom_A(\mbf U_1,\mbf U_2)=0$. Indeed, choose any $T\in\Hom_A(\mbf U_1,\mbf U_2)=0$. If $T\neq 0$, then $\Ker(T)$ is a  $A$-submodule of $\mbf U_1$ not equal to $\mbf U_1$. So $\Ker(T)=0$. Also, since the range of $T$ is a non-trivial $A$-submodule of $\mbf U_2$, $T$ must be surjective. This is impossible.

We say that $\mbf U$ is \textbf{strongly irreducible}, if $\mbf U$ is irreducible, and $\End_A(\mbf U)=\Cbb\id_{\mbf U}$. We say that $\mbf U$ is \textbf{strongly and completely reducible} if $\mbf U\simeq\bigoplus_i\mbf U_i$ where each $\mbf U_i$ is strongly irreducible. For instance, this is so if $\mbf U$ is a direct sum of irreducible finite dimensional representations. 

We have seen in Proposition \ref{lb131} that any irreducible (ordinary) $\Vbb$-module is strongly irreducible. More generally, we have:

\begin{thm}\label{lb108}
Let $\Wbb$ be an irreducible admissible $\Vbb$-module with finite-dimensional $\Omega(\Wbb)$. Then $\Wbb$ is strongly irreducible (as an $\fk U(\Vbb)$-module), and $L_0$ differs $\wtd L_0$ by a constant. Moreover, $\Omega(\Wbb)$ equals the lowest non-trivial eigenspaces of both $L_0$ and $\wtd L_0$, and is a (strongly) irreducible $\fk U_0(\Vbb)$-module.
\end{thm}

Note that by Corollary \ref{lb130}, if $\Wbb$ is an irreducible (ordinary) $\Vbb$-module, then it automatically has finite-dimensional $\Omega(\Wbb)$.

\begin{proof}
Choose any $T\in\End_\Vbb(\Wbb)$. Since $[T,Y_\Wbb(v)_n]=0$ for any homogeneous $v\in\Vbb$ and $n\geq\wt(v)$, we have $T\Omega(\Wbb)\subset\Omega(\Wbb)$. Since $\Omega(\Wbb)$ is non-trivial (since it contains the lowest $\wtd L_0$-weight space) and finite-dimensional, $T|\Wbb_{(s)}$ has an eigenvalue $\lambda\in\Cbb$. It follows that $\Ker(T-\lambda\id_\Wbb)$ is a non-trivial weak $\Vbb$-submodule of $\Wbb$, which must be $\Wbb$. So $T=\lambda\id_\Wbb$. In particular, $\wtd L_0-L_0$ is a constant. The rest of the statements follows from Proposition \ref{lb98}.
\end{proof}

\begin{co}\label{lb99}
Assume that $\Vbb$ is $C_2$-cofinite. Let $\Wbb$ be an irreducible  admissible $\Vbb$-module. Then $\Wbb$ is an irreducible (ordinary) $\Vbb$-module.
\end{co}

\begin{proof}
By Theorem \ref{lb122} and the description \eqref{eq156}, it is clear that each $\wtd L_0$-weight space is finite-dimensional. In particular, this is true for the lowest non-trivial $\wtd L_0$-weight space, which by Proposition \ref{lb98} is $\Omega(\Wbb)$. So $\Omega(\Wbb)$ is finite-dimensional. Thus, by Theorem \ref{lb108}, $L_0$ is diagonalizable with finite dimensional eigenspaces, and the eigenvalues are in $\lambda+\Nbb$ for some $\lambda\in\Cbb$.
\end{proof}

\begin{pp}\label{lb115}
Let $\mbf U$ be a strongly irreducible $A$-module. Let $V,W$ be  vector spaces. Consider the $A$-modules $\mbf U\otimes V,\mbf U\otimes W$ where the actions of $A$ are on the $\mbf U$-component. Define a linear  map
\begin{align*}
\Phi:\Hom(V,W)\rightarrow \Hom_A(\mbf U\otimes V,\mbf U\otimes W),\qquad T\mapsto \id_{\mbf U}\otimes T
\end{align*}
Then $\Phi$ is an isomorphism of vector spaces.
\end{pp}
\begin{proof}	$\Phi$ is clearly injective. Let us prove that $\Phi$ is surjective. Choose any $S\in\Hom_A(\mbf U\otimes V,\mbf U\otimes W)$. We shall show that for each $v\in V$, there exists a (necessarily unique) $w\in W$ such that $S(u\otimes v)=u\otimes w$ for each $u\in U$. Then $S=\id_U\otimes T$ where $T$ sends each $v$ to $w$. By the uniqueness of $w$ with respect to $v$, the map $T$ is linear.

Let us fix any $v\in V$, and let $\id\otimes v$ denote the homomorphism $\mbf U\rightarrow\mbf U\otimes V$ sending each $u_0\in\mbf U$ to $u_0\otimes v$. For each $w'\in W^*$, let $\id\otimes w'$ denote the homomorphism from $\mbf U\otimes V$ to $\mbf U$ sending each $u_0\otimes v_0$ to $v'(v_0)\cdot u_0$. Then $(\id\otimes w')S(\id\otimes v)$ is an endomorphism of $\mbf U$, which is of the form $\lambda_{w'}\id_{\mbf U}$ for some $\lambda_{w'}\in \Cbb$. 
	
Choose a basis $\{e_i\}$ of $\mbf U_i$. Fix a basis element $e_i$. Then we can find a set of vectors $\{w_j\}$ in $W$ such that $S(e_i\otimes v)=\sum_j e_j\otimes w_j$. Choose any $w'\in W^*$. Then
	\begin{align*}
	\lambda_{w'}e_i=(\id\otimes w')S(\id\otimes v)e_i=(\id\otimes w')S(e_i\otimes v)=\sum_j w'(w_j)\cdot e_j.
	\end{align*}
	Thus, whenever $j\neq i$, we have $w'(w_j)=0$ for any $w'$, and hence $w_j=0$. Therefore $S(e_i\otimes v)=e_i\otimes w_i$, and hence $S(u\otimes v)=u\otimes w_i$ for each $u\in\mbf U$. 
\end{proof}

Set $V=\Cbb$. We obtain:

\begin{co}\label{lb107}
Let $\mbf U$ be a strongly irreducible $A$-module. Let $V$ be a vector space. Consider the $A$-module $\mbf U\otimes W$ where the action of $A$ is on the $\mbf U$-component. Define a linear  map
\begin{align*}
\Phi:W\rightarrow \Hom_A(\mbf U,\mbf U\otimes W),\qquad w\mapsto \Phi(w)
\end{align*}
where $\Phi(w)(u)=u\otimes w$ for each $u\in\mbf U$. Then $\Phi$ is an isomorphism of vector spaces.
\end{co}

Suppose that $A,B$ are associative algebras, and $\mbf W$ is both an $A$-module and a $B$-module. We say that $\Wbb$ is an $A\times B$-module if the actions of $A$ and $B$ on $\mbf W$ commute.

\begin{pp}\label{lb113}
If $\mbf U$ is a strongly irreducible $A$-module, and $\mbf V$ is an irreducible $B$-module, then $\mbf U\otimes \mbf V$ is an irreducible $A\times B$-module. If $\mbf V$ is moreover strongly irreducible, then so is the $A\times B$-module $\mbf U\otimes \mbf V$.
\end{pp}

\begin{proof}
Choose any non-zero vector $w\in\mbf U\otimes\mbf V$. Since $\mbf U\otimes\mbf V$ is clearly completely reducible as an $A$-module, by basic representation theory, the submodule $(A\otimes 1)w$ is also completely reducible,\footnote{A submodule of a completely reducible module is completely reducible; see \cite{Lang} section XVII.2.} hence contains an irreducible submodule, which must be equivalent to $\mbf U$. By Corollary \ref{lb107}, this submodule must be of the form $\mbf U\otimes v$ for a non-zero vector $v\in\mbf V$. Thus, as $(A\otimes B)(\mbf U\otimes v)=\mbf U\otimes \mbf V$, we have $(A\otimes B)w=\mbf U\otimes \mbf V$, which proves that $\mbf U\otimes\mbf V$ is irreducible.

Choose any $S\in\End_{A\times B}(\mbf U\otimes\mbf V)$. Then $S$ commutes with the actions of $A$. Thus, by Proposition \ref{lb115}, $S=\id_{\mbf U}\otimes T$ where $T\in\End(\mbf V)$. Since $S$ commutes with the actions of $B$, so does $T$. Thus $T\in\End_B(\mbf V)$, which must be a scalar multiplication if $\mbf V$ is strongly irreducible.
\end{proof}

\begin{pp}\label{lb101}
Let $\mbf W$ be a   representation of $A\times B$. Suppose that $\mbf W$ is strongly and completely reducible as an $A$-module, and (resp. strongly and) completely reducible  as a $B$-module. Then, there exist  strongly irreducible $A$-modules $\{\mbf U_i\}_{i\in\mc I}$ and (resp. strongly) irreducible $B$-modules $\{\mbf V_i\}_{i\in\mc I}$, such that
\begin{align*}
\mbf W\simeq\bigoplus_{i\in\mc I} \mbf U_i\otimes\mbf V_i.
\end{align*}
\end{pp}

\begin{proof}
Since $A\curvearrowright\mbf W$ is strongly and completely reducible, by Corollary \ref{lb107}, we can find a collection of mutually inequivalent strongly irreducible $A$-modules $\{\mbf U_i\}_i$ such that the equivalence
\begin{align}
\mbf W\simeq\bigoplus_i\mbf U_i\otimes \Hom_A(\mbf U_i,\mbf W)\label{eq206}
\end{align}
holds for $A$-modules. Here, $A$ is acting on $\mbf U_i\otimes \Hom_A(\mbf U_i,\mbf W)$ by acting on the $\mbf U_i$-component. Moreover, the natural embedding $\mbf U_i\otimes \Hom_A(\mbf U_i,\mbf W)\hookrightarrow\mbf W$ is given by $u\otimes T\mapsto Tu$ if $u\in\mbf U_i,T\in \Hom_A(\mbf U_i,\mbf W)$. 
	
Since the actions of $A$ and $B$ commute, each $b\in B$ can be viewed as an element of $\End_A(\mbf W)$. Thus, $\Hom_A(\mbf U_i,\mbf W)$ is naturally a $B$-module where each $b\in B$ acts on $T\in \Hom_A(\mbf U_i,\mbf W)$ as $bT$. It is easy to see that \eqref{eq206} is an equivalence of $A\times B$-module. 

By the irreducibility of $\mbf U_i$, for each nonzero $u\in\mbf U_i$, the map $\Hom_A(\mbf U_i,\mbf W)\rightarrow\mbf W,T\mapsto Tu$ is an injective homomorphism of $B$-modules. It follows that each $\Hom_A(\mbf U_i,\mbf W)$ is  equivalent to a $B$-submodule of $\mbf W$. Thus, it is a direct sum of  irreducible $B$-modules. We can thus write $\Hom_A(\mbf U_i,\mbf W)$ as a  direct sum $\bigoplus_j\mbf V_{i,j}$ of  irreducible $B$-modules. It follows that $\mbf W\simeq\bigoplus_{i,j}\mbf U_i\otimes\mbf V_{i,j}$. It is easy to see that any irreducible submodule of a direct sum of strongly irreducible modules is strongly irreducible (since it is isomorphic to one of the strongly irreducible summand). So $\mbf V_{i,j}$ is strongly irreducible if $\mbf W$ is strongly and completely reducible as a $B$-module.
\end{proof}

\section{Dual tensor products}\label{lb110}

Let
\begin{align}
\fk X=(C;x_1,\dots,x_N;x';x'';\eta_1,\dots,\eta_N;\xi;\varpi)\label{eq196}
\end{align}
be an $(N+2)$-pointed compact Riemann surface with local coordinates. Throughout this section, we fix mutually disjoint neighborhoods $W_1,\dots,W_N,W',W''$ of $x_1,\dots,x_N,x',x''$ respectively, on which the local coordinates are defined. We assume that each connected component of $C$ contains at least one of $x_1,\dots,x_N$, and call such $\fk X$ an \textbf{$N$-pointed compact Riemann surface with local coordiates and $2$ outputs}. We let \index{SX@$\SX$, $\SX(b)$} \index{DX@$\DX$}
\begin{gather*}
\SX=x_1+\cdots+x_N,\qquad \DX=x'+x''.
\end{gather*}
For each $a,b\in\Nbb$, define \index{VX@$\scr V_{\fk X,a,b}^{\leq n},\scr V_{\fk X,a,b}$}
\begin{gather*}
\scr V_{\fk X,a,b}^{\leq n}\equiv \scr V_{\fk X}^{\leq n}(-(L_0\DX+ax'+bx''))\qquad(\forall n\in\Nbb),\\
\scr V_{\fk X,a,b}=\varinjlim_{n\in\Nbb}\scr V_{\fk X,a,b}^{\leq n}.
\end{gather*}
Here, $\scr V_{\fk X,a,b}^{\leq n}$ is a locally free $\scr O_C$-submodule of $\scr V_{\fk X}^{\leq n}$ described as follows: Outside $x'$ and $x''$, $\scr V_{\fk X,a,b}^{\leq n}$ equals $\scr V_{\fk X}^{\leq n}$; $\scr V_{\fk X,a,b}^{\leq n}|_{W'}$ and $\scr V_{\fk X,a,b}^{\leq n}|_{W''}$ are generated by 
\begin{gather*}
\mc U_\varrho(\xi)^{-1}\xi^{a+L_0}v\qquad \text{resp.} \qquad \mc U_\varrho(\varpi)^{-1}\varpi^{b+L_0}v
\end{gather*}
where $v$ is any homogeneous vector of $\Vbb^{\leq n}$. It is easy to check that this definition is independent of the local coordinates $\xi,\varpi$.

Let $\Vbb$ be a  VOA. We do not assume $\Vbb$ to be $C_2$-cofinite. Let $\Wbb_1,\dots,\Wbb_N$ be $\Vbb$-modules associated to $x_1,\dots,x_N$. We define a \textbf{truncated  $\fk X$-tensor product} of $\Wbb_1,\dots\Wbb_N$  to be the vector space \index{T@$\scr T_{\fk X,a,b}(\Wbb_\blt),\scr T_{\fk X,a,b}^*(\Wbb_\blt)$}
\begin{align}
\scr T_{\fk X,a,b}(\Wbb_\blt)=\frac{\Wbb_\blt}{H^0(C,\scr V_{\fk X,a,b}\otimes\omega_C(\blt\SX))\cdot\Wbb_\blt}.
\end{align}
Its dual space is denoted by $\scr T_{\fk X,a,b}^*(\Wbb_\blt)$ and called a \textbf{truncated dual $\fk X$-tensor product}. Note that when $a'>a$ and $b'>b$, we have a natural injective linear map $\scr T_{\fk X,a,b}^*\hookrightarrow \scr T_{\fk X,a',b'}^*$. We can thus define \index{W@$\boxbackslash_{\fk X}(\Wbb_\blt),\boxbackslash_{\fk X}^\low(\Wbb_\blt)$}
\begin{gather*}
\boxbackslash_{\fk X}(\Wbb_\blt)=\varinjlim_{a,b\in\Nbb}\scr T_{\fk X,a,b}^*(\Wbb_\blt),
\end{gather*}
called the \textbf{dual $\fk X$-tensor product} of $\Wbb_1,\dots,\Wbb_N$.

In a similar way, one can define $\scr V_{\fk X,a,b}^{\leq n},\scr V_{\fk X,a,b}$, and (the sheaves of) (truncated) (dual) $\fk X$-tensor products  when $\fk X$ is a family of $N$-pointed compact Riemann surfaces with local coordinates and $M$-outputs. In the case that $M=0$, we obtain the spaces/sheaves of covacua and conformal blocks. All the results in chapters \ref{lb90} and \ref{lb91} can be generalized to these sheaves/vector spaces using almost the same idea. For instance, notice that in the setting \eqref{eq196}, we have
\begin{align}
\scr V_{\fk X,a,b}^{\leq n}/\scr V_{\fk X,a,b}^{\leq n-1}\simeq\Vbb(n)\otimes_\Cbb\Theta_C^{\otimes n}(-n\DX-ax'-bx'').
\end{align}
Thus, we have the vanishing Theorem \ref{lb21} with $\scr V_C^{\leq n}$ replaced by $\scr V_{\fk X,a,b}^{\leq n}$. (Of course, the integer $k_0$ in that theorem should now also depend on $a,b$.) In the following, we will directly claim and use the generalizations of those results without  proving them again.

\begin{eg}\label{lb112}
Let $\Mbb_1,\Mbb_2$ be irreducible $\Vbb$-modules. By Convention \ref{lb36},  their lowest $\wtd L_0$-weights  (with non-trivial weight subspaces $\Mbb_1(0),\Mbb_2(0)$) are both $0$. Choose any $\upphi\in\scr T_{\fk X}^*(\Wbb_\blt\otimes\Mbb_1\otimes\Mbb_2)$. Then there is a natural linear map
\begin{gather*}
\Psi_\upphi:\Mbb_1^{\leq a}\otimes\Mbb_2^{\leq b}\rightarrow\scr T_{\fk X,a,b}^*(\Wbb_\blt)
\end{gather*}
defined such that for any $w_\blt\in\Wbb_\blt,m_1\in\Mbb_1^{\leq a},m_2\in\Mbb_2^{\leq b}$,
\begin{gather}
\Psi_\upphi(m_1\otimes m_2)(w_\blt)=\upphi(w_\blt\otimes m_1\otimes m_2).\label{eq216}
\end{gather}
Here, $\Mbb_1^{\leq a}$ is the subspace spanned by all $\wtd L_0$-homogeneous vectors with weight $\leq a$, and $\Mbb_2^{\leq b}$ is understood in a similar way. By taking the direct limit over $(a,b)$, we obtain
\begin{gather*}
\Psi_\upphi:\Mbb_1\otimes\Mbb_2\rightarrow \boxbackslash_{\fk X}(\Wbb_\blt).
\end{gather*}
Our next goal is to define a weak $\Vbb\times\Vbb$-module structure on $\boxbackslash_{\fk X}(\Wbb_\blt)$ such that $\Psi_\upphi$ is a homomorphism. This will imply, by the irreducibility of $\Mbb_1,\Mbb_2$, that $\Psi_\upphi$ is injective when $\upphi\neq 0$.
\end{eg}

\subsection*{Actions of $\Vbb$}

Similar to the treatment for conformal blocks, one can prove that the formal sewing of a dual tensor product element is a formal dual tensor product element. Moreover, if one is sewing $\fk X$ with  $(\Pbb^1;0,1,\infty)$ (which has $0$ outputs), the a.l.u. of sewing can be proved with the help of the strong residue theorem. Thus, we are able to prove the \textbf{propagation of dual tensor product elements}.  For each $a,b,n\in\Nbb,\upphi\in\scr T_{\fk X,a,b}^*(\Wbb_\blt),w_\blt\in\Wbb_\blt$, we have a homomorphism of $\scr O_C^{\boxtimes n}$-modules
\begin{gather*}
\wr^n\upphi(w_\blt):\scr V_{\fk X,a,b}^{\boxtimes n}\rightarrow\scr O_{\Conf_n(C-\SX)},
\end{gather*}
i.e., for each open subsets $U_1,\dots,U_n\subset C$, we have a homomorphism of $\scr O(U_1)\otimes\cdots\otimes\scr O(U_n)$-modules
\begin{gather*}
\wr^n\upphi(w_\blt):\scr V_{\fk X,a,b}(U_1)\otimes_\Cbb\cdots\otimes_\Cbb\scr V_{\fk X,a,b}(U_n)\rightarrow\scr O_{\Conf_n(C-\SX)}(U_1\times\cdots\times U_n),
\end{gather*}
(Recall that $\scr O_{\Conf_n(C-\SX)}(U_1\times\cdots\times U_n)=\scr O(\Conf(U_1-\SX,\dots,U_n-\SX))$.) and these maps are compatible under restrictions to open subsets. We have $\wr^0\upphi=\upphi$. Moreover, for the $\eta_1,\dots,\eta_N,x_1,\dots,x_N$ chosen at the beginning of this section, Theorem \ref{lb92} holds verbatim if we replace $\scr V_C$ with $\scr V_{\fk X,a,b}$. Since $\scr V_C$ equals $\scr V_{\fk X,a,b}$ outside the output points $x',x''$, Theorem \ref{lb92} also holds if (following the notations of that theorem) for each $1\leq k\leq n$, we still choose $v_k\in\scr V_k(U_k)$, but assume in addition that $x',x''\notin U_k$. In particular, $y_1,\dots,y_n$ cannot be $x',x''$. We will use this theorem only for $n=1,2$.

\begin{rem}
Write $C-\{x',x''\}$ as $C-\DX$ for simplicity. Note that $(\scr V_{\fk X,a,b}|_{C-\DX})^{\boxtimes n}$ equals $\scr V_{C-\DX}^{\boxtimes n}$. The restricted homomorphism
\begin{gather}
\wr^n\upphi(w_\blt):\scr V_{C-\DX}^{\boxtimes n}\rightarrow\scr O_{\Conf_n(C-\SX-\DX)}\label{eq197}
\end{gather} 
is independent of the numbers $a,b$ making $\upphi$ belonging to $\scr T_{\fk X,a,b}^*(\Wbb_\blt)$. Indeed, the case $n=0$ is obvious. Suppose the case for $n-1$ is true.  By  Theorem \ref{lb92}-(1), $\wr^n\upphi(v_1,\dots,v_n,w_\blt)$ is independent of $a,b$ when $v_1$ is a section of $\scr V_{C-\DX}$ defined near $x_1,\dots,x_N$.  Thus, by the argument in the proof of Proposition \ref{lb41}, the independence of $a,b$ is true for any $v_1$. To summarize, we have a well defined $\wr^n\upphi$ in \eqref{eq197} for any $\upphi\in\boxbackslash_{\fk X}(\Wbb_\blt)$. 
\end{rem}

We also regard
\begin{gather*}
\wr^n\upphi(w_\blt):\scr V_C^{\boxtimes n}\rightarrow\scr O_{\Conf_n(C-\SX-\DX)}
\end{gather*}
sending each $v_1\in \scr V_C(U_1),\dots,v_n\in\scr V_C(U_n)$ to $\wr^n\upphi(v_1|_{U_1-\DX},\dots,v_n|_{U_n-\DX},w_\blt)$. In particular, for each homogeneous $v\in\Vbb$, considered as a constant section of $\Vbb\otimes_\Cbb\scr O(W')$, we have $\wr\upphi(\mc U_\varrho(\xi)^{-1} v,w_\blt)\in\scr O(W'-\{x'\})$. Moreover, choose $a,b$ such that $\upphi\in\scr T_{\fk X,a,b}^*(\Wbb_\blt)$. Then $\xi^{\wt(v)+a}v$ is an element of $\scr V_{\fk X,a,b}(W')$. So 
\begin{gather}
\wr\upphi(\mc U_\varrho(\xi)^{-1} v,w_\blt)\in \scr O_C((\wt(v)+a)x')(W').\label{eq198}
\end{gather}

Now, for each $v\in\Vbb,n\in\Zbb$, we define
\begin{gather*}
Y_+(v)_n:\boxbackslash_{\fk X}(\Wbb_\blt)\rightarrow\Wbb_\blt^*
\end{gather*}
as follows. Identify $W'$ with $\xi(W')$ via $\xi$. (So $\xi$ is identified with the standard coordinate $z$.) Identify  $\scr V_{W'}\simeq\Vbb\otimes_\Cbb\scr O_{W'}$ via $\mc U_\varrho(\xi)$.  Then, for any $\upphi\in\boxbackslash_{\fk X}(\Wbb_\blt)$, the evaluation of $Y_+(v)_n\upphi$ with any $w_\blt\in\Wbb_\blt$ is
\begin{gather}
\boxed{~~Y_+(v)_n\upphi(w_\blt)=\Res_{z=0}\wr\upphi( v,w_\blt)z^ndz~~}\label{eq199}
\end{gather}
Then, by \eqref{eq198}, we have
\begin{align}
Y_+(v)_n\upphi=0\qquad (\text{if }n\geq \wt(v)+a).\label{eq200}
\end{align}
Set $Y_+(v,z)=\sum_{n\in\Zbb}Y_+(v)_nz^{-n-1}$.  Then the lower truncation property is satisfied: $Y_+(v,z)\upphi\in \Wbb_\blt^*[[z]]$.  Thus, $(\boxbackslash_{\fk X}(\Wbb_\blt),Y_+)$ becomes a linear representation of $\Vbb$ if we can show that $Y_+(v)_n\upphi\in \boxbackslash_{\fk X}(\Wbb_\blt)$.

\begin{lm}\label{lb93}
For any homogeneous vector $v\in\Vbb$, $n\in\Zbb$, and $\upphi\in\scr T_{\fk X,a,b}^*(\Wbb_\blt)$, we have $Y_+(v)_n\upphi\in \scr T_{\fk X,a',b}^*(\Wbb_\blt)$ where $a'=a+\max\{0,\wt(v)-n-1\}$.
\end{lm}

\begin{proof}
Let $\zeta$ be another standard coordinate of $\Cbb$. So both $z$ and $\zeta$ are identified with $\xi$ as coordinates. (But they are independent as variables.) In the following, for a two-variable holomorphic function, we will let $\zeta$ (resp. $z$) denote the first (resp. second) complex variable.

Identify $W'\simeq\xi(W')$ via $\xi$ and  $\scr V_{W'}\simeq\Vbb\otimes_\Cbb\scr O_{W'}$ via $\mc U_\varrho(\xi)$ as above. One can define $\uppsi\in H^0(C-\SX-\DX,\scr V_C^*)$ such that for any section $u$ of $\scr V_C$ defined in an open subset $W$ of $C-\SX-\DX$,
\begin{align*}
\uppsi(u)=\Res_{z=0}\wr\wr\upphi(u, v,w_\blt)z^ndz,
\end{align*}
or more precisely, if we also identify $W$ with an open subset of $\Cbb$ so that $\zeta$ is a complex variable on $W$, then $\uppsi(u)(\zeta)=\Res_{z=0}\wr\wr\upphi(u, v,w_\blt)(\zeta,z)z^ndz$. If $u$ is defined near $x''$, then $\wr\wr\upphi(\xi^{b+L_0}u, \varpi^{a+L_0}v,w_\blt)$ is holomorphic (with no poles) near $\zeta=x'',z=x'$. Thus, $\uppsi(\xi^{b+L_0}u)$ has no  pole at $x''$. So $\uppsi\in H^0(C-\SX-x',\scr V_{\fk X,a',b}^*)$. We shall show $\uppsi\in H^0(C-\SX,\scr V_{\fk X,a',b}^*)$. Suppose this can be proved. By Theorem \ref{lb92}-(1), if we identify $W_i\simeq \eta_i(W_i)$ via $\eta_i$ and identify $\scr V_{W_i}\simeq\Vbb\otimes_\Cbb\scr O_{W_i}$ via $\mc U_\varrho(\eta_i)$, then for any section $u$ of $\scr V_{\fk X,a',b}(W_i)=\scr V_C(W_i)$  (which restricts to a section on $W_i-\{x_i\}$), we have
\begin{align*}
&\uppsi(u)(\zeta)=\Res_{z=0}\wr\upphi(v,w_1\otimes\cdots\otimes Y_{\Wbb_i}(u,\zeta)w_i\otimes\cdots\otimes w_N)z^ndz\\
=&Y_+(v)_n\upphi(w_1\otimes\cdots\otimes Y_{\Wbb_i}(u,\zeta)w_i\otimes\cdots\otimes w_N).
\end{align*}
So $\uppsi$ restricts to $(Y_+(v)_n\upphi)_{x_i}$ (defined similarly as in \eqref{eq125}) near each $x_i$. Thus, as in the proof of Theorem \ref{lb32}, the Residue theorem implies that $Y_+(v)_n\upphi\in \scr T_{\fk X,a',b}^*(\Wbb_\blt)$.

Choose any homogeneous vector $u\in\Vbb$ with weight $\wt(u)$, considered as a constant section of $\scr V_C(W')$.  Consider $\uppsi(u)$ as a holomorphic function with variable $\zeta$. We shall show that $\uppsi(u)\in\scr O_{W'}((\wt(v)+a')x')(W')$. Set $f=f(\zeta,z)$ to be
\begin{align*}
f=\wr\wr\upphi(u,v,w_\blt)\qquad\in\scr O(\Conf_2(W'-\{x'\})).
\end{align*}
Then, as $\zeta^{a+\wt(u)}u,z^{a+\wt(b)}\in\scr V_{\fk X,a,b}(W')$, we have
\begin{align*}
\zeta^{a+\wt(u)}z^{a+\wt(v)}f(\zeta,z)\in\scr O(\Conf_2(W')).
\end{align*}
Choose circles $C_1,C_2,C_3$ in $W'$ surrounding $x'$ with radii $r_1<r_2<r_3$. For each $z\in C_2$, choose a circle $C(z)$ with center $z$ and radius less than $r_2-r_1$ and $r_3-r_2$. Let $m\in\Zbb$.  Then
\begin{align*}
&\Res_{\zeta=0}~\zeta^m\uppsi(u)d\zeta=\oint_{C_3}\oint_{C_2}\zeta^mz^nfdzd\zeta=\oint_{C_2}\oint_{C_3}\zeta^mz^nfd\zeta dz\\
=&\oint_{C_2}\oint_{C_1}\zeta^mz^nfd\zeta dz+\oint_{C_2}\oint_{C(z)}\zeta^mz^nfd\zeta dz.
\end{align*}
When $z\in C_2$, since $\zeta^{a+\wt(u)}f$ has no pole at $\zeta=0$, we have $\oint_{C_1}\zeta^mz^nfd\zeta=0$ whenever $m\geq a+\wt(u)$. 

Apply Theorem \ref{lb92}-(2) (by choosing $U_1=U_2$ to be $W'-\{x'\}$), we have
\begin{align*}
&\oint_{C_2}\oint_{C(z)}\zeta^mz^nfd\zeta dz=\oint_{C_2}\oint_{C(z)}\zeta^mz^n\wr\wr\upphi(u,v,w_\blt)d\zeta dz\\
=&\oint_{C_2}\oint_{C(z)}\zeta^mz^n\wr\upphi(Y(u,\zeta-z)v,w_\blt)d\zeta dz\\
=&\sum_{l\in\Nbb}{m\choose  l}\oint_{C_2}\oint_{C(z)}(\zeta-z)^lz^{m+n-l}\wr\upphi(Y(u,\zeta-z)v,w_\blt)d\zeta dz\\
=&\sum_{l\in\Nbb}{m\choose  l}\oint_{C_2}z^{m+n-l}\wr\upphi(Y(u)_lv,w_\blt)dz\\
=&\sum_{l\in\Nbb}{m\choose  l}Y_+(Y(u)_lv)_{m+n-l}\upphi(w_\blt),
\end{align*}
where we have used \eqref{eq199} in the last step. By \eqref{eq200}, the above expression equals $0$ when
\begin{align*}
m+n-l\geq\wt(Y(u)_lv)+a=\wt(u)+\wt(v)-l-1+a,
\end{align*}
and hence when
\begin{align*}
m\geq \wt(u)+\wt(v)+a-1-n.
\end{align*}
Thus, we conclude that $\Res_{\zeta=0}~\zeta^m\uppsi(u)d\zeta$ equals $0$ when $m\geq a+\wt(u)+\max\{0,\wt(v)-n-1\}$, i.e., when $m\geq \wt(u)+a'$. This finishes the proof.
\end{proof}

\subsection*{The weak $\Vbb$-module $\boxbackslash_{\fk X}(\Wbb_\blt)$}

We now show that the linear representation  $(\boxbackslash_{\fk X}(\Wbb_\blt),Y_+)$ is indeed a weak $\Vbb$-module.

\begin{lm}\label{lb95}
Choose any $\upphi\in\boxbackslash_{\fk X}(\Wbb_\blt)$, $m,n\in\Zbb$. Identify $W'\simeq\xi(W')$ via $\xi$. Choose $u,v\in\Vbb$, considered as constant sections of $\scr V_C(W')$ defined by $\mc U_\varrho(\xi)$. Then for any $w_\blt\in\Wbb_\blt$,
\begin{align}
Y_+(u)_m Y_+(v)_n\upphi(w_\blt)=\Res_{\zeta=0}\Res_{z=0}\wr\wr\upphi(u,v,w_\blt)\zeta^m z^n dzd\zeta.\label{eq201}
\end{align}
\end{lm}

As previously, $\wr\wr\upphi(u,v,w_\blt)$ is short for $\wr\wr\upphi(u,v,w_\blt)(\zeta,z)$ where $\zeta,z$ are both standard complex variables of $\Cbb$.

\begin{proof}
By \eqref{eq199}, we have
\begin{align*}
Y_+(u)_m Y_+(v)_n\upphi(w_\blt)=\Res_{\zeta=0}\wr(Y_+(v)_n\upphi)(u,w_\blt)\zeta^md\zeta.
\end{align*}
Thus, \eqref{eq201} will follow if we can show
\begin{align}
\wr(Y_+(v)_n\upphi)(u,w_\blt)=\Res_{z=0}\wr\wr\upphi(u,v,w_\blt) z^n dz\label{eq202}
\end{align}
for any section $u$ of $\scr V_{C-\SX}$. If $u$ is defined on an open subset of $W_i$, then, by Theorem \ref{lb92}, under the identification $\scr V_{W_i}\simeq\Vbb\otimes_\Cbb\scr O_{W_i}$ defined by $\mc U_\varrho(\eta_i)$, we have
\begin{align*}
&\wr(Y_+(v)_n\upphi)(u,w_\blt)=(Y_+(v)_n\upphi)(w_1\otimes\cdots\otimes Y(u,\eta_i)w_i,\otimes\cdots\otimes w_N)\\
=&\Res_{z=0}\wr\upphi(v,w_1\otimes\cdots\otimes Y(u,\eta_i)w_i,\otimes\cdots\otimes w_N)z^ndz\\
=& \Res_{z=0}\wr\wr\upphi(u,v,w_\blt)z^ndz.
\end{align*}
Thus, \eqref{eq202} holds when $u$ is near $x_1,\dots,x_N$.  Using the argument in the proof of Proposition \ref{lb41} together with the fact that each connected component of $C$ contains at least one of $x_1,\dots,x_N$, it is easy to see that \eqref{eq202} holds for any $u$.
\end{proof}

\begin{pp}
$(\boxbackslash_{\fk X}(\Wbb_\blt),Y_+)$ is a weak $\Vbb$-module.
\end{pp}

\begin{proof}
We shall show that $(\boxbackslash_{\fk X}(\Wbb_\blt),Y_+)$ satisfies the criteria in Proposition \ref{lb94}. It is clear  that $\Wbb_\blt$ projects to a dense subspace of the dual space of $\boxbackslash_{\fk X}(\Wbb_\blt)$. Also, we have $Y_+(\id,z)=\id_{\boxbackslash_{\fk X}(\Wbb_\blt)}$ by Theorem \ref{lb92}-(3). Choose any $u,v\in\Vbb$, $\upphi\in\boxbackslash_{\fk X}(\Wbb_\blt)$, $w_\blt\in\Wbb_\blt$. $\upphi$ and $w_\blt$ play the role of $w,w'$ in Proposition \ref{lb94}. Identify $W'$ with $\xi(W')$ via $\xi$, and choose $r>0$ such that $\mc D_r\subset W'$. Identify $\scr V_{W'}$ with $\Vbb\otimes_\Cbb\scr O_{W'}$ via $\mc U_\varrho(\xi)$ as usual. Set $f=f(\zeta,z)\in\scr O(\Conf_2(\mc D_r^\times))$ to be
\begin{align*}
f(\zeta,z)=\wr\wr\upphi(u,v,w_\blt)(\zeta,z).
\end{align*}
By Lemma \ref{lb95} and Theorem \ref{lb92}-(4), we clearly have
\begin{gather*}
Y_+(v)_nY_+(u)_m\upphi(w_\blt)=\Res_{z=0}\Res_{\zeta=0}~f(\zeta,z)\zeta^m z^nd\zeta dz,\\
Y_+(u)_mY_+(v)_n\upphi(w_\blt)=\Res_{\zeta=0}\Res_{z=0}~f(\zeta,z)\zeta^mz^ndz d\zeta,
\end{gather*}
which verify \eqref{eq193}, \eqref{eq195}. With the help of Theorem \ref{lb92}, we compute
\begin{align*}
&\Res_{z=0}\Res_{\zeta-z=0}~f(\zeta,z)(\zeta-z)^mz^nd\zeta dz\\
=&\Res_{z=0}\Res_{\zeta-z=0}\wr\wr\upphi(u,v,w_\blt)(\zeta,z)\cdot (\zeta-z)^m z^nd\zeta dz\\
=&\Res_{z=0}\Res_{\zeta-z=0}\wr\upphi(Y(u,\zeta-z)v,w_\blt)(\zeta,z)\cdot(\zeta-z)^m z^nd\zeta dz\\
=&\Res_{z=0}\wr\upphi(Y(u)_mv,w_\blt)(z)\cdot z^ndz,
\end{align*}
which, by \eqref{eq199}, equals $Y(Y(u)_mv)_n\upphi(w_\blt)$. This verifies \eqref{eq194}.
\end{proof}

\subsection*{The weak $\Vbb\times\Vbb$-module $\boxbackslash_{\fk X}(\Wbb_\blt)$}

One can define similarly a weak module structure $Y_-$ on $\boxbackslash_{\fk X}(\Wbb_\blt)$ by using the sections near $x''$. To be more precise, for each $v\in\Vbb,n\in\Zbb$, we define
\begin{gather*}
Y_-(v)_n:\boxbackslash_{\fk X}(\Wbb_\blt)\rightarrow\Wbb_\blt^*
\end{gather*}
as follows. Identify $W''$ with $\varpi(W'')$ via $\varpi$. (So $\varpi$ is identified with the standard coordinate $z$.) Identify  $\scr V_{W''}\simeq\Vbb\otimes_\Cbb\scr O_{W''}$ via $\mc U_\varrho(\varpi)$.  Then, for any $\upphi\in\boxbackslash_{\fk X}(\Wbb_\blt)$, the evaluation of $Y_-(v)_n\upphi$ with any $w_\blt\in\Wbb_\blt$ is
\begin{gather}
\boxed{~~Y_-(v)_n\upphi(w_\blt)=\Res_{z=0}\wr\upphi( v,w_\blt)z^ndz~~}\label{eq213}
\end{gather}
Then $(\boxbackslash_{\fk X}(\Wbb_\blt),Y_-)$ is also a weak $\Vbb$-module.

Recall Definition \ref{lb100}.

\begin{thm}
$Y_+$ and $Y_-$ commute. So $(\boxbackslash_{\fk X}(\Wbb_\blt),Y_+,Y_-)$ is a weak $\Vbb\times\Vbb$-module.
\end{thm}

\begin{proof}
Identify $W'\simeq \xi(W')$ via $\xi$ and $W''\simeq\varpi(W'')$ via $\varpi$. Let $z,\zeta$ be the standard complex variables of $W',W''$ respectively.  Choose any $v\in\Vbb$, considered as a constant section of $\scr V_C(W')$ defined by $\mc U_\varrho(\xi)$.  In the proof of Lemma \ref{lb95}, we have shown that \eqref{eq202} is true for any section $u$ of $\scr V_{C-\SX}$. In particular, this is true if we take $u\in\Vbb$ and consider it as a constant section of $\scr V_C(W''-\{x''\})$ defined by $\mc U_\varrho(\varpi)$.  Thus, we may apply $\Res_{\zeta=0}(\cdots)\zeta^md\zeta$ to obtain
\begin{align*}
Y_-(u)_m Y_+(v)_n\upphi(w_\blt)=\Res_{\zeta=0}\Res_{z=0}\wr\wr\upphi(u,v,w_\blt)(\zeta,z)\cdot \zeta^m z^n dzd\zeta.
\end{align*}
A similar argument shows
\begin{align*}
Y_+(v)_nY_-(u)_m\upphi(w_\blt)=\Res_{z=0}\Res_{\zeta=0}\wr\wr\upphi(v,u,w_\blt)(z,\zeta)\cdot\zeta^m z^n d\zeta dz.
\end{align*}
By Theorem \ref{lb92}-(4), we have $\wr\wr\upphi(u,v,w_\blt)(\zeta,z)=\wr\wr\upphi(v,u,w_\blt)(z,\zeta)$ (when $z\in W',\zeta\in W''$). The commutativity of  $Y_-(u)_m$ and $Y_+(v)_n$ follows.
\end{proof}

To determine the lowest weight subspace of $\boxbackslash_{\fk X}(\Wbb_\blt)$, we first need a lemma.

\begin{lm}\label{lb103}
Choose any $a,b\in\Nbb$ and $\upphi\in\boxbackslash_{\fk X}(\Wbb_\blt)$. Then $\upphi\in\scr T_{\fk X,a,b}^*(\Wbb_\blt)$ if and only if $Y_+(v)_n\upphi=Y_-(u)_m \upphi=0$ whenever $u,v$ are homogeneous, $n\geq\wt(v)+a$, and $m\geq\wt(u)+b$.
\end{lm}

\begin{proof}
The ``only if" part follows from \eqref{eq200} and a similar equation for $Y_-$. We now prove the ``if" part. Suppose  $Y_+(v)_n\upphi=Y_-(u)_m \upphi=0$ for the $u,v,m,n$ described above. Choose any $w_\blt\in\Wbb_\blt$. Consider $\wr\upphi(w_\blt)\in H^0(C-\SX-\DX,\scr V_{\fk X,a,b}^*)$, whose expression near each $x_i$ is given by \eqref{eq125}. If $\nu$ is a section of $\scr V_{\fk X,a,b}$ defined near $x'$, then $\mc U_\varrho(\xi)\nu$ is an $\scr O_{W'}$-linear sum of elements of the form $v\xi^nd\xi$,   where $v$ is homogeneous and $n\geq\wt(v)+a$. By \eqref{eq199}, $\Res_{\xi=0}\wr\upphi(v,w_\blt)\xi^nd\xi=0$. So $\wr\upphi(\nu,w_\blt)$ has no pole near $\xi=0$. Thus $\wr\upphi(w_\blt)\in H^0(C-\SX-\{x''\},\scr V_{\fk X,a,b}^*)$. A similar argument shows $\wr\upphi(w_\blt)\in H^0(C-\SX,\scr V_{\fk X,a,b}^*)$. Thus, as in the proof of Theorem \ref{lb32}, we may use Residue theorem to deduce $\upphi\in\scr T_{\fk X,a,b}^*(\Wbb_\blt)$.
\end{proof}

\begin{co}\label{lb105}
We have
\begin{align}
\Omega_{+-}\big(\boxbackslash_{\fk X}(\Wbb_\blt)\big)=\scr T_{\fk X,0,0}^*(\Wbb_\blt).
\end{align}
Moreover, if $\Vbb$ is $C_2$-cofinite, then $\Omega_{+-}\big(\boxbackslash_{\fk X}(\Wbb_\blt)\big)$ is finite dimensional.
\end{co}

\begin{proof}
The equation follows directly from Lemma \ref{lb103} and the definition of $\Omega_{+-}$ in \eqref{eq207}. If $\Vbb$ is $C_2$-cofinite, we may show that $\scr T_{\fk X,0,0}(\Wbb_\blt)$ is finite dimensional using the idea in the proof of Theorem \ref{lb67}. (In particular, the vanishing Theorem \ref{lb21}, in which  $\scr V_C^{\leq n}$ is replaced by $\scr V_{\fk X,0,0}^{\leq n}$, is used.)
\end{proof}

The following result is claimed in Remark \ref{lb111}. We are now ready to prove it. Note that neither $C_2$-cofiniteness nor rationality is assumed here.

\begin{pp}\label{lb117}
Let $\Mbb,\wht\Mbb$ be irreducible (ordinary) $\Vbb$-modules. Suppose that the $\fk U_0(\Vbb)$-modules $\Omega(\Mbb),\Omega(\wht\Mbb)$ are equivalent. Then $\Mbb\simeq\wht\Mbb$.	
\end{pp}

\begin{proof}
We set $\fk X=(\Pbb^1;1;0;\infty;z-1;z;1/z)$. Namely, we choose $C=\Pbb^1,N=1,x_1=1,x'=0,x''=\infty,\eta_1=z-1,\xi=z,\varpi=1/z$. We set $\Wbb_\blt=\Wbb_1=\Vbb$. Define $\upphi\in(\Vbb\otimes\Mbb\otimes\Mbb')^*$ to be $\upphi(v\otimes m\otimes m')=\bk{Y_\Wbb(v,1)m,m'}$. Then, as shown in example \ref{lb44}, $\upphi$ is a (non-zero) element of $\scr T_{\fk X}^*(\Vbb\otimes\Mbb\otimes\Mbb')$. Define $\Psi_\upphi:\Mbb\otimes\Mbb'\rightarrow\boxbackslash_{\fk X}(\Vbb)$ as in example \ref{lb112}, i.e.,  by
\begin{align*}
\Psi_\upphi(m\otimes m')(v)=\bk{Y_\Wbb(v,1)m,m'}=\bk{Y_\Wbb(v)_{\wt(v-1)}m,m'}.
\end{align*}
$\Psi_\upphi$ is clearly a homomorphism of weak $\Vbb\times\Vbb$-modules. Since $\Mbb$ and (hence) $\Mbb'$ are irreducible ordinary $\Vbb$-modules, by Proposition \ref{lb113}, $\Mbb\otimes\Mbb'$ is an irreducible weak $\Vbb\times\Vbb$-module. Thus $\Psi_\upphi$ must be injective. So $\Mbb\otimes\Mbb'$ is isomorphic to an irreducible weak $\Vbb\times\Vbb$-module $\mbb K:=\Psi_\upphi(\Mbb\otimes\Mbb')$.

In a similar way, we can define $\wht\upphi,\Psi_{\wht\upphi}$ using $\wht\Mbb$, and $\wht\Mbb\otimes\wht\Mbb'$ is isomorphic to $\wht{\mbb K}:=\Psi_{\wht\upphi}(\wht\Mbb\otimes\wht\Mbb')$. Now, let $T:\Omega(\Mbb)\rightarrow\Omega(\wht\Mbb)$ be an isomorphism of $\fk U_0(\Vbb)$-module. Then, by Remark \ref{lb111}, its transpose $T^t:\Omega(\wht\Mbb')\rightarrow\Omega(\Mbb')$ is also an isomorphism. Choose any $m\in\Mbb,m'\in\Mbb'$. Then
\begin{align*}
&\Psi_{\wht\upphi}(T^{-1}m\otimes T^\tr m')(v)=\bk{Y_{\wht\Mbb}(v)_{\wt(v)-1}T^{-1}m,T^\tr m'}\\
=&\bk{Y_\Mbb(v)_{\wt(v)-1}m,m'}=\Psi_\upphi(m\otimes m')(v).
\end{align*}
This shows that $\mbb K$ and $\wht{\mbb K}$ have at least one non-zero common element. So $\mbb K\cap\wht{\mbb K}$ is a non-trivial weak $\Vbb\times\Vbb$-submodule of $\mbb K$, which must be $\mbb K$ by the irreducibility of $\mbb K$. Thus $\mbb K\subset\wht{\mbb K}$ and, similarly, $\mbb K=\wht{\mbb K}$. Therefore, $\Mbb\otimes\Mbb'$ and $\wht\Mbb\otimes\wht\Mbb'$ are both isomorphic to $\mbb K=\wht{\mbb K}$ as weak $\Vbb\times \Vbb$-modules. In particular, they are equivalent as weak $\Vbb\times\id$-modules. So $\Mbb\simeq\wht\Mbb$.
\end{proof}

The following proposition can be thought of as a converse of example \ref{lb112}.

\begin{pp}\label{lb109}
Let $\Mbb,\wht\Mbb$ be $\Vbb$-modules, and let $\Phi:\Mbb\otimes\wht\Mbb\rightarrow\boxbackslash_{\fk X}(\Wbb_\blt)$ be a homomorphism of weak $\Vbb\times\Vbb$-modules. Then there exists $\uppsi\in\scr T_{\fk X}^*(\Wbb_\blt\otimes\Mbb\otimes\wht\Mbb)$ such that for any $m\in\Mbb,\wht m\in\wht\Mbb,w_\blt\in\Wbb_\blt$,
\begin{align*}
\Phi(m\otimes\wht m)(w_\blt)=\uppsi(w_\blt\otimes m\otimes\wht m).
\end{align*}
\end{pp}

Thus, using the notation of example \ref{lb112}, we have $\Psi_\uppsi=\Phi$.

\begin{proof}
Define $\uppsi$ to be a linear functional on $\Wbb_\blt\otimes\Mbb\otimes\wht\Mbb$ whose value at $w_\blt\otimes m\otimes\wht m$ is $\Phi(m\otimes\wht m)(w_\blt)$. Consider $\wr\Phi(m\otimes\wht m)(w_\blt)$, which is an element of $H^0(C-\SX-\DX,\scr V_C^*)$. By Theorem \ref{lb92}, its series expansion near $x_i$ is of the form
\begin{align*}
&\Phi(m\otimes\wht m)(w_1\otimes\cdots\otimes Y_{\Wbb_i}(v,\eta_i)w_i\otimes\cdots\otimes w_N)\\
=&\uppsi(w_1\otimes\cdots\otimes Y_{\Wbb_i}(v,\eta_i)w_i\otimes\cdots\otimes w_N\otimes m\otimes\wht m)
\end{align*}
when evaluated with $\mc U_\varrho(\eta_i)^{-1}v$ ($v\in\Vbb$). When evaluated with $\mc U_\varrho(\xi)^{-1}v$ (considered as a constant section of $\scr V_C(W')$), it becomes, by \eqref{eq199}, 
\begin{align*}
&\wr\Phi(m\otimes\wht m)(v,w_\blt)=(Y_+(v,\xi)\Phi(m\otimes\wht m))(w_\blt)=\Phi(Y_\Mbb(v,\xi)m\otimes\wht m)(w_\blt)\\
=&\uppsi(w_\blt\otimes Y_\Mbb(v,\xi)m\otimes\wht m).
\end{align*}
Similarly, its evaluation with $\mc U_\varrho(\varpi)^{-1}v$ is $\uppsi(w_\blt\otimes m\otimes Y_{\wht\Mbb}(v,\varpi)\wht m)$. Thus, by Theorem \ref{lb32}, $\uppsi$ is a conformal block.
\end{proof}

\subsection*{The weak $\Vbb\times\Vbb$-module $\boxbackslash_{\fk X}^\low(\Wbb_\blt)$}

In application, it would be more suitable to consider   $\boxbackslash_{\fk X}^\low(\Wbb_\blt)$, \index{W@$\boxbackslash_{\fk X}(\Wbb_\blt),\boxbackslash_{\fk X}^\low(\Wbb_\blt)$} the $\fk U(\Vbb)\times\fk U(\Vbb)$-submodule of $\boxbackslash_{\fk X}(\Wbb_\blt)$ generated by $\Omega_{+-}\big(\boxbackslash_{\fk X}(\Wbb_\blt)\big)$. Recall in Definition \ref{lb133} the meaning of rationality.

\begin{lm}\label{lb106}
Assume that $\Vbb$ is $C_2$-cofinite and rational. Then we have the following equivalence of weak $\Vbb\times\Vbb$-modules
\begin{align}
\boxbackslash_{\fk X}^\low(\Wbb_\blt)\simeq\bigoplus_i \Mbb_i\otimes\wht\Mbb_i\label{eq212}
\end{align}
where each $\Mbb_i$ and $\wht\Mbb_i$ are irreducible (ordinary) $\Vbb$-modules.
\end{lm}

\begin{proof}
By proposition  \ref{lb101} and Theorem \ref{lb108}, it suffices to check that $\boxbackslash_{\fk X}^\low(\Wbb_\blt)$ is  a direct sum of irreducible admissible (and hence ordinary by Corollary \ref{lb99}) $\Vbb\times \id$-modules and also a direct sum of irreducible admissible $\id\times\Vbb$-modules. Indeed, suppose this is true. Then we have \eqref{eq212} where each $\Mbb_i$ and $\wht\Mbb_i$ are irreducible weak $\Vbb$-modules. $\Mbb_i$  must be isomorphic to an irreducible weak $\Vbb\times\id$-submodule of  $\boxbackslash_{\fk X}^\low(\Wbb_\blt)$, which is therefore ordinary. Similary, $\wht\Mbb_i$ is ordinary.

For each $m\in\Zbb$, we let $\fk U_{-m}(\Vbb)$ be the elements of $\fk U(\Vbb)$ raising the $\wtd L_0$-weights by $m$. Namely, it is spanned by $(v_1,n_1)\cdots(v_k,n_k)$ satisfying $\sum_{i=1}^k (\wt(v_i)-n_i-1)=m$. For each $b\in\Nbb$, we consider the  $\Vbb\times\id$-module 
\begin{align*}
\mbb X_b=(\fk U(\Vbb)\times\fk U_{-b}(\Vbb))\scr T_{\fk X,0,0}^*(\Wbb_\blt)
\end{align*}
whose lowest weight subspace is denoted by $\Omega_+(\mbb X_b)$. Using the commutator formula \eqref{eq242}, it is easy to see that $\mbb X_b$ is annihilated by  $Y_-(u)_m$ where $v$ is homogeneous and $\wt (u)-m-1+b<0$ (i.e., $(u,m)$ lowers the $\wtd L_0$-weight by more than $b$). Thus, by Lemma \ref{lb103}, we obtain
\begin{align*}
(\id\times\fk U_{-b}(\Vbb))\scr T_{\fk X,0,0}^*(\Wbb_\blt)\subset\Omega_+(\mbb X_b)\subset \scr T_{\fk X,0,b}^*(\Wbb_\blt).
\end{align*}
The first inclusion shows that $\Omega_+(\mbb X_b)$ generates $\mbb X_b$. The second one shows that $\Omega_+(\mbb X_b)$ is finite-dimensional since $\scr T_{\fk X,0,b}^*(\Wbb_\blt)$ is so by the proof of Theorem \ref{lb67}. Thus,  by the rationality of $\Vbb$, $\mbb X_b$ is a direct sum of irreducible admissible $\Vbb\times\id$-modules. Since $\boxbackslash_{\fk X}^\low(\Wbb_\blt)=\bigvee_{b\in\Nbb}\mbb X_b$, we conclude that $\boxbackslash_{\fk X}^\low(\Wbb_\blt)$ is a sum, and hence a direct sum, of irreducible admissible $\Vbb\times\id$-modules. The claim for $\id\times\Vbb$ is proved similarly.
\end{proof}

Let $\mc E$ be a complete list of mutually inequivalent irreducible (ordinary) $\Vbb$-modules, which is finite by Corollary \ref{lb118}. The word ``complete" means that any irreducible $\Vbb$-module is equivalent to one object in $\mc E$. The following theorem gives a complete characterization of $\boxbackslash_{\fk X}^\low(\Wbb_\blt)$ when $\Vbb$ is $C_2$-cofinite and rational.

\begin{thm}\label{lb134}
Define a homomorphism of weak $\Vbb\times\Vbb$-modules
\begin{gather}
\Psi:\bigoplus_{\Mbb,\wht\Mbb\in\mc E}\Mbb\otimes\wht\Mbb\otimes\scr T_{\fk X}^*(\Wbb_\blt\otimes\Mbb\otimes\wht\Mbb)\rightarrow \boxbackslash_{\fk X}^\low(\Wbb_\blt),\\
m\otimes\wht m\otimes\upphi\mapsto \Psi_\upphi(m\otimes \wht m)\nonumber
\end{gather}
where $\Psi_\upphi$ is defined in example \ref{lb112}. Then $\Psi$ is injective. If $\Vbb$ is $C_2$-cofinite and rational, then $\Psi$ is an isomorphism.
\end{thm}

Note that that the image of each $\Psi_\upphi$ is in $\boxbackslash_{\fk X}^\low(\Wbb_\blt)$ follows from the obvious fact that $\Mbb\otimes\wht\Mbb$ is generated by $\Omega(\Mbb)\otimes\Omega(\wht\Mbb)$.

\begin{proof}
If $\Ker(\Psi)$ is non-trivial, then it is a non-trivial $\fk U(\Vbb)\times\fk U(\Vbb)$-submodule of the domain $\Dom(\Psi)$ of $\Psi$. Since $\Dom(\Psi)$	 is clearly completely reducible, by basic representation theory, so is $\Ker(\Psi)$. Thus, $\Ker(\Psi)$ contains an irreducible $\fk U(\Vbb)\times\fk U(\Vbb)$-submodule $\mbf W$. The projection of $\mbf W$ onto one of the irreducible component of $\Dom(\Psi)$ is non-trivial. Thus, $\mbf W\simeq \Mbb\otimes\wht\Mbb$ for some $\Mbb,\wht\Mbb\in\mc E$. By Corollary \ref{lb107}, there exists a nonzero $\upphi\in\scr T_{\fk X}^*(\Wbb_\blt\otimes\Mbb\otimes\wht\Mbb)$ such that $\mbf W=\Mbb\otimes\wht\Mbb\otimes\upphi$. As $\mbf W\subset\Ker(\Psi)$, we have $\Psi_\upphi=0$. By \eqref{eq216}, for any $w_\blt\in\Wbb_\blt,m\in\Mbb,\wht m\in\wht\Mbb$, we have $\upphi(w_\blt\otimes m\otimes \wht m)=0$. So $\upphi=0$, which gives a contradiction.

When $\Vbb$ is $C_2$-cofinite and rational, the surjectivity of $\Psi$ follows from Proposition \ref{lb109} and Lemma \ref{lb106}.
\end{proof}

\section{Factorization}

We assume the setting of Section \ref{lb104}. Thus, we recall that $\fk X$ is a family over $\mc D_{r\rho}$ obtained by sewing an $N$-pointed compact Riemann surface  $\wtd{\fk X}=(\wtd C;x_1,\dots,x_N;x';x'';\eta_1,\dots,\eta_N;\xi,\varpi)$. Recall $\SX=x_1+\cdots+x_N$ and $\DX=x'+x''$.  Note that for each $q\in\mc D_{r\rho}$, the fiber $\mc C_q$ is nodal (with one node) if and only if $q=0$. Moreover, $\wtd C$ is the normalization of the nodal curve $C:=\mc C_0$. In particular, we have  $\nu:\wtd C\rightarrow C$ defined by gluing $x',x''$. Also, $\fk X_0=(C;x_1,\dots,x_N;\eta_1,\dots,\eta_N)$.

We assume  that $\Vbb$ is both $C_2$-cofinite and rational. By Convention \ref{lb36}, for each irreducible $\Vbb$-module $\Mbb$,  its lowest $\wtd L_0$-weight (with nontrivial weight space) is $0$. Thus $\Omega(\Mbb)=\Mbb(0)$.

Choose any $\upphi\in\scr T_{\fk X_0}^*(\Wbb_\blt)$. Then $\upphi$ is a linear functional on $\Wbb_\blt$.

\begin{lm}
$\upphi$ is an element of $\scr T_{\wtd{\fk X},0,0}^*(\Wbb_\blt)$.
\end{lm}

\begin{proof}
Recall that $W',W''$ are open discs in $\wtd C$ centered at $x',x''$  respectively and disjoint from $x_1,\dots,x_N$. A section $\upnu$ of $\scr V_{\wtd{\fk X},0,0}\otimes\omega_{\wtd C}$ defined on $W'\cup W''$	is of the form
\begin{gather*}
\mc U_\varrho(\xi)\upnu|_{W'}=\xi^{L_0}ud\xi,\qquad \mc U_\varrho(\varpi)\upnu|_{W''}=\varpi^{L_0}vd\varpi
\end{gather*}
where $u\in\Vbb\otimes_\Cbb\scr O(W'),v\in\Vbb\otimes_\Cbb\scr O(W'')$.   $\upnu$ can be viewed as 	an element of $\scr V_C\otimes\omega_C((W'-x')\cup (W''-x''))$. By the description of $\scr V_C$ and $\omega_C$ near the node (see \eqref{eq225} and \eqref{eq159}), it is easy to see that $\upnu$ belongs to $\scr V_C\otimes\omega_C(\nu(W'\cup W''))$. (Indeed, one can check that the value of $\upnu$ at the node $y'=\nu(x')=\nu(x'')$ (as an element of $\scr V_C\otimes\omega_C|y'$) is $0$.) So $\scr V_{\wtd{\fk X},0,0}\otimes\omega_{\wtd C}$ is naturally a subsheaf of $\scr V_C\otimes\omega_C$. Thus, as $\upphi$ vanishes on $H^0(C,\scr V_C\otimes\omega_C(\blt\SX))\cdot \Wbb_\blt$, it vanishes on the subspace $H^0(\wtd C,\scr V_{\wtd{\fk X},0,0}\otimes\omega_{\wtd C}(\blt\SX))\cdot \Wbb_\blt$. This proves $\upphi\in\scr T_{\fk X,0,0}^*(\Wbb_\blt)$.
\end{proof}	

\begin{lm}\label{lb116}
For any homogeneous $v\in\Vbb$,
\begin{align}
Y_+(v)_{\wt(v)-1}\upphi=\sum_{l\in\Nbb}\frac{(-1)^{\wt(v)}}{l!}Y_-(L_1^lv)_{\wt(v)-l-1}\upphi.\label{eq214}
\end{align}
\end{lm}
\begin{proof}
Consider $\upphi$ as a conformal block on $C$. Choose any homogeneous $v\in\Vbb^{\leq n}$.  As argued in the proof of Theorem \ref{lb114}, one can construct $\upnu\in H^0(\wtd C,\scr V_{\wtd C}^{\leq n}\otimes\omega_{\wtd C}(\blt\SX+\DX))$ such that
\begin{gather*}
\mc U_\varrho(\xi)\upnu|_{W'}=\xi^{L_0-1}vd\xi +\xi^n(\mathrm{elements~of~}\Vbb^{\leq n}\otimes_\Cbb\scr O(W'))d\xi,\\
\mc U_\varrho(\varpi)\upnu|_{W''}=-\varpi^{L_0-1}\mc U(\upgamma_1)vd\varpi+\varpi^n(\mathrm{elements~of~}\Vbb^{\leq n}\otimes_\Cbb\scr O(W''))d\varpi.
\end{gather*}
By \eqref{eq159} and \eqref{eq54}, it is clear that $\upnu$ can be viewed as an element of $H^0(C,\scr V_C\otimes\omega_C(\blt\SX))$. Thus $\upphi(\upnu\cdot w_\blt)=0$. 

Consider $\wr\upphi(w_\blt)\in H^0(\wtd C-\SX-\DX,\scr V_{\wtd C})$. Then $\wr\upphi(\upnu,w_\blt)\in H^0(\wtd C-\SX-\DX,\omega_{\wtd C})$. By residue theorem, we have
\begin{align*}
\Res_{x'}\wr\upphi(\upnu,w_\blt)+\Res_{x''}\wr\upphi(\upnu,w_\blt)=-\sum_{i=1}^N \Res_{x_i}\wr\upphi(\upnu,w_\blt),
\end{align*}
which, by Theorem \ref{lb92}, equals $-\upphi(\upnu\cdot w_\blt)$ and hence is zero. By \eqref{eq199}, we have, under the identification $\scr V_{W'}\simeq\Vbb\otimes_\Cbb\scr O_{W'}$  via $\mc U_\varrho(\xi)$,
\begin{align*}
\Res_{x'}\wr\upphi(\upnu,w_\blt)=\Res_{\xi=0}\wr\upphi(\xi^{\wt(v)-1}v+\cdots,w_\blt)d\xi=Y_+(v)_{\wt(v)-1}\upphi(w_\blt).
\end{align*}
By \eqref{eq73},
\begin{align*}
\varpi^{L_0-1}\mc U(\upgamma_1)v=\varpi^{L_0-1}e^{L_1}(-1)^{L_0}v=\sum_{l\in\Nbb}\frac{(-1)^{\wt(v)}}{l!}L_1^lv\cdot \varpi^{\wt(v)-1-l}d\varpi.
\end{align*}
Thus, by \eqref{eq213},
\begin{align*}
\Res_{x''}\wr\upphi(\upnu,w_\blt)=-\sum_{l\in\Nbb}\frac{(-1)^{\wt(v)}}{l!}Y_-(L_1^lv)_{\wt(v)-l-1}\upphi(w_\blt).
\end{align*}
This proves \eqref{eq214}.
\end{proof}

Let $\mc E$ be a complete list of mutually inequivalent irreducible $\Vbb$-modules.  By Corollary \ref{lb118}, $\mc E$ is a finite set. Recall the map
\begin{gather*}
\wtd{\fk S}_q:\bigoplus_{\Mbb\in\mc E}\scr T_{\wtd{\fk X}}^*(\Wbb_\blt\otimes\Mbb\otimes\Mbb')\rightarrow\scr T_{\fk X_q}^*(\Wbb_\blt),\\
\bigoplus_{\Mbb\in\mc E}\uppsi_\Mbb\mapsto\sum_{\Mbb}\wtd{\mc S}\uppsi_\Mbb(q)
\end{gather*}
defined for any $q\in\mc D_{r\rho}^\times$ in Section \ref{lb104}. When $q=0$, an element in the image of $\wtd{\fk S}_q$ is a linear functional on $\Wbb_\blt$, which must be also in $\scr T_{\fk X_0}^*(\Wbb_\blt)$ by Proposition \ref{lb41}.  Recall our assumption on $\wtd L_0$, which implies that
\begin{align*}
\wtd{\mc S}\uppsi_\Mbb(w_\blt)(0)=\uppsi_\Mbb(w_\blt\otimes P(0)\btr\otimes_\Mbb~\btl).
\end{align*}
Here, $P(n)$ is the projection onto the $\wtd L_0$-weight $n$ subspace, and hence $P(0)$ is the projection onto the lowest weight subspace.

\begin{pp}[Nodal factorization]\label{lb120}
$\wtd{\fk S}_0$ is a surjective linear map.
\end{pp}

Indeed, one can use the same method for Theorem \ref{lb114} to prove that $\wtd{\fk S}_0$ is also injective.

\begin{proof}
Choose $\upphi\in\scr T_{\fk X_0}^*(\Wbb_\blt)$. By Corollary \ref{lb105}, $\upphi$ is in the finite-dimensional space $\Omega_{+-}\big(\boxbackslash_{\wtd{\fk X}}(\Wbb_\blt)\big)$. By Theorem \ref{lb106}, we can find finitely many irreducible $\Vbb$-modules $\Mbb_1,\dots,\Mbb_l,\wht\Mbb_1,\dots,\wht\Mbb_l$ such that $\bigoplus_{1\leq i\leq l}\Mbb_i\otimes\wht\Mbb_i'$ is a weak $\Vbb\times\Vbb$-submodule of $\boxbackslash_{\wtd{\fk X}}(\Wbb_\blt)$, that $\upphi\in \bigoplus_{1\leq i\leq l}\Mbb_i(0)\otimes\wht\Mbb_i'(0)$, and that the projection of $\uppsi$ to each $\Mbb_i(0)\otimes\wht\Mbb_i'(0)$, which we denote by
\begin{align*}
\phi_i\in\Mbb_i(0)\otimes\wht\Mbb_i'(0),
\end{align*}
is non-zero. By Lemma \ref{lb116} and equation \eqref{eq215}, we have for each $i$ that 
\begin{align*}
(Y_{\Mbb_i}(v)_{\wt(v)-1}\otimes\id)\phi_i=(\id\otimes Y_{\wht\Mbb_i}(v)_{\wt(v)-1}^\tr)\phi_i.
\end{align*}
Identify $\Mbb_i(0)\otimes\wht\Mbb_i'(0)$ naturally with $\Hom_\Cbb(\Mbb_i(0),\wht\Mbb_i(0))$. Then the above equation shows that $\phi_i\in\Hom_{\fk U_0(\Vbb)}(\Mbb_i(0),\wht\Mbb_i(0))$. Since $\phi_i$ is assumed to be non-zero, we have $\Mbb_i(0)\simeq\wht\Mbb_i(0)$. Thus, by Proposition \ref{lb117}, we have the equivalence of $\Vbb$-modules $\Mbb_i\simeq\wht\Mbb_i$. We may thus assume that $\wht\Mbb_i=\Mbb_i$. As $\upphi_i\in\End_{\fk U_0(\Vbb)}(\Mbb_i(0))$, it is a scalar multiplication. We may thus find  a (non-zero) $\lambda_i\in\Cbb$ such that $\phi_i=\lambda_i P(0)\btr\otimes_i\btl$.

By Theorem \ref{lb134} or Proposition \ref{lb109}, we can find for each $i$ a non-zero conformal block $\uppsi_i\in\scr T_{\wtd{\fk X}}^*(\Wbb_\blt\otimes\Mbb_i\otimes\Mbb_i')$ such that
\begin{align*}
(m_i\otimes\wht m_i')(w_\blt)=\uppsi_i(w_\blt\otimes m_i\otimes\wht m_i')
\end{align*}
for each $w_\blt\in\Wbb_\blt,m_i\in\Mbb_i,\wht m_i'\in\wht\Mbb_i'$. In particular,
\begin{align*}
\phi_i(w_\blt)=\uppsi_i(w_\blt\otimes\phi_i)=\lambda_i\upphi_i(w_\blt\otimes P(0)\btr\otimes_i\btl)=\lambda_i\wtd{\mc S}\uppsi_i(w_\blt)(0)
\end{align*}
Thus $\upphi=\sum_i\lambda_i\wtd{\mc S}\uppsi_i(0)$. It is now clear that $\upphi$ is in the image of $\wtd{\fk S}_0$.
\end{proof}

\begin{thm}[Factorization]\label{lb135}
Assume that $\Vbb$ is $C_2$-cofinite and rational. Then $\wtd{\fk S}_q$ is an isomorphism for each $q\in\mc D_{r\rho}$, and $\fk S_q$ is an isomorphism for each $q\in\mc D_{r\rho}^\times$.
\end{thm}

\begin{proof}
As explained before Theorem \ref{lb114}, it suffices to prove that $\wtd{\fk S}_q$ is an isomorphism for each $q\in\mc D_{r\rho}$.  Let $D\in\Nbb$ be the dimension of the domain of $\wtd{\fk S}_q$, which is finite by Theorem \ref{lb67} and Corollary \ref{lb118}. Let $K_q$ be the dimension of the image of $\wtd{\fk S}_q$. By Corollary \ref{lb119}, if $q\neq 0$ then $K_q$ is independent of $q$. We fix $q\in\mc D_{r\rho}^\times$. By Proposition \ref{lb60} (Nakayama's lemma), $K_q\leq K_0$. By Theorem \ref{lb114}, $D\leq K_q$. By Proposition \ref{lb120}, $K_0\leq D$. Thus $D=K_q=K_0$.
\end{proof}

\begin{rem}\label{lb136}
The above two theorems show that if a nodal curve $C$ (with normalization $\wtd C$) has one node, then its dimensions of spaces of conformal blocks can be calculated from those of $\wtd C$. Note that the results in Section \ref{lb110} can be generalized to $N$-pointed compact Riemann surfaces with local coordinates and arbitrary numbers  of outputs. Accordingly, we may prove the nodal factorization for an $N$-pointed nodal curve $(C;x_1,\dots,x_N)$ with an arbitrary number $M$ of nodes in the same way: the dimension of the space of conformal blocks associated to the modules $\Wbb_1,\dots,\Wbb_N$ is
\begin{align}
\sum_{\Mbb_1,\dots,\Mbb_M\in\mc E} \dim\scr T_{\wtd{\fk X}}^*(\Wbb_\blt\otimes\Mbb_\blt\otimes\Mbb_\blt')\label{eq243}
\end{align}
where $\wtd{\fk X}$ is the compact Riemann surface $\wtd C$ with $N+2M$ marked points, and $\Wbb_\blt=\Wbb_1\otimes\cdots\otimes\Wbb_N$, $\Mbb_\blt=\Mbb_1\otimes\cdots\otimes\Mbb_M$, $\Mbb_\blt'=\Mbb_1'\otimes\cdots\otimes\Mbb_M'$, as usual.
\end{rem}

\begin{thm}
Assume that $\Vbb$ is $C_2$-cofinite and rational. Let $\fk X=(\pi:\mc C\rightarrow\mc B;\sgm_1,\dots,\sgm_N)$ be a family of $N$-pointed complex curves. Then $\scr T_{\fk X}(\Wbb_\blt)$ and (hence) $\scr T_{\fk X}^*(\Wbb_\blt)$ are locally free. Consequently, the function $b\in\mc B\mapsto \dim\scr T_{\fk X_b}^*(\Wbb_\blt)$ is locally constant.
\end{thm}

\begin{proof}
Recall that by Theorem \ref{lb69}, the dimensions of  the fibers of $\scr T_{\fk X}(\Wbb_\blt)$ are given by the dimensions of the spaces of conformal blocks.   Outside the discriminant locus $\Delta$, the dimensions equal \eqref{eq243} by applying Theorem \ref{lb135} several times. In $\Delta$, the same is true by nodal factorization (Remark \ref{lb136}). Thus, the rank function of $\scr T_{\fk X}(\Wbb_\blt)$ is locally constant. By Theorem \ref{lb67}, $\scr T_{\fk X}(\Wbb_\blt)$ is finite-type. Thus, it is locally free by Theorem \ref{lb59}.
\end{proof}

\chapter{More on the connections}

\section{Connections and local coordinates}\label{lb142}

Consider a family of $N$-pointed compact Riemann surfaces $\fk X=(\pi:\mc C\rightarrow\mc B;\sgm_1,\dots,\sgm_N)$. Assume $\mc B$ is biholomorphic to a Stein open subset of $\Cbb^m$, and $\fk X$ admits a set of local coordinates $\eta_1,\dots,\eta_N$. For each $\yk\in\Theta_{\mc B}(\mc B)$, we have defined in Section \ref{lb84} a differential operator $\nabla_\yk$ depending on $\eta_\blt$ and a lift $\wtd\yk\in H^0(\mc C,\Theta_{\mc C}(\blt\SX))$ (satisfying $d\pi(\wtd\yk)=\pi^*\yk$). We have also seen that if $\eta_\blt$ is fixed, then $\nabla_\yk$ is determined up to an $\scr O(\mc B)$-scalar addition by $\wtd\yk$. In this section, we show the same is true for the dependence of $\nabla_\yk$ on $\eta_\blt$ \emph{if $\Wbb_1,\dots,\Wbb_N$ are simple $\Vbb$-modules}.

Let $\tau_\blt=(\tau_1,\dots,\tau_m)$ be coordinates of $\mc B$, and  write $\tau_\blt\circ\pi$ also as $\tau_\blt$ for simplicity. Then we can write
\begin{align*}
\yk=\sum_j g_j(\tau_\blt)\partial_{\tau_j}
\end{align*}
for some $g_1,\dots,g_m\in\scr O(\tau_\blt(\mc B))$. Choose mutually disjoint neighborhoods $W_1,\dots,W_N$ of $\sgm_1(\mc B),\dots,\sgm_N(\mc B)$ on which $\eta_1,\dots,\eta_N$ are defined respectively. Then $(\eta_i,\tau_\blt)$ is a coordinate of $W_i$. So we can find $h_i\in\scr O((\eta_i,\tau_\blt)(W_i\setminus\sgm_i(\mathcal B)))$ such that
\begin{align}
\wtd\yk|_{W_i}=h_i(\eta_i,\tau_\blt)\partial_{\eta_i}+\sum_j g_j(\tau_\blt)\partial_{\tau_j/\eta_i},\label{eq259}
\end{align}
where $\partial_{\tau_j/\eta_i}$ means the partial derivative $\partial_{\tau_j}$ defined with respect to the coordinate $(\eta_i,\tau_\blt)$.

Let $\nabla_\yk^{(\eta_\blt)}$ denote the differential operator defined by $\eta_\blt$ and $\wtd\yk$. Let $\mu_\blt$ be another set of local coordinates of $\fk X$, and let $\nabla_\yk^{(\mu_\blt)}$ denote the differential operator defined by $\mu_\blt$ and $\wtd\yk$. Let 
\begin{align*}
\alpha_i=(\mu_i|\eta_i):\mc B\rightarrow\Gbb
\end{align*}
be a holomorphic family of transformations whose value at each $b\in\mc B$ is defined similar to \eqref{eq244}, namely  $\alpha_{i,b}(z)=\mu_i\circ(\eta_i,\pi)^{-1}(z,b)$ for each $b\in\mc B$. Equivalently, consider $\alpha_i$ as a holomorphic function defined near $0\times\mc B\subset\Cbb\times\mc B$,
\begin{align}
\alpha_i\circ(\eta_i,\pi)=(\mu_i,\pi).\label{eq249}
\end{align}
Let $\alpha_i'(0)\in\scr O(\mc B)$ whose value at each $b\in\mc B$ is  $\partial_z\alpha_{i,b}(z)|_{z=0}$. Assume $\Vbb$ has central charge $c$, and $\Wbb_1,\dots,\Wbb_N$ are simple $\Vbb$-modules. For each $\Wbb_i$, we let $\Delta_{\Wbb_i}\in\Cbb$ be the unique number such that $\Wbb_i=\bigoplus_{n\in\Nbb}(\Wbb_i)_{(\Delta_{\Wbb_i}+n)}$ and $(\Wbb_i)_{(\Delta_{\Wbb_i})}$ is non-trivial.   Thus, according to Convention \ref{lb36}, $L_0-\wtd L_0=\Delta_{\Wbb_i}\id$ on $\Wbb_i$. $\Delta_{\Wbb_i}$ is called the \textbf{conformal weight} of $\Wbb_i$. \index{zz@$\Delta_{\Wbb}$}

\begin{thm}\label{lb144}
Let $f\in\scr O(\mc B)$ be 
\begin{align}\label{eq245}
f=\sum_{i=1}^N\Big(\Res_{\eta_i=0} ~\frac c{12}h_i(\eta_i,\tau_\blt)\Sbf_{\eta_i}\mu_i d\eta_i-\Delta_{\Wbb_i}\sum_{j=1}^m g_j(\tau_\blt)\partial_{\tau_j}\log\alpha_i'(0)  \Big).
\end{align}
Then for each section  $s$ of $\scr T_{\fk X}(\Wbb_\blt)$,
\begin{align*}
\nabla_\yk^{(\eta_\blt)}s-\nabla_\yk^{(\mu_\blt)}s=-fs.
\end{align*}
Consequently, for each section $\upphi$ of $\scr T_{\fk X}^*(\Wbb_\blt)$,
\begin{align*}
\nabla_\yk^{(\eta_\blt)}\upphi-\nabla_\yk^{(\mu_\blt)}\upphi=f\upphi.
\end{align*}
\end{thm}

Note that the first part on the right hand side of \eqref{eq245} is similar to $\frac c{12}C_i$ where $C_i$ is defined in Remark \ref{lb138}. Also, the residue $\Res_{\eta_i=0}$ is taken with respect to the coordinate $(\eta_i,\tau_\blt)$.

\begin{proof}
	
	
Choose $w_1\in\Wbb_1,\dots,w_N\in\Wbb_N$. Recall $w_\blt=w_1\otimes\cdots\otimes w_N$. Let $s=\mc U(\eta_\blt)^{-1}w_\blt\in\scr W_{\fk X}(\Wbb_\blt)(\mc B)$. So
\begin{align}
\mc U(\mu_\blt)s=\mc U(\alpha_1)w_1\otimes\cdots\otimes \mc U(\alpha_N)w_N.\label{eq250}
\end{align}
By \eqref{eq185} and \eqref{eq247}, 
\begin{align*}
\mc U(\eta_\blt)\nabla_\yk^{(\eta_\blt)} s=\sum_i w_1\otimes w_2\otimes\cdots\otimes \wtd w_i\otimes\cdots\otimes w_N
\end{align*}
where
\begin{align*}
\wtd w_i=-\Res_{\eta_i=0}~h_i(\eta_i,\tau_\blt)Y_{\Wbb_i}(\cbf,\eta_i)w_i d\eta_i.
\end{align*}
Thus
\begin{align}
\mc U(\mu_\blt)\nabla_\yk^{(\eta_\blt)} s=\sum_i \mc U(\alpha_1)w_1\otimes\cdots\otimes \mc U(\alpha_i)\wtd w_i\otimes\cdots\otimes\mc U(\alpha_N) w_N.\label{eq251}
\end{align}
By \eqref{eq169},
\begin{align*}
\mc U(\varrho(\alpha_i|\id_\Cbb))\cbf=\mc U(\varrho(\mu_i|\eta_i))\cbf=(\partial_{\eta_i}\mu_i)^2\cbf+\frac c{12}\Sbf_{\eta_i}\mu_i\cdot\id.
\end{align*}
Using Theorem \ref{lb31}, we compute
\begin{align}
&\mc U(\alpha_i)\wtd w_i=-\Res_{\eta_i=0~}h_i\cdot\mc  U(\alpha_i)Y_{\Wbb_i}(\cbf,\eta_i)w_i d\eta_i\nonumber\\
=&-\Res_{\eta_i=0}~h_i\cdot  Y_{\Wbb_i}(\mc U(\varrho(\alpha_i|\id_\Cbb))\cbf,\mu_i)\mc U(\alpha_i)w_i d\eta_i\nonumber\\
=&-\Res_{\mu_i=0}~h_i \partial_{\eta_i}\mu_i\cdot Y_{\Wbb_i}(\cbf,\mu_i)\mc U(\alpha_i)w_i d\mu_i-\Res_{\eta_i=0}~ \frac c{12}h_i\Sbf_{\eta_i}\mu_i\cdot \mc U(\alpha_i)w_i d\eta_i.\label{eq252}
\end{align}

We write $\wtd\yk$ in the $(\mu_i,\tau_\blt)$-coordinate:
\begin{align*}
\wtd\yk|_{W_i}=\Big(h_i\partial_{\eta_i}\mu_i+\sum_j g_j\partial_{\tau_j/\eta_i}\mu_i\Big)\partial_{\mu_i}+\sum_j g_j\partial_{\tau_j/\mu_j}.
\end{align*}
Recall \eqref{eq250}, apply \eqref{eq185} and \eqref{eq247} again, and use the above expression of $\wtd\yk$ in the $(\mu_i,\tau_\blt)$-coordinate, we have
\begin{align}
\mc U(\mu_\blt)\nabla_\yk^{(\mu_\blt)} s=\sum_i\mc U(\alpha_1)w_1\otimes\cdots\otimes \upupsilon_i\otimes\cdots\otimes \mc U(\alpha_N)w_N\label{eq253}
\end{align}
where (defining $\wtd\upnu$ using $\mu_\blt$)
\begin{align*}
&\upupsilon_i=\sum_j g_j\partial_{\tau_j}\mc U(\alpha_i)w_i-\wtd\upnu(\wtd\yk)\mc U(\alpha_i)w_i.\nonumber\\
=&\sum_j g_j\partial_{\tau_j}\mc U(\alpha_i)w_i-\Res_{\mu_i=0}~\Big(h_i\partial_{\eta_i}\mu_i+\sum_j g_j\partial_{\tau_j/\eta_i}\mu_i\Big)Y_{\Wbb_i}(\cbf,\mu_i)\mc U(\alpha_i)w_id\mu_i.
\end{align*}
From \eqref{eq249}, it is easy to see
\begin{align*}
\partial_{\tau_j/\eta_i}\mu_i=(\partial_{\tau_j}\alpha_i)(\eta_i,\pi)=(\partial_{\tau_j}\alpha_i)(\alpha_i^{-1}(\mu_i,\pi),\pi)
\end{align*}
where $\alpha_i^{-1}$ is the fiberwise inverse of $\alpha_i$. Identify $W_i$ with a neighborhood of $0\times\mc B$ via $(\mu_i,\pi)$ so that $(\mu_i,\pi)$ is identified with $(z,\id_{\mc B})$, and think of $\alpha_i$ as a family of transformation and write the parameter of $\mc B$ as the subscript of $\alpha_i$, we have
\begin{align*}
\partial_{\tau_j/\eta_i}\mu_i(z,b)=(\partial_{\tau_j}\alpha_i)(\alpha_i^{-1}(z,b),b)=(\partial_{\tau_j}\alpha_i)_b(\alpha_{i,b}^{-1}(z))
\end{align*}
or simply
\begin{align*}
\partial_{\tau_j/\eta_i}\mu_i=(\partial_{\tau_j}\alpha_i)\circ\alpha_i^{-1}.
\end{align*}
Use this relation and apply Lemma \ref{lb141} to the family $\alpha_i$, we have
\begin{align*}
&\partial_{\tau_j}\mc U(\alpha_i)w_i=\Res_{z=0}~(\partial_{\tau_j}\alpha_i)(\alpha_i^{-1}(z)) Y_\Wbb(\cbf,z) \mc U(\alpha_i)w_i dz-\Delta_{\Wbb_i}\partial_{\tau_j}\log\alpha_i'(0)\mc U(\alpha_i)w_i\nonumber\\
=&\Res_{\mu_i=0}~\partial_{\tau_j/\eta_i}\mu_i\cdot Y_\Wbb(\cbf,\mu_i) \mc U(\alpha_i)w_i d\mu_i-\Delta_{\Wbb_i}\partial_{\tau_j}\log\alpha_i'(0)\mc U(\alpha_i)w_i.
\end{align*}
Thus
\begin{align}
\upupsilon_i=-\Res_{\mu_i=0}~h_i\partial_{\eta_i}\mu_i Y_{\Wbb_i}(\cbf,\mu_i)\mc U(\alpha_i)w_id\mu_i-\Delta_{\Wbb_i}\sum_j g_j\partial_{\tau_j}\log\alpha_i'(0)\mc U(\alpha_i)w_i.\label{eq254}
\end{align}
Combine \eqref{eq251}, \eqref{eq252}, \eqref{eq253}, \eqref{eq254} together, and notice \eqref{eq250}, we obtain $\mc U(\mu_\blt)\nabla_\yk^{(\eta_\blt)}s-\mc U(\mu_\blt)\nabla_\yk^{(\mu_\blt)}s=-f\mc U(\mu_\blt)s$.
\end{proof}

\begin{lm}\label{lb141}
Let $T$ be an open subset of $\Cbb$. Let $\rho:T\rightarrow\Gbb,\zeta\mapsto\rho_\zeta$ be a holomorphic family of transformations. Then for any $\Vbb$-module $\Wbb$, if we let $A=L_0-\wtd L_0$, then
\begin{align}
\partial_\zeta\mc U(\rho_\zeta)w=\Res_{z=0}~(\partial_\zeta\rho_\zeta)(\rho_\zeta^{-1}(z)) Y_\Wbb(\cbf,z) \mc U(\varrho_\zeta)w dz-\partial_\zeta\log\rho_\zeta'(0)A\mc U(\rho_\zeta)w.
\end{align}
\end{lm}

\begin{proof}

Choose any $\zeta_0\in T$ and apply Lemma \ref{lb29} and Remark \ref{lb140} to the family $\zeta\mapsto\rho_\zeta\circ\rho_{\zeta_0}^{-1}$, we have
\begin{align*}
&\partial_\zeta\mc U(\rho_\zeta)w\big|_{\zeta=\zeta_0}=\partial_\zeta\mc U(\rho_\zeta\circ\rho_{\zeta_0}^{-1})\mc U(\rho_{\zeta_0})w\big|_{\zeta=\zeta_0}\\
=&\Res_{z=0}~\partial_\zeta(\rho_\zeta\circ\rho_{\zeta_0}^{-1})(z) Y_\Wbb(\cbf,z)\mc U(\rho_{\zeta_0})w dz\big|_{\zeta=\zeta_0}-\partial_\zeta(\rho_\zeta\circ\rho_{\zeta_0}^{-1})'(0)A\mc U(\rho_{\zeta_0})w\big|_{\zeta=\zeta_0}.
\end{align*}
$\partial_\zeta(\rho_\zeta\circ\rho_{\zeta_0}^{-1})(z)$ is just $(\partial_\zeta\rho_\zeta)(\rho_{\zeta_0}^{-1}(z))$. Note that $\rho_{\zeta_0}^{-1}$ is the inverse function of $\rho_{\zeta_0}$, whose derivative is $1/(\rho_{\zeta_0}'\circ \rho_{\zeta_0}^{-1})$. Thus
\begin{align*}
(\rho_\zeta\circ\rho_{\zeta_0}^{-1})'=(\rho_\zeta'\circ\rho_{\zeta_0}^{-1})\cdot(\rho_{\zeta_0}^{-1})'=\frac{\rho_\zeta'\circ\rho_{\zeta_0}^{-1}}{\rho_{\zeta_0}'\circ \rho_{\zeta_0}^{-1}}
\end{align*}
whose value at $z=0$ (noticing $\rho_{\zeta_0}^{-1}(0)=0$) is $\rho_\zeta'(0)/\rho_{\zeta_0}'(0)$. Therefore, 
\begin{align*}
\partial_\zeta(\rho_\zeta\circ\rho_{\zeta_0}^{-1})'(0)\big|_{\zeta=\zeta_0}=\partial_\zeta\rho_\zeta'(0)/\rho_{\zeta_0}'(0)\big|_{\zeta=\zeta_0}=\partial_\zeta\log\rho_\zeta'(0)\big|_{\zeta=\zeta_0}.
\end{align*}
This proves the desired equation at $\zeta=\zeta_0$.
\end{proof}

\section{Projective flatness of connections}

Our  goal of this section is to calculate the curvature of the connection associated to a family of $N$-pointed compact Riemann surfaces with local coordinates $\fk X=(\pi:\mc C\rightarrow\mc B;\sgm_1,\dots,\sgm_N;\eta_1,\dots,\eta_N)$. Choose sections $\yk,\zk$ of $\Theta_{\mc B}$ defined on Stein open subsets of $\mc B$. Choose lifts $\wtd\yk,\wtd\zk$ as in the previous section as in Section \ref{lb142} or \ref{lb84}. We write their local expressions at $W_i$ as

\begin{gather*}
\wtd\yk|_{W_i}=h_i(\eta_i,\tau_\blt)\partial_{\eta_i}+\sum_j g_j(\tau_\blt)\partial_{\tau_j},\\
\wtd\zk|_{W_i}=k_i(\eta_i,\tau_\blt)\partial_{\eta_i}+\sum_j l_j(\tau_\blt)\partial_{\tau_j}.
\end{gather*}
For brevity, in the above expressions, we set
\begin{gather*}
Y=\sum_j g_j(\tau_\blt)\partial_{\tau_j},\qquad Z=\sum_j l_j(\tau_\blt)\partial_{\tau_j},
\end{gather*}
which have the same expressions as $\yk,\zk$, although their meanings are slightly different. We let $[\wtd\yk,\wtd\zk]$ be the lift of $[\yk,\zk]$. Define $\nabla_\yk,\nabla_\zk,\nabla_{[\yk,\zk]}$ using these lifts and the local coordinates $\eta_\blt$. Let $R(\yk,\zk)=\nabla_\yk\nabla_\zk-\nabla_\zk\nabla_\yk-\nabla_{[\yk,\zk]}$. Choose  $\Vbb$ with central charge $c$, and $\Vbb$-modules $\Wbb_1,\dots,\Wbb_N$.

\begin{thm}\label{lb145}
Let $f\in\scr O(\mc B)$ be
\begin{align}
f=-\sum_{i=1}^N \Big(\Res_{\eta_i=0} ~\frac c{12}\partial_{\eta_i}^3h_i(\eta_i,\tau_\blt)\cdot k_i(\eta_i,\tau_\blt) d\eta_i\Big).
\end{align}
Then for each section  $s$ of $\scr T_{\fk X}(\Wbb_\blt)$,
\begin{align*}
R(\yk,\zk)s=-fs.
\end{align*}
Consequently, for each section $\upphi$ of $\scr T_{\fk X}^*(\Wbb_\blt)$,
\begin{align*}
R(\yk,\zk)\upphi=f\upphi.
\end{align*}
\end{thm}

Thus, the (local) connections defined by these differential operators are projectively flat, and the curvatures depend only on the central charge $c$ of $\Vbb$, but not on $\Vbb$ or its modules.

\begin{proof}
We have
\begin{align*}
\nabla_\zk s=Zs-\sum_i (k_i\cbf d\eta_i)\cdot s
\end{align*}
and hence
\begin{align*}
\nabla_\yk\nabla_\zk s=YZs-Y\sum_i (k_i\cbf d\eta_i)\cdot s-\sum_i(h_i\cbf d\eta_i) Zs+\sum_{i,j} (h_i\cbf d\eta_i)(k_j\cbf d\eta_j) s.
\end{align*}
Similarly,
\begin{align*}
\nabla_\zk\nabla_\yk s=ZYs-Z\sum_i (h_i\cbf d\eta_i)\cdot s-\sum_i(k_i\cbf d\eta_i) Ys+\sum_{i,j} (k_j\cbf d\eta_j)(h_i\cbf d\eta_i) s.
\end{align*}

Note that $[Y,k_i\cbf d\eta_i]=(Yk_i)\cbf d\eta_i$ and $[Z,k_i\cbf d\eta_i]=(Zk_i)\cbf d\eta_i$ since $Y,Z$ are orthogonal to $d\eta_i$. Also, if $i\neq j$ then $h_i\cbf d\eta_i$ and $k_j\cbf d\eta_j$ are acting on different tensor-components of $s$. So they commute. Using Proposition \ref{lb51}, we see that in the case that $i=j$, the action of $[h_i\cbf d\eta_i,k_i\cbf d\eta_i]$ on $\scr T_{\fk X}(\Wbb_\blt)$ equals
\begin{align*}
&[h_i\cbf d\eta_i,k_i\cbf d\eta_i]=\mbf L_{h_i\cbf d\eta_i}k_i\cbf d\eta_i \xlongequal{\eqref{eq149}} \sum_{n\geq 0}\frac 1{n!}(\partial_{\eta_i}^nh_i)k_i Y(\cbf)_n\cbf d\eta_i\\
=&\sum_{n\geq 0}\frac 1{n!}(\partial_{\eta_i}^nh_i) k_i L_{n-1}\cbf d\eta_i \xlongequal{\text{Rem. }\ref{lb58}} h_ik_i L_{-1}\cbf d\eta_i+ 2(\partial_{\eta_i}h_i)k_i\cbf d\eta_i+\frac c{12} (\partial_{\eta_i}^3h_i)k_i\id d\eta_i
\end{align*}
where the three summands in the last expression correspond respectively to $n=0,1,3$, the only cases that $L_{n-1}\cbf\neq 0$. By Lemma \ref{lb55}, 
\begin{align*}
&[h_i\cbf d\eta_i,k_i\cbf d\eta_i]=-\partial_{\eta_i}(h_ik_i) \cbf d\eta_i+ 2(\partial_{\eta_i}h_i)k_i\cbf d\eta_i+\frac c{12} (\partial_{\eta_i}^3h_i)k_i\id d\eta_i\\
=&(\partial_{\eta_i}h_i)k_i\cbf d\eta_i-h_i(\partial_{\eta_i}k_i)\cbf d\eta_i+\frac c{12} (\partial_{\eta_i}^3h_i)k_i\id d\eta_i
\end{align*}
when acting on $\scr T_{\fk X}(\Wbb_\blt)$. Thus, 
\begin{align}
[\nabla_\yk,\nabla_\zk]s=&[Y,Z]s-\sum_i (Yk_i\cbf d\eta_i) s+\sum_i(Zh_i\cbf d\eta_i) s\nonumber\\
&+((\partial_{\eta_i}h_i)k_i\cbf d\eta_i)s-(h_i(\partial_{\eta_i}k_i)\cbf d\eta_i)s-fs.\label{eq255}
\end{align}
On the other hand,
\begin{align*}
[\wtd\yk,\wtd\zk]|_{W_i}=(h_i\partial_{\eta_i}k_i-k_i\partial_{\eta_i}h_i+Yk_i-Zh_i)\partial_{\eta_i}+[Y,Z],
\end{align*}
which shows $\nabla_{[\yk,\zk]}s$ equals the sum of all the terms on the right hand side of \eqref{eq255} except $-fs$. This proves the desired relation.
\end{proof}

\section{Constructing flat connections}\label{lb146}

The goal of this section is to define flat connections on sheaves of conformal blocks depending on as few parameters as possible. We adopt the following notation: If $\scr L$ is a line bundle on a complex manifold $X$, then for any sections $s_1,s_2$ of $\scr L$ on an open $U\subset X$, if $s_2$ is nowhere zero, then $\frac {s_1}{s_2}$ is the unique element of $\scr O(U)$ whose multiplication with $s_2$ is $s_1$.

Assume $\Wbb_1,\dots,\Wbb_N$ are simple $\Vbb$-modules. We explain how to obtain a flat connection associated to sheaves of covacua and conformal blocks of $\Vbb$. Let
\begin{align*}
\fk X=(\pi:\mc C\rightarrow\mc B;\sgm_1,\dots,\sgm_N;\nu_1,\dots,\nu_N)
\end{align*}
be a \textbf{family of $N$-pointed compact Riemann surfaces with jets}. This means that $(\pi:\mc C\rightarrow\mc B;\sgm_1,\dots,\sgm_N)$ is $N$-pointed, and the jet
\begin{align*}
\nu_i\in\sgm_i^*\omega_{\mc C/\mc B}(\mc B)
\end{align*}
is nowhere zero for each $1\leq i\leq N$. Thus, for each $b\in\mc B$, $\nu_i(b)$ can be regarded as a (holomorphic) cotangent vector of $\mc C_b$ at $\sgm_i(b)$.

\begin{eg}
If $\fk X$ is a family of $N$-pointed compact Riemann surfaces with local coordinates $\eta_\blt$, then $\fk X$ has a natural choice of jets: let $\nu_i=\sgm_i^*d\eta_i$.
\end{eg}

\begin{eg}
Assume $\fk X=(\pi:\mc C\rightarrow\mc B;\sgm_\blt)$ is $N$-pointed. For each $i$, we let $\mc B^i$ be the open subset of non-zero vectors of the line bundle $\sgm_i^*\omega_{\mc C/\mc B}$.  Let $p_i:\mc B^i\rightarrow\mc B$ be the projection sending the vectors to their initial points. Using these projections, we define the relative product $\mc B^\triangle=\mc B^1\times_{\mc B}\mc B^2\times_{\mc B}\cdots \times_{\mc B}\mc B^N$, i.e., the closed submanifold of all $(\gamma^1,\dots,\gamma^N)\in \mc B^1\times\cdots\times \mc B^N$ satisfying $p_1(\gamma^1)=\cdots= p_N(\gamma^N)$. Let $p:\mc B^\triangle\rightarrow\mc B$ be the natural projection defined by $p_1,\dots,p_N$. Then we may pull back $\fk X$ along $p:\mc B^\triangle\rightarrow\mc B$ to obtain an $N$-pointed $\fk X^\triangle=(\pi:\mc C^\triangle\rightarrow\mc B^\triangle;\sgm_\blt^\triangle)$. More precisely: we let $\mc C^\triangle=\mc C\times_{\mc B}\mc B^\triangle$ which can be considered as a submanifold of $\mc C\times\mc B^\triangle$. $\sgm_i^\triangle$ is determined by $\sgm_i^\triangle(b^\triangle)=(\sgm_i(p(b^\triangle)),b^\triangle)$ for every $b^\triangle\in\mc B^\triangle$. Then $\fk X^\triangle$ has natural jets $\nu_\blt$ such that for each $b^\triangle\in\mc B^\triangle$, if we consider $b^\triangle=(\gamma^1,\dots,\gamma^N)$  as an element of $\mc B^1\times\cdots\times\mc B^N=\sgm_1^*\omega_{\mc C/\mc B}\times\cdots\times\sgm_N^*\omega_{\mc C/\mc B}$ and set $b=p(b^\triangle)$, then $\nu_i(b^\triangle)$, a cotangent vector of $\mc C^\triangle_{b^\triangle}=\mc C_b\times b^\triangle\simeq \mc C_b$ at $\sgm_i^\triangle(b^\triangle)=(\sgm_i(b),b^\triangle)\simeq \sgm_i(b)$, is $\gamma^i$.
\end{eg}

Let $\fk X$ be $N$-pointed with jets $\nu_\blt$ as above. Fix a $C_2$-cofinite rational VOA $\Ubb$ with non-zero central charge $c_\Ubb$. We assume that $\Ubb$ is holomorphic, i.e., $\Ubb$ has only one simple module which is $\Ubb$ itself. For instance, one can take $\Ubb$ to be the VOA associated to an even self-dual lattice, or the moonshine VOA. By factorization, any space of conformal block associated to $\Ubb$ and a pointed curve has dimension one. Thus the sheaves of conformal blocks of $\Ubb$ are line bundles. We fix the sheaf of conformal blocks of $\Ubb$ associated to $\fk X$ and the $\Ubb$-modules \index{LXU@$\scr L_{\fk X}^\Ubb$} $\Ubb,\dots,\Ubb$:
\begin{align*}
\scr L_{\fk X}^\Ubb=\scr T_{\fk X}^*(\Ubb\otimes\cdots\otimes\Ubb)
\end{align*}
and consider it as a line bundle on $\mc B$.

\subsection*{Flat connections depending on $\nu_\blt$ and a nowhere zero $\uptheta\in\scr L_{\fk X}^\Ubb(\mc B)$}

We assume that there is a nowhere zero section  $\uptheta\in\scr L_{\fk X}^\Ubb(\mc B)$. Then we shall define a flat connection $\nabla^\uptheta$ independent of local coordinates and lifts of tangent vectors. It suffices to define such connection locally. So we assume temporarily that $\mc B$ is Stein and small  enough so that we can choose local coordinates $\eta_\blt$. For each $i$, $\sgm_i^*d\eta_i\in\sgm_i^*\omega_{\mc C/\mc B}(\mc B)$ is nowhere zero. Thus $\sgm_i^*d\eta_i/\nu_i\in\scr O(\mc B)$. 
If $\mu_\blt$ is another set of local coordinates, we define $\alpha_i$ and hence $\alpha_i'(0)\in\scr O(\mc B)$ as in \eqref{eq249}. Then it is easy to see
\begin{align}
\sgm_i^*d\mu_i=\alpha_i'(0)\cdot \sgm_i^*d\eta_i.\label{eq258}
\end{align}
For each section $\yk$ of $\Theta_{\mc B}$, choose a lift $\wtd\yk$. Define $\nabla^{(\eta_\blt)}$ using $\eta_\blt$ and $\wtd\yk$ as in Section \ref{lb84}.  For each section $\upphi$ of $\scr T_{\fk X}^*(\Wbb_\blt)$, let 
\begin{align}
\nabla^\uptheta_\yk \upphi=\nabla_\yk^{(\eta_\blt)}\upphi-\frac c{c_\Ubb}\cdot \frac{\nabla_\yk^{(\eta_\blt)}\uptheta}{\uptheta}\cdot \upphi-\sum_{i=1}^N \Delta_{\Wbb_i}\cdot \yk\Big(\log\frac{\sgm_i^*d\eta_i}{\nu_i}\Big)\cdot \upphi.\label{eq264}
\end{align}

\begin{thm}
$\nabla^\uptheta$ is a flat connection of $\scr T_{\fk X}^*(\Wbb_\blt)$. Moreover,  $\nabla^\uptheta$ depends on the jets $\nu_\blt$ and the nowhere zero section $\uptheta\in\scr L_{\fk X}^\Ubb(\mc B)$ but not on  $\eta_\blt$ or the lift $\wtd\yk$ of $\yk$. Therefore, $\nabla^\uptheta$ can be defined globally without assuming $\mc B$ is Stein or local coordinates exist.
\end{thm}

\begin{proof}
	By Proposition \ref{lb73} and Remark \ref{lb143}, $\nabla^{\uptheta}_\yk$ is independent of $\wtd\yk$.  That $\nabla_{\yk^\vartheta}$ is independent of the local coordinates $\eta_\blt$ follows from \eqref{eq258} and Theorem \ref{lb144}. Note that if $f\uptheta=\nabla_\zk^{(\mu_\blt)}\uptheta$, then $\nabla_\yk^{(\mu_\blt)}\nabla_\zk^{(\mu_\blt)}\uptheta=\nabla_\yk^{(\mu_\blt)}(f\uptheta)=\yk(f)\uptheta+f\nabla_\yk^{(\mu_\blt)}\uptheta$. From this we see
	\begin{align*}
	\yk\Big(\frac{\nabla_\zk^{(\mu_\blt)}\uptheta}{\uptheta} \Big)=\frac{\nabla_\yk^{(\mu_\blt)}\nabla_\zk^{(\mu_\blt)}\uptheta}{\uptheta}-\frac{\nabla_\yk^{(\mu_\blt)}\uptheta}{\uptheta}\cdot \frac{\nabla_\zk^{(\mu_\blt)}\uptheta}{\uptheta} .
	\end{align*}
	With help of this relation and Theorem \ref{lb145}, it is straightforward to check that $\nabla$ is has zero curvature.
\end{proof}

\begin{eg}
In the case that one cannot find a nowhere zero $\uptheta\in\scr L_{\fk X}^\Ubb(\mc B)$, one can consider a ``central extension" of $\fk X$ as follows. Regard $\mc B$  as a closed submanifold of $\scr L_{\fk X}^\Ubb$ consisting of zero vectors. Note that we have a natural projection 
\begin{align*}
p:\scr L_{\fk X}^\Ubb\rightarrow\mc B
\end{align*}
sending each vector to its initial point. We can pull back $\fk X$ along $p:\scr L_{\fk X}^\Ubb-\mc B\rightarrow\mc B$ and obtain a new family $\fk Y$ (with base manifold $\scr L_{\fk X}^\Ubb-\mc B$). One can also pullback the jets of $\fk X$. Then $\fk Y$  has a natural global nowhere section $\uptheta$ of $\scr L_{\fk Y}^\Ubb$.
\end{eg}

\subsection*{Connections depending on $\nu_\blt$ and a projective structure $\fk P$}

Suppose $\fk X$ has a projective structure $\fk P$ and jets $\nu_\blt$,  one can define a connection $\nabla^{\fk P}$ as follows. Choose a lift $\wtd\yk$ of the tangent field $\yk$, and let $h_i$ be as in \eqref{eq259}. Then
\begin{align}
\nabla^{\fk P}_\yk \upphi=\nabla_\yk^{(\eta_\blt)}\upphi-\sum_{i=1}^N\Res_{\eta_i=0} ~\frac c{12}h_i(\eta_i,\tau_\blt)\Sbf_{\eta_i}\fk P d\eta_i-\sum_{i=1}^N \Delta_{\Wbb_i}\cdot \yk\Big(\log\frac{\sgm_i^*d\eta_i}{\nu_i}\Big)\cdot \upphi.\label{eq263}
\end{align}

\begin{thm}
$\nabla^{\fk P}$ is independent of the choice of $\eta_\blt$ and the lift $\wtd\yk$ of $\yk$.
\end{thm}

\begin{proof}
To compare the definition of $\nabla$ using two sets of local coordinates $\eta_\blt$ and $\mu_\blt$, it suffices to assume $\mu_1,\dots,\mu_N$ belong to $\fk P$. Then the coincidence  follows from Theorem \ref{lb144}. When the local coordinates belong to $\fk P$, $\nabla$ is independent of the choice of lift by Remark \ref{lb143} and Lemma \ref{lb74} (or equation \eqref{eq179}). Thus, for a general $\eta_\blt$, the independence on $\wtd\yk$ is also true.
\end{proof}

Note that unlike the previous connection, $\nabla^{\fk P}$ might not be flat.

\section{Functoriality}

Assume for simplicity that $\Vbb$ is $C_2$-cofinite so that the sheaves of conformal blocks are holomorphic vector bundles. If $F_i:\Wbb_i\rightarrow\Mbb_i$ is a homomorphism of $\Vbb$-modules for each $1\leq i\leq N$, then we clearly have an $\scr O_{\mc B}$-module homomorphism $F_\blt^*:\scr T_{\fk X}^*(\Mbb_\blt)\rightarrow \scr T_{\fk X}^*(\Wbb_\blt)$ defined by sending each $\upphi$ to $\upphi\circ(F_1\otimes\cdots\otimes F_n)$. If we have $G_i:\Mbb_i\rightarrow \mathbb P_i$ where each $\mbb P_i$ is also an $\Vbb$-module, then by setting $(GF)_i=G_iF_i$, we have $(GF)_\blt^*=F_\blt^* G_\blt^*$. Moreover, if each $F_i$ is identity, then so is $F_\blt^*$. Thus, we have a contravariant functor $\Wbb_\blt\rightarrow\scr T_{\fk X}^*(\Wbb_\blt)$.

We may also fix $\Vbb$-modules $\Wbb_1,\dots,\Wbb_N$, and consider morphisms between two families of compact Riemann surfaces. To be more precise, if $\fk X^j=(\pi^j:\mc C^j\rightarrow\mc B^j;\sgm_1^j,\dots,\sgm_N^j)$ ($j=1,2$) are families of $N$-pointed compact Riemann surfaces, then a \textbf{morphism} $F:\fk X^1\rightarrow\fk X^2$ is a pair $F=(F_\Crm,F_\Brm)$ where $F_\Crm:\mc C^1\rightarrow\mc C^2$ and $F_\Brm:\mc B^1\rightarrow \mc B^2$ are holomorphic maps, $\pi^2\circ F_\Crm=F_\Brm\circ \pi^1$, $F_\Crm\circ \sgm_i^1=\sgm_i^2$ for each $1\leq i\leq N$, and $F_\Crm$ restricts to an isomorphism of compact Riemann surfaces $\mc C^1_{b}\rightarrow \mc C^2_{F_\Brm(b)}$ for each $b\in\mc B^1$. We will write both $F_\Crm$ and $F_\Brm$ as $F$ for short when no confusion arises. 

We can pull back $\scr T_{\fk X^2}^*(\Wbb_\blt)$ along $F_\Brm$ to get an $\scr O_{\mc B^1}$-module $F^*\scr T_{\fk X^2}^*(\Wbb_\blt)\equiv F_\Brm^*\scr T_{\fk X^2}^*(\Wbb_\blt)$. Thus we have $F^*:\scr T_{\fk X^2}^*(\Wbb_\blt)(V)\rightarrow F^*\scr T_{\fk X^2}^*(\Wbb_\blt)(F^{-1}(V))$ for each open $V\in\mc B^2$. We can define a similar \index{F@$F^\diamond$} map 
\begin{align}
F^\diamond:\scr T_{\fk X^2}^*(\Wbb_\blt)(V)\rightarrow \scr T_{\fk X^1}^*(\Wbb_\blt)(F^{-1}(V))
\end{align}
as follows.  For each $b\in F^{-1}(V)$, $F_\Crm$ restricts to an isomorphism of $N$-pointed fibers $\mc C^1_b\rightarrow\mc C^2_{F(b)}$. This gives a natural isomorphism
\begin{align*}
F_b^\diamond:\scr T_{\fk X^2}^*(\Wbb_\blt)|F(b)\rightarrow \scr T_{\fk X^1}^*(\Wbb_\blt)|b.
\end{align*}
Then $F^\diamond$ is defined such that for each $\upphi\in\scr T_{\fk X^2}^*(\Wbb_\blt)(V)$ and $b\in F^{-1}(V)$, $(F^\diamond\upphi)(b)=F_b^\diamond\upphi(F(b))$. One can write down the explicit formula: Assume the restriction $\fk X^2_V$ admits local coordinates $\eta_\blt^2$. Then one can define local coordinates $\eta_\blt^1$ of the restricted family $\fk X^1_{F^{-1}(V)}$ such that 
\begin{align}
\eta_i^1=\eta_i^2\circ F_\Crm
\end{align}
for each $1\leq i\leq N$. Choose $w\in\Wbb_\blt$. Then
\begin{align}
(F^\diamond\upphi)(\mc U(\eta_\blt^1)^{-1}w)=\upphi(\mc U(\eta_\blt^2)^{-1}w).\label{eq261}
\end{align}
If $G:\fk X^2\rightarrow\fk X^3$ is another morphism, we clearly have
\begin{align}
(GF)^\diamond=F^\diamond G^\diamond.\label{eq262}
\end{align}
Also, if $F$ is the identity map, then $F^\diamond$ is clearly also the identity.

$F^\diamond$ and $F^*$ can be related in the following way. Define an $\scr O_{\mc B^1}$-module \index{zz@$\Phi_F$} isomorphism
\begin{gather}
\Phi_F:F^*\scr T_{\fk X^2}^*(\Wbb_\blt)\xrightarrow{\simeq} \scr T_{\fk X^1}^*(\Wbb_\blt),\nonumber\\
\Phi_F F^*\upphi=F^\diamond\upphi
\end{gather}
for each section $\upphi$ of $\scr T_{\fk X^2}^*(\Wbb_\blt)$. Namely, it is defined by $\Phi_F (f\cdot F^*\upphi)=fF^\diamond\upphi$ for any holomorphic function $f$ of $\mc B^1$. To check that $\Phi_F$ is well-defined and is an $\scr O_{\mc B^1}$-isomorphism, note that for each $b\in \mc B^1$ we define an isomorphism of vector spaces $F_b^*\equiv (F_B)_b^*:\scr T_{\fk X^2}^*(\Wbb_\blt)|F(b)\rightarrow F^*\scr T_{\fk X^2}^*(\Wbb_\blt)|b$ by pullback. Then we can define an isomorphism
\begin{gather*}
\Phi_{F,b}:F^*\scr T_{\fk X^2}^*(\Wbb_\blt)|b\xrightarrow{\simeq} \scr T_{\fk X^1}^*(\Wbb_\blt)|b,\\
\Phi_{F,b}F^*_b=F_b^\diamond.
\end{gather*}
Then it is clear that $\Phi_F(f\cdot  F^*\upphi)(b)=f(b)\cdot (F^\diamond\upphi)(b)=f(b)\cdot F_b^\diamond\cdot \upphi(b)=f(b)\cdot \Phi_{F,b}F^*_b\cdot \upphi(b)$ which depends only on $f(b)\upphi(b)$. So $\Phi_F$ is a well-defined isomorphism of vector bundles whose restriction to each fiber over $b$ is $\Phi_{F,b}$. To summarize, we have
\begin{thm}
For each morphism $F:\fk X^1\rightarrow\fk X^2$ of families of $N$-pointed compact Riemann surfaces, there is an isomorphism of (holomorphic) vector bundles $\Phi_F:F^*\scr T_{\fk X^2}^*(\Wbb_\blt)\xrightarrow{\simeq} \scr T_{\fk X^1}^*(\Wbb_\blt)$ such that $\Phi_F F^*\upphi=F^\diamond\upphi$ for each section $\upphi$ of $\scr T_{\fk X^2}^*(\Wbb_\blt)$, and $F^\diamond\upphi$ is described by \eqref{eq261}. If $F$ is the identity morphism, then $\Phi_F$ is the identity map. If $G:\fk X^2\rightarrow\fk X^3$ is also a morphism, then for each open $W\subset\mc B^3$, the following maps from $\scr T_{\fk X^3}^*(\Wbb_\blt)(W)\rightarrow \scr T_{\fk X^1}^*(\Wbb_\blt)(F^{-1}G^{-1}(W))$ are equal.
\begin{align*}
\Phi_{GF}\cdot (GF)^*=\Phi_F \cdot F^*\cdot  \Phi_G\cdot  G^*.
\end{align*}
\end{thm}
The last equation is due to \eqref{eq262}.

Recall that if $\nabla$ is a connection on $\scr T_{\fk X^2}^*(\Wbb_\blt)$, then its pullback $F^*\nabla=F_\Brm^*\nabla$ is a connection on $F^*\scr T_{\fk X^2}^*(\Wbb_\blt)$ defined by $(F^*\nabla)_{\yk}(F^*\upphi)=F^*(\nabla_{dF(\yk)}\upphi)$ for each section $\upphi$ of $\scr T_{\fk X^2}^*(\Wbb_\blt)$ and each tangent vector $\yk$ of $\mc B^1$. 

Suppose that $\fk X^2$ admits jets $\nu_\blt^2$. Then one can define jets $\nu_\blt^1$ of $\fk X^1$ such that for each $b\in\mc B^1$, the cotangent vector $\nu^1_i(b)$ of $\mc C^1_b$ at $\sgm^1_i(b)$ is $F^*d\nu_i^2(b)$. If $\fk X^2$ also admits a projective structure $\fk P^2$, then one can define a projective chart (and hence a projective structure) $\fk P^1$ consisting of all $(F_\Crm^{-1}(U),\eta\circ F_\Crm)$ where $(U,\eta)$ belongs to $\fk P^1$. (In particular, $U$ is an open subset of $\mc C^2$ and $\eta\in\scr O(U)$ is univalent on each fiber.) Then for the connections defined by \eqref{eq263} (assume $\Wbb_1,\dots,\Wbb_N$ are simple), it is not hard to check that for any tangent field $\yk$ of $\mc B^1$,
\begin{align}
\nabla^{\fk P^1}_\yk=\Phi_F\cdot (F^*\nabla^{\fk P^2})_\yk \cdot \Phi_F^{-1}\label{eq265}
\end{align}
when acting on sections of $\scr T_{\fk X^1}^*(\Wbb_\blt)$.

Alternatively, suppose that, instead of projective structures, we have a global section $\uptheta^2$ of $\scr L_{\fk X^2}^\Ubb$. We can define a global section of $\scr L_{\fk X^1}^\Ubb$ to be $\uptheta^1=F^\diamond\uptheta^2=\Phi_F F^*\uptheta^2$. Then for the flat connections defined by \eqref{eq264}, we also have
\begin{align}
\nabla^{\uptheta^1}_\yk=\Phi_F\cdot (F^*\nabla^{\uptheta^2})_\yk \cdot \Phi_F^{-1}.
\end{align}

\begin{eg}\label{lb148}
Let $G$ be a group of automorphisms of $\fk X$, i.e., we have a homomorphism $G\rightarrow\mathrm{Aut}(\fk X)$. For each $g\in G$, we have an action $g^\diamond=\Phi_g\cdot g^*:\scr T_{\fk X}^*(\Wbb_\blt)(V)\rightarrow \scr T_{\fk X}^*(\Wbb_\blt)(g^{-1}V)$ for each open $V\subset\mc B$, and we have $(gh)^\diamond=h^\diamond g^\diamond$ for each $g,h\in G$. Thus, we have a right action of $G$ on $\scr T_{\fk X}^*(\Wbb_\blt)(\mc B)$.

Suppose that $\fk X$ admits jets $\nu_\blt$ and a projective structure $\fk P$, and both are invariant under the action of $G$. Suppose also that $\mc B$ is simply-connected and $\nabla^{\fk P}$ is flat. We can define the vector space  $\scr C_{\fk X}(\Wbb_\blt)$ of all $\upphi\in\scr T_{\fk X}^*(\Wbb_\blt)(\mc B)$ which are parallel under $\nabla^{\fk P}$, i.e., annihilated by $\nabla^{\fk P}_\yk$ for each tangent field $\yk$ of $\mc B$. Then $\dim\scr C_{\fk X}(\Wbb_\blt)$ is equal to the rank of the vector bundle $\scr T_{\fk X}^*(\Wbb_\blt)$. Moreover, $g\in G\rightarrow (g^{-1})^\diamond$ defines a (left) action of $G$ on  the vector space $\scr C_{\fk X}(\Wbb_\blt)$. This is also true when we not assume the existence of a set of $G$-invariant jets, but assume $\Wbb_1,\dots,\Wbb_N$ are all $\Vbb$ so that $\Delta_{\Wbb_1}=\dots=\Delta_{\Wbb_N}=0$.
\end{eg}

\section{Modular invariance}

Let $\Vbb$ be $C_2$-cofinite and rational. Let $\wtd{\fk Y}=(\Pbb^1;1,0,\infty)$. We associate local coordinates $z$ to $0$ and $z^{-1}$ to $\infty$. Then we can sew $\wtd{\fk Y}$ along $0,\infty$ to get a family $\fk Y=(\mc R\rightarrow\mc D_1^\times;\sigma)$ of $1$-pointed elliptic curves. Recall $\mc D_1^\times$ is the punctured unit open disc. Associate to $1,0,\infty$ simple $\Vbb$-modules $\Wbb,\Mbb,\Mbb'$ where $\Mbb'$ is contragredient to $\Mbb$. Then for each $\uppsi\in\scr T_{\wtd{\fk Y}}^*(\Wbb\otimes\Mbb\otimes\Mbb')$, we have the sewn conformal block $\wtd{\mc S}\uppsi\in \scr T_{\fk Y}^*(\Wbb)(\mc D_1^\times)$ and $\mc S\uppsi=q^{\Delta_\Mbb}\wtd{\mc S}\uppsi$. Let $\wtd{\fk Q}$ be unique projective structure of $\Pbb^1$, i.e., the one containing $(\Cbb,z)$. This in turn gives a projective structure $\fk Q$ of $\fk Y$. More precisely: the local coordinate $z-1$ of $\mbb P^1$ at $1$ extends constantly (with respect to sewing) to a local coordinate $\mu$ of $\fk Y$. $\fk Q$ is the projective structure containing $\mu$. 

Let $\gamma$ be a jet of $\fk Y$ whose value at each $q\in\mc D_1^\times$ is the cotangent vector $d\mu$. Then the definition of the connection $\nabla^{(\mu)}$ of $\scr T_{\fk Y}^*(\Wbb)$ using the local coordinate $\mu$ (as in Section \ref{lb84})  is the same as the connection $\nabla^{\fk Q,\lambda}$ defined by $\fk Q$ and the jet $\gamma$ (as in \eqref{eq263}), and is independent of the choice of lifts. Since the local coordinates of $0,\infty$ belong to $\fk Q$, by Theorem \ref{lb85} and Remark \ref{lb138}, $\mc S\uppsi$ is parallel under $\nabla^{\fk Q,\mu}$.

Let $\mbb H$ be the (open) upper half plane of $\Cbb$. Define
\begin{align*}
\mc B=\mbb H\times\Cbb^\times.
\end{align*}
Define an action of $\Zbb^2$ on $\Cbb\times\mc B$ such that for each $a,b\in\Zbb$ and $(z,\tau,\zeta)\in\Cbb\times\mbb H\times\Cbb^\times$, $(a,b)(z,\tau,\zeta)=(z+a+b\tau,\tau,\zeta)$. Then we have a (universal) family of $1$-pointed elliptic curves
\begin{align*}
\fk X=(\pi:\mc C\rightarrow\mc B;\sgm)
\end{align*}
where 
\begin{align*}
\mc C=(\Cbb\times\mc B)/\Zbb^2,
\end{align*}
the projection $\pi$ is defined by the standard one $\Cbb\times\mc B\rightarrow\mc B$, and $\sgm$ comes from the section $\mc B\rightarrow\Cbb\times \mc B$, $(\tau,\zeta)\mapsto (0,\tau,\zeta)$. Let $\Gamma$ be the modular group, i.e.,
\begin{align*}
\Gamma=SL_2(\Zbb)=\bigg\{\begin{pmatrix}
a & b\\
c & d
\end{pmatrix}
:a,b,c,d\in\Zbb,~ad-bc=1
\bigg\}.
\end{align*}
Then we have a group action of $\Gamma$ on $\fk X$ such that for any $g=\begin{pmatrix}
a & b\\
c & d
\end{pmatrix}\in \Gamma$, 
\begin{gather*}
g:\mc B\rightarrow\mc B,\qquad g(\tau,\zeta)=\Big(\frac{a\tau+b}{c\tau+d},(c\tau+d)\zeta  \Big),
\end{gather*}
and the action of $g$ on $\mc C$ descends from the one on $\Cbb\times\mc B$ determined by
\begin{gather*}
g:\Cbb\times\mc B\rightarrow \Cbb\times\mc B,\qquad g(z,\tau,\zeta)=\Big(\frac z{c\tau+d},\frac{a\tau+b}{c\tau+d},(c\tau+d)\zeta  \Big).
\end{gather*}
Then we have  $\fk X/\Gamma=(\mc C/\Gamma\rightarrow \mc B/\Gamma;\sgm)$ where $\mc B/\Gamma$ is the (fine) moduli space of $1$-pointed elliptic curves with jet.\footnote{One reason to work with $\mc B$ instead of $\mbb H$ and to consider jets is that pointed elliptic curves with jet have trivial automorphism groups. So we can have a fine moduli space.} The jet of $\fk X/\Gamma$ come from $\nu$ of $\fk X$ which will be described later. By example \ref{lb148}, $\Gamma$ acts on $\scr T_{\fk X}^*(\Wbb)(\mc B)$.

There are two natural choices of flat connections on $\scr T_{\fk X}^*(\Wbb)$. We have a morphism $F:\fk X\rightarrow\fk Y$ described as follows. As a holomorphic map between base manifolds, we have
\begin{gather*}
F:\mc B=\mbb H\times\Cbb^\times\rightarrow \mc D_1^\times,\qquad (\tau,\zeta)\mapsto \exp(2\im\pi\tau).
\end{gather*}
The map $F:\mc C\rightarrow\mc R$ is defined such that  for each $(\tau,\zeta)\in\mc B$, the map 
\begin{align*}
\Cbb\rightarrow\Pbb^1,\qquad z\mapsto \exp(2\im\pi z)
\end{align*}
descends to $\Cbb\rightarrow \mc R_{\exp(2\im\pi\tau)}$, and furthermore descends to $\mc C_{(\tau,\zeta)}\xrightarrow{\simeq}\mc R_{\exp(2\im\pi\tau)}$. Then we can pullback the projective structure $\fk Q$ and the jet $\lambda$ of $\fk Y$ to $\fk P',\nu'$ of $\fk X$ as described above \eqref{eq265}. Then by \eqref{eq265}, $\nabla^{\fk P',\nu'}$ is equivalent to $F^*\nabla^{\fk Q,\lambda}$ via $\Phi_F:F^*\scr T_{\fk Y}^*(\Wbb)\xrightarrow{\simeq}\scr T_{\fk X}^*(\Wbb)$, and is therefore flat since $\nabla^{\fk Q,\lambda}$ is acting on a $1$-dimensional complex manifold. So 
\begin{align*}
F^\diamond\mc S\uppsi=\Phi_F F^*\mc S\uppsi=e^{2\im\pi\tau\Delta_\Mbb}\cdot \Phi_F F^*\wtd{\mc S}\uppsi
\end{align*}
is a global section of $\scr T_{\fk X}^*(\Wbb)$ parallel under $\nabla^{\fk P',\nu'}$. (Recall that $F^\diamond$ is described by \eqref{eq261}.) By factorization, $\scr T_{\fk X}^*(\Wbb)$ is $\scr O_{\mc B}$-generated by, and hence the  vector space of $\nabla^{\fk P',\nu'}$-parallel sections in $\scr T_{\fk X}^*(\Wbb)(\mc B)$ is spanned by all $F^\diamond\mc S\uppsi$ where $\uppsi\in\scr T_{\wtd{\fk Y}}^*(\Wbb\otimes\Mbb\otimes\Mbb')$ and $\Mbb$ is a simple $\Vbb$-module.

Unfortunately, neither $\fk P'$ nor $\nu'$ is modular invariant (i.e. $\Gamma$-invariant). As a consequence, $\nabla^{\fk P',\nu'}$ is not modular invariant. To get modular invariant ones, we let $\fk P$ be the projective structure of $\fk X$ whose pull back to the family $\Cbb\times\mc B\rightarrow\mc B$ is the standard one, i.e., its restriction to each fiber $\Cbb$ is the one containing $(\Cbb,z)$. We let $\nu$ be the jet of $\fk X$ to be 
\begin{align*}
\nu=\zeta\cdot dz
\end{align*}
i.e., for each $(\tau,\zeta)\in\mc B$, $\nu(\tau,\zeta)$ is the cotangent vector $\zeta dz$ of the fiber $\mc C_{(\tau,\zeta)}$ at $0$ (when lifted to $\Cbb$) where $dz$ is defined by the standard coordinate $z$ of $\Cbb$. Then both $\fk P$ and $\nu$ are modular invariant. So is $\nabla^{\fk P,\nu}$.

\begin{thm}
Let $\tau,\zeta$ also denote the standard coordinates of $\mbb H,\mbb C^\times$ respectively. Then, when acting on sections of $\scr T_{\fk X}^*(\Wbb)$, we have
\begin{gather*}
\nabla_{\partial_\tau}^{\fk P,\nu}=\nabla_{\partial_\tau}^{\fk P',\nu'}+\frac{\im c\pi}{12}\id,\\
\nabla_{\partial_\zeta}^{\fk P,\nu}=\nabla_{\partial_\zeta}^{\fk P',\nu'}+\frac{\Delta_\Wbb}\zeta\id.
\end{gather*}
\end{thm}
As an immediate consequence,  $\nabla^{\fk P,\nu}$ is also flat.

\begin{proof}
$\nu'$ is the differential of $\exp(2\im\pi z)$ at $z=0$. So $\nu'=2\im\pi dz$, and hence $\nu=\frac{\zeta}{2\im\pi}\nu'$. Hence  $\partial_\zeta(\log(\nu/\nu'))=\zeta^{-1}$. This, together with \eqref{eq263}, shows the second the identity. Since $\partial_\tau(\log(\nu/\nu'))=0$, we have $\nabla_{\partial_\tau}^{\fk P',\nu}=\nabla_{\partial_\tau}^{\fk P',\nu'}$. Thus, it suffices to prove $\nabla_{\partial_\tau}^{\fk P,\nu}=\nabla_{\partial_\tau}^{\fk P',\nu}+\frac{\im c\pi}{12}\id$. Let $\eta$ be the local coordinate of $\fk X$ defined by the standard coordinate $dz$ of $\Cbb$. Let $\eta'$ be the pullback of $\mu$ along $F$, i.e., $\eta'=\exp(2\im\pi z)-1$. Then $\eta$ belongs to $\fk P$ and $\eta'$ belongs to $\fk P'$. Let $\wtd\yk$ be a lift of $\partial_\tau$, and assume its expression near $\sgm(\mc B)$ is
\begin{align*}
h(\eta,\tau,\zeta)\partial_\eta+\partial_\tau
\end{align*}
where the partial derivatives are defined by the coordinates $(\eta,\tau,\zeta)$. Then by Theorem \ref{lb144} and relation \eqref{eq263}, we have $\nabla_{\partial_\tau}^{\fk P,\nu}-\nabla_{\partial_\tau}^{\fk P',\nu}=f\id$ where
\begin{align*}
f=\Res_{\eta=0}~\frac c{12}h(\eta,\tau,\zeta)\Sbf_\eta \eta' d\eta.
\end{align*}

It is easy to calculate that $\Sbf_\eta\eta'=\Sbf_z(\exp(2\im\pi z)-1)=(2\im\pi)^2-\frac 32 (2\im\pi)^2=2\pi^2$. We can pullback $\wtd\yk$ to a global meromorphic tangent field of $\Cbb\times\mc B$ whose poles are in $\Zbb^2\cdot (\{0\}\times\mc B)$. Denote this pullback also by $\wtd\yk$, and notice $\eta$ is just the standard coordinate $z$ of $\Cbb$, we have $\wtd\yk=h(z,\tau,\zeta)\partial_z+\partial_\tau$ and $h$ is a meromorphic function on $\Cbb\times\mc B$ with poles in $\Zbb^2\cdot (\{0\}\times\mc B)$. Moreover, $\wtd\yk$ is invariant under the action of $\Zbb^2$. From this it is easy to see that
\begin{align*}
h(z+1,\tau,\zeta)=h(z,\tau,\zeta),\qquad h(z+\tau,\tau,\zeta)=h(z,\tau,\zeta)+1.
\end{align*}
Let $\gamma_\tau$ be an anticlockwise parallelogram of $\Cbb$ around $0$ described by $A_\tau\rightarrow B_\tau\rightarrow C_\tau\rightarrow D_\tau\rightarrow A_\tau$, where $A_\tau=-0.5-0.5\tau$, $B_\tau=0.5-0.5\tau$, $C_\tau=0.5+0.5\tau$, $D_\tau=-0.5+0.5\tau$. Then by the above relation,
\begin{gather*}
\int_{A_\tau B_\tau}hdz+\int_{C_\tau D_\tau}hdz=\int_{A_\tau B_\tau}hdz-\int_{A_\tau B_\tau}(h+1)dz=-\int_{A_\tau B_\tau}dz=-1,\\
\int_{B_\tau C_\tau}hdz+\int_{D_\tau A_\tau}hdz=\int_{A_\tau D_\tau}hdz-\int_{A_\tau D_\tau}hdz=0.
\end{gather*}
So
\begin{align*}
\Res_{z=0}~hdz=\frac 1{2\im\pi}\oint_{\gamma_\tau}hdz=-\frac 1{2\im\pi}.
\end{align*}
Thus $f=\Res_{z=0}~\frac c{12} h\cdot 2\pi^2dz=\frac{\im c\pi}{12}$, which completes the proof.
\end{proof}

\begin{co}
For any simple $\Mbb$ and any $\uppsi\in\scr T_{\wtd{\fk Y}}^*(\Wbb\otimes\Mbb\otimes\Mbb')$, 
\begin{align}
\zeta^{-\Delta_\Wbb} \exp\Big(-\frac{\im c\pi}{12}\tau\Big)F^\diamond\mc S\uppsi\label{eq266}
\end{align}
is a multivalued (with respect to $\zeta$) global section of $\scr T_{\fk X}^*(\Wbb)$ parallel under the modular invariant flat connection $\nabla^{\fk P,\nu}$. Moreover, any such $\nabla^{\fk P,\nu}$-parallel section of $\scr T_{\fk X}^*(\Wbb)$ is a $\Cbb$-linear combination of sections of this form.
\end{co}
If we let $q_\tau=F(\tau)=\exp(2\im\pi\tau)$, then the projective factor $\exp(-\frac{\im c\pi}{12}\tau)$ becomes the celebrated $q_\tau^{-\frac c{24}}$.

\begin{co}[Modular invariance]
For any $\upphi$ in the form \eqref{eq266}, and for any $g\in\Gamma=SL_2(\Zbb)$, $g^\diamond\upphi$ is also a $\Cbb$-linear combination of sections of the form \eqref{eq266}.
\end{co}

\printindex

\noindent {\small \sc Yau Mathematical Sciences Center, Tsinghua University, Beijing, China.}

\noindent {\textit{E-mail}}: binguimath@gmail.com
\end{document}